\documentclass{amsart}
\usepackage[utf8]{inputenc} 
\usepackage{amsmath, amssymb, amsthm, amsfonts,bigints}
\usepackage{csquotes} 
\usepackage[american]{babel}
\usepackage{color}
\usepackage{tikz-cd}
\usepackage[isbn=false, doi=false, url=false] {biblatex}
\renewbibmacro{in:}{} 
\bibliography{kdve} 
\usepackage{enumitem} 
\usepackage{hyperref}

\newtheorem{theorem}{Theorem}
\newtheorem{lemma}[theorem]{Lemma}
\newtheorem{definition}[theorem]{Definition}
\newtheorem{proposition}[theorem]{Proposition}
\newtheorem{corollary}[theorem]{Corollary}
\newtheorem{remark}[theorem]{Remark}

\numberwithin{theorem}{section}
\numberwithin{equation}{section}
\numberwithin{claim}{section}

\newcommand{\N}{\mathbb{N}}
\newcommand{\Z}{\mathbb{Z}}

\newcommand{\R}{\mathbb{R}}
\newcommand{\X}{\mathbb{X}}

\newcommand{\C}{\mathbb{C}}
\newcommand{\Sc}{\mathcal{S}}
\newcommand{\T}{\mathcal{T}}
\newcommand{\I}{\mathfrak{I}}
\newcommand{\flux}{\operatorname{Fl}} 
\newcommand{\sech}{\operatorname{sech}}
\newcommand{\KdV}{{\operatorname{KdV}}}
\newcommand{\mKdV}{{\operatorname{mKdV}}}

\newcommand{\NLS}{\operatorname{NLS}}
\newcommand{\AKNS}{\operatorname{AKNS}}
\newcommand{\Gardner}{{\operatorname{Gardner}}}
\newcommand{\Wadati}{{\operatorname{Wadati}}}
\newcommand{\Magri}{{\operatorname{Magri}}}

\newcommand{\real}{\operatorname{Re}}
\newcommand{\im}{\operatorname{Im}}
\DeclareMathOperator{\trace}{\rm tr\,} 
 
\newcommand{\supp}{\operatorname{supp}}
\newcommand{\id}{\operatorname{Id}}
\newcommand{\extra}[1]{}

\title{Well-posedness for the KdV hierarchy}
\author{Friedrich Klaus} 
\address{Karlsruhe Institute of Technology, 76131 Karlsruhe, Germany}
\email{friedrich.klaus@kit.edu} 
\author{Herbert Koch}
\address{Mathematisches Institut der Universität Bonn, Endicher Allee 60, 53115 Bonn, Germany}
\email{koch@math.uni-bonn.de}
\author{Baoping Liu}
\address{Beijing International Center for Mathematical Research,
 Peking University,
 Beijing,
  China}
\email{baoping@math.pku.edu.cn}

\begin{document}

\begin{abstract} We prove a version of well-posedness for all equations of the KdV hierarchy in $H^{-1}$. Ingredients are 
\begin{enumerate}
    \item The Miura map which allows to define the Gardner hierarchy through the generating function of the energies so that the $N$th Gardner equation is equivalent to the $N$th KdV equation. 
    \item A rigorous relation between the generating functions of the energies and the KdV resp. Gardner Hamiltonians. 
    \item Kato smoothing estimates for weak solutions and approximate flows. 
\end{enumerate}
\end{abstract}
\maketitle
\tableofcontents

\section{Introduction}

The  Korteweg-de Vries (KdV) equation
\begin{equation} \label{eq:KdV}  u_t+ u_{xxx}- 6 u u_x = 0 \end{equation} 
is a fascinating and basic object in diverse areas: It is a generic asymptotic equation for propagating waves, and it has a deep algebraic structure visible in the existence of an infinite sequence of formally conserved energies $H_n^{\KdV}$.
The KdV equation is a member of the KdV hierarchy: It is the Hamiltonian equation with 
the Hamiltonian function 
\[ H_1^{\KdV}(u) = \frac12  \int u_x^2+ 2u^3 dx  \]  
with respect to the Gardner Poisson structure (defined in  \eqref{eq:defGardner}).

The study of rough initial data is related to weaker assumptions on the frequency localization for the validity of the asymptotic equation, whereas the study of higher order equations is related, or more precisely a necessary ingredient, to larger time scales.

Well-posedness for the KdV equation itself has been an active and stimulating area for the last three decades.  Well-posedness results on $H^s$ spaces have been proven by Kato~\cite{MR535697} (local in $H^s(\R)$, $ s >3/2$ and global in $H^2(\R)$) by energy methods.  Bourgain~\cite{MR1215780} introduced the $X^{s,b}$ space and showed global well-posedness in $L^2(\R)$ and $L^2(\mathbb{S})$. Zhou~\cite{Zhou} proved uniqueness of weak solutions in $L^2(\mathbb{R})$. Kenig, Ponce and Vega~\cite{MR1329387} proved sharp bilinear estimate in $X^{s,b}$ space and thus showed local well-posedness of KdV in $H^{s}(\R), s>-\frac34$ and $H^s(\mathbb{S}), s>-\frac12$. 
The results were extended to global well-posedness by Colliander, Keel, Staffilani, Takaoka and Tao~\cite{MR1969209}, via the construction of almost conserved quantities.  The local existence in $H^{-\frac34}(\R)$ was shown by Christ, Colliander, Tao~\cite{MR2018661}, and global existence  in $H^{-3/4}(\R)$ was proved independently by Guo~\cite{MR2531556} and Kishimoto~\cite{MR2501679}.
Liu~\cite{MR3292346} used modified energies to establish local in time apriori bounds in $H^s(\R)$ for $ s \ge -\frac45$. Buckmaster and Koch~\cite{MR3400442} and later Killip, Visan and Zhang~\cite{MR3820439}% 
and
Koch and Tataru~\cite{MR3874652} proved uniform in time apriori estimates in $H^{-1}$. The apriori estimates remain true for all higher flows. 
The estimates in \cite{MR3400442}  are strong enough to construct weak solutions (which was done in \cite{MR3400442} for KdV in the nonperiodic case, with a weaker notion of weak solution  through a compactness argument). Finally, Killip and Visan~\cite{MR3990604} proved global well-posedness for initial data in $H^{-1}(\R)$ and $H^{-1}(\mathbb{S})$, in the sense that the solution map extends uniquely from Schwartz space to a jointly continuous flow map $\Phi: \R\times H^{-1}\rightarrow H^{-1}$.  This also provides a new proof of the well-posedness in the periodic case, which was first shown by Kappeler and Topalov~\cite{MR2267286}. The results in $H^{-1}$ are sharp since Molinet~\cite{MR2830706, MR2927357} showed that the solution map can not be continuously extended to $H^s, s<-1$ in both periodic and nonperiodic case.

For higher order equations in the KdV hierarchy, Saut~\cite{Saut-KdV} proved  global existence of persistent solutions of the $k$-th KdV equation for initial data in $H^k(\R)$, but the uniqueness was left open. Kenig, Ponce and Vega~\cite{KPV-higher1, KPV-higher2} studied generalized higher order equations, not necessarily integrable, and showed local well-posedness for initial data in weighted spaces.   Pilod~\cite{Pilod} showed that all higher equations are 
ill-posed in any $H^s(\R), s\in \R$, in the sense that data-to-solution map is not $C^2$ at the origin.   Grünrock~\cite{Grunrock} proved local well-posedness for the KdV hierarchy in  Fourier–Lebesgue spaces. Kenig and Pilod~\cite{Kenig-Pilod}, Guo, Kwak and Kwon~\cite{GuoKwakKwon}  showed global well-posedness for the general 5th order KdV in energy space $H^2(\R)$. 
Bringmann, Killip and Visan~\cite{BKV} proved global well-posedness for 5th order KdV for initial data in $H^{-1}(\R)$ by a new strategy that integrates dispersive effects into the method of commuting flows. 
In the periodic case, Kappeler and Molnar~\cite{Kappeler-Molnar}
showed that the 5th KdV equation is $C^0$ well-posed
in $H^s(\mathbb{S})$ if $s\geq 0$, and strongly ill-posed if $s<0$, in the sense that data-to-solution map does not admit a continuous extension to $H^s(\mathbb{S}), s<0$.  

Beyond the KdV hierarchy 
there have been striking new developments at the interface of PDE-techniques
and integrable structures: The work of Killip and Visan on the KdV equation~\cite{MR3990604} introduced a new perspective and powerful technique which has motivated the  study of a number of integrable problems: sharp global well-posedness for cubic NLS and mKdV 
in $H^s(\R)
, s>-\frac12$~\cite{harropgriffiths2020sharp}, for the derivative NLS~\cite{harropgriffiths2023global} (which uses  crucially the work of Bahouri and Perelman  \cite{MR4448993} as well as the equicontinuity of
Harrop-Griffiths,
Killip and Visan~\cite{MR4565673}) and for the Benjamin-Ono 
equation
by Killip, Laurens and Visan~\cite{killip2023sharp}. 

G\'erard and coworkers introduced new integrable
PDEs,
the cubic Sz\"ego~\cite{MR3652066} being only the first, with striking new ideas like an explicit 
formula for solutions to the Benjamin-Ono equation~\cite{gerard2022explicit}, which in turn become a crucial element in~\cite{killip2023sharp}.  They constructed global Birkhoff coordinates for the Benjamin-Ono equation on the torus~\cite{BObirkhoff}, which was used to prove sharp well-posedness for Benjamin-Ono equation in $H^s(\mathbb{S}), s>-\frac12$~\cite{BOActa}. See also the work of Gassot~\cite{Gassot} on sharp well-posedness for the third equation of the integrable Benjamin-Ono hierarchy on the torus.

Clearly this list of results is incomplete and there are many omissions. Similarly we omit a presentation of consequences, possibly the most important being on random initial data. 

We conceive these developments as evidence that 
beyond the single results a new picture of integrable PDEs seems to be emerging, which is still incomplete and to which we hope to contribute.

 \bigskip

The Schrödinger operator 
\begin{equation}  L^{\KdV} \phi = -\phi'' + u \phi \label{eq:schroedinger} \end{equation}
is the Lax operator of the KdV hierarchy. The KdV Hamiltonians are defined as the coefficients of the asymptotic series for the logarithm of the transmission coefficient (see \eqref{eq:asymptotic})  
and they all Poisson  commute.    

The main result of this paper is  well-posedness, more precisely 
\begin{enumerate}
    
    \item Existence and uniqueness for the $N$th Gardner equation \eqref{eq:GardnerN} in $ C(\R, H^{N}(\mathbb{X}))$, $ \mathbb{X}= \R$ or $ \mathbb{X}= \mathbb{S}^1$ 
    for initial data $u_0 \in H^{N}(\mathbb{X}))$ which has a continuous extension to $ H^{N-1} \ni w_0 \to w \in C(\R; H^{N-1})$.  Analogous results  for the $N$th KdV equation \eqref{eq:HN} in spaces with one derivative less are an immediate consequence.
     See Theorem~\ref{thm:KdV} for $N=1,$ and Theorem~\ref{thm:high}  for $N\geq 2$.  
    \item In the case $\mathbb{X}= \R$, Kato smoothing estimates and tightness of weak solutions with initial data in $L^2(\R)$.
    \item \label{item:uniqueness} Uniqueness for a class of weak solutions  to the $N$th KdV  equation with  initial data in $L^2(\R)$.
\end{enumerate} 
The following theorem collects some of these statements. 
\begin{theorem} \label{th:main} 
Let $N \geq 1$. Suppose that 
\begin{equation}\label{eq:kdvregweak}  u\in L^\infty (\R;L^2),  \quad  \partial_x^{N-1} u \in L^2_{loc}( \R \times \R) \end{equation} 
and that it is a weak solution to the $N$th KdV equation \eqref{eq:HN}. Suppose that in addition for every $t_0 \in \R$  
\begin{equation}\label{eq:kdvlimit}   \limsup_{x_0 \to \pm \infty}  \Vert  u^{(N-1)} \Vert_{L^2( (t_0,t_0+1)  \times (x_0-1,x_0+1))} = 0. \end{equation}  Then $u \in C([0,\infty); L^2(\R))$ and we  denote the initial trace by $u_0:=u(0) \in L^2$. The Kato smoothing estimate 
\begin{equation} \label{eq:kdvKato}  \sup_t \Vert u(t) \Vert_{H^{-1}} +   \sup_{x_0} \Vert \sech(   x- \kappa^{2N}  t-x_0) u^{(N-1)} \Vert_{L^2(\R\times \R)} \le  c( \Vert u_0 \Vert_{H^{-1}})  \Vert u_0 \Vert_{H^{-1}} \end{equation}  
holds for all $ \kappa \ge \kappa_0( \Vert u_0 \Vert_{H^{-1}})   $. 

Given $u_0 \in L^2(\R)$ there is a unique weak solution to the $N$th KdV equation which satisfies \eqref{eq:kdvregweak} and \eqref{eq:kdvlimit}. 
The map $L^2(\R) \ni u_0 \to u\in C(\R; H^{-1}(\R))   $ extends to a continuous map  \[  H^{-1} \ni u_0 \to  u \in  C(\R; H^{-1} )\]  
to weak solutions which satisfy in addition \eqref{eq:kdvlimit} and \eqref{eq:kdvKato}. 
\end{theorem}
While the regularity assumption in \eqref{eq:kdvregweak} looks inconsistent with the Kato smoothing estimate in \eqref{eq:kdvKato} we need $L^2$ regularity to get equivalence of weak solutions to the KdV and the Gardner hierarchy, see Theorem \ref{thm:equivalenceweaksolutionsuw}, and the condition $u \in L^\infty L^2$ is used only to prove the equivalence. We prove the theorem by studying the analogous statement for the Gardner hierarchy (Theorem \ref{thm:Gardneruniqueness}).

Why should one be interested in higher KdV equations? From applications one would like to explain why the Korteweg-de Vries equation provides a good description of nonlinear waves in the KdV regime. Typical results are consistency results up to a certain time scale for localized well-prepared initial data (see \cite{MR1780702, MR1780703} and the references in these papers)  - whereas the KdV solitons seem to be relevant in many situations like tsunamis  despite the interaction with other waves and despite large time scales and amplitudes. 

It is striking that the theory of integrable systems provides a very detailed and geometric picture of the simultaneous dynamics of all the KdV flows. Inverse scattering allows a  linearization of the evolutions \cite{MR512420,MR874343} and hence implies  well-posedness of the hierarchy on the Schwartz space. The geometric contents is clearly visible in the relation of the Korteweg-de Vries hierarchy to the diffeomorphism group on $\R^1$ and the torus $\mathbb{S}^1$, more precisely on a central extension, the Virasoro-Bott group (see \cite{MR2456522, MR2326419,MR2216268}). 
Ignoring the topology, the tangent space at the identity of the Virasoro-Bott group is given by the pairs $(v\partial , g)$ of vector fields times $\R$. It is a Lie algebra.
The set of Lax operators 
\[ 2a \partial^2 + u \] 
can be understood as the dual  space of the Virasoro algebra, the Lie algebra which is the tangent space of the Virasoro-Bott group at the identity, on which the Virasoro-Bott group acts by the coadjoint representation. The orbit structure is well understood in the torus case and can be classified in terms of the spectrum of the Lax operator. The Korteweg-de Vries equation can be realized through moment maps and the natural biHamiltonian structure allows to construct a countable sequence of Poisson commuting Hamiltonians.  
It is tempting to ask whether larger time scales for asymptotic equations can be understood in terms of this striking symmetry of the KdV hierarchy.

Another geometric interpretation of the KdV hierarchy is as flow on restricted Grassmanians \cite{MR783348, MR900587}. This is the origin of the ubiquitous $\tau$ function \cite{MR730247,MR1736222}  and the bilinear relation of Hirota  \cite{MR2085332}. 
Let $(t_j)_{j\in \N}$ be a sequence with only finitely many nonzero components. We denote by 
$u(\cdot,t_1,t_2. \dots)\in H^{-1}$ the function resp. distribution obtained from $u_0 \in H^{-1}$ by moving the times $t_j$ along the $j$th KdV flow. This is well defined since the flows commute. The $\tau$ function satisfies 
\[ \partial_x^2 \ln \tau = u \] 
and $u$ (and hence the second derivative of $\ln \tau$) is defined for such sequences $(t_j)$ provided $u_0 \in H^{-1}(\R)$. The expression  
\[ \frac{\partial^2}{\partial t_i \partial t_j} \ln \tau \] 
is  a differential polynomial (see \cite{MR4222602}) and can be evaluated for $u\in H^{\infty}$. It is not difficult to see that it can be integrated and hence a $\tau$ function  exists in this situation. It seems a natural question  whether a $\tau$ function can be defined for $u \in H^{-1}$ (as unique continuous extension of a $\tau$ function on $H^\infty$) and whether related objects  like  vertex operators  can be defined for $u \in H^{\infty}$ or even for $u \in H^{-1}$.

In this paper we study simpler and more basic questions. Crucial points are 

\begin{enumerate} 
\item 
Rigorous estimates for the difference between the generating function of the KdV Hamiltonians and the partial sums for Sobolev functions.
\item A study of the Miura map 
\[    w \to u=w_x + 2\tau w + w^2 \] 
resp. the operator factorization  
\[    (\partial +  \tau + w)(-\partial +  \tau +w ) = -\partial^2 + u + \tau^2. \]
We completely characterize the global mapping properties. The hardest part is a bound of $\Vert w \Vert_{L^2}$ in terms of $ \Vert u \Vert_{H^{-1}}$ and the distance of the ground state energy of $-\partial^2+ u $ and $-\tau^2$.
\item The Miura map allows to translate well-posedness questions for the KdV hierarchy to the Gardner hierarchy for $w$, which has better properties. The Miura map itself and its inverse enter at a number of points. 
\item We prove uniqueness of rough weak solutions for the Gardner equations, and not only continuous extensions of flows on more regular function spaces.  
\end{enumerate}

\subsubsection*{Acknowledgements} F.K. was supported by the Deutsche Forschungsgemeinschaft (DFG, German Research Foundation) – Project- ID 258734477 – SFB 1173. H.K. was supported by the DFG through the Hausdorff Center for Mathematics under Germany's Excellence Strategy - GZ 2047/1, Projekt-ID 390685813 Hausdorff Center for Mathematics and   Project ID 211504053 - SFB 1060. B.L. was supported by the NSF of China (No. 12571254, 12341102).
Work for this project was done during a Research in Pairs stay of the first and second author in Oberwolfach.

\section{Proof of Theorem \ref{th:main}}

\label{sec:proof} 

This section explains the structure and the strategy. We will prove the main theorems using various results from later sections. This will make the structure highly modular, it weakens the dependencies of the later sections and allows for a coherent presentation of the proof. 
A large part is the same on the line and on the circle, and we denote  $ \mathbb  X= \mathbb S^1$ or $\mathbb \R$ if statements are true for both spaces.

\subsection{The generating function  of the KdV hierarchy} 
 We formulate the setting  of  the Korteweg-de Vries hierarchy and describe our approach. More details and proofs about the structure can be found in Section \ref{sec:AKNS}.
A central quantity of the scattering theory of $L^{\KdV}$ is the transmission coefficient. To avoid technical complications in the discussion, we take Schwartz functions as potentials in this paragraph whenever needed. 
For $z$ in the upper half plane, the equation 
\[ L^{\KdV} \phi = z^2 \phi \] 
has a unique solution $\phi_l$ called left Jost function, normalized by 
\begin{equation}  \lim_{x\to - \infty} \phi_l e^{izx} = 1 \end{equation} 
The transmission coefficient is the meromorphic function defined by 
\begin{equation}\label{def:TKDV}
(T^{\KdV}(z))^{-1} = \lim_{x\to \infty} \phi_l e^{izx} \end{equation}
on the upper half plane with the poles given by the square root of the eigenvalues of the Schrödinger operator. It satisfies for Schwartz potentials 
\begin{equation}\label{eq:tkdv}  |T^\KdV(\xi) | \le 1   \qquad \text{ for } \xi \in \R \end{equation} 
Its logarithm can be related to an asymptotic series 

\begin{equation}\label{eq:asymptotic}  \frac{i}2 \ln T^{KdV}(z) \sim \sum_{n=-1}^\infty H_n^{\KdV} (2z)^{-2n-3}  \end{equation}
 where 
 \begin{equation}\label{eq:kdvhamiltonian}  
\begin{split} 
	H_{-1}^{\KdV}\, & =\frac{1}{2}\int u d x, \qquad H_{0}^{\KdV}\, =\frac{1}{2}\int u^{2} d x,
	\\
	H_{1}^{\KdV}\, & =\frac{1}{2}\int u_{x}^{2}+2u^{3} d x, \qquad H_2^{\KdV}\, = \frac{1}{2}\int u_{xx}^2 + 10u_x^2 u + 5 u^4 dx,
\end{split}
\end{equation} 
see \cite{MR783348}   and \cite{MR3874652}.  The Hamiltonians are integrals over differential polynomials (polynomials in $u$ and its derivatives) with the structure
 \[ H_n^{\KdV}(u)  = \frac12 \int |u^{(n)}|^2 dx +  O\big( (1+ \Vert u \Vert_{H^{-1}})^{n-1} \Vert u \Vert_{H^{-1}} \Vert u \Vert_{H^n}^2\big). \]
 The precise meaning of \eqref{eq:asymptotic} is the content of Proposition \ref{eq:kdv:Lipschitz} below:  We define the difference Hamiltonian as $(2z)^{2N+3}$ times the difference to the partial sum  
  \begin{equation} \label{eq:tauN} \mathcal{T}^{\KdV} _N(z,u) :=  (2z)^{2N+3}  \left(\frac{i}2 \ln T^\KdV(z) - \sum_{n={-1}}^N H^{\KdV}_n( 2z)^{-2n-3} \right).  
 \end{equation} 
We make the meaning of the asymptotic expansion precise by estimating the difference Hamiltonians in Proposition \ref{eq:kdv:Lipschitz}. It is classical that for Schwartz functions 
\[  |\T^\KdV_N(i\tau,u)| \le c \tau^{-2} \qquad \text{ for } \tau \ \text{ large }   \]
and 
\[ \lim_{\tau\to \infty}   -(2\tau)^2 \T^\KdV_N(i\tau,u) = H_{N+1}^\KdV. \]
We define the  generating function  
\begin{equation} \label{def:T-1} \T^{\KdV}_{-1}(z,u) =   iz  \ln T^\KdV(z) -  \frac12 \int u dx  \sim \sum_{n=0}^\infty (2z)^{-2n-2} 
H_n^\KdV 
\end{equation} 
\extra{
\bigskip 

{\color{blue}
Here is what determines the sign. By the inequality in \eqref{eq:TKdVpos}
\[ \real (-z^2 \ln T^\KdV_r) \ge 0.  \]
We construct $ \T^\Gardner_{-1} $ so that 
\[  \lim_{\tau\to \infty} (2i\tau)^2\T^\KdV_{-1}(i\tau, u) = \frac12 \Vert u \Vert_{L^2}^2 \]
hence, if $u \ne 0 $ and $ \tau $ is large then $\T^\KdV_{-1}(i\tau, u) < 0 $.
Thus 
\[ \T^\KdV_{-1}(z,u) = iz \Big(\ln T^\KdV(z,u) - \frac1{2iz} \int u dx\Big)  \]
so that 
\[   \tau \T^\Gardner_{-1}(i\tau,u) =-iz \T^\Gardner(z)\Big|_{z=i\tau} = z^2 \Big(\ln  T^\KdV(z) - \frac1{2iz} \int u dx \Big) \Big|_{z=i\tau} < 0  \]
} 

\bigskip 
}
   for the KdV energies which plays a central role in this paper
  for several reasons:
  \begin{itemize} 
   \item It gives uniform control over the $H^s$ norm of solutions, see  \cite{MR3820439} and \cite{MR3874652}. 
   \item It is (locally) a smooth function of the potential $u$ for fixed $z$ and it defines Hamiltonian flows with respect to the Gardner Poisson structure. These flows all commute.
   \item The KdV Hamiltonian $H^\KdV_{N+1}$ is the  limit of $-(2\tau)^2\mathcal{T}^{\KdV}_N(i\tau,u)$ as $ \tau \to \infty$. We will use a linear combination of such Hamiltonians  so that   the corresponding flows converge to the flow of $H^\KdV_{N+1}$. However we will do that at the level of the Gardner hierarchy.  
 \item We will use a diffeomorphism on Sobolev spaces to translate the problem into a problem for the Gardner equation in the variable $w$ depending on $\tau$. Then $\frac12 \Vert w \Vert_{L^2}^2 =  -\T^{\KdV}_{-1}(i\tau, u)$, see \eqref{eq:rhologT}, we define a generating function for the Gardner hierarchy, so that the $N$th Gardner equation for $w$ is essentially equivalent to the $N$th KdV equation for $u$.
  \end{itemize} 
  By definition the transmission coefficient is a meromorphic function on the upper half plane with simple poles at the square roots of the eigenvalues of the Lax operator in the upper half plane. Let $-\kappa_j^2$ be the eigenvalues of the Laplace operator with $\kappa_{j+1} < \kappa_j$.
  Then, ignoring the question about convergence (for example assuming that $ u$ is a Schwarz function so that there is a finite product),   the product 
  \[   \prod_{j} \frac{ z- i\kappa_j}{z+i\kappa_j}   T^{\KdV}(z) \]
  is holomorphic in the upper half plane. If the potential $u$ is a Schwartz function then $T^\KdV(\xi) $ is defined for $ \xi \in \R $ as the limit from the upper half plane. We recall that $ |T^\KdV(\xi)| \le 1  $ (see \cite{MR3874652} for example). $T^\KdV(i\tau)$ is real if $ \tau$ is sufficiently large. It is holomorphic in the upper half plane outside $i[0, \kappa_0] $. We choose $\ln \T^\KdV(z)$ so that the logarithm is real at $i\tau$ for $\tau $ large. 

     We choose the standard  complex logarithm for the Blaschke factors $\frac{z-i\kappa_j}{z+i\kappa_j} $.
  
  Then  
  \[  \real ( - z^2 \ln \Big(  \prod_{j} \frac{ z- i\kappa_j}{z+i\kappa_j}   T^{\KdV}(z))\Big) \ge 0    \]
on the real axis, it is holomorphic on the upper half plane, and nonnegative on the upper half plane by the maximum principle.  Moreover 
\[   \real  \ln \frac{z+i \kappa_j}{z-i\kappa_j} = \ln \Big| \frac{z+i\kappa_j}{z-i\kappa_j} \Big|   \]
 is harmonic in the upper half plane outside $i[0,\kappa_j]$. It  is continuous away from  $ i \kappa_j$, vanishes on the real axis, is positive on $i(0,\kappa_j) $ and $ \sim -\log |z-i\kappa_j| $ near $i\kappa_j$.
Hence  it is positive on the upper half plane outside  $\kappa_j$ and 
  \begin{equation}\label{eq:TKdVpos}  
  \real \big[ iz  \T^\KdV_{-1}(z,u)\big] =\real \Big(- z^2\ln T^\KdV(z) - \frac{iz}2 \int u dx  \Big) \ge 0.
  \end{equation}
 on the upper half plane away from the $ \kappa_j$.

  On the real axis the absolute values of the factors 
  $\frac{ z- i\kappa_j}{z+i\kappa_j}$
  are identically $1$. We choose the standard branch of the complex logarithm,  and compute
  \begin{equation}
      \begin{split}  -\Delta \real \Big( z^2 \ln \frac{ z +i \kappa}{z-i\kappa} \Big) 
 \, &      = \real z^2 (-\Delta \ln \frac{z+i\kappa}{z-i\kappa}) - \real 2iy \partial_x \log \frac{ix-y-\kappa}{ix-y+\kappa}
 \\ & = 2\pi \kappa^2 \delta_{i\kappa} - 4\pi i z \delta_{i (0,\kappa)}=:\mu_{\kappa}
 \end{split} 
   \end{equation} 
Here $\delta_{i(0,\kappa)}$ is the one dimensional Hausdorff measure restricted to this line segment, $ \delta_{i\kappa} $ is the Dirac measure at $ i \kappa$
and  $\mu_{\kappa}$ is a signed Radon measure in the upper half plane supported in $i(0,\kappa]$. 
The contribution at
$i \kappa$ comes from
the singularity of 
the logarithm and the contribution on the line section from the jump of the $x$ derivative across.
    As a consequence - independent of the convergence question above, and choosing the obvious branch of the logarithm - $ \real z^2 \ln ( \frac{z +i\kappa}{z-i\kappa} ) $ is continuous in the complement of $ i \kappa$ in the upper half plane. In particular, $ \Delta z^2 \ln ( \frac{z +i\kappa}{z-i\kappa} ) $ is supported in $i [0,\kappa]$. We obtain on the upper half plane 
     \begin{equation}\label{eq:super} - \Delta \im z \T^{\KdV}_{-1} 
  = - \Delta \real z^2 \ln T^\KdV =\sum_j  -\Delta \real z^2 \ln \frac{z+i\kappa_j}{z-i\kappa_j} = \sum_j     \mu_{\kappa_j}  
  \end{equation}
We may also reverse the steps and obtain
$\im z \T^{\KdV}_{-1}(z,u)=  \real \big[-iz \T^\KdV_{-1}(z,u)\big]  $ 
from \eqref{eq:super} and the values of $ -2\xi \T^\KdV_{-1}(\xi, u)  $ on the real axis. 
The map $z \to \im 2z\T^\KdV_{-1}(z,u)$ is harmonic in the upper half plane without $i [0, \kappa_0]$ where $ -\kappa_0^2 $ is the ground state energy.  We compute
  using the fundamental solution $\frac1{2\pi} \ln |z| $
  by integration by parts for the second identity and 
  \[ \frac{s^2}{s-iz}= s+iz - \frac{z^2}{s-iz},    \]
\begin{equation} \begin{split}  \kappa^2 \ln\frac{|z+i\kappa|}{|z-i\kappa|} 
- \int_0^\kappa  2s \ln \frac{|z+is|}{|z-is|} ds \,  & =  
 \real\Big(  \int_0^\kappa s^2  \partial_s\Big(\ln (s-iz) - \ln(s+iz) \Big) ds\Big)
\\=  & \real  \int_0^\kappa s^2  \Big( \frac1{s-iz} - \frac1{s+iz} \Big) ds 
\\ = & \real \Big(  2i\kappa z + z^2 \int_0^\kappa  -\frac1{s-iz} +\frac1{s+iz} ds  \Big) 
\\ =&  \real \Big( 2i\kappa z + z^2 \big( \ln(1+\frac{\kappa}{iz}) - \ln( 1- \frac{\kappa}{iz} )\big) \Big)  
\end{split}
\end{equation}
and hence, using the Poisson formula and summing over the eigenvalues, 
\begin{equation} \label{eq:realpart}  \begin{split}  0\le \real\big( 2iz \T^{\KdV}_{-1}(z,u)\big) =\, & \real \Big[ \frac{1}{\pi} \int_{\R} \frac{\im z }{| z-\xi|^2} \real (2i\xi  \T^{\KdV}_{-1}(\xi,u)) d\xi \\ &  -2z^2 \sum_j  \ln \Big(\frac{z+i\kappa_j}{z-i\kappa_j}\Big)+    4iz \kappa_j \Big]. \end{split}   \end{equation}  
The right hand side of the following identity 
\begin{equation}\label{eq:trace.h-1}
    \T^{\KdV}_{-1}(z,u) = \frac{1} {2\pi z} \int_{\R} \frac{-2\xi}{z-\xi} \im \T^{\KdV}_{-1}(\xi,u) d\xi   +iz \sum_j  \ln (\frac{z+i\kappa_j}{z-i\kappa_j})+ \sum_j 2 \kappa_j 
\end{equation}
is holomorphic on the upper half plane away from $i[0,\kappa_0]$.  If we multiply both sides by $2\pi z$ and take the real part we obtain the identity \eqref{eq:realpart}. The difference $U$ of the two sides is holomorphic, with 
$\real (-2iz U) = 0 $, hence $U$ is an imaginary constant. But $ \lim_{\tau \to \infty} U( i\tau)=0$ and we obtain the identity \eqref{eq:trace.h-1}.

 \begin{lemma}\label{lem:kdvestT1} Let $ d(z^2, S(-\partial^2+u)) $ be the distance to the spectrum. Then for $ \im z > 0 $
 \begin{equation}\label{eq:estT1} |\T^{\KdV}_{-1}(z,u) | \lesssim   \frac{1}{\im z |z| }  \max\Big\{ 1,-\ln\big(2 d(z^2, S(-\partial^2 +u))|z|^{-2}\big)\Big\} \Vert u \Vert^{2}_{L^2}.\end{equation}
 Let $ \tau_0>0$ be such that $ -\partial^2+ u +\tau_0^2$ is positive semidefinit. Then 
 \begin{equation}\label{eq:monoton} (\tau_0,\infty) \ni \tau \to \T^\KdV_{-1}(i\tau, u) \end{equation}
 is monotonically increasing. 
\end{lemma}
 The lemma provides bounds on the transmission coefficients away from the spectrum, and the estimate blows up mildly as we approach the spectrum.% 

 \begin{proof}
By construction
$T^\KdV(i\tau) \in \R$
for $\tau $ large (so that we avoid the eigenvalues) and, using   the trace identity \eqref{eq:trace.h-1} we obtain the $L^2$ trace identity
\begin{equation}\label{eq:trace.l2}  
\frac12 \Vert u \Vert^2_{L^2} = \lim_{\tau \to \infty} (2i\tau)^2 \T^{\KdV}_{-1}(i\tau) 
 = \frac2{\pi}   \int_{\R} -2\xi \im \T^{\KdV}_{-1}(\xi) d\xi  + 
\sum_j \frac83 \kappa_j^3.
\end{equation}
%The trace identity \eqref{eq:trace.h-1} now gives, first for nice $u$,
%\begin{equation}\label{eq:trace}  \T^\KdV_{-1}(z,u) = 
%  \frac{1} {2\pi z  } \int_{\R} \frac{1} {\xi-z} (-2\xi) \im \T^{\KdV}_{-1}(\xi,u) d\xi   +z \sum_j  \ln (\frac{z+i\kappa_j}{z-i\kappa_j})-  \sum_j  2\kappa_j.   \end{equation}
We estimate  $  \frac{1}{2\pi |z| |\xi-z|} \le    \frac{1}{2\pi |z|\im z } $
and, for an absolute constant $C$
\[ 
\Big|z \ln ( \frac{z+i\kappa}{z-i\kappa} )-2i\kappa\Big|   = \Big|z \ln \frac{ 1+i \kappa/z}{1-i\kappa/z}-2i\kappa\Big| 
\le C \frac{\kappa^3}{|z|^2}   (1- \min\{ \log |1-i\kappa/ z|,0 \} ). 
\]
We compare with \eqref{eq:trace.l2} and obtain \eqref{eq:estT1} for nice (Schwartz) functions $u$.
In compact sets of the upper halfplane without $i[0,\kappa_0]$ the map $L^2\ni u \to \T^\KdV_{-1}(z,u)$
is continuous and we have to understand $-2\xi \T^\KdV_{-1}(\xi,u) dx $ as nonnegative measure which in general is not absolutely continuous with respect to the Lebesgue measure. 
This implies \eqref{eq:estT1}. 

To check the monotonicity we observe that (using $ \im \T^\KdV_{-1}(\xi, u) = \im \T^\KdV_{-1}(-\xi,u)$)
\[ (\tau_0,\infty)\ni \tau \to    \frac1{2\pi i\tau } (\frac1{i\tau -\xi}+ \frac{1}{i\tau +\xi})
= -\frac1{\pi} \frac1{\tau^2+\xi^2}
\]
is monotonically increasing. Similarly we check that 
\[ (1,\infty) \ni \tau \to \rho(\tau) = -\tau \ln \frac{\tau+1}{\tau-1}+2
= -2 \sum_{n=0}^\infty \frac1{2n+1} \tau^{-2n-1} 
\]
is monotonically increasing. Hence the integrands in \eqref{eq:trace.h-1} are monotinically increasing which implies the monotonicity of $ (\tau_0, \infty) \to \T^\KdV_{-1}(i\tau,w)$.
\end{proof}

The difference
Hamiltonian $\mathcal{T}^\KdV_{N}$ plays
a central role,  not only for  giving  \eqref{eq:asymptotic} a precise meaning. It is a function of $u$ and $z$. We introduce some notation for its domain of definition.
\begin{definition}\label{def:sigmau}
Let $ \tau >0$. We define
\[ \sigma(u) = \left\{ \begin{array}{l}  0  \text{ if } -\partial^2+ u \text{
is positive definite.}  \\ 
\text{ the square root of the negative of the lowest eigenvalue otherwise, }\end{array} \right. 
\]
the set of functions 
\[ U(\tau_0) = \{ u \in H^{-1}: \sigma(u) < \tau_0\}  \]
and $\C_{\tau_0} = \{ z \in \C: \im z >0 \} \backslash i[0 , \tau_0] $.
\end{definition} 

The space $H^{-1}$ will play a prominent role. For scale invariant estimates we introduce the norm 
\begin{equation}\label{eq:XX} \Vert u \Vert_{H^{-1}_\tau} =   \Vert  (\xi^2 + \tau^2)^{-1/2}  \hat u \Vert_{L^2}. \end{equation}

\begin{proposition} \label{eq:kdv:Lipschitz} Let $N \ge -1 $.
The difference Hamiltonians are defined on $ \C_{\tau_0}  \times    U(\tau_0)$. They are holomorphic in the first component and analytic in the second component and satisfy 
\begin{equation}  \label{eq:tauNdec}
          | \T^{\KdV}_N(z,u) |\le  C_N  
          \left( \frac{|z|}{\im z} \right)^{2N+3}  
          (\im z)^{-2} \Big( \Vert u \Vert_{L^{N+2}}^{N+2} + 
          \Vert u \Vert^2_{H^{N+1}}\Big),
\end{equation}           
 if $ (\im z)^{-1/2} \Vert u \Vert_{H^{-1}_{\im z}}  \le \delta < 1$ for an absolute constant $\delta$.         
\end{proposition}
Observe that for all $ u \in H^{-1} $
\[ \lim_{\tau \to  \infty}  \tau^{-1/2}\Vert u \Vert_{H^{-1}_{2\tau}} =0.\]
 Proposition  \ref{lem:TNest} in Section \ref{subsec:diffHamiltonian} provides an analogous estimate  for $\T^\Gardner_N(iz,w,\tau)$. It is likely that one may combine the estimate for the Gardner case and the relation between KdV and Gardner 
 (see \eqref{eq:Gardner-KdV} and Subsection \ref{subsec:diffeo}) to deduce this proposition from the Gardner case. 
  Theorem 9.3 in \cite{MR3874652} provides a weaker estimate than Proposition \ref{eq:kdv:Lipschitz}. However the bounds and convergence results for the multilinear integrals there are similar to the ones in this paper and are strong enough to imply Proposition \ref{eq:kdv:Lipschitz}. We do not rely on this proposition and state it only for completeness. 
  
The structure of the coefficients $H^{\KdV}_N$ is described in Theorem \ref{thm:formofkdv}.

\begin{remark}
Proposition \ref{eq:kdv:Lipschitz} has  the  immediate consequence that  for $ u \in H^{N+1}$
\begin{equation} \label{eq:approximate}   \lim_{\tau\to \infty} (2i\tau)^2 \mathcal{T}^{\KdV}_{N-1}(i\tau,u) = H_N^{\KdV}(u)  \end{equation}  
by applying the estimate to $\mathcal{T}^{KdV}_N(i\tau)$.
It then implies 
\[ \limsup_{\tau \to \infty} |\T^\KdV_N (i \tau ,u) |+ |H^\KdV_N(u)|  \le c  \Big(\Vert u \Vert_{L^{N+1}}^{N+1} + \Vert u \Vert_{H^N}^2\Big) .   \]
The convergence above \eqref{eq:approximate} holds for $ u \in H^N$ by an approximation argument.  Since we use this techniques only for the Gardner case 
we omit detailed proofs here. 
\end{remark}

The $N$th equation of the KdV hierarchy is the Hamiltonian equation of the Hamiltonian $H_N$ with respect to the  Gardner Poisson bracket
\begin{equation} \label{eq:defGardner}   \{ F, G \}^{\Gardner} = \int \frac{\delta}{\delta u} F(u) \partial  \frac{\delta}{\delta u} G(u) dx. \end{equation}
where $ \frac{\delta}{\delta u}$  denotes the variational derivative defined by 
\[ \int \frac{\delta}{\delta u} F v dx = \frac{d}{dt} F(u+tv) |_{t=0} \] 
for $v \in C^\infty_c$, assuming that the right hand side is defined. 
More explicitly, the $N$th KdV equation is 
\begin{equation} \label{eq:HN}u_t = \partial_x \frac{\delta}{\delta u} H^\KdV_N. \end{equation}

\subsection{Poisson commutation and the theorem of Frobenius}
\label{subsec:Frobenius} 
We begin by considering the finite dimensional case. Let $M^n$ be a smooth $n$ dimensional manifold. A Poisson bracket is a bilinear map from $C^{\infty}(M) \times C^{\infty}(M) \to C^{\infty}(M) $ which we denote by $\{f,g\}$ with the properties: 
\[   \{ f, g \} = -\{g,f \}  \qquad \text{ antisymmetry, } \] 
\[ \{f,\{g,h\}\}+ \{ g, \{h,f\}\} + \{ h, \{f,g\} \} = 0 \qquad \text{ Jacobi identity,}   \]
\[ \{ f,gh\} = \{f, g\} h + \{f,h\} g \qquad \text{ Leibnitz rule.} \]
The Leibniz rule says  that $f \in C^\infty(M) $ defines a derivation, and hence a vector field $X_f$, the Hamiltonian vector field of $f$, so that 
%$ \{g,f\} = - \{f,g\}= X_f g$.
$ \{f,g\} =  X_f g$.
Then 
%\[\begin{split}  X_{\{f,g\} } h \, &  =\{ h, \{ f,g \} \}
%\\ & = \{\{ h,g \},f \} - \{ \{h,f\} ,g \} 
%= (X_f X_g - X_g X_f ) h 
%\\ & = [ X_f, X_g] h 
%\end{split}
%\] 
\[\begin{split}  X_{\{f,g\} } h \, &  =\{ \{ f,g \} , h\}
\\ & = \{f, \{g ,h \} \} - \{g, \{f,h\}  \} 
= (X_f X_g - X_g X_f ) h 
\\ & = [ X_f, X_g] h 
\end{split}
\] 
and 
\[ X_{\{f , g\}} = [X_f,X_g]. \]
In particular the functions $f,g$ Poisson commute iff the Hamiltonian vector fields commute. 
At $ x \in M$, $\{ f,g \}$ is linear in the derivative $Dg$ of $g$ at the point  $x$. By antisymmetry, it is also linear in the derivatives of $f$ at $x$, and there is a bounded antisymmetric bilinear form 
$\Omega$ on the cotangent space of $M$ at $x$ such that $\{f,g\}  = \Omega( Df(x), Dg(x))$. 
This defines a linear map $J: T^* _x M \to T_x M$ so that $ \Omega(Df,Dg) = Df(J Dg)$. Then 
\[  J Df = X_f. \]
Smooth vector fields on a manifold define a Lie algebra. Similarly a Poisson product defines a Lie algebra of smooth functions on $M$.
A symplectic manifold is an even dimensional manifold with a nondegenerate closed two form $\omega$. Then $ \Omega $ is the induced antisymmetric bilinear form on $T^*M$, $J$ is invertible  and 
\[ \omega_x(v_x,w_x)= J^{-1} v_x (w_x). \]
 Let $X$ and $Y$ be two $C^2$ vector fields and $\phi^X(t,\cdot)$ resp $\phi^Y(t,\cdot)$ the corresponding local flows.  One version of the Frobenius theorem is (whenever the local flow is defined) 
  \[     [X,Y] = 0 \quad  \Longleftrightarrow   \phi^X(s,\phi^Y(t,\cdot)) = \phi^Y(t,\phi^X(s,\cdot)). \] 
  The implication $\Longleftarrow$ requires differentiation of the identity on the right with respect to $\partial_s \partial_t $ follows by an evaluation at $s=t=0$. The opposite direction is a bit more involved. 
  Assume $[X,Y]=0$. We differentiate 
\begin{equation}\label{eq:flowcommute}  \begin{split} \hspace{1cm} & \hspace{-1cm}  \frac{d}{ds} (\phi^X(s,\phi^Y(t,\cdot)) -\phi^Y(t,\phi^X(s,\cdot)) ) \\ & = X(\phi^X(s,\phi^Y(t,\cdot)) - (D\phi^Y)(t,\phi^X(s,\cdot)) X(\phi^X(s,\cdot)) \end{split} \end{equation}  
We claim 
\begin{equation} \label{eq:claimcom} (D\phi^Y)(t,\cdot) X(\cdot) = X(\phi^Y(t,\cdot)). \end{equation}  
Assume the claim. Then the right hand side of \eqref{eq:flowcommute} is 
\[ X(\phi^X(s,\phi^Y(t,\cdot)) - X(\phi^Y(t,\phi^X(s,\cdot)) \] 
and by Grönwall's inequality, the identity for $s=0$, and local Lipschitz continuity of $X$, we find $\phi^X(s,\phi^Y(t,\cdot)) = \phi^Y(t,\phi^X(s,\cdot))$.
To prove the claim we differentiate \eqref{eq:claimcom} with respect to $t$, 
\[
\begin{split} 
\frac{d}{dt} ((D\phi^Y)(t,\cdot) X(\cdot) - X(\phi^Y(t,\cdot)))\, &    \\
&\hspace{-4.5cm} = D(Y(\phi^Y(t,\cdot))) X(\cdot)  - DX(\phi^Y(t,\cdot))Y (\phi^Y(t,\cdot)) \\ 
&\hspace{-4.5cm} = DY(\phi^Y(t,\cdot))D\phi^Y(t,\cdot) X(\cdot) - DX(\phi^Y(t,\cdot))Y (\phi^Y(t,\cdot)) \\ 
&\hspace{-4.5cm}  = ((DY)X - (DX)Y)(\phi^Y(t,\cdot)) +  DY(\phi^Y(t,\cdot))(D\phi^Y(t,\cdot)X(\cdot)-X(\phi^Y(t,\cdot))) \\
&\hspace{-4.5cm}  = [X,Y](\phi^Y(t,\cdot)) +  DY(\phi^Y(t,\cdot))(D\phi^Y(t,\cdot)X(\cdot)-X(\phi^Y(t,\cdot)))
\end{split} 
\] 
and the claim follows again by $[X,Y] = 0$, Grönwall's inequality, the identity for $t=0$ and local Lipschitz continuity of $DY$. 

We specialize to Hamiltonian vector fields. The Hamiltonians $f$, $g$  Poisson commute if and only if 
$[X_f, X_g] = 0 $ which holds if and only if 
\[ \exp(sJ \nabla f) \exp(tJ g ) = \exp(tJ \nabla f) \exp(sJ \nabla g) \] 
which holds if and only if $g$ is conserved on orbits of $\exp(tJ \nabla f)$.
To see this we calculate 
\[ \frac{d}{dt} g( \exp(t J \nabla  f) = \{ f, g\} \circ \exp(tJ \nabla f)  \]
and the left hand side vanishes iff the right hand side vanishes.

 We are interested in the infinite dimensional setting. The Gardner Poisson structure on suitable smooth functionals (the most important functionals  being integrals over polynomials in $u$ and its derivatives and $\T^{\KdV}_{-1}$) is 
 \[ \{F, G\}^\Gardner= \int_{\R}  \frac{\delta}{\delta u}F \partial \frac{\delta}{\delta u}  G dx. \]
It is straight forward to verify that this defines a Poisson bracket. 
The Hamiltonian vector field is 
 \[ X_F = \partial \frac{\delta}{\delta u} F.  \]

The core of integrability  is the fact the Hamiltonians  $\T^{\KdV}_{-1}$ Poisson commute (see Lemma \ref{kdvpoisson}). 

\begin{proposition} \label{prop:commutation}
The Hamiltonian $\T^{\KdV}_{-1}(z, u) $  Hamiltonian Poisson commutes with respect to the Gardner bracket, 
\[ \{ \T^{\KdV}_{-1}(z_1,\cdot ) , \T^{\KdV}_{-1}(z_2,\cdot ) \}^\Gardner = 0. \]
\end{proposition}

 The Hamiltonian vector field of $\T^\KdV_{-1}(i\tau,.) $
is $ \partial \frac{\delta}{\delta u} \T^\KdV_{-1}(i\tau, .) $. We will see that it defines a smooth flow on the subset of $H^{-1}$ of potentials for which  $-\partial^2+ u + \tau^2$ is positive definite. Moreover, by an adaptation of the finite dimensional arguments,  
\[\begin{split}  
0\, & =\Big\{\T^\KdV_{-1}(i\tau_1),  \T^\KdV_{-1}(i\tau_2, .) \Big\} (u) 
\\ & =  \frac12\frac{d}{dt}   \T^{\KdV}_{-1}\big(i\tau_1, u + t \partial \frac{\delta}{\delta u} \T^\KdV_{-1}(i\tau_2,u)\big)\big|_{t=0}
 \\ & \qquad - 
 \frac12\frac{d}{dt}   \T^{\KdV}_{-1}\big(i\tau_2, u + t \partial \frac{\delta}{\delta u} \T^\KdV_{-1}(i\tau_1,u)\big)\big|_{t=0} .
\end{split}
\]

The KdV Hamiltonians define unbounded vector fields. The  Poisson commutator 
$\{H_n^\KdV, H_N^\KdV\} $  is  defined for $u \in H^{2n+2N+2}$ and it can be extended to $u \in H^{n+N+1} $ and even slightly beyond that, but certainly not to open subsets of $H^{-1}$.

The KdV Hamiltonians are defined as limits of Poisson commuting functions, hence they Poisson commute with one another and with $\T^\KdV_{-1}(i\tau)$, at least for sufficiently regular functions.  As a consequence 
$T^\KdV_{-1}(i\tau, u) $ and $H^\KdV_n$ are conserved along any KdV flow. Of course this extension of results in finite dimensions has to be proven. Nevertheless the finite dimensional view  gives a good orientation.

\subsection{Miura map and Gardner hierarchy}

We find   it is easier to study well-posedness   questions  for the Gardner hierarchy, which we define and study here and in Section \ref{sec:AKNS}. The Hamiltonians $H_N^{\Gardner}(w,\tau_0)$ with 
 \[ H_0^{\Gardner} = \frac12\int w^2 dx , \quad H_{1}^{\Gardner}(w, \tau_0) = \frac12\int w_x^2+ w^4+ 4\tau_0 w^3 dx \] 
 and the Gardner equations  (the Hamiltonian 
  $H_N^{\Gardner}$ is defined in \eqref{eq:Gardnerexpansion})
 \begin{equation} \label{eq:GardnerN}  w_t = \partial \frac{\delta}{\delta w} H_N^{\Gardner}(w,\tau_0)  \end{equation} 
 depend on a spectral parameter $\tau_0$. They are connected to KdV by the remarkable modified Miura map (with $z = i \tau_0$)
 \begin{equation} \label{eq:Miura}  M(-iz, w):=   w_x + w^2 - 2iz w \end{equation} 
 which has been used by Miura, Gardner and Kruskal  to formally derive the Hamiltonians of the KdV hierarchy \cite{MR252826}. A short calculation shows that the inverse is given by% 
 (recall that $ \phi_l$ is the  Jost function)
 \[  u \to \partial_x \ln \phi_l + iz. \] 
The modified Miura map defines an analytic diffeomorphism  for $N \ge 0$ and $z = i \tau_0$ (see Definition \ref{def:sigmau} for $U(\tau_0)$) 
  \begin{equation}
 M(\tau_0,.): H^{N} \to   U(\tau_0)  \cap H^{N-1},
 \end{equation} 
 which is contained in Lemma \ref{lem:uw} and Proposition  \ref{prop:equivalencewu} in Section \ref{sec:mmiura}.  
 The function  $ w\in C(I; H^{N})$  is a weak solution to the $N$th Gardner equation  \eqref{eq:GardnerN}  
 if and only if  $u= M(\tau_0,w)\in C(I, H^{N-1})$ is a weak solution to the $N$th KdV equation (see Theorem \ref{thm:equivalenceweaksolutionsuw}).
 
A short calculation  \eqref{eq:rhologT}   shows that for $M(-iz,w(z)) = u$,
\begin{equation}
     \T^\KdV_{-1}(z, u) =- \frac12 \int w^2(z) dx \label{eq:T-1}
\end{equation} 
and as a consequence, using the chain rule, we compute
\[w= - (-\partial-2iz+2w)\frac{\delta }{\delta u} \T^\KdV_{-1}(z,u) \]
and we can write the equation for%
 the flow defined by $u_t= \partial_x \frac{\delta}{\delta u} \T^{\KdV}_{-1}(i\tau, u)=-\partial_x (-\partial+2 \tau +2w)^{-1} w$ as a system of differential equations
\begin{equation} \label{eq:tauflow}  
\begin{split} u_t \, & =  \partial_x F(u)  \\  w_x +2\tau w + w^2\, & = u 
\\-\partial_x F + 2\tau F + 2w F &= -w. 
\end{split}
\end{equation}
Observe that $F(u)$ gains two derivatives over $u$! The map $ w \to u $ in the second equation and its inverse is studied in detail in Section \ref{sec:mmiura}.

We call the flow defined by the Hamiltonian% 
$ \T^\KdV_{-1}(i\tau)$
$\tau $-flow.  The formulas above easily imply the following
(see also Section \ref{sec:mmiura}). 

\begin{proposition} \label{prop:kdvLipschitz}
The Hamiltonian $\T^{\KdV}_{-1}(i\tau, u) $ is real for  
$\tau > \sigma(u)$ (see Definition \ref{def:sigmau}).
The map to the variational derivative  
\[   H^{N-1} \ni u \to  \frac{\delta}{\delta u } \T_{-1}^{\KdV} (z, u) \in H^{N+1} \]
for $N \in \N$ 
is locally Lipschitz and smooth, and also the Hamiltonian vector field 
\[ H^{N-1}\ni u \to \partial_x \frac{\delta }{\delta u} \T_{-1}^{\KdV}(z,u)  \in H^N \] 
is locally Lipschitz and smooth. 
\end{proposition}

\begin{proof} 
Let $ \tau > \sigma(u)$. Then the second equation of \eqref{eq:tauflow} can be solved for $w$ in terms of $u$ (see also \cite{MR2189502} ) with  bounds for  $w$ in terms, see Section \ref{sec:mmiura}, Lemma \ref{lem:uw}. Similarly the third equation can be solved (see Subsection \ref{subsec:goodvariable} with a somewhat different notation). It is not hard to see that the map $u \to F$ gains two derivatives and it is smooth - it is arbitrarily often Frech\'et differentiable.\end{proof}

By the Cauchy-Lipschitz theorem the Hamiltonians $ \T^{\KdV}_{-1}(i\tau,.)$ define a local flow on 
$H^N$ for $N \ge -1$ by 
\[ u_t = \partial \frac{\delta}{\delta u} \T^{\KdV}_{-1} (u). \]
 The Hamiltonians $\T_{-1}^{\KdV}(i\tau)$ for different $\tau$ are conserved under the $\tau$ flow as a consequence of Proposition \ref{prop:commutation} and 
  the flows  commute with themselves.
The KdV Hamiltonians are conserved under the $\tau$ flows. The KdV Hamiltonians control the Sobolev norms $H^N$. This is a consequence of  Theorem 9.3 in \cite{MR3874652}. We prove the analogous results for the Gardner hierarchy and its generating function. The two hierarchies are closely related and one could deduce the claims here from the claims for Gardner. Hence  $\tau$ flows are global in time and preserve higher regularity. 

\medskip

 Killip and Visan \cite{MR3990604} (see also \cite{MR1995460}, Chapter 11)   introduced the diagonal Green's function (the diagonal of the Green's function of the Lax operator) into the well-posedness question. 
 It is related to $w$  through the factorization 
 \[ -\partial^2 + u+\tau_0^2 = (\partial + w + \tau_0) ( -\partial+ w + \tau_0). \]
 The linear operators on the right hand side can be inverted and one obtain the formula for the  integral kernel of the resolvent \[
 G(x,y)= \int_{\max\{y,x\}}^{\infty}   \exp\Big( \tau_0(x+y-2t) - 2\int_{\max\{x,y\}}^t  w ds- \int_{\min\{x,y\}}^{\max\{x,y\}} w ds ) \Big)       dt. 
 \]
 Slightly deviating from the notation of Killip and Visan we define the good variable $v$ via the diagonal Green's function (see Lemma \ref{lem:relationuvw}, formula \eqref{eq:fdT} expressing the variational derivative in terms of Jost solutions, and \eqref{eq:Gxy} which expresses the Green's function in terms of the Jost solutions.  Compare also to \cite{MR1995460}, Chapter 11) 
 \[ \beta:
 = G(i\tau_0 , x,x)= \frac{\delta \ln T^{\KdV}(i\tau_0)}{\delta u}, \] 
 \[ v := \frac{1}{2\tau_0 \beta} -1.  \]

Now $u, v, w$ are related by the following relations (see Lemma \ref{lem:relationuvw}) 
\begin{align}
\label{eq:wv}w&= W(\tau_0, v):= \tau_0 v - \frac12 \partial \ln( 1+ v) \\
\label{eq:uw}u&=M(\tau_0,w) =  w_x+w^2+2\tau_0 w\\
\label{eq:uv}u&= -\frac{1}{2}\frac{v_{xx}}{v+1}+\frac34\frac{v_x^2}{(v+1)^2}+\tau_0^2v^2+2\tau_0^2 v
    \end{align}
For $ s > -\frac12$ 
\begin{equation} \label{eq:V}  W : \mathcal{V}^s:= \{ v \in H^{s+1}: v > -1\}\ni v \to \tau_0 v- \frac12 \partial_x \ln (1+v)  \in H^{s}  \end{equation}  
is an analytic diffeomorphism (Theorem~\ref{thm:equivalence}). With these definitions  $u \in C(I,H^{N-1}))$ 
 is a  solution to $N$th KdV equation \eqref{eq:HN} if and only if $v \in C(I, \mathcal{V}^{N+1})$ is a weak solution to 
\begin{equation}\label{eq:v}   v_t =  2\partial \Big( (v+1)\sum_{n=-1}^{N-1} \frac{\delta H_n^{\KdV}}{\delta u}(2i\tau_0)^{2(N-1-n)}(u) \Big).   \end{equation}  
if and only if $ w= W(v) $ satisfies 
the $N$th Gardner equation \eqref{eq:GardnerN}. This is the content of 
Theorem \ref{thm:HNflows} for smooth solutions and 
Theorem \ref{thm:equivalenceweaksolutionsuw}  and Theorem \ref{thm:vwweak} for weak solutions under weaker resp. different regularity assumptions. 
We write down the equation for $v$ for $N=1$
and recall that $H_{-1}^\KdV=\frac12\int u dx$ and $H_0^\KdV = \frac12 \int u^2 dx $
\begin{equation} \label{eq:goodvariablekdv} 
  v_t =2 \partial \Big[ (2i\tau_0^2) \frac{\delta H^\KdV_{-1}}{\delta u}
  + \frac{\delta H_0^\KdV}{\delta u} \Big]   
   = 
   2\partial\big[ (v+1)(u- 2\tau_0^2) \big]. 
\end{equation}
It is remarkable that all the equations \eqref{eq:v} are differential equations. 

The relation between the $N$th KdV equation \eqref{eq:HN}, the $N$th Gardner equation  \eqref{eq:GardnerN}  and the $N$th 'good variable' equation \eqref{eq:v} via the diffeomorphisms 
\eqref{eq:uw} (Miura map $M$), \eqref{eq:wv} (the map $W$) and the composition \eqref{eq:uv} extends to more general weak solutions  (Theorem \ref{thm:equivalenceweaksolutionsuw} and Theorem \ref{thm:vwweak}). This reduces the proof of Theorem \ref{th:main} to a similar statement for the Gardner hierarchy, Theorem \ref{thm:Gardneruniqueness} below.

\subsection{Well-posedness for the  KdV equation in \texorpdfstring{$H^{-1}(\X)$}{H-1}.} 
\label{subsec:KdVwell}

We find it instructive to follow the proof of Killip and Visan \cite{MR3990604} in our setup to prepare for the case of the higher KdV equations. 
We  follow the strategy of Killip and Visan   and prove well-posedness and the commutation property simultaneously, by approximating by flows defined by the Hamiltonians $ -(2\tau)^2\mathcal{T}^\KdV_{0}(i\tau)$, motivated by \eqref{eq:approximate}. Alternatively one can deduce  that the flows defined on smooth functions commute, and approximate the initial data to verify commutation of the flows. 

Using the previous section we study well-posedness for the Gardner equation for $ N \ge 2$ resp. well-posedness for the good variable equation for $N=1$, similar to Killip and Visan.

There is a difference in the case $N=1$, for which we cannot define weak solutions to the Gardner equation assuming only $ w\in L^\infty(\R;L^2)$ since the nonlinearity contains the term $w^3$. For $N \ge 2$ this is no issue since $L^\infty H^{N-1} \subset L^\infty L^\infty$. Also when $N = 1$ and $\mathbb{X} = \R$ the local smoothing estimates allow to make sense of weak solutions, see Theorem \ref{thm:Gardneruniqueness}. In this case $N=1$ for general geometry $\mathbb{X}$ we can use the equation for the good variable $v$ instead, as Killip and Visan do. We give this argument now.

The approximate KdV flow is defined by the Hamiltonian 
\[ 
\begin{split} (2i\tau)^2\mathcal{T}^\KdV_0 (i \tau)\, &  = (2\tau)^4 \T^{\KdV}_{-1}(i\tau,u)+ 2\tau^2 \int u^2 dx
\\ & =      \frac{i(2i\tau)^5}2 \ln T^{\KdV}(i\tau) - \frac{(2i\tau )^4}2 \int u dx - \frac{(2i\tau )^2}2 \int u^2 dx .   
\end{split} 
\]
The functional  $ \frac12 \int u^2 dx $ generates the translations,
hence (see Lemma~\ref{lem:poissonbrackets}) 
\begin{equation}\label{eq:hap} \begin{split} \hspace{1cm}& \hspace{-1cm} \{ v(\tau)  , (2i\tau_1)^2 \mathcal{T}^{\KdV}_0(i \tau_1)\} = \frac{4\tau_1^4}{\tau^2_1-\tau^2} \partial_x \frac{v(\tau)-v(\tau_1)}{v(\tau_1) +1} + 4\tau_1^2 \partial_x v(\tau)
\\ & =  \frac{4\tau_1^2}{\tau_1^2-\tau^2} \partial_x \Big[\frac1{1+v(\tau_1)}  \big( \tau_1^2( v(\tau_1) - v(\tau)) + (\tau^2_1-\tau^2) v(\tau)(v(\tau_1)+1  )   \big) \Big] 
\\ & =   \frac{4\tau^2_1}{\tau^2_1- \tau^2} \partial_x \Big[ 
\frac{ \tau_1^2  v(\tau_1)(v(\tau)+1)        }{ v(\tau_1)+1}- \tau^2 v(\tau)\Big]. 
\end{split} 
\end{equation}

We recall the good variable version of KdV \eqref{eq:goodvariablekdv} 
\begin{equation} \label{eq:KdVgood} \begin{split} 
v_t \, & = 2 \partial_x \Big[ (v+1) \Big( -(2  \tau)^{2} \frac{\delta H_{-1}}{\delta u} + \frac{\delta H_0}{\delta u }\Big) \Big] 
\\ & = \partial_x \Big[ -4\tau^2 v+ 2(v+1)u \Big].
\end{split}
\end{equation}

The building blocks for the proof of well-posedness for KdV - following Killip and Visan - are also building blocks for the higher order KdV equations.  

\noindent
{\bf  Step 1: Well-posedness of the approximate flow.} The approximate flow is the flow of the Hamiltonian vector field 
$(2i\tau)^2 \T^{\KdV}_0(i\tau, u)$
which we discussed above. The term $\frac{(2\tau)^2}2 \int u^2 dx$ generates translations
and we have seen that the Hamiltonian $ \T^\KdV_{-1}(i\tau,u)$ defines a global in time smooth flow, hence the same is true for the Hamiltonian $\T^\KdV_{0}(i\tau, u)$.

\noindent
{\bf Step 2: Equicontinuity.} Equicontinuity of a bounded set $Q \subset H^{-1}$  can be characterized as 
\begin{equation} \label{eq:equicon}  Q \subset H^{-1} \text{ equicontinuous} \Longleftrightarrow \lim_{\tau\rightarrow \infty}\sup_{u\in Q}\|u\|_{H^{-1}_{\tau}}=0 \Longleftrightarrow  \lim_{\tau \to \infty} 
\sup_{u \in Q}   |\T_{-1}^{\KdV} (i\tau, u)| = 0. \end{equation} 
This fact is an immediate consequence of \eqref{eq:T-1},
\[ \T^{\KdV}_{-1}(i\tau ,u) = -\frac12 \int w(\tau)^2 dx \]
with $w$ satisfying \eqref{eq:uw}. Notice that $\|w\|_{L^2}$ is controlled by $ \|u\|_{H^{-1}_{\tau}}$ and the distance between the square root of the negative of the lowest eigenvalue $-\kappa_0^2$ and $ \tau$ (see  Lemma~\ref{lem:uw}).
Since $ \T_{-1}^{\KdV}( i\tau,\cdot ) $ is preserved under the $ \T_{-1}^{\KdV}(i\tau_1,\dot) $ flow also equicontinuity is preserved along the flow. 

\noindent{\bf Step 3: Convergence of the  difference vector field  in $H^{-2}$ to zero uniformly on equicontinuous sets.}  
Let $Q\subset H^{-1}$ be an equicontinuous bounded set. We claim 
\begin{equation}\label{uniformconv}    \lim_{\tau_1 \to \infty}  \sup_{u \in Q}  \Big\Vert \big\{ v(\tau,u) , \T^{\KdV}_1(i\tau_1,\cdot)  \big\}^\Gardner  \Big\Vert_{H^{-2}} = 0. 
\end{equation} 
The vector field here is the $\T^\KdV_1(i\tau_1,u)= (2i\tau_1)^2 \T^\KdV_0(i\tau_1,u)- H^\KdV_1(u)$  vector field in the $v$ coordinates, given by \eqref{eq:KdVgood} (see also \eqref{eq:hap}: 
\[ \begin{split} 
\big\{ v(\tau), \T^\KdV_1(i\tau_1, . )\big\}\, 
\hspace{-2cm} &\hspace{2cm}   = 
\big\{ v(\tau),  (2i\tau_1)^2\T^\KdV_0(i\tau_1, .) - H_1^{\KdV}(.) \big\} 
\\  &  = \partial_x \Big[   \frac{4\tau_1^2}{\tau_1^2-\tau^2} \big(  \frac{\tau_1^2 v(\tau_1)(v(\tau)+1)}{v(\tau_1)+1} - \tau^2  v(\tau)\big)  + 4\tau^2 v(\tau) - 2(v(\tau)+1) u )\Big] 
\\ & = \partial_x \Big[\Big( \frac{4\tau_1^2 v(\tau_1)}{v(\tau_1)+1}-  2u \Big) (v(\tau)+1) +\frac{4\tau^2\tau_1^2}{\tau_1^2-\tau^2} \frac{  v(\tau_1)(v(\tau)+1)}{v(\tau_1)+1} - \frac{4\tau^2}{\tau_1^2-\tau^2}v(\tau)   \Big].
\end{split} 
\]
The linearization of 
 \[  w(\tau_1) = -\frac12 \partial_x \ln(v+1) + \tau_1 v,    \qquad u=w_x + 2\tau_1 w + w^2 \] 
 is  (with a dot denoting the linearized variables)
\[ -\frac12 \dot v_{xx} +2\tau_1^2 \dot v          = \dot u. \]
It is not hard to see that the linearization dominates when $\tau_1 $ is large. Thus 
\[  \lim\limits_{\tau_1\to \infty} \tau_1^2 v(\tau_1) \to \frac{u}2 \text{ in }  H^{-1}, \qquad  \lim_{\tau_1\to \infty}  \Vert v(\tau_1) \Vert_{H^1} = 0. \]  
uniformly on bounded equicontinuous sets and  
\[   \frac{v(\tau_1)}{v(\tau_1)+1} \to 0  \qquad \text{ in } L^2, \qquad   \frac{4\tau_1^2 v(\tau_1)}{v(\tau_1)+1} -2u \to 0 \qquad \text{ in } H^{-1}     \]
as $ \tau_1 \to \infty$ uniformly on equicontinuous sets. 
This implies the claimed uniform convergence to zero in $H^{-2} $ on equicontinuous sets.

\noindent{\bf Step 4: The difference flow.}
Let $ \tau_1,\tau_2 \ge 1$ and consider the difference  flow
\[  u_t =   \partial \frac{\delta}{\delta u } \Big[  - (2\tau_2)^2 \T^{\KdV}_0( i\tau_2, \cdot ) + (2\tau_1)^2\T^{\KdV}_0 (i\tau_1,\cdot )\Big] . \]
In order to keep the notation brief we introduce the formal notation 
\[     \exp \Big(  t\partial \frac{\delta}{\delta u} H \Big)  u_0 \]
for the solution to the Hamiltonian equations with Hamiltonian $H$ and initial data $u_0$. 
The Hamiltonians $\T^\KdV_0(i\tau_1) $ and $\T^\KdV_0(i\tau_2)$ Poisson commute by Proposition \ref{prop:commutation}. They are continuously differentiable and the generated flows commute, which is expressed as 
\[ 
\begin{split}  
u ( \tau_1,\tau_2,s,t) \, & := \exp\Big(- t (2\tau_2)^2 \partial \frac{\delta}{\delta u} \T^{\KdV}_0(i\tau_2,\cdot) - s(2 \tau_1)^2 \partial \frac{\delta}{\delta u} \T^{\KdV}_0 (i\tau_1,\cdot) \big\}  \Big)u_0 
\\ & = \exp\Big(-t (2\tau_2)^2 \partial \frac{\delta}{\delta u} \T^{\KdV}_0(i\tau_2,\cdot)\Big) \exp\Big( - s  (2\tau_1)^2 \partial \frac{\delta}{\delta u} \T^{\KdV}_0 (i\tau_1,\cdot) \Big)u_0  
\\ & = \exp\Big(-s (2\tau_2)^2 \partial \frac{\delta}{\delta u} \T^{\KdV}_0(i\tau_1,\cdot)\Big) \exp\Big( - t  (2\tau_1)^2 \partial \frac{\delta}{\delta u} \T^{\KdV}_0 (i\tau_2,\cdot) \Big)u_0  
\end{split} 
\]

The set 
\begin{equation} \label{eq:Qdef}  Q= \Big\{ \exp(t_1 \partial \frac{\delta}{\delta u}\T^{\KdV}_0(i\tau_1,\cdot)) \exp( t_2 \partial \frac{\delta}{\delta u} \T^\KdV_0(i\tau_2,\cdot) ) u_0 :    \tau_1,\tau_2 > \tau_0, t_1,t_2 \in \R  \Big\} \end{equation} 
is bounded and equicontinuous in $H^{-1}$ by Step 2 -   the functionals  $\T^\KdV_{-1}(i\tau, .) $ are constant on $Q$. Let (by an abuse of notation) $v(t)=v(\tau_1,\tau_2,t) $ be the $v$ function corresponding to $u( \tau_1,\tau_2,t,t) \in Q$. Then (we consider $u \to v$ as a change of coordinates)  
\[ \begin{split}  v_t \, & = \{ v , -(2\tau_2)^2 \T^{\KdV}_0(i\tau_2,\cdot)  + (2\tau_1)^2 \T^{\KdV}_0(i\tau_1,\cdot) \} 
\\ & = \{ v, \T^{\KdV}_1(i\tau_2,\cdot) \} - \{ v, \T^{\KdV}_1(i\tau_1,\cdot) \} 
\end{split} 
\]
and (suppressing some arguments) 
\[\begin{split} \hspace{1cm} & \hspace{-1cm}  \Vert v(t) - v_0 \Vert_{H^{-2}}    \le \int_0^t  \Vert \{ v , \T^{\KdV}_1(i\tau_2,\cdot)\} - \{v, \T^{\KdV}_1(i\tau_2,\cdot) \}\Vert_{H^{-2}} ds 
\\ & \le |t| \Big(    \sup_{u \in Q} \Vert \{ v, \T^{\KdV}_1(i\tau_2,\cdot) \}(v(u)) \Vert_{H^{-2}} 
+ \sup_{u \in Q}  \Vert \{ v, \T^{\KdV}_1(i\tau_1,\cdot) \}(v(u)) \Vert_{H^{-2}}\Big) 
\\ & \to 0 
\end{split} 
\]
as $ \tau_1,\tau_2 \to \infty$  by Step 3, uniformly on compact time intervals.

\noindent{\bf Step 5: Convergence of the approximate flow.}
We want to prove that $ e^{t (2\tau)^2\partial \frac{\delta}{\delta u} \T^{\KdV}_0(i\tau,\cdot)}u_0\in L^2$ is a Cauchy sequence in $\tau$ for $t$ in compact intervals for $ \tau_1, \tau_2 > \tau_0$. By commutativity of the flows   (and a suggestive abuse of notation) 
\[\begin{split}  \hspace{1cm} & \hspace{-1cm} u(\tau_1,\tau_2,t) :=  \exp\Big(  t  (2\tau_2)^2\partial \frac{\delta}{\delta u} \T^{\KdV}_0( i\tau_2,\cdot) \Big) u_0 
- \exp\Big(  t (2\tau_1)^2 \partial \frac{\delta}{\delta u} \T^{\KdV}_0( i\tau_1,\cdot) \Big) u_0
\\ & =
    \Big\{\exp\Big(  t (2\tau_2)^2 \partial \frac{\delta}{\delta u} \T^{\KdV}_0( i\tau_2,\cdot) \Big) 
  \exp\Big( - t (2\tau_1)^2\partial \frac{\delta}{\delta u} \T^{\KdV}_0( i\tau_1,\cdot) \Big)-\id
  \Big\}\\ & \quad  \times  \exp\Big(t (2\tau_1)^2 \partial \frac{\delta}{\delta u} \T^{\KdV}_0(i\tau_1,.) \Big) u_0.
\end{split} 
\]
Let $\tau_0 \le  \tau_1,\tau_2$  and $Q\subset H^{-1}$ be the equicontinuous set as in \eqref{eq:Qdef}. 
Let $ v(\tau_1,\tau_2,t) $  the $v$ function corresponding to $u(\tau_1,\tau_2,t)$ and $ v(\tau_1,t) $ the one corresponding to 
\[  \exp\Big( - t (2\tau_1)^2\partial \frac{\delta}{\delta u} \T^{\KdV}_0( i\tau_1,.) \Big) u_0 \in Q.\]

By Step 4 and equicontinuity 
\[ \lim_{\tau_1,\tau_2 \to \infty} \Vert v(\tau_1,\tau_2,t) -v_0 \Vert_{H^{-2} } = 0. \]
However all functions $u(\tau_1,\tau_2,t)$ are in the  fixed equicontinuous set $Q\subset H^{-1}$ (and the corresponding functions $v$ are equicontinuous in $H^1$), hence 
\[ \lim_{\tau_1,\tau_2\to \infty} \Vert v(\tau_1,\tau_2,t)-v(\tau_1,t)  \Vert_{H^{1} } = 0.\] 
Thus   $ v(t,\tau_1)\in H^1 $ is  Cauchy  in $\tau_1$, uniformly for $t$ in compact sets. It  extends to a continuous map (with a small abuse of notation) 
\[  H^{-1} \times \R \times  (\tau, \infty] \ni (u_0,t,\tau_1) \to u(t,\tau_1) \in H^{-1}     \] 
resp.
\[  H^{-1} \times \R \times  (\tau, \infty] \ni (u_0,t,\tau_1) \to v(t,\tau_1) \in H^{1}.     \] 
We have proven   a slightly stronger version of the seminal theorem of Killip and Visan \cite{MR3990604}. 
  \begin{theorem}   \label{thm:KdV}
Let $\mathbb{X}=\R$ or $\mathbb{X}= \mathbb{S}^1$ and $Q \subset H^{-1}(\mathbb{X}) $ by an equicontinuous bounded subset of $H^{-1}(\mathbb{X}) $ and let $ \tau $ be sufficiently large and $ \tau_1 >\tau$. Then the approximate flow  defined by the Hamiltonian equation
\[ u_t = \partial\big(  -(2\tau_1)^2 \frac{\delta}{\delta u} \mathcal{T}^\KdV_0 (i  \tau_1,\cdot) \big) 
\] 
with initial data in $Q$ has a unique global solution $u(t,\tau_1)$ in $C(\R,H^{-1}(\mathbb{X}))$.  The set $ \{ u(t,\tau_1) : u_0 \in Q, t \in \R, \tau_1 > \tau \}\subset H^{-1}  $ is bounded and equicontinuous. 
The good variable $v(t,\tau_1) $ converges  in $H^1(\mathbb{X})$ uniformly on compact time intervals  to a weak solution of  
\eqref{eq:KdVgood}, the good variable KdV equation, as $ \tau_1 \to \infty$.  The flow commutes with the $ \tau $ flows. Higher Hamiltonians are preserved.
\end{theorem} 

The claim on commutation and preservation of Hamiltonians are an immediate consequence of the construction and Proposition \ref{prop:commutation}.
We can use Theorem \ref{thm:vwweak} and Theorem \ref{thm:equivalenceweaksolutionsuw} to translate this result to the Gardner  equation.

\subsection{The generating function of the Gardner hierarchy}

For $N \geq 2$ we consider the Gardner equations instead of the good variables equation. 
 The starting point is the $\tau$ flow for the Gardner hierarchy defined by the generating function of the Gardner hierarchy (recall \eqref{eq:T-1},
also see \eqref{eq:Gardner0}):
\begin{equation}\label{eq:gardnergenerating}
\begin{split} 
\T_{-1}^{\Gardner}(z,w,\tau_0)\, & = \frac{1}{4\tau_0^2+4z^2  } \Big(\frac12 \int w^2\, dx +\T^{\KdV}_{-1}(z,w_x+2\tau_0 w+w^2)  \Big)\\ &  =  \frac1{8(\tau_0^2+ z^2)}  \int w^2 -w^2(z) dx  . 
\end{split} 
\end{equation}
 where we recall that $w(z)$ is defined by the left Jost solution $\phi_l$ for $ -\partial^2 + u - z^2$, 
\begin{equation}\label{eq:w(z)}   w(z)  = \partial_x \ln \phi_l + i z \end{equation}
or equivalently, as unique solution to 
\begin{equation} 
    w_x(z) -2izw(z) + w^2(z) = w_x +2\tau_0 w +w^2. 
\end{equation} 
The  relevance of the generating function $\T^\Gardner_{-1}( z, w, \tau)$ comes from the fact 
that the  map $w \to u = w_x+2\tau w + w^2$ is a diffeomorphism in a large class of function spaces, see Proposition  \ref{prop:equivalencewu} and allows to relate solution to KdV equations in $H^{-1}$ to solutions to equations of the Gardner hierarchy. 
In this Subsection we give bounds for $ \T^\Gardner(z, w, \tau_0)$ for $ w \in L^2 $ and $\im z \ge \tau_0$, Proposition \ref{prop:T2} similar to Lemma \ref{lem:kdvestT1},
bound Frech\'et derivatives of its variational derivative which implies wellposedness of the Hamiltonian flow of $\T^\Gardner( i \tau, w, \tau_0)$ for $ \tau > \tau_0$
in $L^2$ together with a number of bounds, Theorem \ref{thm:flow}.  Moreover the Hamilonians 
Poisson commute for different $ \tau$ and hence  the flows commute by Proposition \ref{prop:commutation}. We define approximate Hamiltonians $\T^\Gardner_{app,N}(i\tau, w, \tau_0)$ for $ H^\Gardner_N( w, \tau_0)$, give bounds for them and explain their relation to equicontinuity. They will play an important role in the sequel.  The next proposition gives bounds for $ \im z > \tau_0$, which we will not need, but they complement the bounds for $ \im z$ large.

\begin{proposition} \label{prop:T2} 
Let $ 0 <  \tau_0 $. $\T^{\Gardner}_{-1}(z, w, \tau_0)$ is holomorphic in $z$ for 
\[  \{ z: \im z >0 \text{ and either } \real z \ne 0 \text{ or } \real z =0 \text{ and } 
\im z > \tau_0 \}.  
\]
It satisfies for $ \im z > \tau_0$ 
\[  |\T^{\Gardner}_{-1}(z, w, \tau_0) | \lesssim   \frac1{ \im z |z-i\tau_0|} \Vert w \Vert^2_{L^2}\]
\end{proposition}
\begin{proof} 
It suffices to consider $\tau_0=1$. By \eqref{eq:trace.h-1} 
\begin{equation} 
\begin{split} \hspace{1cm}& \hspace{-1cm} 
\T^{\Gardner}_{-1}(z, w, 1) =  \frac1{8+ 8 z^2} \left( \Vert w \Vert_{L^2}^2-  \int w^2(z) dx  \right)\\ &   = \frac1{4+ 4z^2} \Big(-\T^\KdV_{-1}(i,u)+ \T^{\KdV}_{-1}( z, u)\Big)
\\ & =  \frac1{4+  4z^2}  \Big[ \frac1{2\pi}  \int_{\R}\Big( \frac1{i\xi+1} -\frac1{z(\xi -z)}\Big)     \big[\real 2i \xi \T^{\KdV}_{-1} ( \xi, u)\big]  d\xi  
\\ & \quad     +   \sum_j   \ln  \frac{1+\kappa_j}{1-\kappa_j}   +iz \ln \frac{z+i\kappa_j}{z-i\kappa_j} \Big]. 
\end{split} \label{eq:Gardnertrace}
\end{equation} 
By \eqref{eq:monoton}  
\[ 0\le  -\T^\Gardner_{-1}(i\tau,w,1) =  \frac{-1}{4(\tau^2-1)} (\T_{-1}^\KdV(i\tau,u)- \T_{-1}^\KdV(i,u)) \le \frac1{8(\tau^2-1)} \Vert w \Vert_{L^2}^2.
\]
We obtain 
\[ |\T^\Gardner_{-1}(2i,w,1) | \le \frac1{6} \Vert w \Vert_{L^2}^2. \]
We apply formula \eqref{eq:Gardnertrace} first with $z= 2i$ to obtain bounds for the integral and the sum below, using 
$ \kappa_j <1$ and that $  \real (2iz  \T^\KdV_{-1}(\xi, u))$ is even in $\xi$ and nonnegative
\[\begin{split} \hspace{0.5cm}& \hspace{-0.5cm} \frac18\frac1{2\pi} \int \frac{2}{1+\xi^2} \real [ 2i \xi \T^{\KdV}_{-1}(\xi, u)] d\xi  
+ \frac1{24} \sum_j \Big(\ln \frac{1+\kappa_j}{1-\kappa_j}- 2\kappa_j\Big) 
\\ & \le  \frac1{6} \Big[  \frac1{2\pi} \int \real ( \frac1{i\xi+1}-\frac1{2i\xi +4})  \real [ 2i\xi \T^{\KdV}_{-1}(\xi,u) ] d\xi \\ &\qquad +  \sum_j \ln \frac{1+\kappa_j}{1-\kappa_j} - 2 \ln \frac{1+\kappa_j/2}{1-\kappa_j/2}  \Big] 
\\ & =  -\real \T^\Gardner_{-1}(2i,w,1)\le \frac16 \Vert w \Vert_{L^2}^2.
\end{split}\] 
We complete the estimate for $\im z \ge 1$ by 
\[ \left|\frac{1+\xi^2}{z(\xi-z)}+ \frac{1+\xi^2}{z(-\xi -z)} \right|  \le      \frac{2(1+\xi^2)} { |z^2- \xi^2|} 
\lesssim   \frac{|z|}{\im z}
\]
and 
\[ 
\Big|z \ln \frac{1-i\kappa_j/z}{1+i\kappa_j/z} -2i\kappa_j\Big| \le \left\{ \begin{array}{rl} \kappa_j^3/|z|^2 & \text{ if } |\kappa_j|\le \frac12 |z| 
\\   \kappa_j  \ln |1-i\kappa_j/z| & \text{  if } |z| \le 2\kappa_j  \end{array} \right\} \lesssim \ln\frac{1+\kappa_j}{1-\kappa_j}- 2\kappa_j.
\]
A comparison with the formula for $ \real \T^\Gardner(2i, w, 1) $ and the bound above imply for $\im z \ge \tau_0=1 $ 
\[ |\T_{-1}^\Gardner(z, w, i )| \lesssim   \frac{1}{|1+ z^2|} \frac{|z|}{\im z}    |\T_{-1}^\Gardner(2i,w,i) |
\lesssim   \frac{1}{ \im z |z-i| }       \Vert w \Vert^2_{L^2}.  
\]

\end{proof}

We will need various bounds on the variational derivatives, including  
 weighted estimates in later parts with slowly varying weights. 
 
 \begin{definition} We say $\gamma$  is $\tau$-slowly varying (see also Definition \ref{def:slowlyvarying}) if 
\begin{equation} |\eta_x| \le \tau \eta, \qquad |\eta^{(j)}|\le c_j \tau^j \eta. \end{equation}
\end{definition} 

The generating function  $\T^\Gardner_{-1}(z,w,\tau_0)$ generates a Hamiltonian flow which satisfies very strong estimates, which are central for the proof of wellposedness of the Gardner hierarchy in $L^2$.

\begin{theorem}\label{thm:flow}  The generating function $ \T^\Gardner(i\tau , w , \tau_0)$ is real and smooth in $w \in L^2$. Let $ \gamma$, $ \gamma_1$ and $\gamma_2$ be $\tau/2$ slowly varying. The  variational derivative of the generating function $\T^\Gardner_{-1}(i\tau, w , \tau_0)$ satisfies 
\begin{equation} \label{eq:T1smooth0}
\Big\Vert \gamma \frac{\delta}{\delta w} \T^\Gardner_{-1}(i\tau, w, \tau_0) \Big\Vert_{H^{N+1}}
\le c  \Vert \gamma w \Vert_{H^N} 
\end{equation} 
and 
\begin{equation} \label{eq:T1smooth}
\big\Vert \gamma_1 \gamma_2    D  \frac{\delta}{\delta w}\T^\Gardner_{-1}(i\tau, w , \tau_0) [ w_1]\big\Vert_{H^{N+1}}  \le c \big( \Vert \gamma_1 w \Vert_{H^N} \Vert \gamma_2 w_1 \Vert_{L^2} 
+\Vert \gamma_1\gamma_2  w_1 \Vert_{H^N} \big)   
\end{equation}
with a constant depending on $\tau_0^{-1} \Vert w \Vert_{L^2}$ and $ \frac{\tau}{\tau-\tau_0}$ and $N$(here  $D$ denotes the total derivatives resp. the Frech\'et derivative). 
It is continuous in $w$ in the corresponding norms. 
It generates a smooth global Hamiltonian flow on $H^N$ which we denote by 
\[ \exp\Big( t \partial \frac{\delta}{\delta w} \T^\Gardner_{-1}( i \tau,., \tau_0\Big) \]
on $H^N$ for $N \ge 0$. It satisfies 
\[  \Big\Vert \gamma  \exp\Big( t \partial \frac{\delta}{\delta w} \T^\Gardner_{-1}( i \tau,., \tau_0)\Big)(w) \Big\Vert_{H^{N}} \le c \Vert \gamma w \Vert_{H^N}.  \]
Its  Frech\'et derivative satisfies 
\begin{equation}  \label{eq:weight}
\Big\Vert \gamma_1  \gamma_2  D \exp\Big( t \partial \frac{\delta}{\delta w} \T^\Gardner_{-1}( i \tau,w, \tau_0\Big) \big[ \dot w\big] \Big\Vert_{H^{N}} 
 \le c  \Big(  \Vert \gamma_1 w \Vert_{H^N} \Vert \gamma_2 \dot w \Vert_{L^2} + \Vert \gamma_1 \gamma_2 \dot w \Vert_{H^N} \Big).  
\end{equation} 
 where $c$ depends on $t$,   $\frac{\tau}{\tau- \tau_0}$ and $ \tau_0^{-1/2} \Vert w \Vert_{L^2}$.  
The map is continuous in $w$ and $\dot w$ in the appropriate norms.  
\end{theorem}

We compute the variational derivative of $\T^\Gardner_{-1}$.  
\begin{proposition}
Let $ \tau> \tau_0$. Then $ \T_{-1}^{\Gardner}(i\tau, w,\tau_0) \in \R$ and 
\begin{equation}\label{eq:variationalder} 
\begin{split} \hspace{.5cm} & \hspace{-.5cm} 
\frac{\delta}{\delta w} \T^{\Gardner}_{-1}(i\tau, .,\tau_0)(w) 
\\ &  =\frac1{4\tau^2 -4 \tau_0^2}\big[(-\partial + 2\tau_0 +2w)(-\partial + 2\tau + 2w(i\tau))^{-1} w(i\tau)-w\Big]
\\ &  =  \frac1{4\tau^2 -4 \tau_0^2}\big[w(i\tau)-w - 2 (\tau-\tau_0 + w(i\tau)- w) (-\partial + 2 \tau + 2w(i\tau))^{-1} w(i\tau) \big]. 
\end{split} 
\end{equation} 
\end{proposition}

\begin{proof} The functional $\T^\Gardner_{-1}(i\tau, w, \tau_0)$ is given in \eqref{eq:gardnergenerating}.
 Let  $ \tau > \tau_0$. By Lemma \ref{lem:uw} the solution  $w(i\tau) \in L^2$  to 
\begin{equation} \label{eq:wit}    \partial w(i\tau) + 2 \tau w(i\tau) + w^2(i\tau)= w_x+ 2\tau_0 w + w^2 \end{equation}
is uniquely determined and satisfies by the monotonicity statement in Lemma  \ref{lem:kdvestT1}
$\Vert w(i\tau) \Vert_{L^2} \le \Vert w \Vert_{L^2}$
. 
We compute with $ w(s) = w+ s \phi$ and the corresponding $w(i\tau,s)$ the Fr\'echet derivative  \[ Dw(i\tau)[\phi]=  \frac{d}{ds} w(i\tau,s) \Big|_{s=0} = ( \partial + 2 \tau + 2 w(i\tau))^{-1} ( \partial + 2\tau_0 + 2 w) \phi  \]
\[   \frac{d}{ds} \frac12 \int  w^2(i\tau,s) - (w+s\phi)^2 dx \Big|_{s=0} 
= \int  \Big[ (\partial + 2 \tau_0 + 2 w) \phi \Big]( - \partial + 2 \tau + 2 w(i\tau))^{-1} w( i\tau) 
- \phi w dx 
 \]
which implies the first identity in \eqref{eq:variationalder}. The second identity follows from 
\[ \begin{split} \hspace{1cm} & \hspace{-1cm}   (-\partial + 2\tau_0 +2w)(-\partial + 2\tau + 2w(i\tau))^{-1} w(i\tau)-w  
\\ & = (w(i\tau)-w) +(2(\tau-\tau_0 -( w(i\tau)-w) )( -\partial + 2\tau + 2w(i\tau))^{-1} w(i\tau). 
\end{split} 
\] 
\end{proof}

\begin{proof}[Proof of Theorem \ref{thm:flow}]

The function $ \T^{\Gardner}_{-1}(i\tau, w, \tau_0) $ is real valued by construction. 
Suppose that the estimates \eqref{eq:T1smooth0} and \eqref{eq:T1smooth} hold.
The Hamiltonian equations of the Hamiltonian $\T^\Gardner_{-1}(\tau, w, \tau_0)$ are 
\[ 
  w_t  = \partial \frac{\delta}{\delta w } \T^\Gardner_{-1}(i\tau, w , \tau_0). 
\]
The estimates \eqref{eq:T1smooth0} and \eqref{eq:T1smooth} with $\gamma=1$ resp. $ \gamma_1=\gamma_2=1$ imply that the right hand side of the equation is bounded and Lipschitz continuous on $H^N$.
Local existence in $L^2$ in $H^N$  follows then the Cauchy-Lipschitz theorem.  The $L^2$ norm is conserved: 
\[ \begin{split} \frac{d}{dt} \frac12 \Vert w(t) \Vert_{L^2}^2\, &  = -\int \frac{\delta}{\delta w} \T^\Gardner_{-1} ( i\tau, . ,\tau_0) \partial_x w dx \\ & =  \frac{d}{ds} \T^\Gardner_{-1}(i\tau , w( x+ s) ,\tau_0) |_{s=0} = 0 \end{split}   \]
since   $\T_{-1}^\Gardner$ does not change when one translates $w$. Thus the flow is global in $L^2$.  Global wellposedness in $H^N$ is in the same fashion a consequence of Proposition \ref{prop:commuting} and Proposition \ref{prop:diffham}: The $N$ Hamiltonian $H^\Gardner_N(w,\tau_0)$ is conserved and controls (together with the $L^2$ norm) $ \Vert w \Vert_{H^N}$. We  complete the proof in the weighted norms after the proof of Proposition \ref{prop:diffham}.

We turn to the proof of \eqref{eq:T1smooth0} and begin with the linear equation 
\[    (\partial_x +  2\tau+ \tilde w) \phi = f.   \]
If $ \tilde w \in L^2$, $ f \in L^2 $ there exists a unique solution $ \phi \in H^1$. If moreover  $ \gamma_2 \tilde w \in H^N $ and $ \gamma_1 \gamma_2 f \in H^N$  then 
\begin{equation}\label{eq:linesttw}  \Vert \gamma \phi \Vert_{H^{N+1}_\tau} \le c ( \tau^{-1/2} \Vert \tilde w \Vert_{L^2}) 
( \Vert \gamma_1\gamma_2 f \Vert_{H^N} + \Vert \gamma_1 \tilde w \Vert_{H^N} \Vert \gamma_2 f \Vert_{L^2} ) 
\end{equation} 
by Lemma \ref{lem:linearesteq}, \eqref{eq:linwinverse}. Lemma  \ref{lem:kdvestT1}  gives the bound 
$\Vert w(i\tau) \Vert_{L^2} \le  \Vert w \Vert_{L^2}$. 
 All constants in the sequel may depend on $ \tau$, $\tau_0$ and $ \tau_0^{-1/2}\Vert w \Vert_{L^2} $.
 With $ \tilde w = w(i\tau)- w $, 
\[  (\partial_x + 2 \tau +  (2w+ \tilde w) )  \tilde w = -2(\tau -\tau_0)w. \]
and we apply \eqref{eq:linesttw}  possibly changing from line to line, by induction on $N$, for $n \ge 0$,
\begin{equation}\label{eq:tildew} \Vert \gamma ( w(i\tau) -w)\Vert_{H^{N+1}}=   \Vert  \gamma  \tilde w \Vert_{H^{N+1}_\tau} \le  c 
( \Vert \gamma  w \Vert_{H^N} + \Vert \gamma \tilde w \Vert_{H^N})
\le c  \Vert \gamma w \Vert_{H^N} 
\end{equation}
The last inequality follows from iterating the estimate and notice $\Vert \gamma \tilde w \Vert_{L^2}\lesssim \|\gamma w\|_{L^2}$.
We estimate the variational derivative $\frac{\delta}{\delta w} \T^\Gardner_{-1}$ using the second identity  \eqref{eq:variationalder}).
The first summand   $w(i\tau)-w$ is bounded by \eqref{eq:tildew}. Next we bound
\[  \begin{split}  \Vert \gamma  (-\partial+2\tau + 2 w(i\tau))^{-1} w(i\tau)\Vert_{H^{N+1}} 
\, & \le c  \Vert \gamma  w(i\tau) \Vert_{H^N} 
\\ & \le c\big(  \Vert \gamma w \Vert_{H^N} 
+ \Vert \gamma( w(i\tau) - w) \Vert_{H^N}\big) 
\\ & \le c \Vert \gamma w \Vert_{H^N}
\end{split} 
\]
by induction and \eqref{eq:tildew}. 
\[ \begin{split} \Vert \gamma ( \tau-\tau_0 + w(i\tau)-w  ) \phi \Vert_{H^{N}_\tau} \, & \le c \big(\Vert \gamma \phi \Vert_{H^N} \Vert \tilde w \Vert_{H^1} + \Vert \gamma \tilde w \Vert_{H^{N+1}} \Vert \phi \Vert_{L^2}\big).
\\ & \le c ( \Vert  \gamma  \phi \Vert_{H^N} + \Vert \gamma w \Vert_{H^N} \Vert \phi \Vert_{L^2}).
\end{split} 
\]
Together we obtain 
\[ \Big\Vert \gamma \frac{\delta}{\delta w} \T^\Gardner_{-1}(z, w, ,\tau_0) \Big\Vert_{H^{N+1}} \le c \Vert \gamma w \Vert_{H^N}. \]
It remains to estimate Frech\'et derivatives. We linearize the various operators and functions and denote linearized variables with a dot. Then  $ D\frac{\delta}{\delta w} \T^\Gardner_{-1}(i\tau, w, \tau_0) [\dot w] $ is given by $ \frac1{4\tau^2-4\tau_0^2} $ times 
\[ 
\begin{split} \hspace{1cm} & \hspace{-1cm} 
\dot w(i\tau) - \dot w  - 2 ( \dot w(i\tau)-\dot w ) ( -\partial+2\tau + 2 w(i\tau))^{-1} w
\\ & + 4(\tau-\tau_0+ w(i\tau)-w)(-\partial+ 2\tau w + 2 w(i\tau))^{-1} \dot w(i\tau) (-\partial+2\tau+ 2 w(i\tau))^{-1} w(i\tau) 
\\ & - 2 (\tau-\tau_0+ w(i\tau)-w)(-\partial+ 2\tau w + 2 w(i\tau))^{-1}  \dot w(i\tau) 
\end{split} 
\]
and 
\[ \dot w(i\tau) = (\partial +2 \tau + 2w(i\tau) )^{-1} (\partial + 2\tau_0 + 2 w)  \dot w,  \]
\[ \dot w(i\tau) - \dot w =  (\partial+2\tau + 2 w(i\tau))^{-1} (  2(\tau_0-\tau  + w - w(i\tau))\dot w       \]
and we obtain the estimates as above 
\[ \Vert \gamma_1 \gamma_2 ( \dot w(i\tau) - \dot w ) \Vert_{H^{N+1}} \le   c (\Vert \gamma_1 \dot w \Vert_{L^2}  \Vert \gamma_2 w \Vert_{H^N} 
+ \Vert \gamma_1 \gamma_2 \dot w \Vert_{H^N} ). 
\]
For notational simplicity we restrict to $\gamma= \gamma_2$ and $\gamma_1=1$ in the sequel.
\[ 
\begin{split} \hspace{1cm} & \hspace{-1cm} 
\Vert \gamma ( \dot w(i\tau)-\dot w ) ( -\partial+2\tau + 2 w(i\tau))^{-1} w \Vert_{H^{N+1} }
\le \Vert \gamma (\dot w(i\tau) - \dot w \Vert_{H^{N+1} }
\Vert ( -\partial+2\tau + 2 w(i\tau))^{-1} w \Vert_{H^1} 
\\ & \qquad + \Vert \dot w(i\tau) + \dot w \Vert_{H^1}  \Vert \gamma (( -\partial+2\tau + 2 w(i\tau))^{-1} w) \Vert_{H^1} \Vert_{H^{N+1}} 
\\  \le  & c \Vert \gamma \dot w \Vert_{H^N} + \Vert \dot w \Vert_{L^2} \Vert \gamma w \Vert_{H^N}),
\end{split} 
\] 
\[ \begin{split} \hspace{1cm} & \hspace{-1cm} 
\Vert \gamma ( (\tau-\tau_0+ w(i\tau)-w)(-\partial+ 2\tau w + 2 w(i\tau))^{-1} \dot w(i\tau) (-\partial+2\tau+ 2 w(i\tau))^{-1} w(i\tau)) \Vert_{H^{N+1} } 
\\ & \le c\big(  \Vert w(i\tau) - w \Vert_{H^1} \Vert \gamma  (-\partial+ 2\tau w + 2 w(i\tau))^{-1} \dot w(i\tau) (-\partial+2\tau+ 2 w(i\tau))^{-1} w(i\tau)) \Vert_{H^{N+1}}
\\ & \qquad +  \Vert \gamma ( w(i\tau)- w ) \Vert_{H^{N+1} } \Vert(-\partial+ 2\tau w + 2 w(i\tau))^{-1} \dot w(i\tau) (-\partial+2\tau+ 2 w(i\tau))^{-1} w(i\tau)) \Vert_{H^{1}} \big) 
\\ & \le c \big(  \Vert \gamma w  \Vert_{H^N} \Vert \dot w(i\tau) (-\partial+2\tau+ 2 w(i\tau))^{-1} w(i\tau)) \Vert_{L^2}  
\\ & \qquad + \Vert \gamma  \dot w(i\tau) (-\partial+2\tau+ 2 w(i\tau))^{-1} w(i\tau)) \Vert_{H^N} + \Vert \gamma w \Vert_{H^N} \Vert \dot w \Vert_{L^2}\big) 
\\ & \le c \big( \Vert \gamma w \Vert_{H^N}\Vert \dot w \Vert_{L^2} +  \Vert \gamma \dot w(i\tau) \Vert_{H^N} ) 
\\ & \le  c\big(  \Vert \gamma w \Vert_{H^N} \Vert \dot w \Vert_{L^2} + \Vert \gamma \dot w \Vert_{H^N}\big) 
\end{split} 
\] 
and finally 
\[ \begin{split} \hspace{1cm} & \hspace{-1cm} 
\Vert \gamma  (\tau-\tau_0+ w(i\tau)-w)(-\partial+ 2\tau w + 2 w(i\tau))^{-1}  \dot w(i\tau) \Vert_{H^{N+1}} 
\\ & 
\le c \big( \Vert \gamma ( w(i\tau)-w ) \Vert_{H^{N+1}} \Vert (-\partial+ 2\tau w + 2 w(i\tau))^{-1}  \dot w(i\tau) \Vert_{H^1}
\\ & \qquad +  \Vert w(i\tau)-w \Vert_{H^1} \Vert \gamma (-\partial+ 2\tau w + 2 w(i\tau))^{-1}  \dot w(i\tau) \Vert_{H^{N+1}} 
\\ & \le c \big( \Vert \gamma w \Vert_{H^N} \Vert \dot w \Vert_{L^2} 
+ \Vert \gamma w \Vert_{H^N} \Vert \gamma \dot w \Vert_{L^2}+ \Vert \gamma \dot w \Vert_{H^N}.
\end{split} 
\] 
This implies the important estimate \eqref{eq:T1smooth}. As stated above we postpone the proof of the weighted estimates for the evolution 
 until after Proposition \ref{prop:diffham} below.
\end{proof} 

\begin{proposition}\label{prop:commuting} 
Let $ 0 < \tau_0 < \tau_1, \tau_2$. The Hamiltonians $ \T^\Gardner_{-1} (i\tau_2, w , \tau_0) $ and $\T^\Gardner_{-1} ( i\tau_1, w , \tau_0) $ Poisson commute, i.e. 
\[ \int \frac{\delta}{\delta w} \T^\Gardner_{-1}( i\tau_2, .,\tau_0) \partial \frac\delta{\delta w} \T^\Gardner_{-1}( i\tau_1, w, \tau_0) dx = 0. \]
The Hamiltonian $ \T^\Gardner_{-1} ( i \tau_2,w, \tau_0)$ is conserved along the Hamiltonian flow of $ \T^\Gardner_{-1}(i \tau_1, w , \tau_0)$.
The flows commute, .i.e
\[\begin{split} \hspace{1cm} & \hspace{-1cm}  \exp\Big( s \partial \frac{\delta}{\delta w} \T^\Gardner_{-1}(i\tau_2, ., \tau_0)\Big) \circ \exp\Big( t \partial \frac{\delta}{\delta w} \T^\Gardner_{-1}(i\tau_1, . , \tau_0)\Big) 
\\ & =  \exp\Big( t \partial \frac{\delta}{\delta w} \T^\Gardner_{-1}(i\tau_1, ., \tau_0)\Big) \circ \exp\Big( s \partial \frac{\delta}{\delta w} \T^\Gardner_{-1}(i\tau_1, . , \tau_2)\Big). 
\end{split} 
\]
\end{proposition}

The statement about Poisson commutation is contained in  Lemma \ref{lem:GardnerPoisson}.
As explained in in Subsection \ref{subsec:Frobenius} this is equivalent to 
the commutation of the flows in $L^2$, using Theorem \ref{thm:flow}.

We define the difference Hamiltonian for $ N \ge 0$ similarly to the difference  Hamiltonians for KdV by 
\[  \mathcal{T}_N^{\Gardner}(z,w,\tau_0) =    
  \frac{(2z)^{2N+2}}{ 8 ( \tau_0^2 + z^2) }  \int w^2-w^2(z) dx- \sum_{n=0}^N (2z)^{2(N-n)} H_n^{\Gardner}  . \] 
Its bilinear part is  $-\frac12\Vert w^{(N+1)}\Vert_{H^{-1}_{2\tau}}$. We denote 
\begin{equation} \label{eq:TNL}  \mathcal{T}_N^{NL}(i\tau,w, \tau_0) = \mathcal{T}_N^\Gardner(i\tau,w, \tau_0) + \frac12 \Vert w^{(N+1)}\Vert^2_{H^{-1}_{2\tau}}   \end{equation}

To shorten the notation we often write $\T^{\Gardner}_N(z)$ instead of $\T^{\Gardner}_N(z,w,\tau_0)$.
\begin{proposition}\label{prop:diffham}  The Gardner difference Hamiltonians $ \T^\Gardner_N$ are defined in $ \{ z \in C: \im z > \tau_0\} \times H^N$. They are holomorphic in $z$ and analytic in $w$.  On the imaginary axis the difference Hamiltonians are real: $\mathcal{T}_N^{\Gardner}( i  \tau,w,\tau_0 ) \in \R $. 
They satisfy with 
\[   c = c (N, (\tau_0/\im z)^{\frac12} \Vert w \Vert_{L^2}) \]
 \begin{equation} \label{eq:TnGardner}  (2\im z)^2 |\mathcal{T}^{\Gardner}_{N}( z,w, \tau_0) | \le c \Big(\frac{|z|}{\im z} \Big)^{2N+3} \Big(\Vert w \Vert^2_{H^{N+1}_{\tau_0}}+ \Vert w \Vert_{L^{2(N+2)}}^{2(N+2)} + \tau_0^{N+2} \Vert w \Vert^{N+2}_{L^{N+2}}\Big),  \end{equation}
%\begin{equation} \big| \mathcal{T}_N^{NL}(i\tau,w)
%\mathcal{T}^{\Gardner}_N(i\tau,w, \tau_0)  - \frac12 \Vert w^{(N+1)} \Vert^2_{H^{-1}_{2\tau}} \big| 
%\le     c( 1+ \Vert w \Vert_{L^2} )^{2N-1}  \Vert w \Vert_{L^\infty} \Vert w \Vert^2_{H^{N-1}}.   
%\end{equation} 
A set $Q \subset H^N $ is equicontinuous if and only if 
\begin{equation}\label{eq:propequi}  \lim_{\tau \to \infty} \sup_{w \in Q} |\T^{\Gardner}_N(i\tau, w,\tau_0)| = 0. \end{equation} 
 If $Q\subset H^{N-1}$ is equicontinuous  then 
\begin{equation} \label{eq:propuni}  \Big\Vert   \frac{\delta}{\delta w} \T^{ \Gardner}_N(i\tau, w,\tau_0)   \Big\Vert_{H^{-N-2}} \to 0 \quad \text{ as } \tau \to \infty \end{equation} 
uniformly in $ w \in Q$.
\end{proposition} 
The estimates are proven in Section \ref{subsec:diffHamiltonian}. The statement about equicontinuity is in Corollary \ref{cor:Gardnerequi}.
We obtain as an immediate consequence for $ w \in H^N$
\[ \lim_{\tau \to \infty} (2i\tau)^2  \mathcal{T}^{\Gardner}_{N-1}(i\tau,w, \tau_0) = H_N^{\Gardner}(w, \tau_0). \] 
The recursion relation  follows:
\[ \mathcal{T}^{\Gardner}_N(z, w,\tau_0) = -H_N^{\Gardner} + (2 z)^2 \mathcal{T}_{N-1}^{\Gardner}(z, w, \tau_0) .\]

With this information we complete the proof of Theorem \ref{thm:flow} about global wellposedness of the Hamiltonian flow of $ \T^\Gardner_{-1}(i\tau, ., \tau_0)$ in weighted $H^N$. We begin with the unweighted case.  
The flows for different values of $ \tau$ commute by Proposition \ref{prop:commuting}. The same is true for limits and hence $H^\Gardner_N(w,\tau_0)$ is conserved. Since 
\[\Vert w \Vert^2_{H^N} \le c \Big(    H^\Gardner_N(w,\tau_0) +   \Vert w \Vert_{L^2}^{4N+2} \Big) \]
(the quadratic term of $H^\Gardner$ is $ \frac12 \Vert w^{(N)} \Vert_{L^2}^2$, the other terms are described in  Theorem \ref{thm:formofkdv} and estimated by 
\[     \Vert w^{(N)} \Vert^{1+\frac{m}{N}}_{H^N} \Vert w \Vert^{(2+\frac1N)(N-m)}_{L^2}
\le \varepsilon \Vert w^{(N)} \Vert_{L^2}^2 + c(\varepsilon) \Vert w \Vert^{4N+2}_{L^2},       \]
where $ m\le N- 2$, depending on the homogeneity in $w$ and $ \tau_0$ of the term, see estimate \ref{est:GardnerHamEst} in Subsection \ref{subsec:unique}). 
Thus $\Vert w \Vert_{H^N} $ is uniformly bounded in terms of the initial data and there is a global in time solution in $H^N$.  
Next we consider $N=1$ with $\gamma_1=1$. In this case the we can put the weight always on $\Vert \gamma_2 \dot w \Vert_{L^2} $ and the dependence on this term is linear. If $N \ge 1$
we use the bound on $ \gamma_1 \dot w$ and observe that the estimates a linear in 
\[ \Vert \gamma_1 w(t) \Vert_{H^N} \Vert \gamma_2 \dot w(t) \Vert_{L^2}. \]
Thus Grönwall's inequality gives global bounds for the weighted norms which grow exponentially in time. 
\bigskip

\subsubsection{The approximate Gardner Hamiltonian.}
Next we approximate $H^\Gardner_N(i\tau,w,\tau_0)$ by a linear combination of $\T^\Gardner_{-1}(i2^j \tau, w, \tau_0)$.

\begin{lemma} \label{lem:approximate}
Let $ \lambda >1$.
There exist coefficients $ a_j$  so that 
\begin{equation} \label{eq:appdef}   \sum_{j=0}^N   a_j (2i \lambda^j \tau)^{2N+2} \T^{\Gardner}_{-1}(\lambda^j \tau) - H^\Gardner_N = 
\sum_{j=0}^N  a_j  \T^\Gardner_N( i\lambda^j \tau )  
\end{equation} 
with the terms quadratic in $w$ given by the Fourier multiplier
\[  M(\xi)= \Big( \prod_{k=0}^N \frac{\xi^2} {\xi^2 +  (\lambda^k \tau)^2}\Big)  \frac{\xi^{2(N+1)}-(\tau \mu)^{2(N+1)}}{\xi^2 -\mu^2 \tau^2}   \]  
where 
\[ \mu = \frac{\lambda^{2(N+1)}-1}{\lambda^2-1}. \]
Moreover 
\begin{equation}\label{eq:locappsmoothing}  \frac{d}{d\xi} ( \xi M(\xi)) \sim   \frac{\xi^{2N+2}}{  \xi^2 + \tau^2}.  \end{equation} 
\end{lemma}

\begin{proof} 
By the expansion of $\T^{\Gardner}_{-1}$  the coefficients $a_j$ are determined by the linear Vandermonde equation 
\[ \left( \begin{matrix}  1 & 1& \dots & 1 \\ 
                          1 & \lambda^2 & \dots & \lambda^{2N} \\
                          \vdots & \vdots & \ddots & \vdots \\
                          1 & \lambda^{2N} & \dots & \lambda^{2N^2} \end{matrix} \right) 
                          \left( \begin{matrix} a_0 \\ a_1 \\ \vdots \\ a_n \end{matrix} \right) =
                          \left( \begin{matrix} 1 \\ 0 \\ \vdots \\ 0 \end{matrix} \right) 
                                                  \]
which has a unique solution.  The left hand side of \eqref{eq:appdef} is then 
\[ \begin{split}\hspace{1cm} & \hspace{-1cm}  \sum_{j=0}^N a_j \Big( (2i\lambda^j \tau)^{2(N+1)}   \T^\Gardner_{-1}( i\lambda^j \tau,w, \tau_0)- H_N^\Gardner(w,\tau_0) \Big)\\ & = \sum_{j=0}^{N} a_j \Big( \T_{N}^{\Gardner}(i\lambda^j \tau , w, \tau_0)  + \sum_{k=0}^{N-1} (2i \lambda^j \tau)^{2(N-k)} H_k^\Gardner(w,\tau_0)   \Big)
\\ & =  \sum_{j=0}^N a_j \T^\Gardner_N( i \lambda^j \tau, w, \tau_0). 
\end{split} 
 \]
This proves \eqref{eq:appdef}.
The quadratic part is the Fourier multiplier (recall that the quadratic part of $ \T^\Gardner_N( \tau, w , \tau_0)$
is $ -\frac12 \Vert w^{(N+1)} \Vert_{H^{-1}_{2\tau}}$) 
\[ \sum_{j=0}^N \frac{a_j\xi^{2N+2} }{\xi^2+ (\lambda^j \tau)^2} 
= \xi^{2N+2} \Big(  \prod_{k=0}^N \frac1{\xi^2 + (\lambda^k \tau)^2}\Big)
\sum_{j=0}^N a_j \prod_{k\ne j} (\xi^2 + (\lambda^k \tau)^2). 
\]
We recall $ \sum_{j=0}^N a_j = 0 $, 
\[ \sum_{j=0}^N a_j  \sum_{k\ne j } \lambda^{2k} = \sum_{k=0}^N \lambda^{2k} - \sum_{j=0}^N a_j  \lambda^{2j} = \sum_{k=0}^N \lambda^{2k},   \]
\[\begin{split} \sum_{j=0}^N a_j \sum_{k_1\ne k_2\ne j \ne k_1}  \lambda^{2(k_1+k_2)} 
\, & = \sum_{k_1\ne k_2} \lambda^{2(k_1+k_2)} 
- \sum_{j=0}^N a_j \sum_{k=0}^N  \lambda^{2(k+j)} + \sum_{j=0}^N a_j \lambda^{4j}  
\\ & = \sum_{k_1\ne k_2} \lambda^{2(k_1+k_2)}. 
\end{split}
\]
Recursively we obtain similar estimates for larger products. 
We expand in even powers of $ \lambda$ and the previous consideration in the second identity
\[
\begin{split} \hspace{1cm} & \hspace{-1cm} 
\sum_{j=0}^N a_j \prod_{k\ne j}(\xi^2 + (\lambda^k \tau)^2)
 =  \sum_{j=0}^N a_j \Big(  \xi^{2N} + \tau^2 \xi^{2(N-1)} \sum_{k\ne j } \lambda^{2k}\\ & \qquad    + 
\tau^4 \xi^{2(N-2)} \sum_{k_1\ne  k_2\ne  j\ne k_1} \lambda^{2(k_1+k_2) } + \dots + \tau^{2N} \prod_{k\ne j} \lambda^{2k}  \Big)  
\\ & =  \sum_{n=0}^N \xi^{2(N-n)} \tau^{2n} \Big( \sum_{j=0}^N \lambda^{2j}\Big)^n 
\\ & = \sum_{n=0}^N      \left( \frac{\lambda^{2(N+1)}-1}{\lambda^2-1} \right)^{N-n}  \xi^{2n}\tau^{2(N-n)} 
\\ & =    \frac{\xi^{2(N+1)}  -  (\tau \mu)^{2(N+1)}}{\xi^2- (\tau\mu)^2}.
\end{split} 
\]
Equation \eqref{eq:locappsmoothing} is an immediate consequence.
\end{proof} 

In the sequel we will fix $\lambda=2$.
\begin{definition}\label{def:approx} Let $ \lambda=2$ and the coefficients as in Lemma \ref{lem:approximate}.
We define 
\[ \T^\Gardner_{app,N}(i\tau, w, \tau_0)  =\sum_{j=0}^N a_j (2^{j+1} i  \tau)^{2N+2} \T^\Gardner_{-1}(i2^j \tau ,w, \tau_0)  \]
\end{definition} 
 By our previous consideration $\T_{app,N}^{\Gardner} $
defines a global Hamiltonian flow on $L^2$.

\subsection{Well-posedness for the \texorpdfstring{$N$}{N}th KdV equation in \texorpdfstring{$H^{N-2}$}{HN-1}.}

Let $N \in \N$. Let  $ \T_{N,app}^\Gardner(\tau,w,\tau_0)$ be the approximate Hamiltonian of Definition \ref{def:approx}. We denote  $  \T^\Gardner_{N, app} ( \infty,w,\tau_0):= H^\Gardner_N(w,\tau_0)$. 

\begin{theorem}\label{thm:high}  
Suppose that $ 2\le N \le M+1$. The flow map 
\[   \R \times (0,\infty) \times H^M \ni (t,\tau, w_0) \to \exp\Big( t \partial \frac{\delta}{\delta w} \T^\Gardner_{N, app}\Big) (w_0) \in H^M 
\]
has a continuous extension to 
\[ \R \times (1,\infty] \times H^M \ni (t, \tau, w_0) \to  w(t)=\exp\Big( t \partial \frac{\delta}{\delta w} \T^\Gardner_{N, app}\Big) (w_0) \in H^M.  \]
If $ \tau= \infty$ then $w(t) $ is a weak solution to the $N$th Gardner equation. The flows defined by the Hamiltonians $\T^\Gardner_{-1}(i\tau) $ and $H^\Gardner_n $ with $n \le N $ commute. 
\end{theorem}

The case $N=1$ is covered in Theorem \ref{thm:KdV}, with a slightly different approximate Hamiltonian.  Theorem \ref{thm:equivalenceweaksolutionsuw} says that $w$ is a weak solution as in the theorem for $M\ge N \ge 2$ iff the Miura map $u= w_x+2\tau_0 w +w^2$ is a solution to the $N$th KdV  in $H^{N-1}(\mathbb{X})$.
Both flows have a continuous extension to $M=N-1$.
As a consequence of the theorem we obtain well-posedness of the $N$th KdV equation with initial data in $H^{N-2}$ for $N \geq 2$, in the sense of a continuous extension.  
The first Gardner equation contains the term $w^3$ and we cannot define weak solutions in $C(\R; L^2(\mathbb{X}))$.
The proof relies on a modifcation of  the strategy we used for the first KdV equation in Subsection \ref{subsec:KdVwell}.

\begin{proof} The argument does not distinguish between the equation on the line or on the circle. Let $ N \ge 2$.

\noindent
{\bf  Step 1: Well-posedness of the approximate flow.} The approximate flow is the flow of the Hamiltonian vector field of the Hamiltonian $ \T^\Gardner_{app,N}(i\tau)$
(defined in Definition  \ref{def:approx}) constructed in Theorem \ref{thm:flow}. By Proposition \ref{prop:commuting} it preserve all Hamiltonians $ \T^\Gardner_{-1}(i \tilde \tau,w,\tau_0)$ and, by taking limits, it preserves  the Hamiltonians $H^\Gardner_n(w,\tau_0) $ for $ n \le N-1$ provided the initial data is in $H^{N-1}(\X)$. The Gardner Hamiltonians $H^\Gardner_n$, with $n=0$ and $ n = N-1$ control the $H^{N-1}$ norm of $w$. 
Thus $\T^\Gardner_{app,N}(i\tau, w, \tau_0)$ defines a global flow on $H^{N-1}(\X)$.

\noindent
{\bf Step 2: Equicontinuity.} Equicontinuity of a set $Q \subset H^{N-1}$  can be characterized as (see Proposition \ref{prop:diffham})
\begin{equation} \label{eq:equiconw}  Q \text{ is equicontinuous in } H^{N-1}(\X) \Longleftrightarrow  \lim_{\tau \to \infty} 
\sup_{w \in Q}   |\T_{N-1}^{\Gardner} (i\tau, w)| = 0. \end{equation} 
Since $ \T_{N}^{\Gardner}( i\tau,\cdot ) $ is preserved under the flow of the approximate Hamiltonian $ \T^{\Gardner}_{app,N}(i\tau_1,\dot) $ flow also equicontinuity is preserved along the flow. 
 
\noindent{\bf Step 3: Convergence of the  difference vector field to $0$ in $H^{-N-2}$ uniformly on equicontinuous sets.}  
Let $Q\subset H^{N-1}$ be an equicontinuous bounded set. By Lemma \ref{lem:diffN1} (and the remark following it)  this is equivalent to 
\begin{equation}\label{eq:uniformconvw}    \lim_{\tau_1 \to \infty}  \sup_{w \in Q}  \Big\Vert \big\{ w , \T^{\Gardner}_{app,N}(i\tau_1,\cdot)  \big\}^\Gardner(w)  \Big\Vert_{H^{-N-2}} = 0. 
\end{equation}

\noindent{\bf Step 4: The approximate  flow.}
Let $ \tau_1,\tau_2 > \tau_0 $ and consider the difference of the approximate  flows
\[  w_t =   \partial \frac{\delta}{\delta w } \Big[  -  \T^{\Gardner}_{app,N}( i\tau_2, \cdot ) + \T^{\Gardner}_{app,N} (i\tau_1,\cdot )\Big]  \]
with initial data $w(0)=w_0\in H^{N-1}(\R)$. We denote 
\[ 
w ( \tau_1,\tau_2,t) = \exp\Big( t  \partial \frac{\delta}{\delta w}  \Big( \T^{\Gardner}_{app,N}(i\tau_1,\cdot) -  \T^{\Gardner}_{app,N} (i\tau_2,\cdot) \Big)  \Big)w_0. 
\]

Let $ Q \subset H^{N-1}(\X)$ be an equicontinuous subset. Then also 
 \begin{equation} \label{eq:Qdefw}  \tilde Q= \Big\{ \exp\Big(t \partial \frac{\delta}{\delta w}\Big(\T^{\Gardner}_{app,N}(i\tau_1,\cdot) \Big)\Big) w_0 :    \tau_1 > \tau_0, t \in \R , w_0 \in Q \Big\} \end{equation} 
is equicontinuous in $H^{N-1}$ by Step 2.  Moreover, from Definition 
\ref{def:approx} and Lemma~\ref{lem:approximate}
\[ \begin{split}  w_t \, & = \{ w , \T^{\Gardner}_{app,N}(i\tau_1,\cdot,\tau_0)  - \T^{\Gardner}_{app,N}(i\tau_2,\cdot,\tau_0) \} 
\\ & = \sum_{j=0}^N a_j \Big( \{ w, \T^{\Gardner}_N(i2^{j}\tau_1,\cdot,\tau_0) \} - \{ w, \T^{\Gardner}_N(i2^j \tau_2,\cdot,\tau_0) \} \Big)
\end{split} 
\]
and (suppressing some arguments of functions), by Step 3,  
\[\begin{split} \hspace{1cm} & \hspace{-1cm}   \Vert w(t) - w_0 \Vert_{H^{-N-2}} \\ &   \le \sum_{j=0}^N |a_j|\int_0^t   \Vert \{ w , \T^{\Gardner}_N(i 2^j\tau_2,\cdot)\} - \{w, \T^{\Gardner}_1(i2^j \tau_1,\cdot) \} \Vert_{H^{-N-2}} ds 
\\ & \le \sum_{j=0}^N |a_j|  |t| \sup_{w \in Q} \Big( \Vert \{ w, \T^{\Gardner}_N(i2^ j\tau_2,\cdot,\tau_0) \} \Vert_{H^{-N-2}} 
+  \Vert \{ w, \T^{\Gardner}_N(i 2^j \tau_1,\cdot) \} \Vert_{H^{-N-2}}\Big) 
\\ & \to 0 \qquad \text{ as }\tau_1,\tau_2 \to \infty 
 \end{split} 
\]
as $ \tau_1,\tau_2 \to \infty$  by Step 3, uniformly on compact time intervals and for $w_0 \in Q$.

\noindent{\bf Step 5: Convergence of the approximate flow.}
We want to prove that 
\[ 
e^{t \partial \frac{\delta}{\delta w}  \T^{\Gardner}_{app,N} (i\tau,\cdot,\tau_0)}u_0
\] is  Cauchy in $H^{N-1}(\X)$ as $\tau\to \infty $ uniformly for $t$ in compact intervals. By commutativity of the flows  (and a suggestive abuse of notation) 
\[\begin{split}  \hspace{1cm} & \hspace{-1cm}   \exp\Big(  t  \partial \frac{\delta}{\delta w} \T^{\Gardner}_{app,N}( i\tau_2,\cdot) \Big) w_0 
- \exp\Big(  t \partial \frac{\delta}{\delta w} \T^{\Gardner}_{app,N}( i\tau_1,\cdot) \Big) w_0
\\ & =
    \Big\{\exp\Big[  t \partial \frac{\delta}{\delta w} \Big(\T^{\Gardner}_{app,N}( i\tau_2,\cdot) - \T^\Gardner_{app,N}(i\tau_1,\cdot ) 
    \Big)\Big] -\id
  \Big\}\\ & \quad  \times  \exp\Big(t  \partial \frac{\delta}{\delta w} \T^{\Gardner}_{app,N}(i\tau_1,.) \Big) w_0.
\end{split} 
\]
The set 
\[  Q = \Big\{ \exp\Big(t  \partial \frac{\delta}{\delta w} \T^{\Gardner}_{app,N}(i\tau,.) \Big) w_0  :  \tau > \tau_0, t \in \R \Big\}\subset H^{N-1}(\X) \]
is equicontinuous. Thus the left hand side of the equation above converges to $0$ in $H^{-N-2}(\X)$, uniformly for $t$ in compact time intervals as $ \tau_1,\tau_2 \to \infty$ by Step 4.
Thus 
\[ 
e^{t (2\tau)^2\partial \frac{\delta}{\delta w}  \T^{\Gardner}_{app,N} (i\tau,\cdot,\tau_0)}u_0
\] 
is Cauchy in $H^{-N-2}(\X)$ as $\tau \to \infty$ uniformly for $t$ in compact time intervals. 
Let $ w$ be the limit. Since 
\[
\sup_{\tau> \tau_0,t} \Vert e^{t (2\tau)^2\partial \frac{\delta}{\delta w}  \T^{\Gardner}_{app,N} (i\tau,\cdot,\tau_0)}w_0 \Vert_{H^{N-1}(\X)} < \infty 
\] 
and since the functions are equicontinuous
the convergence in $H^{-N-2} $ implies convergence in $H^{N-1}(\X)$. This completes the proof. 
\end{proof}

\begin{proposition} \label{prop:uniqueweak}
Weak solutions in $C( \R; H^N)$ to the $N$th Gardner equation are unique.
\end{proposition}

\begin{proof}
Let $w$ be a weak solution to the $N$th Gardner equation in $C(\R;H^N)$ and let 
\begin{equation}  v(\tau, t) =   \exp\big( - t \partial \frac{\delta}{\delta w}  T^{\Gardner}_{app,N}(i\tau,. ,\tau_0)\big)   w(t).  \label{eq:diffapp}  \end{equation}
The orbit $ \{ v(\tau,t): \tau > \tau_0, t\in \R\} $ is again an equicontinuous set in $H^{N}$.
In a distributional sense (see Proposition \ref{lem:com} for a more involved statement)
\begin{equation} \label{eq:commutativity}   \partial_t v(\tau, t) = -\sum_{j=0}^N a_j \partial \frac{\delta}{\delta w} \T^{\Gardner}_N(i2^j \tau,.,\tau_0)(v).      \end{equation} 
The right hand side of \eqref{eq:commutativity} converges uniformly to zero in $H^{-N-1}$ on equicontinuous sets as $\tau \to \infty$. Thus, again using equicontinuity 
\[   v(\tau, t) \to w_0 \]
uniformly in $H^{N}$ on compact time  intervals and $w_0$ in bounded equicontinuous sets in $H^{N}$.
As above this implies 
\[  \lim_{\tau\to \infty} v(\tau,t) = w(t)\]
in $H^N$ uniformly on compact time intervals.  
To complete the proof we apply $\exp\Big( t\partial \frac{\delta}{\delta w} \T^\Gardner_{app,N}(i\tau,. ,\tau_0)\Big) $ to both sides of \eqref{eq:diffapp}.
\end{proof} 
 
Here we need one derivative more for uniqueness than for the extension of the approximate flow. It may be possible to remove this gap. However our main focus is on uniqueness of weak solutions with $L^2$ data on $\R$. We will use Theorem  \ref{thm:high} as a step in the consideration of solutions with initial data in $L^2$ on $\R$.
  
 \subsection{Precompactness of weak solutions to the Gardner equations in \texorpdfstring{$L^2$}{L2}.} 
 In previous subsection, we studied the problem for general geometry $\mathbb{X}$ at high regularity. However, on the line there is a gain of regularity by Kato smoothing. This gain of regularity  will allow us to prove existence and uniqueness of weak solutions to the Gardner hierarchy with $L^2$ data. 

Theorem \ref{th:main} follows by applying the modified Miura map, more precisely Theorem \ref{thm:equivalenceweaksolutionsuw}. The results for the Gardner hierarchy are cleaner and slightly stronger than those for the KdV hierarchy: We will obtain uniqueness of weak solutions.

Let $I$  be an interval and let
\[  \Vert f \Vert_{L^2_u(I \times \R)}
= \sup \{   \Vert f \Vert_{L^2( J \times  K) }: J \subset I, |J|=|K|=1\}.\] 

\begin{definition}\label{def:Katosmoothing} We say that $w$  is an element of the Kato smoothing space $X_N$ if
\[  w \in L^\infty(\R, L^2), \quad w^{(N)} \in L^2_u(I\times \R) \quad \text{ for all compact intervalls }I  \]
and for all $t_0 \in \R$
\begin{equation}\label{eq:limit}   \lim_{x_0 \to \pm \infty}  \Vert  w^{(N)} \Vert_{L^2( (t_0,t_0+1)  \times (x_0-1,x_0+1))} = 0. \end{equation} 
\end{definition}

 \begin{theorem} \label{thm:Gardnerprecompact} 
Suppose that $ w \in X_N$ is a weak solution to the $N$th $\tau$ Gardner equation
\[  w_t = \partial \frac{\delta}{\delta w} H^\Gardner_N(w,\tau_0).\]
Then $w \in C([0,\infty); L^2(\R))$. We  define  the initial trace of $ w$ by $w_0:=w(0) \in L^2(\R)$. Then $ \Vert w(t) \Vert_{L^2} = \Vert w_0 \Vert_{L^2}$ and
the Kato smoothing estimate 
\begin{equation} \label{eq:Kato} \begin{split}\hspace{1cm}&\hspace{-1cm}  \sup_t \Vert ( 1+\tanh( x-\kappa \tau^{2N} t)^{1/2})    w(t) \Vert_{L^2}  +   \Vert \sech(   x- \kappa \tau^{2N}  t) w^{(N)} \Vert_{L^2(\R\times \R)} \\ & \le  c( \tau^{-1/2} \Vert w_0 \Vert_{L^2})  \Vert (1+\tanh(x) )^{1/2} w_0 \Vert_{L^2} 
\end{split} 
\end{equation}  
holds for all $ \kappa \ge \kappa_0$ where $\kappa_0$ depends only on $N$ and  $\tau^{-1/2} \Vert w_0 \Vert_{L^2}$. 
If $ Q \subset L^2$ is precompact, if  $\tilde Q$  is a set of weak solutions so that $w(0)\in Q$, $R>1$ and if $ I\subset \R$ is a compact time interval then the sets 
\[  \{  w(t) ; w\in \tilde Q, t \in I \} \subset L^2(\R) \] 
and 
\[    \{  \sech( R^{-1} x) D^N w: w \in \tilde Q\}  \in L^2(I\times \R)\} \]
are precompact. 
\end{theorem} 

It is not obvious that this regularity suffices to define weak solutions. This point is elaborated in Section \ref{sec:weak}. This theorem and its proof are purely about PDEs. Nevertheless integrability is important for the statement about precomapctness. 

 \begin{proof}
 Lemma \ref{lem:energyflux} states that for weak solution in $X_N$ 
 the energy flux identity 
 \[ \partial_t   w^2  = \partial_x   \flux_N (x) \]
holds in a distributional sense. The structure of $ \flux_N $ is given in Lemma \ref{lem:energyflux}. 
  It will be convenient to split the Hamiltonian into a quadratic and a higher order part, corresponding to splitting the $N$th Gardner equation into a linear and a nonlinear part, 
\begin{equation}\label{eq:HNL}  H^\Gardner_N = \frac12 \Vert w^{(N)} \Vert_{L^2}^2 + H^{NL}_N,  \end{equation} 
\begin{equation}  w_t = (-1)^N \partial_x w^{(2N)} + \partial_x \frac{\delta}{\delta w} H^{NL}_N \end{equation} 
Similarly 
\begin{equation} \label{eq:fluxNL}  \flux_N =  \flux_N^L + \flux^{NL}_N  \end{equation} 
with  linear part 
\[ \flux_N^L = (2N+1 )  |w^{(N)}|^2 + \sum_{j=1}^N f_{N,j} \partial^{2j} |w^{(N-j)}|^2       \]
for some unimportant combinatorical constants $f_{N,j}$.  The nonlinear part is described in more detail in Lemma \ref{lem:energyflux}. 

We multiply the energy flux identity by Schwartz function  
\[(1+ \tanh(\frac{x}{R}-1 ))(   1-\tanh(\frac{x}R+1)) \]
and integrate over $(0,T) \times \R $ to obtain 
\[ 
\begin{split} \hspace{1cm} & \hspace{-1cm}  
\frac12\int (1+ \tanh(\frac{x}R-1))(   1-\tanh(\frac{x}{R}+1)) |w|^2 dx \Big|_0^T 
\\ & = R^{-1} \int_0^T \int_{\R} ((1+\tanh(\frac{x}R-))\sech^2(\frac{x}R+1)\\ & \qquad - \sech^2(\frac{x}R-1)(1-\tanh(\frac{x}R+1)) ) \flux_N(x) dx   \, dt.
\end{split} 
\]
The quadratic part on the right hand side is bounded by (up to a constant)  
\[  
 \int_0^T \int_{\R}  R^{-1} (\sech^2(\frac{x}R-1) + \sech^2(\frac{x}R+1)) (|w^{(N)}|^2 
+ R^{-N} |w|^2 dx dt 
\]
which converges to $0$ with $R \to \infty $ since $ w$ lies in the local smoothing space  $ X_N$. The nonlinear part is estimated in Lemma \ref{lem:L2cont}. 
It also converges to $0$ as $R \to \infty$. Thus $ t \to \Vert w(t) \Vert_{L^2} $ is constant. 
Using the equation and the fact that $ w \in X_N $ we see that 
\[   w \in C( \R; H^{-N-1}). \]
Together this implies $ w \in C(\R;L^2(\R))$.

We multiply  the energy-flux identity by  $1+ \tanh(( x-\kappa \tau^{2N})/R) $ 
(which is allowed by the previous calculation) 
and integrate to obtain the following identity. 

\begin{proposition} \label{prop:kato2} 
Let $w\in X_N$ be a weak solution to the $N$th Gardner equation. Then $w \in C(\R; L^2)$, 
$ \Vert w(t) \Vert_{L^2} = \Vert w_0 \Vert_{L^2}$ and (we omit the argument $(x-\kappa \tau^{2N}t)/R$
of $\sech$)
\begin{equation} \label{eq:Katosmoothing} \begin{split} \hspace{.3cm} & \hspace{-.3cm} \int  (1+ \tanh\big( (x-  \kappa  \tau^{2N} t)/R\big) w^2 dx \Big|_{t=0}^{t=T} 
\\ & + \frac1R \int_0^T \int  \sech^2\Big[ (2N+1)  (w^{(N)})^2 
+ \sum_{j=1}^N  f_{N,j}  |w^{(N-j)}|^2\cosh^2\partial^{2j}\sech^2  + \kappa \tau^{2N}   w^2\Big] dx dt \\ & \qquad = \frac1R \int_0^T \int \sech^2 \flux^{NL}_N dx  dt. 
\end{split} 
\end{equation}
\end{proposition} 
The second line is the linear Kato smoothing term. It is equivalent to (we choose $ \kappa_0$ sufficiently large)
\[ \frac1R \int_0^T \Vert \sech w \Vert_{H^N_\tau}^2 + \kappa \tau^{2N} \Vert \sech w \Vert_{L^2}^2 dt \]
if  $R$ is large, independent of $ \tau$, since the middle terms carry a factor $R^{-2j} $ from the differentiation. We bound  the nonlinear term (see Lemma \ref{lem:L2cont} and its proof) by, ignoring the shift by $-\kappa \tau^{2N}t$,
\begin{equation}\label{eq:nonlinearw} \begin{split}  \left| \int \sech^2(x/R) \flux^{NL}_N dx \right| \, & \le   c( \tau, \Vert w \Vert_{L^2}) \Big( \tau^N \Vert \sech(x/R) w  \Vert_{L^2}\Big)^{\frac2{2N-1}}  \Vert \sech(x/R) w \Vert^{2\frac{2N-2}{2N-1}}_{H^N_\tau} 
\\ & \le \varepsilon \Vert \sech(x/R) w \Vert_{H^N_\tau}^2 + c(\varepsilon)  \tau^{2N} \Vert w \Vert^2_{L^2}.
\end{split} 
\end{equation} 
We choose $ \varepsilon $ small and then $\kappa_0$ large so that we can subsume the second line under right hand side of the first line. Integration with respect to time implies the Kato smoothing estimate \eqref{eq:Kato} for weak solutions.

Precompactness  is contained in Theorem \ref{thm:gardnerprecompactness}. It relies on
\begin{enumerate} 
\item Kato smoothing \eqref{eq:Kato} and its translated versions which gives tightness to the right,
\item  the modified Miura map and the estimates of Proposition  \ref{prop:equivalencewu} applied twice to obtain equicontinuity of the orbit and high frequency Kato smoothing,
\item a backward Kato smoothing with 'bad' terms controlled by the high frequency estimates of the previous step. 
\end{enumerate} 
This is quite intricate and the object of Subsection \ref{subsec:Katosmoothing}, Theorem \ref{thm:gardnerprecompactness}, together with more details on the arguments in this subsection. 
\end{proof}

\subsection{Continuous extension  of the approximate Gardner flow in \texorpdfstring{$L^2$}{L2}}
The continuous extension of the approximate flows  uses a variant of the second commuting vector field method of Bringmann, Killip and Visan 
\cite{bringmann2019global}. 
Recall that for $w \in H^N$
\[ \lim_{\tau \to \infty} \T^\Gardner_{app,N}( i\tau, w, \tau_0) = H^\Gardner_N(w,\tau_0). \]
By a small abuse we write $ H^\Gardner_N(w, \tau_0)= \T^\Gardner_{app,N}(i\infty, w, \tau_0)$.

\begin{theorem}\label{thm:Gardnercontinuous}  
 The flow map 
\[   \R \times (0,\infty) \times L^2 \ni (t,\tau, w_0) \to \exp\Big( t \partial \frac{\delta}{\delta w} T^\Gardner_{N, app}(i\tau, . , \tau_0 \Big) (w_0) \in L^2 
\]
has a continuous extension to 
\[ \R \times (1,\infty] \times L^2 \ni (t, \tau, w_0) \to  w(t)=\exp\Big( t \partial \frac{\delta}{\delta w} T^\Gardner_{N, app}(i\tau, ., \tau_0)\Big) (w_0) \in L^2.  \]
If $ \tau= \infty$ then $w(t) $ is a weak solution to the $N$th Gardner equation. The flows defined by the Hamiltonians $\T^\Gardner_{-1}(i\tau) $ and $H^\Gardner_n $  commute. 
\end{theorem} 

By Theorem \ref{thm:high}  there exists a  weak solution for $ w_0 \in H^{N-1}$. By Proposition \ref{prop:uniqueweak} it is unique if $ w_0 \in H^N$ (in $H^1$ if $N=1$).
The proof has shown more: There is a continuous map 
\[ \Xi_N :    (\tau_0, \infty]\times H^N \times \R  \ni (\tau,w_0,t) 
\to  \exp\big( t  \partial \frac{\delta}{\delta w} \T^\Gardner_{app,N}(i\tau,.,\tau_0)\big) (w_0) 
\in H^{N}
\]
where the Hamiltonian for $ \tau=0$ is the $N$th Gardner Hamiltonian.

It is an easy consequence of compactness that the set of weak solutions in $X_N$ is closed in the following sense: Let $w^n$ be a sequence of weak solutions which converges in $C(I, L^2)$ for every bounded interval $I$. Theorem \ref{thm:Gardnerprecompact}  provides uniform bounds and  precompactness. There exists a subsequence so that $w^{n_j} $ weak* converges in $L^\infty(L^2)$ 
and $ \partial^N w^{n_j}$ converges in $L^2$ on  every compact subset of $ \R \times \R$. Thus the limit is  a weak solution in $X_N$.  

We prove more than that: The uniform convergence (and hence uniform continuity) of the flow 
$\T^\Gardner_{app,N}$ to weak solutions to the $N$th Gardner flow on precompact sets of initial data. A key element is a Kato smoothing estimate for weak solutions $w\in X_N$ to the difference  equation 
\begin{equation}\label{eq:diffeq}
\begin{split} 
w_t\, &  = \partial\left(  \frac{\delta}{\delta w} \big[ H^\Gardner_N (w,\tau_0)- \T^\Gardner_{app,N}(i\tau, w, \tau_0) \big] \right) 
\\ & = -\partial\sum_{j=0}^N a_j \frac{\delta}{\delta w}\T^\Gardner_N(i\lambda^j \tau, w, \tau_0).  
\end{split}
\end{equation} 
 
\begin{proof}[Proof  of Theorem \ref{thm:Gardnercontinuous}:]
We begin with a Kato smoothing estimate for the difference flow. The Kato smoothing is very weak, but it will suffice, together with an estimate of the nonlinearity (Lemma  \ref{lem:diffflow}) to prove convergence to a weak solution at $ \tau \to \infty$.  
\begin{proposition} \label{prop:Kato:diff} 
Suppose that $ w$ in the Kato smoothing space $ X_N$ of Definition \ref{def:Katosmoothing} is a weak solution to \eqref{eq:diffeq}
Then there exists $ \kappa $ depending only on $N$ so that 
\begin{equation} \label{eq:Katodiff} 
 \Vert \sech(x-\kappa  \tau_0^{2N} t) w^{(N+1)} \Vert_{L^2(\R_+; H^{-1}_{2\tau})} +\kappa \tau_0^{2N} \Vert \sech(x-\kappa \tau_0^{2N} t) w \Vert_{L^2(\R_+\times \R)}
  \le c \Vert w(0) \Vert_{L^2} 
\end{equation} 
 and, for $ 0\le t \le 1 $ 
 \begin{equation} \label{eq:appconvergence}  \Vert \sech(x)  ( w(t) - w(0)) \Vert_{H^{-N-3}} \le c \tau^{-\frac{3}{N+1}}   \Vert w(0)  \Vert_{L^2}. \end{equation} 
\end{proposition} 
\begin{proof} 
The proposition relies on natural but  fairly sharp estimates. We decompose
the difference Hamiltonian into a quadratic and a higher order part, 
\[  \T^\Gardner_N(\tau, w, \tau_0)= \frac12  \Vert w^{(N+1)} \Vert_{H^{-1}_{2\tau}} 
+ \T^\Gardner_{N,NL}(\tau, w , \tau_0).
\]
There are two parts: The local smoothing estimate, and the estimate of the vector field in terms of the local smoothing bounds. The vector field and the nonlinear Kato smoothing term are estimated in the following lemma. 
\begin{lemma}\label{lem:diffflow}
The linear  estimate 
\begin{equation} \label{eq:diffflow0}
\Big\Vert \sech^2(x/R) \partial \frac{\delta}{\delta w} \Vert w^{(N+1)} \Vert_{H^{-1}_{2\tau}}   
\Big\Vert_{H^{-N-1}}  \le c \tau^{-2} (\Vert \sech^2 w^{(N)} \Vert_{L^2} +  \Vert \sech^2 w \Vert_{L^2} )  
\end{equation} 
and the estimate of the nonlinear part of the vector field
\begin{equation} 
\label{eq:diffflow}
\begin{split}\hspace{1cm}& \hspace{-1cm} 
\Big\Vert \sech^2(x/R) \partial \frac{\delta}{\delta w} \T^{\Gardner}_{N,NL}(i\tau, w, \tau_0)  \Big\Vert_{H^{-N-3}}
\\  & \le c \tau^{-\frac3{2(N+1)}} ( 1+ \Vert w \Vert_{L^2}) (\Vert \sech(x/R) w^{(N+1)} \Vert^2_{H^{-1}_{2\tau}} + \Vert \sech(x/R) w \Vert^2_{L^2})
\end{split}
\end{equation} 
and the estimate for the Kato smoothing term of the nonlinear part of the vectorfield 
\begin{equation} 
\label{eq:diffflowKato}
\begin{split}\hspace{1cm}& \hspace{-1cm}  
\left|\int \tanh(x/R) w(x) \partial\frac{\delta}{\delta w} \T^\Gardner_{N,NL}(i\tau,w,\tau_0)   \Big)\,dx \right|
\\ & \le  c \tau^{-\frac1{N+1}} ( 1+ \Vert w \Vert^{2N+2}_{L^2}) \Big(\Vert \sech(x/R) w^{(N+1)} \Vert^2_{H^{-1}_{2\tau}} + \Vert \sech(x/R) w \Vert^2_{L^2}\Big).
\end{split} 
\end{equation} 
hold for $R \ge 1$ and $ \tau > \tau_0$ with constants depending on $ \tau$,$ \tau_0$ and $ \tau_0^{-1/2} \Vert w \Vert_{L^2}$.
\end{lemma} 
We postpone  the proof of the lemma andmturn to the  Kato smoothing estimate \eqref{eq:Katodiff} for solution $w \in X_N$ to \eqref{eq:diffeq} with  the Hamiltonian 
\[  \T^\Gardner_{app,N}(\tau, w, \tau_0)  - H^\Gardner_N(w,\tau_0) = \sum_{j=0}^N a_j \T^\Gardner_N( i2^j \tau, w, \tau_0).  \]

The linear part is given by Lemma \ref{lem:approximate} with the Fourier multiplier $m(D)$,
\[    \int \tanh((x- \kappa \tau_0^{2N} t)  /R)  w  \partial_x M(D) w dx  
= - R^{-1}  \int \sech((x-\kappa \tau_0^{2N})/R)^2    |m(D) w|^2   \Vert_{L^2}  + \rho  
\]
where 
\[ m(\xi) =   \sqrt{ \frac{d}{d\xi} (\xi M(\xi)) } \sim \frac{ |\xi|^{N+1}}{ \xi^2 +\tau^2} \]
and the remainder $\rho$ satisfies through commutator estimates
\[  |\rho| \le \sum_{n=2}^{N-1} R^{-n} \Vert \sech( (x-\kappa \tau_0^{2N})/R )   (\tau^2 + D)^{-1}  w^{(N+1-n)} 
\Vert_{L^2} +  R^{-N}  \Vert  w \Vert_{L^2}. 
\]
By Lemma \ref{lem:approximate}
\[  m(\xi) \sim \frac{ |\xi|^{N+1}}{\sqrt{\xi^2 + \tau^2}}   \]
hence, commuting $ 1/\sqrt{\tau^2+ D^2}$ with $\sech$,
\[  \int \sech((x-\kappa \tau_0^{2N})/R)^2    |m(D) w|^2   \Vert_{L^2}
\sim     \Vert \sech( (x-\kappa \tau_0^{2N})/R) w^{(N+1)} \Vert_{H^{-1}_{2\tau}}.   
\]
We interpolate the terms in the remainder term and obtain for a small constant $\mu$ (checking the constants, but there would be no harm in introducing a small factor on the second term on the left hand side)
\[  \begin{split} \hspace{2cm}& \hspace{-2cm}    \int \tanh((x- \kappa \tau_0^{2N} t)  /R)  w  \partial_x M(D) w dx 
+ \mu R^{-1} \Vert \sech((x-\kappa \tau_0^{2N} t) /R) w^{(N+1)} \Vert^2_{H^{-1}_{2\tau}}  
\\ & \le c  \Vert  w \Vert_{L^2} 
\end{split}
\] 
When the time derivative falls on $ \tanh( x- \kappa \tau_0^{2N} t )$ we obtain a good term 
\[ -   \frac12 \kappa \tau_0^{2N}  \Vert \sech^2 ( x-\kappa \tau_0^{2N}t) w  \Vert_{L^2}^2.\]
We estimate the nonlinear terms by \eqref{eq:diffflowKato} and obtain 
\[ 
\begin{split} \hspace{1cm} & \hspace{-1cm} 
\frac{d}{dt} \frac12  \int \tanh( x- \kappa \tau_0^{2N} t) w^2 dx 
+ \kappa R^{-1} \Vert \sech^2(x-\kappa \tau_0^2 t) w^{(N+1)}\Vert_{H^{-1}_{2\tau_0} } 
\\ & \qquad + R^{-1} \kappa \tau_0^{2N} \Vert \sech^2(x-\kappa \tau_0^2t) w \Vert_{L^2}^2 
\\ & \le R^{-1} c  \tau^{-\frac1{N-1}} ( \Vert \sech(x- \kappa \tau_0^{2N} t) w^{(N+1)} \Vert_{H^{-1}_{2\tau}}. 
\end{split} 
\]
For $ \tau $ large the right hand side can be controlled by the left hand side and we obtain
the Kato smoothing estimate\eqref{eq:Katodiff} for solutions to \eqref{eq:diffeq}. 
This estimate is uniform for $w_0$ bounded in $L^2$.
By Lemma 2.23, \eqref{eq:diffflow0}, \eqref{eq:diffflow} we can bound the vector field in $H^{-N-2}$
by the Kato smoothing terms \eqref{eq:diffflowKato} and we obtain \eqref{eq:appconvergence}. 

We turn to estimate \eqref{eq:appconvergence}. By the first two estimates of Lemma \ref{lem:diffflow} 
\[ 
\begin{split} 
\Vert \sech^2 ( x/R) w(t)  \Vert_{H^{-N-3}} 
\le \, & c \tau^{-2} \int_0^t \Vert \sech(x/R) w^{(N+1)} \Vert_{H^{-1}_{2\tau}} 
\\ & \hspace{-1cm} + \tau^{-\frac{2}{3(N+1)}} \big(  \Vert \sech(x/R) w^{(N+1)} \Vert^2_{H^{-1}_{2\tau}} + \Vert \sech(x/R) w \Vert_{L^2}^2  \big). 
\end{split} 
\]
We choose $ \tau $ large so that the right hand side is controlled by \eqref{eq:Katodiff}.
This implies \eqref{eq:appconvergence}.

By equicontinuity this implies $w(t) \to w(0) $ in $L^2$, uniformly for initial data in precompact  sets. 
\end{proof}

\begin{proof} [Proof of Lemma \ref{lem:diffflow}]
For the first estimate we calculate 
\[  \sech^2(x/R) \partial  \frac{\delta}{\delta w} \Vert  w^{(N+1)} \Vert_{H^{-1}_{2\tau}} 
=   2 (-1)^{N}       \sech^2 (x/R) \partial \Big[ ((1- (2\tau)^{-2}\partial^2 )^{-1}  \partial^{2N} )       w \Big] 
\]
hence 
\[   \Vert \sech^2(x/R) \partial  \frac{\delta}{\delta w} \Vert  w^{(N+1)} \Vert_{H^{-1}_{2\tau}}  \Vert_{H^{-N-3}} \le c \tau^{-2} \Vert \sech^2(x/R)  w \Vert_{H^N}  \]
and by commuting the derivatives and interpolation. 
The estimates \eqref{eq:diffflow} and \eqref{eq:diffflowKato} are proven after Proposition \ref{prop:weighted}.
Estimate \eqref{eq:diffflowKato} bounds the nonlinear part of the Kato smoothing estimate  for  the Hamiltonian vector fields. 
\end{proof} 

Let $w_0 \in L^2$, $ w^n_0 \in H^N$ so that $w^n_0 \to w_0$ in $L^2$. Let $w^n$  resp. $w^n(\tau)$ be the corresponding solutions to the $n$th Gardner equation resp. the $n$th approximate equation. The flow of the Hamiltonian $ \T^\Gardner_{app,N}(\tau, w , \tau_0)$ is continuous in $L^2$. Thus $ w^n(\tau,t) \to w(\tau,t)  $ in $L^2$, uniformly on compact time intervals.  

By Theorem \ref{thm:high}  
\[    v^n( \tau) = \exp\Big( -t  \partial \frac{\delta}{\delta w} \T^\Gardner_{app,N}(i\tau, . , \tau_0) \Big) w^n(t) \]
satisfies  $v^n\in X_N$ and, since $ H^\Gardner_N(w,\tau_0)  $ and $ \T^\Gardner_{app,N}(\tau, w , \tau_0) $ (and hence also $ \T^\Gardner_{app,N}$ Poisson commute 
on $H^N$
\[  v^n_t(\tau)  = \partial \frac{\delta}{\delta w} (H^\Gardner_{N}(.,\tau_0) - \T^\Gardner_{app,N}(i\tau, . ,\tau_0)  (v^n(\tau)) , \qquad v^n(0,\tau) = w^n_0.  \]
By \eqref{eq:Katosmoothing} and Theorem \ref{thm:flow} the solutions satisfy 
\[     \sup \{   \Vert (v^n(\tau))^{(N)} \Vert_{L^2(I\times K)}: |I|, |K| \le 1 \} \le c(\tau) \Vert w_0 \Vert_{L^2}.   \]
 By Proposition  \ref{prop:Kato:diff}, if $ 0\le t \le 1$, since $w^n_0 \to w_0$, which we assume to be nonzero, 
\[   \Vert \sech( v(t,\tau) - w^n_0) \Vert_{L^2(\R)} \le c (1+t) \tau^{-\frac{3}{2(N+1)}} \Vert w_0 \Vert_{L^2}, \]
at least for large $n$ - by taking a subsequence we may assume that this holds for all $n$.

The $L^2$ norm is conserved, hence  
\[  \Vert v^n(\tau,t) \Vert_{L^2} = \Vert w^n(t) \Vert_{L^2} = \Vert w^n_0 \Vert_{L^2} \]
and we lift the convergence to $L^2$, uniformly for precompact sets of initial data $w_0^n$.

The set $\{ w^n_0\} \subset  L^2$ is precompact since $ w^n_0 \to w_0$ in $L^2$ and hence the convergence in $L^2$ 
\begin{equation}  \lim_{n\to \infty} v^n(t,\tau) \to w_0  \end{equation} 
is uniform for $t \le 1 $. Now 
\[
\begin{split} \hspace{1cm} & \hspace{-1cm} 
 w^n(t) -
\exp\Big( t  \partial \frac{\delta}{\delta w} \T^{\Gardner}_{app,N}(i\tau)\Big)w_0^n  \\ & =
\Big( \id 
 - \exp\Big( t \partial \frac{\delta}{\delta w} ( \T^\Gardner_{app,N}-H^\Gardner_N) \Big)  w^n(t) 
\end{split} 
\]
By Theorem \ref{thm:Gardnerprecompact}  the set 
\[  Q = \{ w^n(t): 0\le t \le 1, n \in \N\} \subset L^2 \] 
is precompact and trivially $ w^n \in  C([0,1];H^N) \subset X_N$.
Thus 
\[ \lim_{\tau\to \infty} \exp\Big( t  \partial \frac{\delta}{\delta w} \T^{\Gardner}_{app,N}(i\tau,.,\tau_0)\Big)w_0^n  
=  \exp\Big(  t \partial \frac{\delta}{\delta w} H^\Gardner_N (w, \tau_0)\Big) w^n_0 
\]
uniformly for $ 0\le t \le 1 $ and $n \in \N$. Thus
\[   \lim_{\tau\to \infty} \exp\Big( t  \partial \frac{\delta}{\delta w} \T^{\Gardner}_{app,N}(i\tau,.,\tau_0)\Big)w_0^n      \]
exists in $L^2$, uniformly for $0 < t < 1$ and the limit is a weak solution as discussed above. 
 The Hamiltonians $ \T^\Gardner_{-1}(i\tau,.,\tau_0) $ and $\T^\Gardner_{app,N}$ Poisson commute and their flows commute. Thus also the limits of the flows commute. \end{proof}

\subsection{Uniqueness of weak solutions}
\label{subsec:unique}

We upgrade the previous construction to uniqueness of weak solutions. 

\begin{theorem} \label{thm:Gardneruniqueness}
Given $w_0 \in L^2(\R) $ there exists a unique weak solutions to the $N$th Gardner equation in the Kato smoothing space $X_N$. 
\end{theorem}

\begin{proof} 
The central part of the argument is contained in the   next proposition. 
\begin{proposition} \label{lem:com} Let $ w \in X_N$  be a weak solution to the $N$th Gardner equation. Then 
\[ v(t)= \exp\Big(-t  \partial \frac{\delta}{\delta w}  \T^{\Gardner}_{app,N}(i\tau,.,\tau_0)\Big) w(t) \] 
lies in $X_N$ and satisfies in a distributional sense 
\begin{equation}\label{eq:vcommute}   v_t =   \partial\frac{\delta}{\delta w} \Big(  H^\Gardner_N  ( v  ,\tau_0) - \T^{\Gardner}_{app,N}(i\tau,v,\tau_0)\Big).  \end{equation} 
\end{proposition} 

We postpone the proof of Proposition \ref{lem:com} and continue with the proof of Theorem \ref{thm:Gardneruniqueness}. Let $w\in X_N$ be a weak solution and $v$ be as above. By Theorem \ref{thm:Gardnerprecompact}
$ w \in C( \R; L^2) $ and hence the set $ \{ w(t) : 0\le t \le1 \}\subset L^2$ is compact. 
Moreover $v$ satisfies \eqref{eq:vcommute} and $ v \in X_N$. By Proposition \ref{prop:Kato:diff} and precompactness we obtain 
\[ \lim_{\tau \to \infty} v(t) \to w_0 \quad \text{ in } L^2 \quad \text{ uniformly for $t$ in compact time intervals}.   \] 
Now  
\[   \exp\Big( t \partial \frac{\delta}{\delta w} \T^\Gardner_{app,N}\Big) w_0 -w(t) 
= \exp\Big(t \partial \frac{\delta}{\delta w} \T^\Gardner_{app,N} \Big) w_0  -  \exp\Big(t \partial \frac{\delta}{\delta w} \T^\Gardner_{app,N} \Big) v(t)  \] 
and uniqueness follows if we prove that the right hand side converges to zero in $L^2$. 
The orbit $v(t)$ stays in a precompact set.   The approximate flow is uniformly continuous on precompact sets and hence 
\[ \lim_{\tau \to \infty}   \exp\Big( t \partial \frac{\delta}{\delta w} \T^\Gardner_{app,N}\Big) w_0 -w(t) = 0 \qquad \text{ in } L^2  \]
uniformly on compact time intervals. 
In  the previous subsection  we have seen that 
\[ \exp\Big( t \partial \frac{\delta}{\delta w} \T^\Gardner_{app,N}(i\tau, . , \tau_0\Big) w_0 \to \exp\Big( t \partial \frac{\delta}{\delta w} H^\Gardner_N(.,\tau_0)\Big) w_0  \qquad \text{ as } \tau \to \infty\]
in $L^2$ uniformly on precompact sets for $t$ in a compact interval.   Thus the weak solution is the one we constructed in Theorem \ref{thm:Gardnercontinuous}. 
\end{proof}

\begin{proof}[Proof of Proposition \ref{lem:com}.]
$ \T^\Gardner_{app,N}$ is a sum of $ \T^\Gardner_{-1}$ functions  at   $ 2^j \tau$, $ 0 \le j \le N$. 
  To simplify the notation and the argument we restrict to a multiple of a single $\T^\Gardner_{-1}$. Since the multiple does not change the argument we fix it to $1$.
  Let $ w \in X_N$ be a weak solution and $ v = \exp\big( -t \frac{\delta}{\delta w} \T^\Gardner_{-1}(i\tau,.,\tau_0)\big) w(t)$. We claim that $v \in X_N$ and 
  \begin{equation} \label{eq:comfund}  v_t= \partial \frac{\delta}{\delta w} \big( H^\Gardner_N(v, \tau_0) - \T^\Gardner_{-1}(i\tau,v,\tau_0)\big)\end{equation}
  in a distributional sense. 
  The argument can easily be iterated. The proof of \eqref{eq:comfund} displays the crucial argument.  The proposition is a consequence of  the  following  three Lemmas. Let $ F,G$ be Poisson commuting Hamiltonians, $\{ F, G\} = 0 $. Let $ X_F$ and $X_G$ be the corresponding commuting vector fields, $ [X_F,X_G]=0 $. Let $\Phi_F(t)$ resp. $ \Phi_G(t)$ be the corresponding flows which commute: $ \Phi_G(t, \Phi_F(s,x)) = \Phi_F(s,\Phi_G(t,x_0)) $. By the chain rule, if we differentiate both sides with respect to $t$ and evaluate at $t=0$
  \begin{equation} \label{eq:comweak}   X_G(\Phi_F(s,x_0) ) = D\Phi_F(s,x_0)X_G(x_0).   \end{equation}

As a consequence 
\begin{equation} \label{eq:comcomplete} \begin{split} \hspace{1cm} & \hspace{-1cm} \frac{d}{dt} \phi_G(ct, \Phi_F(t,x_0)) |_{t=t_0}\\ & =  cX_G( \Phi_G(ct_0,\Phi_F(t_0,x_0)))
+   D\phi_G(ct_0, \Phi_F(t_0,x_0)) X_F(\Phi_F(t_0,x_0)) 
\\ & = cX_G ( \Phi_G(ct_0,\Phi_F(t_0,x_0))) + X_F ( \Phi_G(ct_0,\Phi_F(t_0,x_0))). 
\end{split}
\end{equation} 
Here we used the chain rule for the first equality and \eqref{eq:comweak} for the second.

The  Lemma \ref{lem:chainrule} provides the chain rule  for weak solutions. Lemma \ref{lem:inte} gives the analogue of \eqref{eq:comweak}, which relies on complete integrability.  They provide the formula for almost all $t$.
Lemma \ref{lem:reg} gives a delicate estimate which is crucial for completing the argument. It provides finally integrability in time. 

\begin{lemma} \label{lem:reg}
Let $ w \in L^2$, $ \sech w \in H^N$, $ v(s) = \exp\Big( s \partial \frac{\delta}{\delta w} \T^\Gardner_{-1}(\tau,.,\tau_0)\Big)(w) $ and $ \phi \in C^\infty_c(\R)$. Then 
\[
\begin{split}
\hspace{1cm} & \hspace{-1cm} 
\int  D \exp\Big( s \partial \frac{\delta}{\delta w} \T^\Gardner_{-1} (i\tau, v, \tau_0)\Big) \Big[ \partial \frac{\delta}{\delta w} H^\Gardner_N( w , \tau_0) \Big] \phi dx 
\\ & = -D H^\Gardner_N ( v(s),\tau_0) \Big[ \partial \Big\{ D \exp\Big( s\partial \frac{\delta}{\delta w} \T^\Gardner_{-1} (i\tau, w, \tau_0)\Big)\Big\}^T \phi   \Big]
\end{split} 
\]
and
\[ 
\begin{split} \hspace{1cm} & \hspace{-1cm} 
\Big| \int  D \exp\Big( s \partial \frac{\delta}{\delta w} \T^\Gardner_{-1} (i\tau, v, \tau_0)\Big)\Big[ \partial \frac{\delta}{\delta w} H^\Gardner_N( w , \tau_0) \Big] \phi dx           \Big| 
\\ & \lesssim c \big(1+  \Vert \sech w \Vert_{H^N}^2         \big) \Vert \cosh^4 \phi \Vert_{H^{N+1}}.
\end{split}
\]
Here $D$ denotes the Frech\'et derivative. The expressions are continuous in the appropriate topologies. 

\end{lemma} 
\begin{proof} 
The first formula is a consequence of the definition of the variational derivative. 
For the second formula we have to study the homogeneous terms in $H^\Gardner_N$

The Hamiltonians $H^\Gardner_N$ are integrals over differential polynomials, more precisely sums of integrals over multiples of differential monomials
\[ \prod_{j=0}^N \tau_0^\gamma (w^{((j)})^{\alpha_j}   \]
with homogeneity $H= \gamma+ \sum_{j=0}^N  \alpha_j\ge 2$, either $\gamma=0$ or $H-\gamma \ge 3$ and weight (total number of derivatives) $ M=\sum_{j=1}^N j \alpha_j $ with (Theorem \ref{thm:formofkdv}) 
\[ 2N+2  = \gamma + H+ M. \]
Moreover $M$ is even. 
Let $\gamma=0$  to simplify the notation and  
\[ h_N(w) =   \int \prod_{j=0}^N \{w^{(j)}\}^{\alpha_j} dx.  \]
By integration by parts we may assume that the largest $j$
occurs twice. 
We polarize this expression and obtain a symmetric expression 
\[ g_N( w_1, \dots ,w_H) \]
so that $ h_N(w) = g( w,\dots, w) $. Then 
\[ D h_N(w)[ \psi]  = H  g(\psi, w, \dots ,w).\]
We estimate (now with $ \gamma \ge 0 $, picking a particular term) 
\begin{equation} \label{est:GardnerHamEst}   
\tau^\gamma \left| \int_I   w^{2(N-m)-\gamma}     (w^{(m)})^2  dx  \right| 
 \lesssim  \Vert w\Vert_{L^\infty}^{2(N-m)-\gamma}  \Vert w \Vert^2_{H^m}
\end{equation} 
By interpolation all terms with this homogeneity, weight and power of $\tau$ satisfy the same estimate.   We continue (restricting again to the most important case  $\gamma=0$ for notational simplicity) with $I$ an interval
\[\begin{split} \Vert w \Vert_{L^\infty(I)}^{2(N-m)} \Vert w \Vert^2_{H^m(I)}  \,& \le  \Vert  w \Vert_{L^2}^{N-m} \Vert w \Vert_{H^1(I)}^{N-m} \Vert u \Vert_{H^m(I)}^2
\\ &  \le   \Vert w \Vert_{L^2(I)}^{N-m + \frac{N-1}N (N-m)+ \frac{2(N-m)}N} 
\Vert u \Vert_{H^N(I)}^{\frac1N( N-m)  + 2 \frac{m}N   }
\\& = \Vert w \Vert_{L^2(I)}^{  (2+\frac{1}N)  (N-m)} \Vert w \Vert_{H^N(I)}^{1+\frac{m}{N}}  
\\ & \le   \Vert w \Vert_{L^2(I)}^{H-2}      \Vert w \Vert^2_{H^N(I)}
\end{split}
\]
We polarize and sum using weights
\[  | g( \psi,, w , \dots w)| \lesssim \Vert \cosh^2 \psi \Vert_{H^N} 
 \Vert \sech w \Vert^2_{H^N}  \Vert w \Vert_{L^2}^{H-3}.\] 
We apply this estimate with $\psi  $ given by (with different $\phi$ in the formula) 
\[  \psi= \partial \Big\{ D \exp\Big( s\partial \frac{\delta}{\delta w} \T^\Gardner_{-1} (i\tau, v, \tau_0)\Big)\Big\}^T \phi.      \]

By the Schwarz theorem on symmetry of second derivatives 
\[ \begin{split}  \int \Big[  D \frac{\delta}{\delta w} \T^\Gardner_{-1}(i\tau, w,\tau_0) [ w_1 ] \Big] w_2 dx \, & =  \frac{\delta^2}{\delta w^2} \T^\Gardner_{-1}(i\tau, w, \tau_0) [ w_1,w_2]
\\ & = \int \Big[  D \frac{\delta}{\delta w} \T^\Gardner_{-1}(i\tau, w,\tau_0) [ w_2 ] \Big] w_1 dx
\end{split}
\]
The linearized equation is (the dot indicating linearized variables)
\[ \dot v_t =  \partial   D \frac{\delta}{\delta w} \T^\Gardner_{-1}( i \tau, v , \tau_0) [\dot v] .  \]
Suppose that $\psi$ satisfies 
\[  \psi_t = D \frac{\delta}{\delta w} \T^{\Gardner}_{-1}(i\tau, v(t) , \tau_0) [\partial \psi].  \]
Then 
\[ \frac{d}{dt} \int \dot v(t) \psi(t) dx = 0. \]
and $ \int \dot v(T) \hat \psi(T) dx = \int \dot v(0) \psi(0)dx $.
Let $ \tilde \psi = \partial \psi$. It satisfies 
\[ \tilde \psi_t = \partial D \T^\Gardner_{-1}(i\tau, v , \tau_0) [\tilde \psi(t)] \]
hence, by Theorem \ref{thm:flow} 
\[ \Vert \tilde \psi(t) \Vert_{L^2} \le c \Vert \tilde \psi(T) \Vert_{L^2} \]
and, if $ \tau \ge 4$ (so that $ \cosh^2$ is $ \tau$ slowly varying), by a small variation of Theorem \ref{thm:flow},  
\[ \Vert \cosh^2   \tilde \psi(t) \Vert_{H^N} \le c (1+\Vert \sech w \Vert^2_{H^N} )  \Vert \cosh^4 \tilde \psi(T)  \Vert_{H^N}.  \]
Since $ \int \tilde \psi dx $ is conserved and hence $0$ we also have 
\[ \Vert  \psi(t) \Vert_{L^2} \le c \Vert \psi(T) \Vert_{L^2}, \quad \Vert \cosh^2 \hat \psi(t) \Vert_{H^{N+1}} \le c (1+ \Vert \sech w \Vert^2_{H^N})\Vert \cosh^4 \hat \psi(T) \Vert_{H^{N+1}}. \]
We choose $  \psi(T) = \phi$.
 Thus 
\[  \Vert  \cosh^2 \psi(0)  \Vert_{H^{N+1}} \le c(1+ \Vert \sech w \Vert_{H^N}^2)  \Vert \cosh^4 \phi  \Vert_{H^{N+1}} \]
which completes the proof. 
\end{proof}

 \begin{lemma}[Chain rule] \label{lem:chainrule}
 Let $ \phi $ be  test function and $ w \in X_N$ be a weak solution and $ \sech w(t) \in H^N$, which holds for almost all $t$. Then 
 we have for almost all $t$ (so that $ \sech w^{(N)}(t,.) \in L^2$) 
 \begin{equation}
 \begin{split} \hspace{1cm} & \hspace{-1cm} 
 \frac{d}{dt}  \int    \Big[\exp\Big(s \partial \frac{\delta}{\delta w} \T^{\Gardner}_{-1}(i\tau,.,\tau_0)\Big) w(t)\Big]   \phi   dx  
\\ &  = \int  D \exp\Big(s \frac{\delta}{\delta w} \T^\Gardner_{-1}(i\tau, w(t), \tau_0)\big) \partial_x \Big[ \frac{\delta}{\delta w} H^\Gardner_N ( w(t), \tau_0)  \Big] \phi  dx  
\end{split}
\end{equation}  
 with  the integral over the right hand side  estimated  using $w \in X_N$ and  Lemma \ref{lem:reg}.
\end{lemma} 
\begin{proof} 
The chain rule holds for smooth maps. Let $J_\varepsilon w(t)$ be the convolution (in $x$) with a compactly supported Dirac sequence. Then $ \partial_t J_\varepsilon w \in L^1_{loc} (\R, H^N(\R))$ and 
\[
 \begin{split} \hspace{1cm} & \hspace{-1cm} 
 \frac{d}{dt}  \int    \exp\Big( s \partial \frac{\delta}{\delta w} \T^{\Gardner}_{-1}(i\tau,.,\tau_0)\Big) J_\varepsilon w(t)    \phi   dx  
\\ &  = \int   D \exp\Big(s \frac{\delta}{\delta w} \T^\Gardner_{-1}(i\tau, J_\varepsilon w(t), \tau_0)\Big)\Big[ J_\varepsilon \partial_x \frac{\delta}{\delta w} H^\Gardner_N(w(t),\tau_0) \Big] \phi(t)  dx.  
\end{split}
\]  
The left hand side converges to 
\[  \frac{d}{dt} \int \exp\Big(  s \partial \frac{\delta}{\delta w} \T^\Gardner_{-1} (i\tau, .,\tau_0\Big)[ w(t)]\phi dx\]
as $\varepsilon \to 0 $.
If    $\tau\ge \max\{2,\tau_0\}$ then by Theorem \ref{thm:flow} 
  \[   \Big\Vert \sech   \exp\Big( s \partial \frac{\delta}{\delta w} \T^{\Gardner}_{-1}(i\tau,.,\tau_0)\Big)  w(t)   \Big\Vert_{H^N} \le c(T,\Vert w \Vert_{L^2}) \Vert \sech w(t) \Vert_{H^N}  \]
and 
\[ 
\begin{split} \hspace{1cm} & \hspace{-1cm}  
\Big \Vert  \sech  D \exp\Big( t \partial \frac{\delta}{\delta w} \T^\Gardner_{-1}( i\tau, w, \tau_0)\Big)[\dot w] \Big\Vert_{H^{N+1}} 
\\ & \le  c( \Vert w \Vert_{L^2})  \big( \Vert  \sech \dot w \Vert_{H^{N}} + \Vert \dot w \Vert_{L^2} \Vert \sech w \Vert_{H^N}\big).
\end{split}
\] 
All the  quantities are continuous in all arguments in the norms of the estimates. 
\end{proof} 

The second identity relies on integrability. 

\begin{lemma} \label{lem:inte}
Let $ w \in L^2$, $sech  w \in H^N$, and $ \tau \ge 2 $. Then 
\begin{equation}\label{eq:inte}  \begin{split} \hspace{1cm} & \hspace{-1cm}   \int   D \exp\Big(t \partial\frac{\delta}{\delta w} \T^\Gardner_{-1}(i\tau, ., \tau_0)\Big)\Big[ \partial \frac{\delta}{\delta w} H^\Gardner_N(w,\tau_0)\Big] \phi  dx  
\\ & = -\int \Big[   \frac{\delta}{\delta w}  H^\Gardner_N\Big( \exp\big( t \frac{\delta}{\delta w} \T^\Gardner_{-1}(i\tau, ., \tau_0)\big) w ,\tau_0 \Big)\Big]  \partial \phi dx.  
\end{split}
\end{equation} 
\end{lemma} 

 \begin{proof}[Proof of Lemma \ref{lem:inte}:]
  For $w  \in  H^{2N+1}(\R)$, we have in the $L^2$ sense 
\begin{equation}\label{eq:commutativityequation}
\begin{split} \hspace{1cm} & \hspace{-1cm} 
     D \exp\Big( t\partial \frac{\delta}{\delta w} \T^\Gardner_{-1}(i\tau, w, \tau_0)\Big)   \partial \frac{\delta}{\delta w} H_N^{\Gardner}(w,\tau_0)\\ & = \partial \Big(\frac{\delta}{\delta w} H^\Gardner_N\Big)\Big(\exp( t \partial \frac{\delta}{\delta w} T^\Gardner_{-1}(i\tau, ., \tau_0) w ,\tau_0 \Big)
     \end{split}
\end{equation}
from the considerations above for smooth flows which can be written as weak equation in the form of \eqref{eq:inte} and we have to prove that \eqref{eq:inte} has a continuous extension to 
$w \in L^2$, $ \sech w \in H^N$.

Since (omitting the factor $ \sech$)
\[ \exp\Big( t \frac{\delta}{\delta w} T^\Gardner_{-1} (i\tau,. , \tau_0) \Big)  \] 
is a continuous map from $  \{ w \in L^2: \sech w \in H^N\}$ to itself
 we obtain continuity of the right  hand side. Similarly 
$D \exp\Big( t \frac{\delta}{\delta w} T^\Gardner_{-1} (i\tau,w , \tau_0) \Big)$ 
is a bounded continous map on that space provided $ w \in L^2$, $ \sech w \in H^N$.
 \end{proof} 
We complete the proof of Proposition \ref{lem:com}. By the Leibniz rule, the chain rule Lemma \ref{lem:chainrule} and  integrability (Lemma \ref{lem:inte})
\[ \begin{split}   v_t(t)\, &  =  -\partial \frac{\delta}{\delta w} \T^{\Gardner}_{app,N}(i\tau, .,\tau_0) (v(t)) 
\\ & \qquad +  D \exp( -t \partial \frac{\delta}{\delta w} \T^\Gardner_{app,N}(i\tau, ., \tau_0) [ \partial \frac{\delta}{\delta w} H^\Gardner (w(t),\tau_0)]
\\ & =  -\partial \frac{\delta}{\delta w} \T^{\Gardner}_{app,N}(i\tau, .,\tau_0) (v(t)) 
+ \partial \frac{\delta}{\delta w} H^\Gardner_N( v(t),\tau_0)  
\end{split} 
\]
in the sense of distributions in $x$, provided $ \sech w(t) \in H^N$, hence for almost all $t$.
It remains to verify that the identity holds in the distributional sense in space and time, which amounts to proving integrability when we test by a space time test function.
To prove this we provide bounds for terms in the  equalities above.
 Let $ \phi $ be a test function in $ \R$. By Lemma \ref{lem:reg}
 \[ 
 \begin{split} \hspace{1cm} & \hspace{-1cm} 
 \Big|\int  D \exp( -t \partial \frac{\delta}{\delta w} \T^\Gardner_{app,N}(i\tau, ., \tau_0) [ \partial \frac{\delta}{\delta w} H^\Gardner (w(t),\tau_0)]\phi dx  \Big| 
\\ &  \le  c ( 1+ \Vert \sech w(t) \Vert_{H^N}^2) \Vert \cosh^2 \phi \Vert_{H^{N+1}}  
 \end{split}
 \]
and, using again Theorem \ref{thm:flow} 
\[   
\begin{split}
\Big| \int  \partial \frac{\delta}{\delta w} H^\Gardner_N( v(t),\tau_0) \phi dx \Big| 
& \le  c \Vert \sech v(t) \Vert^2_{H^N} \Vert \cosh^4 \phi \Vert_{H^{N+1}} 
\\ & \le c \Vert \sech w(t) \Vert_{H^N}^2 \Vert \cosh^4 \phi \Vert_{H^{N+1}}. 
\end{split} \]
\end{proof}

\subsection{Outline of the text}
In this section we have proven the main theorems and we postponed most of the fixed time estimates on the Miura map and the Hamiltonian vector fields of Hamiltonian functions, the Kato smoothing estimate and precompactness of weak solutions.   In Section \ref{sec:AKNS} we develop the structure and the formulas for KdV and the Gardner hierarchy.  The modified Miura map and estimates for it are the central object in Section \ref{sec:mmiura}.
Section  \ref{sec:weak} uses the modified Miura map to verify equivalence of weak solutions to the $N$th KdV equation, the $N$th Gardner equation and the $N$th good variable equation. The object of Section \ref{sec:tightness} are properties of weak solutions, Kato smoothing and precompactness of orbits. 
 Section \ref{sec:differenceflow} gives  the asymptotic series a precise meaning. It provides
several  multilinear estimates in particular for the difference flow.  
Some of the formulas and calculations are contained in the appendices. 

\section{The KdV and the Gardner hierarchy}

\label{sec:AKNS} 

The KdV hierarchy is a long studied object  and we refer to Babelon et al \cite{MR1995460}, Faddeev and  Takthajan \cite{MR2348643}, Novikov et al \cite{MR779467},  Dickej \cite{MR1964513} and Gesztesy and Holden \cite{MR1992536}. The KdV Hamiltonians  and related quantities are coefficients of asymptotic series, which are typically treated as formal series. We give the asymptotic series a precise meaning. In contrast to  the expositions above we insist on working in an $H^N$  setting for the potentials and on sharp error estimates for the difference of partials sums and quantities expanded in asymptotic series. 
We deduce consequences on the coefficients from properties of the logarithm of the transmission coefficient and related quantities.  
This allows to use rigorous arguments for potentials in Sobolev spaces, and is central in the approximation of higher flows.  The various connections between the KdV hierarchy and the Gardner hierarchy may be an original contribution, including  a Lenard recursion without antiderivatives. 
In the Subsection \ref{subsec:determinants} on Schatten class operators and determinants, we also used some elegant arguments from
Harrop-Griffiths,
Killip and Visan \cite{harropgriffiths2020sharp}. 

In what follows we will often omit arguments of functions when they are clear from the context.
We also omit the upper index $\KdV$ for the transmission coefficients when the meaning is clear from the context. 
\subsection{The KdV hierarchy} 

We consider the Korteweg-de Vries equation (KdV) \eqref{eq:KdV},
\[
	u_t = -u_{xxx} + 6uu_x
\]
and its hierarchy. The KdV equation \eqref{eq:KdV} has the form of a Hamiltonian equation
\begin{equation*}
    u_t = \partial \frac{\delta H^{\KdV}}{\delta u},
\end{equation*}
where $H^{\KdV} = \frac12 \int u_x^2 + 2 u^3 dx$ and $\frac{\delta H}{\delta u}$ is the functional derivative of $H$ defined to be the unique function such that
\begin{equation}\label{eq:fd}
   \int \phi(x) \frac{\delta H}{\delta u}(x)\, dx = \left.\frac{d}{dt}\right|_{t=0} H(u+t\phi)
\end{equation}
for all test functions $\phi \in \Sc(\R)$. Given a functional of the form
\begin{equation*}
    H(u) = \int_\R h(u,u',u'',...) dx,
\end{equation*}
one can check easily that 
\begin{equation*}
    \frac{\delta H}{\delta u} = \sum_{i=0}^\infty (-\partial)^i \left(\frac{\partial h}{\partial u^{(i)}}(u,u',...)\right).
\end{equation*}
Whenever the functional can be written as an integral over a differential polynomial
 the sum on the right-hand side is finite.

The Lax operator for KdV is the Schr\"odinger operator
\begin{equation}  
	L \phi = (- \partial_{x}^2 + u) \phi, \label{schroedinger}
\end{equation}
with potential $u\in H^{-1}$. We omit the index $\KdV$ for the Lax operator in this part. 
Let $ \im\, z >0 $ and consider the left and right Jost solutions $ \phi_l,\phi_r$ of
	\begin{equation}\label{eq:spectral} (L-z^2) \phi=  -\partial_x^2 \phi + u \phi - z^2 \phi=0, \end{equation} 
with the normalization at $\pm \infty$ \footnote{This deviates slightly from the standard normalization. It has the advantage that it is correct even for $ u \in H^{N}$, $N \ge -1$, in contrast to the standard normalization. It leads to the renormalized transmission coefficient below.}  
\begin{equation}\label{eq:normalJost}  \lim_{x \to -\infty} e^{izx+ \frac1{2iz} \int_0^x u(y) dy   } \phi_l(x) = 1, \end{equation} 
in a local $L^2$ sense, and similarly
\[ \lim_{x\to \infty} e^{-izx-\frac1{2iz} \int_0^x u(y) dy } \phi_r(x) = 1. \]
The integral term in the normalization is needed to obtain a limit for  $u \in H^{-1}$   - it is needed even for $u \in H^N $ for any $N$. For $u \in H^{-1}$ we need an inessential correction to the normalization described below. 
One way to see its origin  is  as follows: We write 
$ u=v_x-2iz  v$ with $ \Vert v \Vert_{L^2} = \Vert e^{2i\real z x} u \Vert_{H^{-1}_{2\im z}}$. Then 
\[ L-z^2 = (\partial - iz +v )(-\partial -iz +v)-v^2,  \]
and  \eqref{eq:spectral}
is equivalent to the system
\begin{equation}\label{eq:renormL}   \psi' = \left( \begin{matrix}   0  & -1\\ 
 v^2 & 2iz - 2v\end{matrix} \right) \psi
\end{equation}
with 
\[ \psi_1 = e^{izx- \int_0^x v dy } \phi,  \quad  \psi_2= e^{izx- \int_0^x v dy}(-\partial - iz  + v) \phi. \] 
The system \eqref{eq:renormL} can easily be solved by a Picard iteration from $-\infty$ starting with the constant function 
$ \left( \begin{matrix} 1 \\0 \end{matrix} \right)$, see Lemma \ref{lem:Jost}. We observe that formally $ \int u dx = -2iz \int v dx $. The limit in \eqref{eq:normalJost} becomes a standard limit  
if we replace $-\frac1{2iz} \int_0^x u $ by  $\int_0^x v$. We define the renormalized transmission coefficent $T_{r} (z)$ on the upper half plane to be the meromorphic function
	\[ T_r(z) = \frac{-2iz}{ W(\phi_l,\phi_r)}= 
\frac{\lim\limits_{x\to -\infty} \exp\Big(izx-\int_0^x v(y) dy\Big)\phi_l(x) }
{ \lim\limits_{x\to \infty} \exp\Big(izx- \int_0^x v(y) dy\Big)  \phi_l(x)}  .\label{eq:Trenorm} \]
Here, $W(f,g) = f'g-fg'$ is the Wronskian and the enumerator is $1$, which we kept to make the expression independent from the normalization. For Schwartz functions the transmission coefficient $T^\KdV$ is defined using the standard normalization of the Jost solutions without the integral in the exponent. Then
\begin{equation} \label{eq:Tren}  T_r =e^{-\frac{1}{2iz} \int_{\R} u dx }    T^\KdV . \end{equation} 
Indeed, for functions $u\in L^1$ the standard left Jost function $\tilde \phi_l$ is related to the Jost function $\phi_l$ defined in \eqref{eq:spectral} and \eqref{eq:normalJost} by $\phi_l = \tilde \phi_l \exp(\frac1{2iz}\int_{-\infty}^0 u \, dy)$.

None of the factors on the right hand side of \eqref{eq:Tren} is defined unless $u \in L^1$ (or $W^{-1,1} \cap H^{-1}$ with a suitable definition of spaces), but the left hand side is defined  for $u \in H^{-1}$.
On the lower half plane we define $T_r(z) = T_r(-z)^{-1}$ which is the correct choice for   \eqref{eq:asymptotickdv} below.
The  inverse of the 
transmission coefficient $(T_r(z,u))^{-1}$ is defined in this way for $ u \in H^{-1}$  meromorphic in $z\in \C\backslash \R$. It is holomorphic in the upper half plane, with the zeros  given by the square roots of the eigenvalues in the upper half plane.

A direct calculation (see the proof  of Lemma \ref{lem:fredholm} and \cite{MR779467}) shows that
\begin{equation} \label{eq:fdTren}\frac{\delta}{\delta u} T_r = \frac{T_r^2}{2iz}\phi_l \phi_r - \frac{T_r}{2iz}, \quad\text{ and hence } \quad  \frac{\delta \ln T_r}{\delta u} = \frac{T_r}{2iz}\phi_l \phi_r - \frac1{2iz},
\end{equation}
respectively for the non-renormalized Transmission coefficient $T^\KdV$ and standard Jost solutions $ \tilde \phi_l$, $\tilde \phi_r$ assuming $u \in L^1$
\begin{equation} \label{eq:fdT}\frac{\delta T^\KdV(z)}{\delta u} = \frac{(T^{\KdV}(z))^2}{2iz} \tilde \phi_l \tilde \phi_r, \quad\text{ and hence } \quad  \frac{\delta \ln T^\KdV}{\delta u} = \frac{T^\KdV}{2iz} \tilde \phi_l \tilde \phi_r.
\end{equation}
We specialize  \eqref{eq:tauN} to  $N=-1$ (See also \eqref{def:T-1} and \eqref{eq:Tren} )
\[ \T^{\KdV}_{-1}(z,u) = iz  \ln T_r(z). \]

An innocent calculation gives a formula which is hard to overestimate: the Lenard recursion.
\begin{lemma}\label{lemma:lnT-recursion}
Let $ u \in H^{-1} \cap L^2_{loc}$. Then the functional derivative  of $\T_{-1}^{\KdV} $ 
 satisfies the ODE
\begin{equation} \label{eq:lenardrecursion}   
    (-\partial^3+4 u\partial + 2u_x) \frac{\delta \T^{\KdV}_{-1}}{\delta u} = 4 z^2 \partial  \frac{\delta  \T^{\KdV}_{-1}}{\delta u} -  u_x .
\end{equation}  
%The same is true for the functional derivative of $T_r$. 
\end{lemma} 
This Lemma is a generating function version of the Lenard recursion. By formally expanding in inverse powers of $z$ one obtains idenitities for the densities of the KdV Hamiltomians which allow a recursive construction of the Hamiltonians. 
In view of \eqref{eq:fdT} and the regularity of the Jost solution 
$\phi_{l,r} \in H^2_{loc}$ the equation \eqref{eq:lenardrecursion}
can be understood in a distributional sense.

\begin{proof} 
Consider $\phi_1, \phi_2$ solutions to $(-\partial^2 + u - z^2)\phi_i = 0$. A short calculation reveals
    \begin{equation*}
        \partial^3 (\phi_1 \phi_2) = 2u_x \phi_1 \phi_2 + 4(u-z^2)\partial(\phi_1 \phi_2).
    \end{equation*}
    Hence \eqref{eq:lenardrecursion} follows from \eqref{eq:fdTren}.
\end{proof} 

There are different formulations of the equations of the Korteweg-de Vries hierarchy: the Gardner hierarchy and the good variable hierarchy, which are equivalent under weak conditions. They correspond to taking different coordinates in a large part of the phase space. The different coordinates are based on relations between the Lax operator, resolvent and Jost solutions. We give connections between these variables in the next Lemma.

\begin{lemma}[Definition of $w$ and $v$]\label{lem:relationuvw} 
Let $ u \in H^{-1} $ and let $\phi_l\in H^1_{loc}$ be the left Jost solution 
with the normalization \eqref{eq:normalJost}. We assume that $ -\partial^2+ u  + \tau_0^2$ is positive semi definite and either  $\real z\ne 0$ or $ z \notin i[0,\tau_0]$. Then $ \phi_l$ never vanishes and we  define $w\in L^2$  (see Lemma \ref{lem:Jost})  as
\begin{equation*}
    w(x,z,u) = (\ln (\phi_l(x,z,u)e^{izx}))'.
\end{equation*}
Then $w$ satisfies the Ricatti equation  
\begin{equation}\label{eq:riccati}
    w' - 2iz w + w^2 = u.
\end{equation}
It allows to factorize the Lax operator
\begin{equation}\label{eq:factor}   L-z^2 = -\partial^2 + u - z^2 = (\partial -iz +w)(-\partial-iz + w). \end{equation}
Moreover
\begin{equation}\label{eq:rhologT}
    \T^\KdV_{-1} (z,u) = -\frac1{2} \int_\R w^2(x,z,u) \, dx.
\end{equation}
Let $G(z,x,y)$ be the integral kernel of the resolvent and $\beta(z,x)= G(z,x,x)$\footnote{Also we note here for comparison that the renormalized perturbation determinant $\alpha$ as defined in \cite{MR3990604} satisfies $\T_{-1}(i\tau,u) = \tau \alpha(\tau,u)$.} its value on the diagonal. We define 
\begin{equation}\label{eq:vdef} v = -\frac{1}{2iz \beta} - 1.\end{equation} 
Then 
\begin{equation}\label{eq:functvw}    w= -\frac12 \partial_x \ln (v+1) -iz  v. \end{equation} 
\end{lemma}
\begin{proof}
     In Lemma \ref{lem:Jost} we will prove that the map 
    $L^2\ni  w \to u \in H^{-1}$ is a diffeomorphism with the natural restrictions on the spectrum. 
 From
    \[\phi_l' = (-iz + w)\phi_l, \qquad \phi_l'' = (w' + (-iz + w)^2)\phi_l= (u-z^2) \phi_l  \]
     we infer \eqref{eq:riccati}. The factorization \eqref{eq:factor} is an immediate consequence.
$L^1 \subset H^{-1}$ is dense.
Then (recall that $\tilde \phi_l$ is the Jost solution with the standard normalization and $ \phi_l$ has the normalization defined in  \eqref{eq:normalJost}), 
 \[   \begin{split}   \T^\KdV_{-1}(u,z) \, &  = iz\Big(\ln T(z) - \frac1{2iz} \int u\Big)
  \\ &      = -iz\Big(\lim_{x\to \infty} \ln \frac{ \tilde \phi_l(x) e^{izx}}{\tilde \phi_l(-x) e^{-izx}} + \frac1{2iz} \int u\Big)
 \\ &  = -iz\int w + \frac1{2iz} u dx\\ & = -iz \int w+\frac1{2iz} ( w_x- 2iz w + w^2) dx
 \\ & = -\frac1{2} \int w^2 dx. \end{split} \]

 The factors in the factorization \eqref{eq:factor} can be inverted independently:
\[
    (L-z^2)^{-1} = (-\partial - iz + w)^{-1}(\partial - iz + w)^{-1},
\]
where we formally write
\begin{align*}
    (\partial - iz + w)^{-1}f(t) &= \int_{-\infty}^t e^{iz(t-y)-\int_y^tw\, ds}f(y)\, dy,\\
    (-\partial - iz + w)^{-1}f(t) &= \int_t^\infty e^{-iz(t-y)+\int_y^tw\, ds}f(y)\, dy.
\end{align*}
Thus,
 \[\begin{split}   (L_z^{-1} f)(x)\, & =       \int_x^\infty e^{iz(t-x)} \exp\Big( \int_x^t -w \Big)                         \int_{\infty}^t e^{-iz(t-y)} \exp \Big( \int_y^t -w \Big)  f(y) dy  dt
 \\ & = \int G(z,x,y) f(y) dy 
 \end{split} 
 \]
where
\begin{equation}  G(z,x,y)= 
\int_{\max\{y,x\}}^{\infty}   \exp\Big(-iz (x+y-2t) - 2\int_{\max\{x,y\}}^t  w ds- \int_{\min\{x,y\}}^{\max\{x,y\}}w ds  \Big)       dt. \label{GreenFunction}\end{equation} 
is the Green's function 
and  the evaluation at the diagonal gives  
\begin{align*}
\beta(z,x) &= G(z,x,x) = \int_x^{\infty} \exp\Big( -2iz(x-t) - 2 \int_x^t w ds \Big) dt\\
&= (-\partial-2iz+2w)^{-1}(1),
\end{align*}
and
\begin{equation} \label{eq:greensfct} 
 \beta' + 2iz  \beta - 2w \beta  = -1.
\end{equation}   
We substitute $ \beta = -\frac1{2iz(v+1)}$ and see
\[ 2iz = - \frac{v'}{(v+1)^2} +\frac{2iz}{v+1}-2 \frac{w}{v+1} \]  
which implies \eqref{eq:functvw}.       
    \end{proof}

We note that \eqref{eq:greensfct} is equivalent to   the ODE for the Lenard recursion, \eqref{eq:lenardrecursion}. Indeed, differentiating \eqref{eq:greensfct} twice and using $\beta'' = -2iz\beta' + 2(w\beta)'$ gives
\[
    \beta''' + 2iz\beta'' - 2(w\beta)'' = \beta''' + 4z^2 \beta' + 4iz(w\beta)' - 2(w\beta)'' = 0.
\]
Now since
\begin{align*}
    4iz(w\beta)' - 2(w\beta)'' &= 4izw\beta' + 4izw'\beta - 2w''\beta - 4w'\beta' - 2w(-2iz\beta' + 2(w\beta)')\\
    &= -2(w''-2izw'+2ww')\beta - 4(w'-2izw+w^2)\beta',
\end{align*}
the equivalence follows by $u = w' - 2iz w + w^2$ and expanding in inverse powers of $z$. 

The KdV Hamiltonians are defined for Schwartz functions $u$ as the coefficients 
\begin{equation}\label{eq:asymptotickdv} 
\T^{\KdV}_{-1}(z,u) = 	iz\ln  T_r(z,u) \sim \sum_{n=0}^{\infty} (2 z)^{-2 n-2} H_{n}^{\KdV}(u)
\end{equation}
of the formal  asymptotic series, or, equivalently 
\begin{align*}
& \T_{-1}^{\KdV}(z) =     iz \ln T(z)- \frac{1}{2}\int u \,dx \\
&\sim  (2z)^{-2}\frac12 \int u^2\,dx + (2z)^{-4} \frac12\int u_x^2 + 2u^3 \, dx + (2z)^{-6}\frac12\int u_{xx}^2 - 5uu_{xx} + 5u^4\, dx...
\end{align*}
We recall  for $N \ge -1$ (see \eqref{eq:tauN}) ,
\begin{equation} \label{eq:partialsum} \T^{\KdV}_N(z,u) := \frac{i}2(2z)^{2N+3} \ln T_r^{\KdV}(z) - \sum_{n=0}^N (2z)^{2(N-n)} H_n^{\KdV}(u), \end{equation} 
where in the case $N = -1$ the sum is empty. By \eqref{eq:asymptotickdv}
 $\T^\KdV_{N+1} = (2z)^2 \T_N - H_{N+1}^{\KdV}$ and in the sense of Proposition  \ref{eq:kdv:Lipschitz},
\[
    \T^\KdV_{N}(z,u) \sim \sum_{n=1}^\infty (2z)^{-2n}H_{N+n}^{\KdV}(u).
\]
\subsection{Poisson structures}

We recall the definition of symplectic forms and Poisson structures on $\R^{2n}$, $ \C^n$ and $ \C^{2n} $. A real symplectic form is a nondegenerate closed real  two form. The most relevant real symplectic form on $\C^n$ is 
\[ \omega(z_1,z_2) = -\im \langle z_1,z_2 \rangle \] 
where we choose the convention that the inner product is complex linear in the first component. 
Equivalently we write $ \R^{2n} = \R^n_x \times \R^n_y $ and define
\[ \omega((x_1,y_1), (x_2,y_2)) = \langle x_1 , y_2 \rangle - \langle y_1,x_2 \rangle = \langle (x_1,y_1), J^{-1} (x_2,y_2)\rangle \] 
where $J =\left( \begin{matrix} 0 & 1_{\R^n} \\ -1_{\R^n} & 0 \end{matrix} \right) $. It defines a bilinear form on the dual space which we denote by 
\[ \omega^{-1} (( \xi_1,\eta_1), (\xi_2,\eta_1) ) = \langle (\xi_1,\eta_1) , J (\xi_2 , \eta_2)\rangle.\]
The Hamiltonian vector field of the function $f$ is  $ J \nabla f$.

A complex symplectic form on $\C^{2n}$ is a holomorphic closed $2$ form, the most important  being (writing 
$\C^{2n} = \C^n_z \times \C^n_\zeta$)
\[ \omega ((z_1,\zeta_1), (z_2 , \zeta_2) ) = \sum_{j=1}^n z_1^j \zeta_2^j - z_2^j \zeta_1^j. \]
Obviously the real part  of the restriction to the real subspace of the  symplectic form is a real symplectic form. 

A real (complex) Poisson bracket is a bilinear map mapping a pair of smooth (holomorphic) functions to smooth (holomorphic) functions which satisfies
\begin{align*} 
&\{ f,g \} =- \{ g, f\} , &   \qquad \text{ (skew symmetry)} \\ 
 &\{ f,gh\} = \{ f, g\} h + \{ f,h\} g,& \qquad \text{ (derivation) } \\ 
  &\{ f, \{g,h\}\} + \{g,\{h,f\}\} + \{h,\{f,g\}\} = 0. & \qquad \text{ (Jacobi identity)} 
\end{align*}

Given a smooth (holomorphic) function $H$ (called Hamiltonian) we define the Hamiltonian vector field by 
\[  X_H f = \{f, H\}. \] 
A real (complex) symplectic form $ \omega$ defines a Poisson structure on smooth (holomorphic) functions by (where in the Hilbert space case by an abuse of notation we identify derivatives with functional derivatives via duality)
\[  \{ f,g \}  = \omega(X_f, X_g) \]
($\omega $ defines an isomorphism between tangent and cotangent space, and $ \omega^{-1} $ is the unique induced two form on the cotangent space). The Jacobi identity is a consequence (and it is equivalent to it) of the closedness of the two form $\omega$.
\[  \{ f,g \}  = Df ( \omega^{-1} Dg). \]
The two form $\omega $ defines an isomorphism between tangent and cotangent space, and $ \omega^{-1} $ is its inverse. $Df$ resp. $Dg$ (recall that $D$ denotes the total derivative) take values in the cotangent space. Then $\{ f, g \} = \omega(X_f,X_g)$ where $X_f$ resp. $X_g$ are the Hamiltonian vector fields of $f$ and $g$. The Jacobi identity is a consequence of (and it is equivalent to)  the closedness of the two form $\omega$.

There is no difficulty
to extend these notions
to  infinite dimensional real and complex Hilbert spaces. 
It is important to note that not every Poisson structure comes from a symplectic form, our most important examples for that being the Gardner Poisson bracket and the Magri Poisson structure \cite{MR286402}:

\begin{definition}\label{def:gardnerstructure}
    The Gardner Poisson structure is defined as
    \begin{equation}
        \{F,G\}^{\Gardner} = \int \frac{\delta F}{\delta u} \partial_x \frac{\delta G}{\delta u} dx.
    \end{equation}
 and the Magri Poisson structure 
 \begin{equation}
        \{F,G\}^{\Magri} = \int \frac{\delta F}{\delta u} (-\partial^3_x+ 2( u \partial_x+ \partial_x u))  \frac{\delta G}{\delta u} dx.
    \end{equation}
\end{definition}
Any constant skew symmetric bracket, and in particular the Gardner bracket,  satisfies the Jacobi identity (see \cite{MR1747916}), 
\[ \{ \{ F,G\},H \} + \{ \{ G, H\}, F\} + \{\{ H, F\} , G\} = 0. \] 
We argue differently for the Magri Poisson bracket.
\begin{lemma} \label{lem:magripoisson} Let   $f(w) = F( w_x-2iz w + w^2)$ and $g(w) = G( w_x-2iz w + w^2)$. Then 
\[ \int \frac{\delta}{\delta w} f \partial_x \frac{\delta}{\delta w} g dx 
= \{F,G\}^{\Magri} + (2iz)^2  \{ F,G\}^{\Gardner}. \] 
In particular every linear combination of the Gardner and the Magri Poisson structure satisfies the Jacobi identity. 
\end{lemma} 
We compute using the chain rule $\frac{\delta}{\delta w} f(w)=(-\partial -2iz +2w) \frac{\delta}{\delta u} F|_{u= w_x-2izw+w^2}$, hence
\[\begin{split}\hspace{1cm}& \hspace{-1cm}   \int \frac{\delta}{\delta w} f \partial_x \frac{\delta}{\delta w} g dx \\
& = \int (-\partial -2iz + 2w) \frac{\delta F}{\delta u}\Big|_{u=w_x -2iz w+ w^2}
\partial (-\partial -2iz + 2w) \frac{\delta G}{\delta u}\Big|_{u=w_x -2iz w+ w^2}
dx\\ &  = \int  \frac{\delta F}{\delta u}(\partial - 2iz + 2w ) \partial 
(-\partial -2iz + 2w) \frac{\delta G}{\delta u} dx
\\ & = \int \frac{\delta F}{\delta u}\Big((- \partial^3 + 2(\partial (w_x-2iz+ w^2) +(w_x-2izw+w^2)  \partial)+ (2iz)^2 \Big) \partial          \frac{\delta G}{\delta u }  dx. 
\end{split} 
\] 

With this at hand we prove
\begin{lemma}\label{kdvpoisson} The transmission coefficients $T_r(z_1)$ and $T_r(z_2)$ and  $\T_{-1}^\KdV(z_1)$ and $\T^\KdV_{-1}(z_2)$ Poisson commute for every linear combination of the Gardner and the Magri Poisson bracket.  
\end{lemma} 
\begin{proof} This is a direct formal calculation, which we do first for the Gardner bracket. $\int u\, dx $ is a Casimir for the Gardner structure ( $\{  f, \int u dx \} =0 $ for all functions $f$) and we may ignore it. Of course this amounts to restricting to integrable $u$ 
and doing the calculations for $ T$ instead of $T_r$. 
We  use \eqref{eq:fdT} to find
\[\begin{split}\hspace{.5cm} & \hspace{-.5cm}   
-8z_1z_2\int \frac{\delta T(z_1)}{\delta u} \partial_x \Big(\frac{\delta T(z_2)}{\delta u}\Big) dx = 2T(z_1)^2 T(z_2)^2 \int (\phi_l \phi_r)(z_1)  \partial_x (\phi_l \phi_r)(z_2)  dx 
\\ & = T(z_1)^2 T(z_2)^2   \int (\phi_l \phi_r)(z_1) \partial_x ( \phi_l \phi_r)(z_2) - (\phi_l \phi_r)(z_2) \partial_x (\phi_l \phi_r)(z_1)  dx
\\ & = \frac{T(z_1)^2 T(z_2)^2}{z_1^2-z_2^2}  \int   \partial_x\Big(W(\phi_l(z_1), \phi_l(z_2)) W(\phi_r(z_1), \phi_r(z_2)) \Big) dx    
\\ & = 0,
\end{split}
\]
where in the second last step the defining equation $-\phi'' +u\phi = z^2\phi$ was used  and that the limits at $\pm \infty$ are the same:
\[ \lim_{x\to \pm \infty} W(\phi_l(z_1),\phi_l(z_2) W( \phi_r(z_1) ,\phi_r(z_2) dx 
=  (z_1-z_2)^2(T(z_1) T(z_2))^{-2}.  
\]
 This implies the same statement for $\ln T$ resp. $\ln T_r$ and also $ \T^{\KdV}_{-1}$.
From the above, and Lemma~\ref{lemma:lnT-recursion}, 
\begin{equation*}
    \int \frac{\delta T(z_1)}{\delta u}  (-\partial^3+4u\partial + 2u_x)\frac{\delta T(z_2)}{\delta u}  dx 
=4z_2^2 \int \frac{\delta T(z_1) }{\delta u}  \partial_x \frac{\delta T(z_2)}{\delta u} dx = 0.
\end{equation*} 
Hence the transmission coefficients Poisson commute also with respect to the Magri structure. We complete the argument for the Magri structure for $\T^\KdV_{-1}$
by 
\[ \begin{split} \{ \ln T(z,u), {\smallint} u dx\}^\Magri\, & = 2\int \frac{\delta \ln T(z,.)}{\delta u} \partial_x u dx 
\\ & = 2\lim_{s\to 0} \frac1s \Big(\ln T(z,u+su_x) -\ln T(z,u)\Big) 
\\ & = 2\lim_{s\to 0} \frac1s \Big(\ln T(z,u(.+s))-\ln T(z,u)\Big) = 0
\end{split}
\]
by translation invariance and regularity of $\T_{-1}^{\KdV}$ on Schwartz functions, which are dense. We obtain Poisson commutation with respect to the Magri structure.
\end{proof}

We want to express Poisson brackets with $ \mathcal{T}^\KdV_{-1}$, which can be read as evolution equations for a number of quantities. Let
\begin{equation} \label{eq:Gardner0}
\begin{split} 
\T_{-1}^{\Gardner}(z,w,\tau)\, & := \frac{1}{4z^2 + 4\tau^2} \Big(\frac12 \int w^2\, dx +\T^{\KdV}_{-1}(z,w_x+2\tau w+w^2)  \Big)\\ &  =  \frac12  \frac1{4z^2 + 4\tau^2}\int w^2 -w^2(z) dx . 
\end{split} 
\end{equation}
Here, $w(z)$ is defined as in \eqref{eq:w(z)} and equality of the first and the second line is ensured by \eqref{eq:rhologT}.

\begin{lemma}\label{lem:GardnerPoisson} Let $ \tau < \tau_1 < \tau_2$. Then $\T^\Gardner_{-1}(i\tau_1,.,\tau)$ and $\T^\Gardner_{-1}(i\tau_2,.,\tau)$ Poisson commute
\[  \{\T^\Gardner_{-1}(i\tau_1,.,\tau),\T^\Gardner_{-1}(i\tau_2,.,\tau)\}^\Gardner = 0.   \]
\end{lemma} 
\begin{proof}
By Lemma \ref{lem:magripoisson}
\[\begin{split} 
\hspace{.3cm} & \hspace{-.3cm} \int \frac{\delta}{\delta w} \T^\Gardner_{-1}(i\tau_1,.,\tau) \partial \frac{\delta}{\delta w} \T^\Gardner_{-1} (i\tau_2,.,\tau) dx \\ & = \frac1{16} \frac1{(\tau^2-\tau_1^2)(\tau^2-\tau_2^2)} 
\Big(  \Big\{ \T^\KdV_{-1}(i\tau,. ) - \T^{\KdV}_{-1}(i\tau_1,. ) ,  \T^\KdV_{-1}(i\tau,. ) - \T^\KdV_{-1}(i\tau_2,. ) \Big\}^\Magri \\ & \quad  -(2\tau)^2    \Big\{ \T^\KdV_{-1}(i\tau,. ) - \T^{\KdV}_{-1}(i\tau_1,. ) ,  \T^\KdV_{-1}(i\tau,. ) - \T^\KdV_{-1}(i\tau_2,. )\Big\}^\Gardner \Big) \Big|_{u = w_x +2\tau w + w^2}.
\end{split} 
\]
The right hand side vanishes by Lemma \ref{kdvpoisson}.
\end{proof}

\begin{lemma}\label{lem:poissonbrackets} The following identities hold for  $u \in H^{-1}$
\begin{equation} \label{eq:poissontau} \{ u, \mathcal{T}^{\KdV}_{-1}(z,u)  \} = \frac12\partial_{x}  \frac{v(z)}{v(z)+1},   \end{equation}
\begin{equation} \label{eq:poissonv}  \{ v(z_1) , \mathcal{T}_{-1}^{\KdV} (z_2,u) \} = \frac{1}{4z_1^2 - 4z_2^2}   \partial_x \frac{ v(z_1)-v(z_2) }{v(z_2)+1}, \end{equation} 
\begin{equation}\label{eq:poissonw} \{ w( z_1), \mathcal{T}_{-1}^{\KdV}(z_2,u)\} = \partial_x \frac{\delta}{\delta w}\T_{-1}^{\Gardner}(z_2,w,z_1). \end{equation}
\end{lemma} 
At first sight the Poisson bracket as defined in Definition \ref{def:gardnerstructure} only makes sense for functionals. We  understand identities like \eqref{eq:poissontau} as follows though: For any test function we consider  
$ \int v(z_1,x)\phi(x) dx$ as a functional of $ u \in H^{-1}$. We ask then that 
\[\left\{  \int v(z_1,x) \phi(x) dx, \mathcal{T}^{\KdV}_{-1}(z_2)  \right\} = -\frac1{4z_1^2-4z_2^2}\int \frac{v(z_1,x)-v(z_2,x)}{v(z_2,x)+1} \partial_x \phi dx \]  
for all test function $ \phi$. In the same way we interpret \eqref{eq:poissontau} and \eqref{eq:poissonw}. Using the chain rule we can rewrite  \eqref{eq:poissonw} for $z_1=i\tau_0$
and $z_2=z$
for $ \im z > \tau_0$ as 
\begin{equation}\label{eq:chainruleGardner}
\partial_x \frac{\delta \T^{\KdV}_{-1}(z)}{\delta u}\Big|_{u= w_x+2\tau_0 w+w^2}
= (\partial + 2\tau_0  + 2w )\partial \frac{\delta \T^{\Gardner}_{-1}(z,.,\tau_0)}{\delta w}. 
\end{equation}
which follows from 
\[  \frac{\delta}{\delta u} \int w \phi dx = \int \phi (\partial +2\tau_0 + 2w)^{-1} u  dx.      \]

\begin{proof}
By definition  of $ \T^\KdV_{-1} $, the  variational derivative of the renormalized transmission coefficient \eqref{eq:fdTren}, the definition of $\beta$ as diagonal Green's function of Lemma \ref{lem:relationuvw}
and the formula for the Green's function 
\begin{equation}\label{eq:Gxy}  G(x,y) =  \frac{T^\KdV}{2iz}  \times \left\{ \begin{array}{ll}   \phi_l(x) \phi_r(y) \quad & \text{ if } x<y \\ \phi_r(x) \phi_l(y) & \text{ if } y<x \end{array} \right.  \end{equation} 

\[ \{ u, \T^{\KdV}_{-1}(z)  \} =  \partial_x \frac{\delta \mathcal{T}^{\KdV}_{-1}(z)}{\delta u} = iz \partial_x (\beta+ \frac1{2iz}) = \frac12\partial_x (\frac{v}{v+1}) \]
which is \eqref{eq:poissontau}.
 Using  the notation $L_z  = -\partial^2 +u -  z^2$, we calculate (see also Killip and Visan \cite{MR3990604} and \cite{MR2348643}) 
\[
\begin{split} \hspace{2mm} & \hspace{-2mm} 
2\left\{ v(z_1) ,  \mathcal{T}^{\KdV}_{-1}( z_2 ) \right\}  = - \frac{2iz_2}{2iz_1}  \left\{ \frac1{\beta(z_1)} , \ln T^{\KdV}_r(z_2) \right\}
\\ & = \frac{z_2}{z_1} \beta^{-2} (z_1) 
\big\{ \beta(z_1) , \ln T^{\KdV}_r(z_2) \big\} 
\\ & =- \frac{z_2}{z_1} \beta^{-2}(z_1) \Big(   L_{z_1}^{-1} \{ u, \ln T^{\KdV}_r(z_2) \} L_{z_1}^{-1} \delta_x\Big)(x) 
\\ & = - \frac{z_2}{z_2}  \beta^{-2} (z_1) 
 \int G(z_1, x,y) (\partial_y \beta(z_2,y)) G(z_1, y,x) dy.
\end{split} 
\]
We use \eqref{eq:lenardrecursion}, and rewrite as operators 
\[
\begin{split}
    \beta'''(z_2)  -2 (u \beta(z_2))' -&2u \beta'(z_2)  + 4 z_1^2 \beta'(z_2) \\
    &= L_{z_1}   \beta'(z_2) + \beta'(z_2)L_{z_1} - 2 L_{z_1}  \beta(z_2) \partial_y + 2 \partial_y \beta(z_2) L_{z_1}.
\end{split}
\]
Thus,
\[
\begin{split} 
\left\{ v(z_1) , \mathcal{T}^{\KdV}_{-1}( z_2 ) \right\} &= -\frac{\frac{z_2}{z_1}}{4 z_1^2 - 4 z_2^2} \beta(z_1)^{-2}
\int G(z_1,x,y)G(z_1,y,x)  
 \Big\{  \beta'''(z_2) \\ & \quad -2(u \beta(z_2))' -2u \beta'(z_2) + 4 z_2^2 \beta'(z_2) + (4z_1^2-4z_2^2) \beta'(z_2)\Big\}  dy 
\\ & =\frac{\frac{z_2}{z_1}}{4 z_1^2 - 4 z_2^2} \beta(z_1)^{-2} ( \beta'(z_1)\beta(z_2) - \beta(z_1) \beta'(z_2) ) 
\\ & =  \frac{\frac{z_2}{z_1}} {4z_2^2-4 z_1^2} \partial \frac{ \beta(z_2)}{\beta(z_1) }\\ 
&   = \frac{1} {4z_2^2- 4 z_1^2} 
\partial \frac{v(z_1)+1}{v(z_2)+1}
\\ & = \frac{1} {4z_2^2- 4z_1^2} 
\partial \frac{v(z_1)-v(z_2)}{v(z_2)+1}
\end{split} 
\]   
which is \eqref{eq:poissonv}.

For the last identity we compute using  $2w = -\partial \ln(1+v) - 2iz v$ and by identifying functions with multiplication operators
\begin{align*}
    2w_t &= -\partial \frac{v_t(z_1)}{v(z_1)+1} - 2iz_1 v_t(z_1) \\
    &= \frac{1}{4(z_2^2-z_1^2)} \partial\Big[(v(z_1) +1)^{-1}\partial(v(z_1)+1) + 2iz_1  (v(z_1)+1)\Big](v(z_2) + 1)^{-1}\\
    & =\frac{1}{4(z_2^2-z_1^2)} \partial \Big[\partial + \big[ \partial \ln( v(z_1)+1)\big]+2iz_1
    + 2iz_1 v(z_1) \Big] (v(z_2) + 1)^{-1}  \\
    & = \frac{1}{4z_1^2-4z^2_2} \partial \big(- \partial + 2w -2iz_1 \big) (v(z_2) + 1)^{-1}\\
    &= \frac{iz_2}{2z_2^2-2z^2_1} \partial \big(- \partial + 2w -2iz_1 \big) \beta(z_2) = \frac{iz_2}{2z_2^2-2z^2_1} \partial \frac{\delta}{\delta w} \ln T(z_2,w_x - 2izw + w^2),
\end{align*}
where in the last step we used the chain rule for $ f(w) = g( w_x -2iz w +w^2)$. We rewrite
\[\begin{split}
    \ln T(z_2,u) &= \ln T_r(z_2,u) - \frac{1}{2iz_2}\int u\, dx = \frac1{iz_2}\Big(\T^{\KdV}_{-1}(z_2,u) - \frac12 \int u\, dx\Big)\\
    &= \frac1{iz_2}\Big(\T^{\KdV}_{-1}(z_2,u) - \frac12 \int w^2\, dx + \text{Casimir}\Big)
\end{split}
\]
and get
\[
    \{w(z_1),\T^{\KdV}_{-1}(z_2)\} = \frac{1}{4z_2^2 - 4z_1^2} \partial \frac{\delta}{\delta w}\Big(\T^{\KdV}_{-1}(z_2,w_x-2iz_1w+w^2) - \frac12 \int w^2\, dx\Big).
\]
which in view of \eqref{eq:Gardner0} is \eqref{eq:poissonw}. \end{proof} 

We obtain the hierarchies for the Gardner variable and good variable by expanding the expressions above into asymptotic power series in $(2z)^{-1}$. The generating function $\T^{\Gardner}_{-1}(z,w,i\tau)$ has an asymptotic expansion in $z$,
\begin{equation} \label{eq:Gardnerexpansion} \T_{-1}^{\Gardner}(z,w,i\tau) \sim \sum_{n=1}^\infty (2z)^{-2n} H^{\Gardner}_{n-1}(w,\tau)  \end{equation}
where we call the coefficients $H_n^{\Gardner}$ the Gardner Hamiltonians.
Then e.g.
\[\begin{split} H_0^{\Gardner}  &= \frac12 \int w^2dx , \quad H_1^{\Gardner}= \frac12 \int w_x^2 + w^4 + 4 \tau w^3 dx,\\
H_2^{\Gardner} &= \frac{1}{2}\int w_{xx}^2 + 10w^2 w_x^2 + 2w^6 + 4\tau(5ww_x^2 + 3w^5) + 24\tau^2 w^4 \, dx.\end{split}\] 

In correspondence with KdV we define
\begin{equation}\label{eq:TNGardner}  
    \T_N^{\Gardner}(z,w,\tau) := (2z)^{2N+2}\T^{\Gardner}_{-1}(z,w,\tau) - \sum_{n=0}^N (2z)^{2(N-n)} H^{\Gardner}_{n}(w,\tau).
\end{equation} 
The relation \eqref{eq:Gardner0} can be written as 
\[    \T^{\KdV}_{-1}(z,w_x+2\tau w +w^2) =(4z^2 + 4 \tau^2)\T^{\Gardner}_{-1}(z,w,\tau)  - \frac12 \Vert w \Vert^2_{L^2}.\]
It  implies a similar identity for all $N$,
\begin{equation} \label{eq:Gardner-KdV}   (2z)^2 \T^{\KdV}_{N}(z,w_x+2\tau w + w^2) = (2z)^2 \T^{\Gardner}_{N+1}(z, w,\tau) + (2 \tau)^2 \T^{\Gardner}_{N}(z,w,\tau) \end{equation} 
and 
\begin{equation} \label{eq:GardnerKdVr}  H_{N}^{\KdV} ( w_x + 2\tau w + w^ 2) = H_{N+1}^{\Gardner}(w,\tau) + 4 \tau^2 H_{N}^{\Gardner}(w,\tau),\end{equation} 
or, equivalently
\begin{equation}\label{eq:generalgardnerhamiltonian} H_N^{\Gardner}(\tau, w) = (-4\tau^2)^N \frac12 \Vert w \Vert_{L^2}^2 + \sum_{n=0}^{N-1}(-4\tau^2)^{N-n-1} H^{\KdV}_n(w_x+2\tau w +w^2).
\end{equation}

From the properties of the generating functions we can derive properties of the Hamiltonians. As a first instance, we note:

\begin{corollary}[Corollary of Lemma \ref{kdvpoisson}]\label{cor:Poissoncommutingenergies} 
The KdV Hamiltonians Poisson commute with respect to the Gardner Poisson structure on the Schwartz space with another and the renormalized transmission coefficient resp. $ \mathcal{T}^{\KdV}_{-1}$.
\end{corollary} 

\begin{proof} This is an immediate
consequence of Lemma \ref{kdvpoisson} 
since the Hamiltonians can be defined by limits of Poisson commuting quantities. By definition \eqref{eq:asymptotickdv} and from \eqref{eq:tauNdec} we see that
\[
    H_{N+1}^{\KdV} = \lim_{\tau \to \infty} (2i\tau)^2 \T_N^{\KdV}(i\tau,u), \quad  \text{if} \quad u \in H^{N+1}.
\]
 and  
\begin{equation} \label{eq:limitpoisson} \big\{H_{N+1}^{\KdV}, T^\KdV_{-1}(z)\big\} = \lim_{\tau \to \infty}\{(2i\tau)^2 \T_N^\KdV(i\tau,u),  T^{\KdV}_{-1}(z)\}
\end{equation} 
because 
\begin{align*}
    \Big|\{H_{N+1}^{\KdV}+&4\tau^2 \T_N^\KdV(i\tau,u),  \T_{-1}^{\KdV}(z)\}\Big|\\
    &= \Big|\int \frac{\delta \T^\KdV_{N+1}(i\tau,u)}{\delta u}\partial \frac{\delta  \T^\KdV_{-1}(z)}{\delta u}\, dx\Big|\\
    &\leq \Big\|\frac{\delta \T^\KdV_{N+1}(i\tau,u)}{\delta u}\Big\|_{H^{-N-2}}\Big\|\frac{\delta T^\KdV_{-1}(i\tau,u) }{\delta u} \Big\|_{H^{N+3}}\\
    &\lesssim \tau^{-2}\big(\|u\|_{H^{N+2}}+\|u\|_{H^{N+2}}^2\big)\|u\|_{H^{N+1}},
\end{align*}
    if $\tau \gtrsim \|u\|_{H^{-1}}^2$, using \eqref{eq:tauNdec}.
We claim that (for $ \tau $ large)
\[   \{(2i\tau)^2 \T_N^\KdV(i\tau,u),  T^\KdV_{-1}(z)\} = 0.   \]
The case $N=-1$ is a consequence of  Lemma \ref{kdvpoisson}.
We proceed by induction. 
Suppose the  inequality holds for $N-1$. Then by \eqref{eq:limitpoisson}
\[ \{ H_N^\KdV, T^\KdV_{-1} \}=0. \]
This implies the claim since $ \T^\KdV_N(i\tau, u) $ is a linear combination of $ \T^\KdV_{-1}$ and $H^\KdV_n$, $ n \le  N$. 
\end{proof}

From Lemma \ref{lem:poissonbrackets} we obtain the evolution of the variables $v(z)$ and $w(z)$ when $u$ evolves according to the $N$th KdV flow, respectively the flows $\T_N$.

\begin{theorem}\label{thm:HNflows}  The following identities hold on Schwartz space 
\[ \{ w(i\tau), H^{\KdV}_N(u) \}= \partial_x \frac{\delta}{\delta w} H_N^{\Gardner}(w,\tau)\]
\[\{ w(i\tau), \T^{\KdV}_N(z,u)\}= \partial_x \frac{\delta}{\delta w} \T^{\Gardner}_N(z,w,\tau)  \] 
\[ \{ v(z), H^{\KdV}_N(u)\} = 2\partial_x \Big[ ( v(z)+1) \frac{\delta}{\delta u} \sum_{j=-1}^{N-1} (2z)^{2(N-1-j)} H^{\KdV}_j(u)
\Big] \]
\[ \{ v(z_1), \T^{\KdV}_N(z_2,u)\} = 2\partial_x \Big[ ( v(z_1)+1) \frac{\delta}{\delta u} \sum_{j=-1}^{N-1}(2z_1)^{2(N-1-j)}  \T^{\KdV}_j(z_2,u)  
\Big] \]
\end{theorem} 

\begin{proof} 
We expand \eqref{eq:poissontau}, \eqref{eq:poissonv} and \eqref{eq:poissonw}
and compare coefficients of $ (2z)^{-3-2n}$. The first two lines follow immediately by definition. For the third and fourth line on the right hand side we write
\[
    \partial_x\frac{v(z_1)-v(z_2)}{v(z_2)+1} = 2\partial_x\Big((v(z_1)+1)\Big(\frac12-\frac{\delta \T_{-1}^{\KdV}(u)}{\delta u}\Big)\Big)
\]
and expand both this and $(4z_1^2-4z_2^2)^{-1}$.
\end{proof}

\subsection{Structure of the Hamiltonians}

The Lenard recursion formula (see \cite{MR2173592} for how it is connected to Lenard) 
\begin{equation}\label{eq:Lenardrecursion}
    \partial \frac{\delta}{\delta u}H^{\KdV}_{N+1}= (-\partial^3+4u\partial + 2u_x)\frac{\delta}{\delta u}H^{\KdV}_{N}
\end{equation}
holds in a  distributional sense  for $ u \in H^{N}$
by \eqref{eq:lenardrecursion},\eqref{eq:asymptotickdv}  
and \eqref{eq:tauNdec}. We will see that the Hamiltonians have a very special structure. For that we introduce the notion of differential polynomials. 

\begin{definition}\label{def:diffpol} A differential polynomial in  $u$ is a polynomial in $u$ and its derivatives.  We say the monomial 
\[  \prod_{j=0}^N (u^{(j)})^{\alpha_j}  \] 
has
\begin{itemize}
    \item homogeneity (total number of factors) $H = \sum_{j=0}^N \alpha_j $,
    \item weight (total number of derivatives) $M=\sum_{j=1}^N j \alpha_j$ and
    \item for KdV, degree $d^{\KdV} = H+M/2$,
    \item for Gardner, degree $d^{\Gardner} = H + M$.
\end{itemize}
We say a differential polynomial has degree $n$ if it is a sum of monomials of degree $n$.
\end{definition}

In Chapter \ref{sec:diffhamiltonian} (see Lemma \ref{lem:tracestructure})
 we prove that the Gardner  Hamiltonians are differential polynomials. By the chain rule we have
    \[
     \partial \frac{\delta H_N^{\KdV}}{\delta u} (u) = (\partial + 2\tau + 2w)\partial\frac{\delta H_N^{\Gardner}}{\delta w}(w,\tau)
    \]
    where $u = w_x + 2\tau w + w^2$. On the left hand side every factor of $\tau$ carries a factor of $w$, hence the homogeneity $\tau$ has to be less or equal than the homogeneity of $w$ in each monomial in $H_N^{\Gardner}(w,\tau)$, since the variational derivative decreases the 
    homogeneity in $w$ by $1$ which can be compensated by the multiplication by $w$. 
   Similarly for each monomial of the integrand on the left hand side of \eqref{eq:GardnerKdVr} the homogeneity in $w$ is at least as high as the homogeneity in $ \tau$. The same is true 
    for the first term on the RHS, hence $H_{N-1}^{\Gardner}$ can be written as integral over multiples of monomials with the homogeneity in $u$ is at least two more than the homogenity 
    in $w$.  The  formulas also show that the sum of the homogeneities in $ \tau $ and  $w$ is always even. 
 
 Together with \eqref{eq:GardnerKdVr} this allows us to  obtain the KdV Hamiltonians directly from the Gardner Hamiltonians 
 by picking the monomials where the homogenity of $w$ is exactly two more than the homogenity of $\tau$, which is positive. 
\[
\begin{split}H^{\KdV}_{N}(u) &= \lim_{\tau\to \infty} H_{N}^{\KdV} ((u/(2\tau))_x + (u/(2\tau))^2 + 2\tau(u /(2\tau))) \\&=  \lim_{\tau \to \infty } 4\tau^2 H_{N}^{\Gardner}(u/(2\tau),\tau)  
\end{split}
\]
We obtain a new recursion formula for the Gardner Hamiltonians
\begin{equation} \label{eq:recursionGardner}   H^{\Gardner}_{N}(w,\tau) = \lim_{\lambda \to \infty } \lambda^2 H_{N-1}^{\Gardner}((w_x+2\tau w + w^2)/\lambda,\lambda)-4\tau^2 H_{N-1}^{\Gardner}(w,\tau).  \end{equation} 
as well as a recursion formula for the KdV Hamiltonians which does not involve taking antiderivatives,
    \begin{align*}
        H^{\KdV}_{N}(u) &= \lim_{\tau \to \infty } 4\tau^2 H_{N}^{\Gardner}\Big(\frac{u}{2\tau},\tau\Big)\\
         &= \lim_{\tau \to \infty } 4\tau^2\Big[H_{N-1}^{\KdV}\Big(\frac{u_x}{2\tau} + u + \frac{u^2}{4\tau^2}\Big) - 4\tau^2 H_{N-2}^{\KdV}\Big(\frac{u_x}{2\tau} + u + \frac{u^2}{4\tau^2}\Big) + \dots\Big].
    \end{align*}

Starting with $H^{\Gardner}_0 = \frac12 \int w^ 2 dx  $ we obtain $H_0^{\KdV} = \frac12 \int u^2 dx$  and 
\[ H_1^{\Gardner} = \frac12  \int  (w_x +2\tau w +w^2)^2 - 4\tau^2 w^2 dx =  \frac12 \int w_x^ 2 +  w^4 + 4\tau w^3 dx.   \]
\[ H_1^{\KdV} =   \frac12  \int u_x^2 + 2 u^3 dx  \]
 \[ 
 \begin{split} 
 H_2^{\Gardner} \, & =   \frac12 \int ( w_{xx} + 2 \tau w_x + 2w w_x)^2 + 2 ( w_x+2\tau w + w^2)^3 
 - 4 \tau^2 (w_x^2 + w^4 + 4 \tau w^3) dx 
 \\ & = \frac12 \int w_{xx}^2 + 10 w^2 w_x^2+   20 \tau w w_x^2     + 2 w^6+ 12\tau w^5  + 24 \tau^2 w^4 dx 
 \end{split} 
 \]
\[ H_2^{\KdV} =   \frac12 \int u_{xx}^2 + 10 u u_x^2   + 6 u^4 dx. \]

We arrive at the following general structure of the Hamiltonians.

\begin{theorem}\label{thm:formofkdv}
A) $H_n^{\KdV}$ can be written as an integral over a sum of  homogeneous differential polynomials  in $u$, 
\[ H_n^{\KdV} = \int e_n dx = \frac12 \int |u^{(n)}|^2 + \sum_{k=3}^{n+2} \int e_{n,k} \, dx  \] 
where 
\[ e_{n,n+2} =    \frac1{n+2}   \binom{2n+2}{n+1}  u^{n+2}. \] 
The degree of $e_n$ is $d^{\KdV} = n+2$. Hence, $e_{n,k}$ is a  sum of products of $k$ factors, each product carrying a total of $2(n+2-k)$ derivatives of order at most $n+2-k$. \newline
B) $H_n^{\Gardner}$  can we written as an integral over a linear combination of differential monomials
\[ H^{\Gardner}_n(\tau,w) =  \frac12 \int |w^{(n)} |^2 dx   + \sum_{m=0}^{n+1}\sum_{k=3}^{2n+2-m}  \int \tau^m  e^{\Gardner}_{n,m,k} dx \]
Here $e^{\Gardner}_{n,m,k}$ are differential polynomials of homogeneity $k \geq 3$ in $w$, of degree $d^{\Gardner} = 2n+2-m$ and a with total number of $2n+2-k-m$ derivatives. No factor contains more than $[n+1 -(k+m)/2]$ derivatives. The term of highest homogeneity is
    \[ e^{\Gardner}_{n,0,2n+2} =  \frac1{2(2n+1)} \binom{2n+2}{n+1}   w^{2n+2}. \]
The homogeneity $k$ of $w$ and homogeneity $m$ of $\tau$ in $e_{n,m,k}$ are related by $k \geq m+2$.
\end{theorem}
\begin{proof}
    A) Since the functional derivative reduces the number of factors by one but keeps the weight the same, we see that the first statement is equivalent to showing that
    \begin{equation*}
        \frac{\delta H_n^{\KdV}}{\delta u} = (-1)^n u^{(2n)} + \dots + \frac{1}{2} \binom{2n+2}{n+1} u^{n+1},
    \end{equation*}
    where the terms in between contain $k$  factors, $1<k<n+1$, and are all of degree $n+1$. This is done inductively. Clearly, for $n=0$ the statement holds. On the other hand, from
    the Lenard recursion, 
    \begin{equation*}
        \partial \frac{\delta H_{n+1}^{\KdV}}{\delta u}= (-\partial^3 + 4u \partial + 2 u_x) \frac{\delta H_n^{\KdV}}{\delta u}
    \end{equation*}
    and the fact that derivatives keep the number of factors the same while they add one derivative, we conclude that the linear term gains two derivatives and the $(n+1)$-linear term becomes
    \begin{equation*}
        \frac{1}{2} \binom{2n+2}{n+1} \partial^{-1}(4u \partial + 2 u_x) u^{n+1} = \frac{1}{2} \binom{2n+2}{n+1}\frac{4n+6}{n+2}u^{n+2} = \frac{1}{2} \binom{2n+4}{n+2}u^{n+2},
    \end{equation*}
    by explicitly calculating the binomial coefficient. The degree of the differential polynomial is increased by the operators $-\partial^2, \partial^{-1}(4u\partial + 2u_x)$ by one, which determines the monomials in between. The bound on the maximal order of derivatives involved in the monomials can then be reached by a finite time of partial integrations.

    B) By \eqref{eq:generalgardnerhamiltonian} we have to analyze what happens when we plug in $u = w_x + 2\tau w + w^2$ into the $N$th KdV Hamiltonians with $N \leq n-1$. We first prove the degree condition. Consider a monomial in the $N$th KdV Hamiltonian. By splitting the sum in $u = w_x + 2\tau w + w^2$, every factor of $u$ becomes a factor of $w$ and gains either a derivative, another factor of $w$ or a factor of $\tau$. Additionally a factor of $(2\tau)^{2(n-1-N)}$ is multiplied. Thus by the degree condition of KdV,
    \[
    \begin{split}
        m + d^{\Gardner} &= 2(n-1-N) + H^{\Gardner} + M^{\Gardner} \\
        &= 2(n-1-N) + 2H^{\KdV} + M^{\KdV} \\
        &= 2(n-1-N) + 2(N + 2) = 2n+2.
    \end{split}
    \]
    The bound on the maximal order of derivatives involved in the monomials can again be reached by a finite time of partial integrations. To obtain the form of the bilinear term we notice that only the part $w_x + 2\tau w$ can contribute to it. Now from the $n$th KdV Hamiltonian respectively we obtain a contribution of
    \[
        \frac12\int \big(w^{(n+1)} + 2\tau w^{(n)}\big)^2 \, dx = \frac12\int \big(w^{(n+1)}\big)^2 + 4\tau^2 \big(w^{(n)}\big)^2 \, dx.
    \]
    This shows that the sum of the bilinear parts in \eqref{eq:generalgardnerhamiltonian} is a telescopic sum in which only the highest order term survives, giving the desired form of the bilinear term in the Gardner Hamiltonian. Likewise, the form of the term of highest homogeneity is reached by setting $u = w^2$ in the $u^{n+1}$ summand of the $n-1$th KdV Hamiltonian. 
   
\end{proof}

Before analyzing the form of the good variable equations we state another result concerning the conservation law for the momentum of the Gardner equations. This will be used in proving the local smoothing properties.

\begin{lemma}[Energy-flux] \label{lem:energyflux} There exist differential polynomials which we call fluxes so that 
\[ \frac{\partial}{\partial t} w^2 = \partial_x \flux_N \] 
if $w$ satisfies the $N$th Gardner equation. 
The flux can be written as 
\[ \flux_N =  \sum_{m=0}^{N+1}\sum_{j=3}^{2N+2-m} \tau_0^m \flux_{m,j,N},\] 
where each $F_{m,j,N}$ has homogeneity $j \geq 3$, degree $2N+2-m$ and weight $2N+2-(m+j)$. 
\end{lemma} 

\begin{proof} 
Let $h_N$ be the density of $H_N^\Gardner$, given as a differential polynomial. Then 
\[ 
\begin{split} 
\partial_x w \frac{\delta}{\delta w} H_N^{\Gardner}\,  \, & = \partial_x w \sum_{j} (-1)^j \partial^j \frac{\partial h_N}{\partial  w^{(j)} } \\ 
\, & = \partial_x  h_N + \sum_{j} (-1)^j\sum_{k_1+k_2=j, k_1 \ge 1}       \partial^{k_1} \Big( w'' \partial^{k_2}\frac{\partial h_N}{\partial  w^{(j)} }  \Big).  
\end{split} 
\] 
We obtain a differential polynomial $\flux_N$ such that 
\begin{equation}  \partial_t w^2 = \partial_x \flux_N. \end{equation}  

We decompose  $\flux_N$ as a finite sum of differential polynomials  $\flux_N = \sum_l \flux_N^l$ where $\flux_N^l$ has degree $l$. Taking derivatives leaves the number of factors invariant and increases the degree by one. On the other hand, we know that
\begin{equation*}
    w\partial \frac{\delta {H}^\Gardner_N}{\delta w} =  \partial_x \sum_{m=0}^{N+1}\sum_{j=3}^{2N+2-m} \tau_0^m\tilde h_{m,j,N}
\end{equation*}
is also a sum over homogeneous differential polynomials of homogeneity $j$ and degree $2N+3-m$. Indeed, the functional derivative leaves the number of derivatives constant while decreasing the homogeneity by one, the derivative increases the degree and the former homogeneity is restored by multiplication by $w$. By identifying terms of like degree we see that
\begin{equation*}
    \flux_N = \sum_{m=0}^{N+1}\sum_{j=3}^{2N+2-m} \tau_0^m\tilde \flux_{m,j,N},
\end{equation*}
where each $\flux_{m,j,N}$ has homogeneity $j \geq 3$, degree $2N+2-m$ and weight $2N+2-(m+j)$. This shows
\begin{equation*}
    \int (1-\tanh(\kappa x))w\partial\frac{\delta H^\Gardner_N}{\delta w}dx = \kappa \int \sech^2(\kappa x) \sum_{m=0}^{N+1}\sum_{j=3}^{2N+2-m} \tau_0^m\tilde \flux_{m,j,N}dx .
\end{equation*}
By partial integration, we can reduce to the situation where each factor has less than $[n+1-(m+j)/2]$ derivatives, by making errors where derivatives fall onto the localization factor $\sech(\kappa x)$ and the degree of the homogenous polynomial is decreased, thus being easier to handle.
\end{proof}

We turn to the good variable equations. We have seen that if $u$ solves the $N$th KdV equation, then by Theorem \ref{thm:HNflows}
\[ \partial_t v(z) = 2\partial_x \Big[ ( v(z)+1) \frac{\delta}{\delta u} \sum_{j=0}^{N-1} (2z)^{2(N-1-j)} H^{\KdV}_j(u)
\Big]
\]
Using the relations between $u, w$, and $v$, we can turn this into a single differential equation. By combining \eqref{eq:riccati} and \eqref{eq:functvw} we see that
\begin{equation}\label{eq:relationuv}
    u = -\frac12\frac{v_{xx}}{v+1} + \frac34\frac{v_x^2}{(v+1)^2} - 2z^2 v - z^2 v^2.
\end{equation}
To turn the above system into a single ODE we plug \eqref{eq:relationuv} into the equation. For $N = 1$ and $z = i\tau$, we find
\begin{equation}\label{eq:gvN=1}
v_t =  \partial_x\Big[-v_{xx} + 6\tau^2 v^2 + 2\tau^2 v^3 + \frac{3}{2} \frac{v_x^2}{v+1}\Big].
\end{equation}
Hence the time evolution equation for the good variable is a deformation of the Gardner equation! For $N = 2$ we obtain the somewhat lengthy equation
\begin{equation}
\begin{split}
    v_t =&  \partial_x\Big[v_{xxxx} - 7\tau^2v_{xx}v^2 - 4\tau^2v v_x^2 -14\tau^2 v_{xx}v-4\tau^2v_x^2 \\
    &\quad + 6\tau^4v^5 +30\tau^4v^4 +40\tau^4v^3\\
    &+ (v+1)^{-1}\big(-\frac52v_{xx}^2 - 5v_{xxx}v_x+18\tau^2v_x^2v + \frac92\tau^2v_x^2v^2 -6\tau^2v_x^2\big)\\
    &+(v+1)^{-2}\Big(\frac{25}{2}v_{xx}v_x^2\Big) + (v+1)^{-3}\Big(-\frac{45}{8}v_x^4\Big)\Big].
\end{split}
\end{equation}

To prove equivalence of weak solutions we need an understanding of the form of the equation for general $N$. This is given in the next theorem whose proof can be found in Appendix \ref{sec:proofofgveq}.
\begin{theorem}\label{thm:formofgoodvariableequation}    
    The $N$th equation for $v$ can be written in the form $v_t = \partial_x F_N$, where
    \begin{equation}
        F_N = \sum_{n,l,L,d}^{2N-1}(v+1)^{-n}\tau^l f_{N,n,k,d}(v),
    \end{equation}
    where $f_{N,n,k,d}$ has homogeneity $k$ in $v$ and a total number of derivatives $d$, and the sum is restricted by
    \[
        \begin{split}
            0 &\leq n \leq 2N-1, \qquad l + d = 2N, \qquad n+1 \leq k \leq 2N+1,\\
            \#&\{\text{factors of } \omega \text{ with at least 1 derivative}\} \geq n+1 \quad \text{if}\quad n\geq 1.
        \end{split}
    \]
    Moreover, the linear part of the equation is $(-1)^N v^{(2N+1)}$, and $\tau^l f_{N,n,k,d}$ contains no term of the form $v^{n+1} v^{(2N)}$. The number of derivatives $d$, and $l$, are always even.
\end{theorem}

In order to estimate later we also need be able to pull out derivatives. Weak solutions will have regularity $L^\infty H^{N}$, and localized one more derivative, which means that the single factors of $v$ in nonlinear terms are not allowed to carry more than $N$ derivatives. If $N = 2$, the only bad term is
\[
    (v+1)^{-1}(v_{xxx}v_x),
\]
and we can pull out derivatives to rewrite it as
\[
\begin{split}
    (v+1)^{-1}&(\partial_x (v_{xx}v_x) - v_{xx}^2) \\
    &= - (v+1)^{-1}(v_{xx}^2) + \partial_x\big((v+1)^{-1}(v_{xx}v_x)\big) + (v+1)^{-2}(v_{xx}v^2_x).
\end{split}
\]
For general $N$ we could have a terms of the form
\[
    (v+1)^{-n}\prod_{i=1}^k v^{(\alpha_i)},
\]
where $\sum \alpha_i = d$ and some of the $\alpha_i$ are larger than $d/2$. Note that $d$ is always even, because the number of derivatives on each differential monomial of $\delta H^{\KdV}_N / \delta u$ is even, as can be seen from the Lenard recursion \eqref{eq:lenardrecursion}. Without loss of generality assume $\alpha_1 > \dots > \alpha_k$. We pull out one derivative from $\alpha_1$. This produces a total derivative of a monomial with $d-1$ derivatives, and a term where we replace $\alpha_1$ by $\alpha_1-1$ one some other $\alpha_i$ by $\alpha_i + 1$. We can iterate this until $\alpha_1 - \alpha_2 = 0$, or $\alpha_1 - \alpha_2 = 1$. In the latter case we simply pull out another derivative using $2u^{(\alpha_1)}(u^{(\alpha_1-1)})^{p-1} = \partial (u^{(\alpha_1-1)})^p$. We arrive at
\[
\begin{split}
    \prod_{i=1}^k v^{(\alpha_i)} &= \sum_{\alpha_1 + \dots + \alpha_k = d, \alpha_i \leq d/2} c_{\alpha_1,\dots,\alpha_k} \prod_{i=1}^k v^{(\alpha_i)} \\
    &\quad + \sum_{l=1}^{\alpha - d/2}\partial^{l}\sum_{\alpha_1 + \dots + \alpha_k = d-l, \alpha_i \leq (d-l)/2} c_{l,\alpha_1,\dots,\alpha_k} \prod_{i=1}^k v^{(\alpha_i)}.
\end{split}
\]
Now we pull the derivatives in front of the factor $(v+1)^{-n}$ as well. The maximal amount of derivatives we have to pull out is $\alpha - d/2 = N-1$, hence we will never create a term with too many derivatives. Moreover, because the differential polynomial with prefactor $(1+v)^{-n}$ has at least $n+1$ factors with derivatives, the most derivatives a single factor can have is $2N-n$, hence at most $N-n$ derivatives have to be pulled out if $n \geq 1$. This may create up to $N-1$ new factors. 

We arrive at the following form. 

\begin{lemma}\label{lem:formofgoodvariableequation}
    The $N$th equation for $v$ can be written in the form $v_t = \partial_x F_N$, where
    \begin{equation}
        F_N = \sum_{j,l,n,k,d}\partial_x^{j}\big((v+1)^{-n} \tau^l F_{N,j,l,n,k,d}(v)\big),
    \end{equation}
    where $F_{N,j,n,k,d}$ has homogeneity $k$ in $v$, a total number of derivatives $d-j$, and no factor of $v$ carries more than $N$ derivatives. The sum is restricted by the conditions
    \[
    \begin{split}
        &0 \leq n \leq 2N-1, \qquad 0 \leq j \leq N-1, \qquad l + d = 2N-j,\\
        &n+1 \leq k \leq 2N+1 \quad \text{if} \quad n \geq N+1,\\
        &n + 1 \leq k \leq 3N \quad \text{if} \quad 1 \leq n \leq N.
    \end{split}
    \]
    The linear part of the equation is $(-1)^N v^{(2N+1)}$.
\end{lemma}

\subsection{Regularised Fredholm determinants and the Wadati Lax operator}
\label{subsec:determinants}
The importance of these objects is that they characterize $\ln T$ and its functional derivatives. We introduce the resolvents $R_{\pm} = (-iz \pm  \partial)^{-1}$ for $ \im z >0$,
\begin{equation} \label{eq:integralkernelsresolvents}
R_+ f(x) = \int_{-\infty}^x   e^{iz(x-y)} f(y) dy, \quad R_- f(x) = 
\int_x^\infty e^{-iz(x-y)} f(y) dy. 
\end{equation} 
For  $q,r \in L^2$ and $z \in \C$ we define the AKNS Lax operator
\[ L(q,r)= i\left( \begin{matrix} \partial & q \\ -r & -\partial  \end{matrix} \right)  \]
so that 
\[ L(q,r)-z 1 = (L(0,0)-z 1)  \left( 1+ \left( \begin{matrix} 0 & R_- q \\ - R_+ r & 0 \end{matrix} \right) \right).\]
Unfortunately the operator in the bracket is only Hilbert-Schmidt for $ q,r \in L^2$, but not trace class, even for Schwartz functions. For trace class operators $K$  one has the expansion 
\begin{equation}\label{eq:logdet}   \ln \det (1 - K ) =  \sum_{n=1}^\infty \frac1n \trace K^n \end{equation} 
where $ \trace K^n $ is defined for $K$ in the $L^n$ Schatten class. In particular only the first term is problematic for the bracket above. 
 On the other hand, formally at least, this trace should be zero due to the off-diagonal block matrix form of the operator. This motivates the use of the renormalized determinant
\[ {\det}_2 ( 1+ K ) = \det( I +K) \exp( - \trace K) \] 
for trace class functions, which has a unique extension to Hilbert Schmidt operators $K$. We refer to Simon \cite{MR3364494} for details.

The Lax operator $-\partial^2 + u $ without potential can be factorized as 
\[ -\partial^2 - z^2 = (\partial + iz)(-\partial +iz). \]
\begin{lemma}\label{lem:fredholm}
Suppose that $ u \in L^1$. Then 
\[ (\partial +iz)^{-1} u (-i\partial +iz)^{-1} \]
is a trace class operator. Moreover, if $ u \in H^{-1}$ and $\im z $ is sufficiently large then
\begin{equation}\label{eq:transmission}  T^r(z,u) {\det}_2( 1+ (\partial +iz)^{-1} u (-\partial +iz)^{-1}) = 1  \end{equation} 
and
\begin{equation} \label{eq:T1}   iz   \ln {\det}_2 ( 1+ (\partial +iz)^{-1} u (-\partial +iz)^{-1}) = - \ln \T^\KdV_{-1}(z,u).  \end{equation}  
\end{lemma} 
\begin{proof}
We factor  
\[(\partial +iz)^{-1} u (-\partial-iz)^{-1}
= \Big( (\partial + iz)^{-1} |u|^{1/2} \Big) \Big( |u|^{-1/2} u (-\partial+iz)^{-1}\Big) 
\]
and verify that the factors are Hilbert-Schmidt operators.
More precisely let $ f\in L^2$. Then $(\partial+ iz)^{-1} f$ has the integral kernel
\[ k(x,y) = \chi_{x<y}  e^{-iz(x-y)}   f(y) \]
which has the $L^2$ norm  $ (\im z)^{-1/2}\Vert f \Vert_{L^2}$. The same argument applies to $f (-\partial+iz)^{-1}$. 
\[ (\partial+iz)^{-1} u (-\partial+iz)^{-1}= \big[ (\partial+iz)^{-1} \sqrt{|u|}\big] \big[|u|^{-1/2} u (-\partial+iz)^{-1} ]  \] 
is the product of two Hilbert-Schmidt operators and thus of trace class. 
Estimate \eqref{eq:T1} is a consequence of \eqref{eq:transmission} in view of the definition of $ \T^\KdV_{-1} $ after \eqref{eq:asymptotickdv}. Let $u\in L^1$ and 
\[ \alpha(u) := \ln {\det} ( 1+ (\partial +iz)^{-1} u (-\partial +iz)^{-1}). \]
Then $ \alpha(0) = 0$ and, if $u,v \in L^1$
\[ \frac{d}{ds} \alpha( u+sv)\Big|_{s=0}
=   \trace \Big( (-\partial^2 -z^2 + u)^{-1} v \Big) 
= \int G(x,x) v(x) dx \]
where $G(x,x)$ is the diagonal Green's function. The diagonal Green's function can be expressed by the nonrenormalized Jost solutions, 
\[   G(x,y) = -\frac{T(z,u)}{2iz}  \left\{ \begin{array}{rl} \phi_l(x) \phi_r(y) & \text{ if } x < y \\
      \phi_r(x) \phi_l(y) & \text{ if } y < x \end{array}\right. \]
and 
\[ \frac{\delta}{\delta u } \alpha = -\frac{T(z)}{2iz} \phi_l(x) \phi_r(x)    \]
We compute
\[ \trace \Big( (\partial+ iz)^{-1} u (-\partial+ iz)^{-1} \Big)  
= \int_{y<x} \exp(2iz (x-y))     u(x) dx dy = \frac1{2iz} \int u dx \] 
Comparison with \eqref{eq:fdTren} implies \eqref{eq:T1}. 
\end{proof}

The Wadati Lax operator is defined by 
\[ L^{\Wadati}(w,\tau)  = i\left( \begin{matrix} \partial & -w \\ w+2\tau &  -\partial \end{matrix} \right)  \] 
for $w\in L^2$. 
The matrix $ \left(\begin{matrix} -iz & 0 \\ 2\tau & iz \end{matrix}\right) $
has the eigenvalues $ \pm iz $ with corresponding eigenvectors 
$ \left( \begin{matrix} 1 \\ - \frac{iz}{2\tau}     \end{matrix} \right) $ and $\left( \begin{matrix}0 \\  1  \end{matrix} \right) $.
 We consider the case $\im z >0$. The Jost solutions are defined by the normalization (if $w \in L^1$, which we assume for simplicity for the moment)
\[ \lim_{x\to -\infty} e^{izx} \phi_l = \left(   \begin{matrix}    
 1 \\ -\frac{iz}{2\tau}   \end{matrix} \right) \quad \text{ resp. } \lim_{x\to \infty} e^{-izx} \phi_r = \left(   \begin{matrix}  0 \\ 1   \end{matrix} \right).      \]
We define the transmission coefficient by 
$ \big(T^{\Wadati}(z,w)\big)=a^{-1}$ where 
\[ a:= \lim_{x\to \infty} e^{iz x}  (\phi_l(x))^1 \]
which is the same as the Wronskian (which does not depend on $x$)
\[ W( \phi_l,\phi_r)= \phi_l^1\phi_r^2 -\phi_l^2\phi_r^1. \]
\begin{lemma} 
The following identity holds for $ w \in L^2 \cap L^1$ 
\begin{equation}\label{eq:wadatikdv}  T^\KdV(z, w_x+2\tau w + w^2) = T^\Wadati(z,w,\tau). \end{equation}
\end{lemma} 
\begin{proof} 
A straight forward calculation gives
 \begin{equation}\label{eq:identity} \begin{split} \hspace{2cm} & \hspace{-2cm}   \left( \begin{matrix} -\partial^2 + w^2+2\tau w +w' & 0  \\ 0 & -\partial^2 +w^2+2\tau  w -w' \end{matrix} \right)  \\ & = -\frac12  \left( \begin{matrix} 1 & 1 \\ -1 & 1 \end{matrix} \right) \left( \begin{matrix}  -\partial & w \\ -(w+2\tau)& \partial \end{matrix} \right)^2\left( \begin{matrix} 1 & -1 \\ 1 & 1 \end{matrix} \right).\end{split}  \end{equation} 
Let $ \phi_l$ be the left Jost function for the Wadati Lax operator. Then, by this calculation, $ \phi_l^1+ \phi_l^2$ is a multiple of the left Jost solution to  the KdV Lax operator. Thus 
we obtain \eqref{eq:wadatikdv}.
\end{proof} 

Next we related $T^\Wadati(z,w,\tau) $ to a Fredholm determinant  of the Wadati Lax operator.

\begin{lemma}
Suppose that $ w \in L^2$ and $\im z >0$. Then 
\begin{equation}\label{eq:Wadatires}
(L^\Wadati(0,\tau)-z )^{-1}\left( \begin{matrix} 0 & -iw \\ i w &0\end{matrix} \right) 
= 
\Big( \begin{matrix}  0 & -(\partial +iz)^{-1} w \\ (-\partial +iz)^{-1} w  & 
2\tau (-\partial^2 -z^2)^{-1} w 
 \end{matrix} \Big) \end{equation}
is a Hilbert-Schmidt operator and  
\begin{equation} \label{identity:wadati}  \ln {\det}_2 \left( 1+ 
\Big( \begin{matrix}  0 & -(\partial +iz)^{-1} w \\ (-\partial +iz)^{-1} w  & 
2\tau (-\partial^2 -z^2)^{-1} w 
 \end{matrix} \Big) \right)  =  \T^\Wadati_{-1}(z,w, \tau).   \end{equation} 
\end{lemma} 

\begin{proof}
By the proof of Lemma \ref{lem:fredholm}
the components of the resolvent of the Wadati Lax operator \eqref{eq:Wadatires} are Hilbert-Schmidt operators, and hence the resolvent operator is Hilbert-Schmidt. 
 The  spectral equation for the Wadati Lax operator can be written as  the system  
 \begin{equation} \label{eq:spectralD} i\partial_x \phi^1 - iw \phi^2 = z \phi^1 \qquad -i\partial_2 \phi^2 + i(w+2\tau) \phi^1 =z \phi^2.      \end{equation}

We assume $ w \in L^2 \cap L^1 $ and calculate the variational derivative of $\ln a=- \ln T^\Wadati(z,w,\tau) $. 
Let  $ \dot \phi_{l} = \frac{d}{ds} \phi_{l}(w+sv) |_{s=0}$. It satisfies (compare \eqref{eq:spectralD})
\[ \begin{split}  i\partial_x \dot \phi^1_l - iw \dot \phi^2_l -z \dot \phi^1_l& = -iv  \phi_l^2, \\ -i\partial_2 \dot \phi^2_l +i(w+2\tau) \dot \phi^1_l -z \dot\phi^2_l & = iv\phi^1_l .   \end{split}   \]

The forward fundamental solution  is, using the notation
 $ \tilde \phi_{l,r}= (-\phi^2_{l,r},\phi_{l,r}^1)$ 
\[   T^{\Wadati} * \left\{ \begin{array}{ll} \phi_r(x)  \tilde \phi_l(y) - \phi_l(x) \tilde \phi_r(y) & \text{ if } x>y  \\ 0 & \text{ if } x<y \end{array}\right.
\]
hence 
\[ \lim_{x\to \infty} e^{izx} \dot \phi_l^1(x) 
= i T^\Wadati(z,w,\tau) \int ( \phi_l^1(x) 
\phi_r^1(x) - \phi^2_l(x) \phi^2_r(x))  v(x) dx.  
\]
and 
\begin{equation}  \frac{\delta }{\delta w} a = iT^\Wadati (\phi_l^1 \phi_r^1 - \phi_l^2 \phi_r^2). \end{equation}
Next we compute the variational derivative of 
\[ \ln \det_2 \left( 1+ (L^\Wadati(w,\tau)-z)^{-1} \left( \begin{matrix} 0 & i w \\ -i w & 0 \end{matrix} \right)\right) \] 
and verify that it agrees with the variational derivative of $ \ln a$.
For the identity we first pretend that  $(L^\Wadati -z )^{-1} $ is trace class to clarify the argument - which it is not.  Then
\[ \begin{split} \hspace{0.2cm} & \hspace{-0.2cm} \frac{d}{ds}  \ln \det \left[\left(  1 +  \left(\begin{matrix} i\partial -z & 0  \\ 2i\tau & -i\partial -z \end{matrix} \right)^{-1}\left( \begin{matrix} 0 & -i (w+sv) \\  i(w +sv)& 0 \end{matrix}\right)\right)^{-1}   \right]_{s=0}   
\\ &  = \trace    \left[\left(  1 +  \left(\begin{matrix} i\partial -z & 0  \\ 2i\tau & -i\partial -z \end{matrix} \right)^{-1}\left( \begin{matrix} 0 & -i w \\  iw& 0 \end{matrix}\right)\right)^{-1} \left( \begin{matrix} i\partial- z & 0 \\ 2i\tau & -i\partial -z \end{matrix} \right)^{-1} \left( \begin{matrix} 0 & -iv \\ iv & 0 \end{matrix}   \right)    \right] 
\\ & = \trace \left[ \left( \begin{matrix} -i\partial -z & -i w \\ i(\tau+w) & i\partial + z \end{matrix} \right)^{-1} \left( \begin{matrix} 0 & -i v \\ iv & 0 \end{matrix} \right)   \right] 
\\ & = iT^\Wadati \int (\phi_l^1(x) \tilde\phi_l^1(x) - \phi^2_l(x) \phi^2_r(x) ) v(x) dx 
\end{split} 
\]
since the fundamental solution is 
\[   k(x,y) = T^\Wadati(z,w,\tau) \left\{\begin{array}{cl}  \phi_r(x) \tilde \phi_l(y)    & \text{ if } x > y \\ 
 \phi_l(x) \tilde \phi_r(y) ) & \text{ if } x < y \\ 
\end{array} \right. 
\]
We have to justify this calculation. Let $0< \sigma<  1$, $w,v \in L^2$ and 
\[  L^\sigma(z,w,\tau) = \left( \begin{matrix}   (i\sigma \partial -i)    (i\partial-z) & -iw \\ i(w+2\tau) &  (-i\sigma \partial -i)   (-i\partial-z) \end{matrix} \right)  \]
and 
\begin{equation}\label{eq:Lsigma}   (L^\sigma(z,w,\tau))^{-1} =  \left( 1 + (L^\sigma(z,0,\tau))^{-1} \left(\begin{matrix} 0 & -i w \\ iw & 0 \end{matrix} \right) \right) (L^\sigma(z,0,\tau))^{-1} \left( \begin{matrix} 0 & -iw \\ iw & 0 \end{matrix} \right).    \end{equation}
The convolution  kernel of $(i\sigma \partial -i )^{-1}  (i \partial-z)^{-1} w  $ is
\[   k_{\sigma}(x)  =  \left\{\begin{array}{cl} 0  & \text{ if } x \ge 0 \\
 \frac1{1 +i\sigma z }   (e^{ x/\sigma} - e^{-izx})   &    \text{ if } x < 0. \end{array} \right. \]
Thus 
\[   \Vert  k_{\sigma} * w - (i\partial -z)^{-1} \Vert_{HS} \to 0 \qquad \text{ as } \sigma \to 0 \]
and, if $ \im z > \tau $ and $\sigma$ is sufficiently small  
\[  \Vert  (L^\sigma(z,w,\tau))^{-1} - (L^\Wadati(w,\tau)-z)^{-1}  \Vert_{HS} \to 0 \qquad \text{ as } \sigma \to 0  \]
where $HS$ denotes the Hilbert-Schmidt norm.
As a consequence 
\[   {\det}_2 \left[  1+ L^\sigma(z,0,\tau)^{-1} \left( \begin{matrix} 0 & -iw\\ iw & 0  \end{matrix}\right)\right]  \to  {\det}_2\left[  1+ (L^\Wadati(0,\tau)-z)^{-1} \left( \begin{matrix} 0 & -iw \\ iw & 0 \end{matrix}\right)\right]. \]
%and 
%\[\begin{split} \hspace{1cm} & \hspace{-1cm} \trace    \left[\left\{ \left(  1 +  \left( %L^\sigma(z,w,\tau) \right)^{-1}\left( \begin{matrix} 0 & -i w \\  iw& 0 %\end{matrix}\right)\right)^{-1}-1 \right\} \left( L^\sigma(z,0,w) \right)^{-1} \left( %\begin{matrix} 0 & -iv \\ iv & 0 \end{matrix}   \right)    \right] 
%\\ & \to \trace    \left[\left\{ \left(  1 +  \left( L^\Wadati(w,\tau) -z\right)^{-1}\left( %\begin{matrix} 0 & -i w \\  iw& 0 \end{matrix}\right)\right)^{-1}-1 \right\} \left( %L^\Wadati(0,w)-z \right)^{-1} \left( \begin{matrix} 0 & -iv \\ iv & 0 \end{matrix}   \right)    %\right]. 
%\end{split}
%\]
Now assume that $v \in L^2 \cap L^1$.
The operators $(\pm i\sigma \partial -i )^{-1}  (\pm i \partial-z)^{-1} v  $ are then trace class operators and we obtain by the  calculation above 
\[  \begin{split} \hspace{1cm} & \hspace{-1cm}  \frac{d}{dt} \det \left[  1+ L^\sigma(z,0,\tau)^{-1}\left(\begin{matrix} 0 & -i(w+tv) \\ i(w+tv) \end{matrix} \right) \right]\\ &  = -i\int (L^{\sigma} (z,w,\tau))^{-1}_{12}- L^\sigma(z,w,\tau)^{-1}_{21}) v dx.   
\end{split} 
\]
The following lemma completes the proof.\end{proof}  

\begin{lemma}
The following convergence holds in the trace class norm
\[ 
\begin{split} \hspace{1cm} & \hspace{-1cm} 
\lim_{ \sigma\to 0 }\big( \Vert (L^\sigma(z,w,\tau)^{-1} - (L^\Wadati(w,\tau)-z)^{-1} )_{12}  v \Vert_{\trace} \\ &  + \Vert (L^\sigma(z,w,\tau)^{-1} - (L^\Wadati(w,\tau)-z)^{-1} )_{21}  v \Vert_{\trace} \big)
= 0. 
\end{split}
\]
\end{lemma}
\begin{proof} 
We calculate
\[ \begin{split}   L^\sigma(z,w,\tau)^{-1}\left( \begin{matrix} 0 & -iv \\ iv & 0 \end{matrix} \right) \,&   = \left(1+ L^\sigma(z,0,\tau)^{-1}\left( \begin{matrix} 0 & -iw \\ iw & 0 \end{matrix} \right)\right)^{-1} 
L^\sigma(z,0,\tau)^{-1} \left( \begin{matrix} 0 & i v \\ iv & 0 \end{matrix} \right) 
\\ & \hspace{-2cm}  = L^\sigma(z,0,\tau)^{-1}  \left( \begin{matrix} 0 & -i v \\ iv & 0 \end{matrix} \right) 
\\ & \hspace{-1cm} -  L^\sigma(z,w,\tau)^{-1} L^\sigma(z,0,\tau) \left( \begin{matrix} 0 & -iw \\ iw & 0 \end{matrix}\right) L^\sigma(z,0,\tau)^{-1} \left( \begin{matrix} 0 & -iv \\ iv & 0 \end{matrix} \right)   
\end{split}
\]
The second term on the right hand side converges in the class of trace class operators. 
The trace of the first summand on the right hand side is 
\[   \left( \begin{matrix} (D^\sigma_+)^{-1}  & 0 \\ -\tau (D^\sigma_+)^{-1} (D^\sigma_-)^{-1} & (D^\sigma_-)^{-1} \end{matrix}\right)_{21} v  =  -\tau   (D^\sigma_+)^{-1} (D^\sigma_-)^{-1} v.      \]
This operator converges in the trace class sense. 
\end{proof}

This connection is used in Chapter \ref{sec:differenceflow} where we derive estimates on the multilinear expansion of the Gardner generating function.

\begin{theorem}\label{thm:gardnermiura}
 Suppose that $ w \in L^1$ and $\im z > \tau$. Then 
\[  T^\KdV(z, w_x+2\tau w + w^2) = T^\Wadati(z,w,\tau) \]
and, if $ w \in L^2$,
\[ \ln T^\KdV_r(w_x+2\tau w + w^2) = \ln {\det}_2 \left[  1 + (L^\Wadati(0,\tau)-z)^{-1} \left( \begin{matrix} 0 & -iw \\ iw & 0 \end{matrix} \right) \right] -    \frac1{2iz}   \Vert w \Vert_{L^2}^2 .   \]
The Gardner generating function $\T_{-1}^{\Gardner}(z,w,\tau)$ defined in \eqref{eq:Gardner0} can be written in the following forms
    \begin{equation}\label{eq:GardnerWadati}
\begin{split} 
        \T_{-1}^{\Gardner}(z,w,\tau)\, &  = (4z^2+4
        \tau^2)^{-1} \Big( \T^{\KdV}_{-1}(z,w_x+2\tau w + w^2 ) + \frac12 \int w^2 dx \Big) 
\\ &         = (4z^2+4
        \tau^2)^{-1}  \frac12 \int w^2- w^2(z) dx
    \\ & \hspace{-2cm} =  
    \frac{iz}{4z^2+4\tau^2}\ln {\det}_2 \Big(1+ (L^{\Wadati}(0,\tau)-z )^{-1}\Big( \begin{matrix} 0 & -iw \\ iw & 0 \end{matrix} \Big)\Big)  + \frac{1}{8z^2+8\tau^2}\int w^2dx
    \end{split} 
    \end{equation}
where $w(z)$ in the second line is defined as the unique solution $w(z) \in L^2$ to 
\[
    w(z)_x -2iz w(z) + w(z)^2 = w_x+2\tau w + w^2.
\]
\end{theorem}

\begin{proof}
The first and second equality in \eqref{eq:GardnerWadati} are \eqref{eq:Gardner0} and the definition of $w(z)$. Recall \eqref{eq:Tren}
    \[
        \T_{-1}^{\KdV}(z,u) = iz\big( \ln T^{\KdV}(z,u) - \frac 1{2iz}\int u \big),
    \]
hence setting $u = w_x + 2\tau w + w^2$ in \eqref{eq:Gardner0} gives
\begin{align*}
        \T_{-1}^{\KdV}(z,u) + \frac12 \int w^2 dx  = iz\ln T^{\Wadati}(z,w,\tau) - \tau  \int w \, dx,
    \end{align*}
 from which the proof of the third equality follows.
\end{proof}

\subsection{Periodic functions} 
Here we assume that $r$ and $q$ are $1$ periodic. The monodromy matrix plays a crucial role now: Let $ \Psi$ be the $\R^{2\times 2}$ valued solution to 
\[ L(z) \Psi = 0 \quad \text{ on } [0,1] \qquad  \Psi(0) = 1_{\R^2} \]  
We define the monodromy matrix $M = \Psi(1)$. Its eigenvalues are the Floquet exponents. 
If $ \im z$ is sufficiently large then one of the Floquet exponents is close to $\exp(-iz)$ and its one dimensional eigenspace $N(z)$ is the initial datum for all exponentially growing solutions (we keep $z$ fixed). Transmission coefficient and resolvent (on $\R$) are defined as in the decaying case, this time without normalizing (we could choose a nonzero element of $N(z)$ for the normalization. The transmission coefficient is then independent of this choice). Of course the entries of the resolvent are now $1$ periodic (since it is unique). There is also the resolvent for the periodic Lax operator in the $1$ periodic problem, which coincides (via natural identifications) with the Lax operator on $\R$. 
We can again express the transmission coefficient in terms of the 2 regularized 
Fredholm determinant for the $1$ periodic Lax operator. The relations between transmission coefficients, resolvent and 
$ \alpha$, $ \beta$ and $ \gamma$ remain the same. 
As a rule of thumb all constructions and formulas remain valid in the periodic case, also in the next section, at least for $|\im z| $ large.

\section{Analytic properties of the Miura map and the map \texorpdfstring{$W$}{W}} 
\label{sec:mmiura} 
The Miura map will play a central role in our analysis. As we have seen, it provides a connection between the KdV and the mKdV hierarchy. It also occurs in the factorization of Schrödinger operators. 
\begin{definition} 
 The Miura map  $M : L^2_{loc}(\R) \to  H^{-1}_{loc}(\R), v \mapsto u  $ is defined by  $u = v_x + v^2$. 
\end{definition} 
The basic questions: 
\begin{enumerate} 
\item Can we characterize the range? 
\item If $u$ is in the range, can we characterize the preimage? 
\item Are there interesting classes of functions resp. function spaces so that we can obtain a complete understanding of the mapping properties? 
\end{enumerate}
are nontrivial and interesting.  Formally 
        \[ L:= - \partial_x^2 + u = (\partial_x + v )(-\partial_x + v) \]
iff $u = M(v) $. Since
\begin{equation}\label{eq:positivity}  \int (-\partial_x^2 \psi + u \psi) \psi dx =  \Vert (-\partial_x  + v ) \psi \Vert^2_{L^2}  \end{equation}  
at least for $ \psi \in C^2_c(\R)$ we see that $u$ can only be in the range of the Miura map if $L$ is positive semidefinite. Kappeler, Perry, Shubin and Topalov 
\cite{MR2189502} have shown that, if $L$ is positive semidefinite in this sense, then $u$ is in the range of the Miura map, and the preimage is either a point, or homeomorphic to an interval. 

If $u= \partial_x v +v^2 \in H^{-1} $ then the spectrum contains $[0,\infty)$.  Kappeler, Perry, Shubin and Topalov characterized the range. Checking $u= t \delta_0$ one easily sees that it is in the range if $ t\ge 0$ but not if $ t <0$, hence the range is not open.  
Given $u \in H^{-1}(\R)$, the Schr\"odinger operator is bounded from below and there exists $ \tau_0$ depending on $\|u\|_{H^{-1}}$ so that 
\[  -\partial_x^2+ u + \tau_0^2 : H^1(\R) \to H^{-1}(\R) \] 
is positive definite, in the sense that the associated quadratic form is strictly positive. Let $ \im z > 0$. If in addition $u \in L^1$ there exist unique left and right Jost solutions 
\begin{equation}  -\psi'' + u \psi  = z^2  \psi \label{schroedingereigen} \end{equation}  
\[\begin{split}   \lim_{x\to -\infty}  e^{ izx} \psi_l(x) = 1  \\ \lim_{x\to \infty} e^{-izx} \psi_r(x) = 1.
\end{split} 
\] 
If $ u \in H^{-1}$ and $\im z=\tau >0$  we write it as $u= -2i z  v+ v' $ with 
\begin{equation}  \Vert v \Vert_{L^2}  \le  \Vert u \Vert_{H^{-1}_{2\tau}}, 
\quad \Vert v \Vert_{L^2}= \Vert u \Vert_{H^{-1}_{2\tau}} \quad \text{ if } z = i \tau. 
\end{equation}
We can replace the normalization  by 
\begin{equation} \label{eq:renormJost} 
\begin{split}   \lim_{x\to -\infty}  e^{ izx- \int_0^x v ds } \psi_l(x) = 1  \\ \lim_{x\to \infty} e^{-izx +  \int_0^x v ds} \psi_r(x) = 1
\end{split} 
\end{equation}
which uniquely defines Jost solutions, see Lemma \ref{lem:Jost}.

 If $ z = i \tau$ then $ \psi_l$ and $ \psi_r$ are real valued, and if in addition $ \tau > \tau_0 $
 then $L + \tau^2 $  is positive definite and $ \psi_l$ and $ \psi_r$ are nonnegative (suppose $\psi_l(x_0)=0$ and $ \psi_l>0 $ on $(-\infty,x_0)$. Then we use $ \chi_{(-\infty,x_0)}\psi_l$ as test function and see 
 \[ 0 = \langle \psi_l, (-\partial^2+ u +\tau^2 )\chi_{(-\infty,x_0)} \psi_l  \rangle,  \] 
 and $L+ \tau^2$ would not be positive definite, a contradiction). 
 Moreover 
  \[ \lim_{x\to \infty}  e^{-\tau x-\int_0^x v  } \psi_l(x)  > 0, \qquad \lim_{x\to- \infty}  e^{\tau x+ \int_0^x v} \psi_r(x) > 0. \] 
 Let $\im z>0$ and $ w = \partial_x \ln \psi_l$.
  A short calculation shows that 
 \[ M( w-iz )  = u -z^2.\] 
%\baoping{In fact, take $\psi=exp(\int_0^x v(s)ds$, we see $\psi$ is positive solution, and $v=\psi_x/\psi$} 
Motivated by this calculation we define the modified Miura map 
\[ M_\tau(w) = w_x + w^2 +2 \tau w \] 
on $L^2(\R)$. An inverse is given by 
\begin{equation} \label{eq:miurainverse} \{ u \in H^{-1}: -\partial^2 +  u + \tau^2 \text{ is } p.d.\} \ni    v \to  w= \frac{d}{dx} \ln \psi_l - \tau. \end{equation} 
The modified Miura map 
\[ H^N_\tau \ni w \to u=w_x + 2\tau w + w^2 \in H^{N-1}_\tau \] 
for $ N \ge 0 $ 
is related to a factorization 
\begin{equation}\label{eq:factorization}   - \partial^2 + u + \tau^2 = ( \partial+ w + \tau)(-\partial + w + \tau) \end{equation}  
and, as a consequence $ - \partial^2+ u + \tau^2$ is positive definite if $u$ is in the range of the Miura map and \eqref{eq:positivity} obtains the form 
\begin{equation}  \langle (- \partial^2+ u + \tau^2 ) \phi , \phi \rangle = \Vert (-\partial + w + \tau) \phi \Vert^2_{L^2}. \end{equation}  
Let 
\begin{equation} \label{eq:utau}  \mathcal{U}_\tau = \{ u \in H^{-1}: -\partial^2+ u + \tau^2 > 0 \}.   \end{equation} 
The modified Miura map is an analytic diffeomorphism from $H^N$ to $\mathcal{U_\tau} \cap H^{N-1}$, see Theorem \ref{thm:equivalence}, it relates (weak) solutions of the KdV hierarchy and the Gardner hierarchy (see Theorem \ref{thm:gardnermiura}   for smooth solutions and Theorem \ref{thm:equivalenceweaksolutionsuw}  for weak solutions)  
and precompactness, equicontinuity and tightness for $u$ and $w$ (see Subsection \ref{subsec:precompactness}). The next key lemma
characterizes the range of the Miura map in a quantitative fashion. 

\begin{lemma} \label{lem:uw} 
 The quantities $ \Vert w \Vert_{L^2}$, $ \Vert w_x + w^2+ 2\tau w \Vert_{H^{-1}_\tau}$
and the ground state energy of $ -\partial^2+ w_x+ w^2+ 2 \tau w  + \tau^2$ are related as follows: 
\begin{equation} \label{eq:equivalent1}
\Vert w_x +w^2 + 2 \tau w \Vert_{H^{-1}_\tau} \le  (2 + \tau^{-1/2} \Vert w \Vert_{L^2}) \Vert w \Vert_{L^2}.
\end{equation} 
For all $ \psi \in H^1$
\begin{equation}\label{eq:equivalent2} 
\langle (-\partial^2+ w_x + 2 \tau w + w^2+\tau^2) \psi, \psi \rangle \ge   \frac{\tau^2}{4}    e^{- \frac{1}{\tau}\Vert w \Vert_{L^2}^2} \Vert \psi \Vert_{L^2}^2,  
\end{equation} 
Let $w \in L^2_{loc}$ and $ -\tau_1^2$ be the infimum of the spectrum of $ -\partial^2+ w_x + w^2 + 2 \tau w$. Then  
there exists an absolute constant $C >0$ so that 
\begin{equation} \label{eq:equivalent3} 
\Vert w \Vert_{L^2} \le C \Big( 1 + \ln( \frac{\tau}{\tau-\tau_1}) \Big)^{1/2} \Vert w_x + 2 \tau w + w^2 \Vert_{H^{-1}_\tau}.
\end{equation}
\end{lemma} 
We will prove Lemma \ref{lem:uw} in Subsection \ref{subsec:uw}

A central element for the diffeomorphism property is the study of the linear equation 
\begin{equation}\label{eq:LW}
    \phi \to L_w \phi : = \phi_x + 2\tau \phi + 2 w \phi=:f .
\end{equation} 
The equation $L_w \phi=f $ can explicitly be solved by 
\[ \phi(x) = \int_{-\infty}^x  \exp\Big(-2\tau(x-y)- 2\int_y^x w dt  \Big) f(y) dy. 
\] 
The  two first order operators on the RHS of \eqref{eq:factorization} can be inverted separately and we can solve $ (-\partial^2+ u + \tau^2 )\phi=f$ by 
\[
\begin{split} 
\phi(x)\, &  =(-\partial+ w+\tau)^{-1} (\partial+w+\tau)^{-1} f 
\\ & = \int_x^{\infty} \exp\Big( -\tau (t-x) -\int_x^t wds\Big) \int_{-\infty}^t  \exp\Big( -\tau(t-y) - \int_y^t w ds' \Big)f(y) dy  dt.
\end{split} 
\] 
The Green's function can be expressed explicitly by \eqref{GreenFunction} 
\[  G(x,y)= \int_{\max\{y,x\}}^{\infty}   \exp\Big( \tau(x+y-2t) - 2\int_{\max\{x,y\}}^t  w ds- \int_{\min\{x,y\}}^{\max\{x,y\}} w ds ) \Big)       dt 
\] 
and  the diagonal Green's function by 
\[ \beta(x) := G(x,x) = \int_x^{\infty} \exp\Big( 2\tau(x-t) - 2 \int_x^t w ds \Big) dt= (-\partial+2\tau+2w)^{-1}(1), \] 
resp. \eqref{eq:greensfct}
\[ 
 -\beta' + 2 \tau \beta +  2w \beta  = 1.
\] 
We define the 'good variable' 
\begin{equation} \label{eq:defv}  v := \frac{1}{2\tau \beta }-  1, \end{equation}   
and using \eqref{eq:greensfct},  we calculate 
\[\partial_x \ln(1+v) = -\frac{\beta'}{\beta}= 2\tau v-2w,\]
%\[\begin{split}  \partial_x ( \ln( 1+ v))\, &  = \frac{v_x}{1+v} \\ & = \frac1{1+v} \Big( -\frac{\beta'} {2\tau \beta^2}\Big)\\ &  = \frac1{1+v} \frac1{2\tau  \beta} 
%\Big(2\tau ( \frac1{2\tau \beta} - 1) -   2w    \Big)\\ &  = 2\tau v -   2w, \end{split}  \] 
and record the final simple formula 
\begin{equation}\label{eq:goodvariables}   2\tau v -  \partial_x \ln( 1+ v  ) =  2w. \end{equation}

We define for $ s >-\frac12$
\begin{equation} \label{eq:Vs}  \mathcal{V}^s = \{  v \in H^{s+1}: v > -1 \}. \end{equation} 
Then
\[ \mathcal{V}^{s} \ni    v \to   \tau v - \frac12 \partial_x \ln(1+v) =: w\in H^s\] is a diffeomorphism (see Theorem \ref{thm:equivalence}) which relates weak solutions to the Gardner hierarchy to weak solution for the good variables of Theorem \ref{thm:formofgoodvariableequation}, and preserves precompactness, tightness and equicontinuity (see Subsection \ref{subsec:precompactness}).

Again linear first order equations resp. operators are central objects: The linearization of \eqref{eq:goodvariables} leads to 
\[  \partial_x ( v/\phi) + 2\tau v = f \] 
as well as the equivalent formulation 
\[  \partial_x \psi + 2 \tau ( \psi \phi ) = f. \]

\begin{remark}\label{rem:complexz} 
As discussed above we can define 
\[   w(z) : = \partial_x \ln \phi_l + i z \] 
for $ \im z > 0 $ and either $ \real z \ne 0 $ or if $ \im z$
so large enough so that 
$-\partial^2+ u - (\im z)^2 $ is positive definite. 
We can even allow complex valued potentials $u$ for which we assume 
\[   \Vert   e^{2i\real z  z x} u  \Vert_{H^{-1}_{2\im z}} <\frac14  \] 
Then, since 
\[  w_x -2iz w + w^2 = u ,\] 
by the triangle inequality and the Sobolev inequality  
\[\begin{split}  \left| \Vert e^{2i\real z x} u  \Vert_{H^{-1}_{2\im z}} 
- \Vert w \Vert_{L^2}\right|\, &  \le \Vert e^{2i\real z x} w^2  \Vert_{H^{-1}_{2\im z }}\\ &  \le (\im z)^{-1/2}  \Vert w^2 \Vert_{L^1}\\ & = (\im z)^{-1/2} \Vert w \Vert_{L^2}^2 \end{split}      \] 
the modified complex Miura map defines a diffeomorphism 
\[ M_z : B^{L^2}_{1/4}(0) \to U \subset H^{-1}\] 
with 
\[ \Big\{ u \in H^{-1}:  \Vert e^{2izx} w \Vert_{H^{-1}_{2\im z}} < \frac14     \Big\}  \subset    U \subset \Big\{ u \in H^{-1}:  \Vert e^{2izx} w \Vert_{H^{-1}_{2\im z}} < \frac34     \Big\}. \] 
It is not difficult to obtain the analogous properties for higher regularity. 
However the proof of last part  of Lemma \ref{lem:uw} below fails: Suppose $z$ is not purely imaginary. Then  $w=-iz (\tanh(-iz x- \zeta)-1) $ for $ \zeta \in \C$ satisfies $ w_x -2iz w + w^2=0$. It is not uniformly bounded on intervals of length $1$ and Claim 2 below fails. We do not know whether there is a bound of $\Vert w(z) \Vert_{L^2} $ in terms of $ \Vert u \Vert_{L^2}$ and the distance of $z$ to the spectrum. 
\end{remark} 

This section is organized as follows. In Subsection \ref{sub:Jost} we define Jost solutions, in  
Subsetion \ref{subsec:uw} we prove Lemma \ref{lem:uw}. 
The diffeomorphism property is made more precise and proven in Subsection \ref{subsec:diffeo}. Finally we sketch the analogous statements for the relation between the good variable hierarchy and the Gardner hierarchy in Subsection \ref{subsec:goodvariable}. 

\subsection{Jost solutions} \label{sub:Jost} 
In this section we define Jost solutions and study some of their properties.

\begin{lemma} \label{lem:Jost} 
There exist unique solutions to \eqref{schroedingereigen} with the
normalization \eqref{eq:renormJost}.  
The left and the right   Jost solutions are linearly dependent iff $z^2$ is an eigenvalue. 
\end{lemma} 
\begin{proof}
Using $ u = -2iz v+ v_x $
we write 
\[ -\partial_{x}^2 -2iz  v + v_x - z^2 = 
(\partial-iz  +v) (-\partial  -iz  + v) -  v^2.   
\]
Let 
\[ \phi_1= e^{izx  - \int_0^x v} \psi,  \phi_2 = e^{izx -\int_0^x v} (-\partial-iz+v ) \psi.\]
Then
\[ \phi' = \left( \begin{matrix}  0   & -1 \\ 
                                    v^2  & 2iz-2v 
\end{matrix} \right) \phi \]
which yields the fixed point identity 
\[ \phi_1(t)=1  -\int_{-\infty}^t \phi_2(y) dy 
= 1- \int_{y <x<t} \exp\Big( 2iz(x-y)  -2 \int_y^x v(s) ds\Big) 
v^2(y) \phi_1(y)dy dx . 
\]
If $\Vert v \Vert_{L^2(-\infty, x_0)} $ is sufficiently small a contraction mapping argument gives existence of a unique solution on $(-\infty,x_0)$. Since smallness can be achieved by chosing $x_0 $, and since we can solve the initial value problem 
we obtain a unique solution. 

It is obvious that $z^2$ is an eigenvalue if the left and the right Jost solutions are linearly dependent.  Suppose that the Jost solutions are  linearly independent. We show that then $z^2$ is not an  eigenvalue.  Equation \eqref{schroedingereigen} is a second order ODE and its space of solutions has dimension $2$. It suffices to prove that there exists a solution which is unbounded on the right. There is an explicit formula in terms of $ \psi_r $,
\[   \psi(x) = \psi_r(x) \int_{x_0}^x (\psi_r(y))^{-2} dy   \]  
where we choose $x_0$ large so that $ \psi_r$ does not vanish on $(x_0,\infty)$.
\end{proof}

It follows from the construction that both functions  $ e^{-\frac1{2iz}\int_0^x v\, ds} \psi_l $ resp. $e^{ \frac1{2iz}\int_0^x v\, ds} \psi_r$  depend analytically on $z$ and $u$.

\subsection{Proof of Lemma \ref{lem:uw}}
\label{subsec:uw}
Surjectivity of
\[ L^2\ni w \to w_x+2\tau w + w^2:=u \in \{ u \in H^{-1}: -\partial^2+ u +\tau^2 > 0\}\] 
follows from Lemma \ref{lem:Jost}. The inverse is given by 
$ w = \partial_x \ln \psi_l - 2\tau $. 
Before beginning with the proof seriously we observe that 
\[ \Vert f_x  +\tau g \Vert_{H^{-1}_\tau} \le \Vert f \Vert_{L^2}+ \Vert g \Vert_{L^2} \]
and by the Fourier transform  
\[ \Vert f \Vert_{L^2} =  \Vert \tau f + f_x \Vert_{H^{-1}_\tau}. \] 
Given an open interval $I$ we define $H^{-1}(I)$ 
as the equivalence classes of distributions in $H^{-1}(\R)$ defining the same distributions on 
$I$ with the norm 
\[ \Vert f \Vert_{H^{-1}_\tau(I)} = \inf \{ \Vert \tilde f \Vert_{H^{-1}_\tau(\R)}: f= \tilde f \text{ on } I \}.  \]
We observe that 
\[ \Vert f_x +\tau g \Vert_{H^{-1}_\tau(I)} \le \Vert f \Vert_{L^2(I)}+ \Vert g \Vert_{L^2(I)} \] 
and every distribution in $H^{-1}(I)$ has a representation of this form, 
\[ \Vert f \Vert_{H^{-1}_\tau} \le  \inf\{ \Vert g \Vert_{L^2(I)}+ \Vert h \Vert_{L^2(I)}:f= \tau g+ h_x\}  \]
and the right hand side is equivalent to
$ \Vert f \Vert_{(H^1_{\tau,0}(I))^*}$, the norm of the dual space of $H^1_{\tau,0}$. Here the index $0$ denotes the space of functions vanishing at the boundary of the interval. 
We fix an extension 
\[ H^{-1}((-\infty,0]) \ni f \to \tilde f \in H^{-1}(\R) \]
so that $\tilde f$ is supported in $(-\infty,1)$ and $\tilde f= f$ if $f$ is supported in $(-\infty,-1)$ so that 
\[ \Vert \tilde f \Vert_{H^{-1}(\R)} \le 2 \Vert f \Vert_{H^{-1}(-\infty,0)} .\]
Given $I=(a,b)$ with $|I|=b-a \ge 2$ we use it to define  an extension $H^{-1}(I)\ni f  \to \tilde f_I \in H^{-1}(\R)$ 
with a uniformly bounded norm so that $\tilde f_I$ is supported 
in $(a-1,b+1)$, $ \tilde f =f $ if $\supp f \subset (a+1,b-1)  $ and $ \Vert \tilde f \Vert_{L^2} \le 2 \Vert f \Vert_{L^2} $, with obvious modifications if $2> b-a \ge 1$. We write 
\[ \tilde f = \tau g + \partial g \]
as above and hence $f = \tau  g + \partial g$ on $(a,b)$.
Then 
\begin{equation} \sup_{I\subset (a,b): |I|=1} \Vert f \Vert_{H^{-1}(I)} \sim  \sup_{I \subset (a,b) : |I|=1} \Vert g \Vert_{L^2(I)}.  \end{equation}

In the same fashion we define Sobolev spaces on open intervals.

\begin{proof} [Proof of Lemma \ref{lem:uw} ]
We begin with some preparations. 
Let $ w \in L^2$. Then
\[ \Vert w_x + 2\tau w + w^2 \Vert_{H^{-1}_\tau} \le 2 \Vert w \Vert_{L^2}
+ \Vert w^2 \Vert_{H^{-1}_\tau}. \]
By duality  and the consequence of the fundamental theorem of calculus 
\[ \Vert f \Vert_{L^\infty}^2 \le \Vert f \Vert_{L^2} \Vert f' \Vert_{L^2}
\le \tau^{-1} \Vert f \Vert^2_{H^1_\tau} 
\] 
we have 
\begin{equation}\label{eq:l1h} 
\Vert f \Vert_{H^{-1}_\tau} \le \tau^{-1/2} \Vert f \Vert_{L^1}.
\end{equation} 
Hence 
\[ \Vert w^2 \Vert_{H^{-1}_\tau} \le \tau^{-1/2} \Vert w^2 \Vert_{L^1}
\le \tau^{-1/2} \Vert w \Vert_{L^2}^2\]
which implies \eqref{eq:equivalent1}.
To prove \eqref{eq:equivalent2}, we recall 
\begin{equation} \label{eq:positive}  \langle (- \partial^2 + w_x +2 \tau w  + w^2+\tau^2) \psi, \psi \rangle = 
\Vert (-\partial +\tau + w ) \psi \Vert_{L^2}^2 \end{equation}  
and  we consider (compare \eqref{eq:LW}) 
\begin{equation}\label{eq:Lw}
     -\psi_x + \tau \psi + w \psi = f.
\end{equation} 
We represent the solution $\psi$ by 
\begin{equation} \label{eq:linpresentation}  \psi(x) = \int_x^{\infty}  \exp\Big( \tau (x-y) - \int_x^y w dt \Big) f(y) dy. \end{equation}  
Denote the integral kernel by $g(x,y)$  (where $g(x,y)=0$ if $ y \le x$) and estimate 
for $ x\le y$
\begin{equation} \label{eq:kernelbound}    \tau(x-y)  -\int_x^y w dt
\le  \tau (x-y) + |x-y|^{1/2} \Vert w \Vert_{L^2}\le \frac12\tau(x-y) + \frac1{2\tau}
\Vert w \Vert_{L^2}^2. \end{equation} 
Then
\[ \max\Big\{ \sup_x \int \exp( g(x,y) ) dy , \sup_y \int \exp(g(x,y)) dx \Big\} \le \exp\Big(\frac1{2\tau}\Vert w \Vert_{L^2}^2 \Big) \frac2\tau \]  
We bound the integral operator \eqref{eq:linpresentation} using Schur's lemma, 
\begin{equation} \label{eq:boundLw}  \Vert  (-\partial + \tau + w)^{-1} \Vert_{L^2\to L^2} \le \frac2\tau \exp\Big( \frac1{2\tau} \Vert w \Vert_{L^2}^2\Big).\end{equation} 
Together with \eqref{eq:positive} we see that  
\begin{equation} -\partial^2+ w_x + 2\tau w + w^2 +   \tau^2 - \frac{\tau^2}4 \exp(-\frac1\tau \Vert w \Vert_{L^2}^2 )    \end{equation} 
is positive semi definite.
%with \friedrich{check constant}
%\begin{equation} \label{eq:deltadef} \delta = \exp\Big(C- \varepsilon %\frac{\|w\|_{L^2}^2}{\|u\|^2_{H^{-1}_\tau}}\Big). \end{equation} 
%where $ \varepsilon$ and $C$ are absolute constants. 
This implies \eqref{eq:equivalent2}. 

\medskip 

Let  $ -\tau^2_1 $ the minimum of the spectrum of $ -\partial^2 + w_x + w^2 + 2\tau w $. We claim that there are absolute constants $C$ and  $\varepsilon > 0$ so that
\begin{equation} \label{eq:spectrum}    \frac{\tau-\tau_1}{\tau} \le  \exp\left(C - \varepsilon \frac{\Vert w \Vert^2_{L^2}}{\Vert w_x +2\tau w + w^2\Vert^2_{H^{-1}_\tau}}    \right). \end{equation}
The estimate \eqref{eq:equivalent3} is an immediate consequence. 

We turn to the proof of \eqref{eq:spectrum} 
 We will construct an approximate eigenfunction $ \phi$ so that 
 \begin{equation}\label{eq:approxeigen}
     \Vert (- \partial_x^2+ w_x +  (w+\tau)^2) \phi \Vert_{L^2} \le \delta^2 \tau^2 \Vert \phi \Vert_{L^2},
 \end{equation}   
 with 
 \[ \delta \le \frac12   \exp\left(C - \varepsilon \frac{\Vert w \Vert^2_{L^2}}{\Vert w_x +2\tau w + w^2\Vert^2_{H^{-1}_\tau}}    \right)    \] 
  which implies \eqref{eq:spectrum}.  
% \[ \tau^2_1 \ge \tau^2(1- \delta^2)\ge \tau^2(1-\delta)^2\]
%and \eqref{eq:boundLw}. \friedrich{From this formulation it looks like $\delta$ is independent of $u,w$. In fact in the proof later get
%\[
%    \delta^2 = \exp\Big(-\frac{\mu^2}{8}\Big(\frac{\|w\|_{L^2}^2}{\|u\|^2_{H^{-1}_\tau}} - c\Big)\Big)
%\]
%which implies the full statement
% }

Recall that with $ \psi_l $ the left Jost function of $ -\partial^2 + w_x+ w^2 + 2\tau w + \tau^2$ one has $ w = \partial_x \ln \psi_l - \tau $. We choose a point $y\in \R$, and $\eta \in C^\infty$, $ \eta(x)=1 $ for $ x \le -1$, $ \eta(x)=0$ for $x\ge 1$ and define
\[ \phi(x) = \eta(\tau (x-y)) \psi_l(x)  \] 
%\baoping{The above discussion is before scaling $\tau=1$, so shouldn't $\eta(\tau(x-y))\psi(x)$?}
The main objective is to choose $y$ so that $ \delta$ above is small if $\|w\|_{L^2}$ is large compared to $\|u\|_{H^{-1}_\tau}$.
By scaling it suffices to consider $ \tau=1$.  The construction depends on the following two claims. 

\noindent{\bf Claim 1: } There exists an absolute constant $C>1$ with the follow property. 
Suppose that for the open interval $J=(a,b) $ (we allow $ a= -\infty$ and $b=\infty$)
\begin{equation} \label{eq:uniformsmall}  \sup_{x\in J=(a,b)}  \Vert   u \Vert_{H^{-1}(x-1,x+1)} \le \frac1{1000}  \end{equation} 
and that $w$ satisfies  
\[ w_x+ 2 w + w^2 =  u  \]
on $(a-1,b+1)$. 
Then  there exists $y \in \R\cup \{ \pm \infty\}$ so that
\begin{equation} \label{eq:connectsol2}
\Vert w - \tanh(\cdot-y) + 1 \Vert_{L^1(a,b)} \le C + 500 \sup_{|I|=1,  I\subset(a-1, b+1)}  \Vert u \Vert_{H^{-1}(I) }(b-a) 
\end{equation} 
and
\begin{equation}\label{eq:connectsol} \Vert w - \tanh(\cdot-y) + 1 \Vert_{L^2(a,b)} \le C  +1000 \Vert u \Vert_{H^{-1}(a-1,b+1)}.  \end{equation}
Here, by an abuse of notation, we denote $ \tanh(\cdot \pm \infty)=\pm 1$.

\noindent{\bf Claim 2: } There exists $C>0$ so that if 
\[ w_x+ 2 w + w^2 = u \]
on $(-2,2)$, then 
\begin{equation} \label{eq:localestp} \Vert w \Vert_{L^2(-1,1)}\le C ( 1+ \Vert   u \Vert_{H^{-1}(-2,2) }) \end{equation}

\bigskip 
We postpone the proof of the claims and deduce  \eqref{eq:spectrum} from the claims by constructing approximate eigenfunctions. On unit size intervals $I$ where  $\Vert u \Vert_{H^{-1}(I)} > \frac{1}{1000} $, we apply the large data estimate \eqref{eq:localestp} and obtain
\[ \Vert w \Vert_{L^2(I)} \le c \Vert u \Vert_{H^{-1}(\tilde I)} \] 
for some enlarged interval $ \tilde I$. Let 
\[ A= \Big\{  x\in \R: \text{ there exists } k \in \mathbb{Z} \text{ with } |x-k| < 3 \text{ and } \Vert u \Vert_{H^{-1} ( k-1,k+1) } > \frac1{1000} \Big\}.
\]
Then $A$ is an open set which can be written as union of of at most $ N= 4 \times 10^6\Vert u \Vert^2_{H^{-1}} $
disjoint open intervals 
$A= \bigcup_j I_j $  since by construction $\Vert u \Vert_{H^{-1}(I_j)} \ge \frac1{1000}$.  Moreover, 
\[ \begin{split} \Vert w \Vert^2_{L^2(A)}
 \, & \le   \frac12 \sum_{k \in A \cap \Z} 
\Vert w \Vert^2_{L^2(k-1,k+1)}
\\  &  \le c  \sum_{k\in A \cap \Z}(1+ \Vert u \Vert^2_{H^{-1}(k-2,k+2)}) 
\\ & \le \tilde c  \sum_{k\in A \cap \Z} \Vert u \Vert^2_{H^{-1}(k-2,k+2)} 
\end{split} 
\] 
where the last inequality holds since for every integer $k \in A$ there is a $k'$ with $|k-k'|< 3$ so that the norm of $u$ is large. By the same reason we may drop the first and the last terms in the summation so that 
\[ \Vert w \Vert_{L^2(A)}^2
\le C \sum_{k: (k-2,k+2) \in A} \Vert u \Vert_{H^{-1}(k-2,k+2)}^2 
\le \tilde{C}  \Vert u \Vert_{H^{-1}(A)}^2. \]

If $J=(a-1,b+1)$ is an interval satisfying \eqref{eq:uniformsmall} then by Claim 1 there exists $y$ such that
\begin{equation} \label{eq:estuniform}   \Vert w-(\tanh(x-y) -1)  \Vert_{L^2((a,b))} \le  c(1+ \Vert u  \Vert_{H^{-1}((a-1,b+1))}) \end{equation}  
hence 
\begin{equation} \label{eq:boundw}\Vert w \Vert_{L^2((a,b))} \le c(1+ \Vert u \Vert_{H^{-1}((a-1,b+1))})+ 2\sqrt{(\min\{y,b\}-a)_+} \end{equation}
which gives a bound on  the length of the interval $[a,y]$ where $ w \sim -2$. Note that in the case $a = -\infty$ \eqref{eq:estuniform} holds with $y = -\infty$ and \eqref{eq:boundw} holds without the second term on the right-hand side.

By construction $ \R = A \cup \bigcup_j J_j$
where $J_j$ are the disjoint intervals decomposing the complement of $ \R \backslash A$ which satisfy \eqref{eq:uniformsmall} and hence \eqref{eq:estuniform} for some $ y= y_j$.

We square and sum the estimate over $A$ and the intervals. The sum of the $H^{-1}$ norms on the right hand sides is bounded by $C\Vert u \Vert_{H^{-1}}$ since each constant from Claim 1 comes in a pair with a large $H^{-1}$ norm from Claim 2. Thus
\begin{equation} \label{eq:sumupw}
\Vert w \Vert^2_{L^2} \le c \Vert u \Vert^2_{H^{-1}} + 8\sum_j  \Big( \min \{y_j,b_j\} - a_j \Big)_+. \end{equation}
Now either 
\[\Vert w \Vert^2_{L^2} \le c \Vert u \Vert^2_{H^{-1}}\]
which immediately implies \eqref{eq:spectrum} and  \eqref{eq:equivalent3}, or
at least for one $j$,
\begin{equation} \label{eq:lowerbound}   (\min \{y_j,b_j\} - a_j)_+ \ge  \frac{1}{8\times 10^6\times\Vert u \Vert_{H^{-1}}^2}   \left( \Vert w \Vert_{L^2}^2 - c \Vert u \Vert_{H^{-1}}^2\right). \end{equation}   
We fix this $j$ in the sequel.

Write $J = (a,b)$ and without loss of generality assume $a < y < b$, otherwise we still get \eqref{eq:spectrum} immediately. Then $y - a$ is bounded from below by \eqref{eq:lowerbound}. With
$ \eta =1 $ on $ (-\infty,-1)$
supported in $(-\infty,1) $
and $ \eta' \in C^\infty_c $ nonnegative we define 
$\phi(x) = \eta(x-(y-2))\psi_l(x)$. We want to estimate
\begin{equation} \label{eq:lowerest} \begin{split} 
    \|(-\partial^2 + u + 1)(\eta \psi_l)\|_{L^2(\R)}\, &  = \|-\eta'' \psi_l - 2\eta' \psi_l'\|_{L^2} 
    \\ & \le c \Vert \psi_l \Vert_{L^2(y-3,y-1)}
    \\ 
    & \leq C \exp\Big(- \frac{\Vert w \Vert_{L^2}^2}{16\times 10^6 \times\Vert u \Vert_{H^{-1}}^2} \Big)  \Vert \psi_l \Vert_{L^2(a,a+1)}
    \\ 
        & \leq C \exp\Big(- \frac{\Vert w \Vert_{L^2}^2}{16 \times 10^6 \times \Vert u \Vert_{H^{-1}}^2} \Big)  
        \|\eta \psi_l\|_{L^2(\R)}.
\end{split} 
\end{equation}
The first inequality is an energy inequality for $ \psi_l$ on a slightly larger interval. Indeed, since $\psi_l$ solves the Schrödinger equation, it satisfies
\[
    -(\eta')^2 \psi \psi'' + (\eta')^2\psi^2 = - u (\eta')^2\psi^2,
\]
and hence
\[
    \|(\eta' \psi)'\|_{L^2}^2 + \|\eta' \psi\|_{L^2}^2 \leq  \|\eta''\psi\|_{L^2}^2 + c\| u\|_{H^{-1}(y-3,y-1)}\|\eta' \psi\|_{H^1}^2\]
%Estimating the right-hand side by a small constant times $\|\eta' \psi\|_{H^1}^2$ plus an $L^2$ term, and 
Using smallness of $\|u\|_{H^{-1}}$ on $J$, shows the first estimate of \eqref{eq:lowerest}. 

%Now
%\[
%    \|\eta'' \psi_l %\|_{L^2(\R)} \lesssim %\|\psi_l\|_{L^2(y-3,y-1)}.
%\]
%On $(y-3,y-1)$ we have up to %a constant
%\[
%\begin{split}
%    \psi_l(x) &= %\exp(\int_a^x w(s)+1\, ds) \\
%    &= \exp(\int_a^x w(s)- %\tanh(s-y)+1\, ds + \int_a^x %\tanh(s-y))\\
%    &\leq %\exp(c_0(1+\|u\|_{H^{-1}})(x-%a)^{1/2} - c_1 (x-a))\\
%    &\leq %\exp(c_2(1+\|u\|_{H^{-1}}^2) %- c_3 (x-a))
%\end{split}
%\]
%In particular
%\[
%    \|\psi_l\|_{L^2(y-3,y-1)} %\lesssim %\exp(c_2(1+\|u\|_{H^{-1}}^2) %- c_3 \sigma).
%\]
%By writing
%\[
%    \psi_l'(x) = %(w(x)+1)\exp(\int_a^x %w(s)+1\, ds),
%\]
%we find the similar bound for %$\eta' \psi_l'$,
%\[
%    \|\eta' %\psi_l'\|_{L^2(\R)} \lesssim %(\|w\|_{L^2} + %1)\exp(c_2(1+\|u\|_{H^{-1}}^2%) - c_3 \sigma).
%\]
For the second inequality, if $a \leq x \leq y-2$ (which we may assume without loss of generality), we use \eqref{eq:connectsol}
\begin{align*}
    \frac{\psi_l(x)}{\psi_l(a)} &= \exp \left(\int_a^x (w(s)+1) ds\right) \\ &=  \exp\left(\int_a^x w(s)- \tanh(s-y)+1\, ds -(x-a)+ \int_a^x \tanh(s-y)+1\,ds\right)\\
    &\leq \exp( \Vert w- \tanh(.-y)+1  \Vert_{L^1(a,x)}  -  (x-a))\\
    &\leq \exp(  C+ (500\sup_{|I|=1, I\subset J} \Vert u \Vert_{H^{-1}(I)} -1)  (x-a)).
\end{align*}
Here we used again that $\|u\|_{H^{-1}}$ is small on $J$. Together with \eqref{eq:lowerbound}, this gives the second inequality of \eqref{eq:lowerest}. 
 
%\[
%    \|\eta \psi_l \|_{L^2(\R)} \gtrsim \|\psi_l\|_{L^2(a,y-2)} \gtrsim \exp(-c_4\|u\|_{H^{-1}}^2)\big(1-\exp(-3(\sigma-2))\big)^{\frac12}.
%\]
%Combining the upper and lower %bounds shows %\eqref{eq:approxeigen} with
%\[
%    \delta^2 = %\frac{(\|w\|_{L^2} + %1)e^{c_5(1+\|u\|_{H^{-1}}^2)}%}{\big(1-e^{(-3(\sigma-2))}\b%ig)^{\frac12}}e^{- c_3 %\sigma}.
%\]
%On this interval %\friedrich{the exponential is %not additive, so the RHS is %not correct}
%\[ \frac{ %\psi_l(y)}{\psi_l(x)}  = %e^{y-x+\int_x^y w(\tau) %d\tau}  \] 
%decays exponentially with %rate $2$ \friedrich{what %exactly means decay here?}. %We argue similarly on all the %intervals $J_j$. 

% The left Jost function is %decaying like $e^{-x}$ on the %interval. We smoothly %truncate in $y_j$ , normalize %so that $ \psi_l(a) = 1$ and %see that 
%\[ \Vert \eta \psi_l %\Vert_{L^2(-\infty, y_j)} \ge %C e^{(1-c \delta) \min\{y_j, %b_j\} -a_j}   \] 
%and 
%\[ \Vert (-\partial^2 + u + %\tau^2 )(\eta\psi_l) \Vert %\le c \exp( - \min\{x_j, %b_j\} - a_j)   \Vert \eta %\psi_l \Vert_{L^2}.   \] 
% Together with \eqref{eq:lowerbound} we arrive at \eqref{eq:spectrum}.   
This completes the proof of Lemma \ref{lem:uw} and it  remains to prove the two claims. Claim 1 relies on Claim 2 which we prove first.

\noindent
{\bf Proof of Claim  2: Large data.} 
We may replace $u$ by (an abuse of notation) $u+v_x$ with $u,v \in L^ 2(I)$ with 
\[ \Vert u \Vert_{L^2(-2,2)}+ \Vert v \Vert_{L^2(-2,2)}\le 2 \Vert u + v_x \Vert_{H^{-1}((-2,2))} .\] 
Suppose that
\begin{equation} w_x + 2 w + w^2= v_x +u \qquad \text{ on } (-2,2). \label{eq:wvu}\end{equation}
We claim 
\begin{equation} \label{eq:wapriori}  \Vert w \Vert_{L^2(-1,1)} \le c \Big(1 + \Vert v \Vert_{L^2(-2,2)} + \Vert u \Vert_{L^2(-2,2)}\Big). \end{equation} 
which implies 
\[ \Vert w \Vert_{L^2(-1,1)} \le c \Big(1 + \Vert v_x+u  \Vert_{H^ {-1}(-2,2)} \Big). \] 

We prove \eqref{eq:wapriori}
with several reductions.  The function $ w_1= w+1$ 
satisfies 
\[  \partial_x w_1 + w_1^2 = v_x + u +1 \] 
and, including $1$ into $u$  it suffices to prove the bound \eqref{eq:wapriori} for solutions to \begin{equation}  w_x + w^2 = v_x + u. \end{equation}  
which we consider from now on. Since 
\[ \Vert w \Vert_{L^2(-1,1)} \le \Vert w_+ \Vert_{L^2(-1,1)} + \Vert w_- \Vert_{L^2(-1,1)} \] 
it suffices to prove the following estimate for the positive part $w_+$ of $w_1$ 
\begin{equation} \label{eq:basicpos}  \Vert w_+ \Vert_{L^2(0,1)} \le 
 c \Big(1 + \Vert v \Vert_{L^2(-1,1)} + \Vert u \Vert_{L^2(-1,1)}\Big). \end{equation}
We apply this and the corresponding shifted estimate on  $(-2,0)$. The argument for $w_{-}$ being similar by reversing the $x$ direction. In this way \eqref{eq:basicpos} implies 
\[ \Vert w \Vert_{L^2(-1,1)} \le \Vert w_+ \Vert_{L^2(-1,1)}
+ \Vert w_- \Vert_{L^2(-1,1)}
\le c\Big( 1+ \Vert v \Vert_{L^2(-2,2)}
+ \Vert u \Vert_{L^2(-2,2)}\Big).
\] 

We prove the estimate first on smaller intervals, and then deduce 
\eqref{eq:basicpos} from the estimates on the smaller intervals. 
Assume that $0<R\le 1$ and $(-R,R)$ is an interval so that \begin{equation}\label{eq:sizev}    \Vert v \Vert_{L^2(-R,R)} \le (2R)^{-1/2} \end{equation}  so that $ \Vert v \Vert_{L^1(-R,R)} \le 1$
%so that
%\begin{equation}\label{eq:boundv}  \max\{1, \Vert v \Vert_{L^2(-R,R)}\}  = (2R)^{-\frac12} \end{equation} 
%so that 
%\[  \Vert v \Vert_{L^1(-R,R)} %\le 1. \] 
%The case $R=\frac12$ is easier, so we assume that we have 
%$ \Vert v \Vert_{L^2(-R,R)}^2 = (2R)^{-1}$. \friedrich{we should assume an inequality for the gluing part later}
The Ansatz $ w_1 = v+w_2$ leads to 
\[ \partial_x w_2 + w_2^2 + 2 v w_2 = u - v^2 \]  
and $ w_3 = \exp\Big(-2\int_0^x v\Big) w_2 $ satisfies 
\[ \partial_x w_3 +  \exp\Big( 2 \int_0^x v\Big) w_3^2 = e^{-2 \int_0^x v}( u-v^2).  \] 
We set
\[ w_4 = w_3 - \int_0^x  \exp\Big(-2\int_0^y v\Big) (u(y)-v^2(y)) dy \] 
which satisfies with $ \kappa = e^2(\Vert u \Vert_{L^1(-R,R)} + \Vert v \Vert^2_{L^2(-R,R)}) $  
\[ \partial_x w_4 + e^{-2}  (|w_4|- \kappa)^2_+ \le 0.  \]   
Thus $w_4$ cannot have a inner local maximum larger than $ \kappa $ or a inner local minimum less than $-\kappa$.  If $J$ is an interval were $w_4 \ge \max(2 \kappa,4e^2/R) $ then 
\[ \partial_x w_4  \le - \frac1{4e^2} w_4^2. \]
Comparison with the general solution to the equation 
\[ \hat w = \frac{4e^2}{x-c} \]  
shows that  the length of the interval is at most $R$.
As a consequence, arguing by contradiction,
\[ w_4 \le \max\Big\{2 \kappa, \frac{4e^2}{R} \Big\}   = \max\Big\{2 e^2  (\Vert u \Vert_{L^1(-R,R)} + \Vert v\Vert^2_{L^2(-R,R)} ), \frac{4e^2}{R} \Big\}  \quad \text{  on } \quad (0,R)\] 
and, by reversing the sign and the $x$ direction 
\[ w_4 \ge -\max\Big\{2 e^2  (\Vert u \Vert_{L^1(-R,R)} + \Vert v\Vert^2_{L^2(-R,R)} ), \frac{4e^2}{R} \Big\}    \quad \text{  on } \quad  (-R,0).\]

Retracing the construction we see (and taking into account \eqref{eq:sizev}) 
\begin{equation}\label{eq:wsmall}   \Vert w \Vert_{L^2(0,R)}  \lesssim  R^{-1/2}+  R \Vert u \Vert_{L^2(-R,R)} + \Vert v \Vert_{L^2(-R,R)}. 
\end{equation}

We want to show \eqref{eq:basicpos}.
If $\Vert v \Vert_{L^2} \lesssim 1 $ we choose $R=1$ and obtain \eqref{eq:basicpos}. Otherwise we choose $0<x_1<1$ so that $  \Vert v \Vert_{L^2(-x_1,x_1)}= (2R_1)^{-1/2} ,  R_1=x_1$  
and obtain 
\[ \Vert w_+ \Vert_{L^2(0,x_1)}
\lesssim R_1 \Vert u\Vert_{L^2(-x_1,x_1)}+ \Vert v \Vert_{L^2(-x_1,x_1)}. 
\]  
We choose recursively the points 
\[x_0< x_1 < \dots < x_N < 1 \leq x_{N+1} \] 
with $ x_0= -x_1$. 
We will prove 
\begin{equation}\label{eq:elementary}  \Vert w_+ \Vert_{L^2(x_{j-1}, x_j)} \le 
c \Big(\Vert v \Vert_{L^2(x_{j-2}, x_j)}+ \Vert u \Vert_{L^2(x_{j-2}, x_j)} \Big)\end{equation}  
and obtain \eqref{eq:basicpos} by squareing \eqref{eq:elementary} and adding over $j$.
We choose the points so that with 
 $R_j = \frac{x_{j}-x_{j-1}}2 $
\begin{equation} \label{eq:normalization}  2R_j \Vert v \Vert_{L^2(x_{j-1},x_j)}^2 = 1 \end{equation}  
or $ R_j = \frac12$. This latter case is easier and we assume the identity \eqref{eq:normalization} . Then 
\[ \Vert w_+ \Vert_{L^2( \frac{x_{j-1}+x_j}2, x_j)} \lesssim 
\Vert u \Vert_{L^2(x_{j-1},x_j )} + \Vert v \Vert_{L^2(x_{j-1},x_j) }. \] 
The estimate on the left half of the intervals $(x_{j-1},x_{j-1}+ R_j)$ is more delicate and we distinguish the cases $ R_j \le R_{j-1}$ and $ R_j \ge R_{j-1}$. To simplify the notation we consider $j=1$ and assume first that $R_2\le R_1$. Then, by the same argument as above  (with slightly worse constants) 
\[ \Vert w_+ \Vert_{L^2(x_1,x_1+R_1)}\lesssim  \Vert u \Vert_{L^2(0,x_1+R_1)}+ \Vert v \Vert_{L^2(x_0,x_2)}. \]

Now assume that $ R_2> R_1$. As above we estimate on $(0,2x_1) $ and on $(x_1,3x_1)$
together 
\[ \Vert  w_+ \Vert_{L^2(x_1,3x_1)} \le c \big( R_1^{-\frac12} + R_1 
 \Vert u \Vert_{L^2(0,3x_1)} + \Vert v \Vert_{L^2(0,3x_1)}\big) \] 
with slightly worse constant due to using the same terms several times. We repeat the argument 
%on $(x_1,5x_1)$ and on $(3x_1,7x_1)$ with $R= 2R_1$, assuming $x_2 \ge 7x_1$.
%Then 
%\[ \Vert  w_+ \Vert_{L^2(3x_1,7x_1)} \le c \big( (2R_1)^{-\frac12} + 2R_1 
% \Vert u \Vert_{L^2(x_1,7x_1)} + \Vert v %\Vert_{L^2(x_1,7x_1)}\big), \] 
 and, as long as $(2^{j+1}-1)x_1 \le x_1 + R_2$, we control $w_+$ on $((2^{j}-1)x_1,  (2^{j+1}-1)x_1)$ by  repeating the argument on $((2^{j-1}-1)x_1,  (3\times 2^{j-1}-1)x_1)$ and 
 $ ((3\times 2^{j-1}-1)x_1), (2^{j+1}-1)x_1)$ and get
 \[\begin{split} \Vert  w_+ \Vert_{L^2((2^j-1)x_1,(2^{j+1}-1) x_1)} \le \, & c \big( (2^{j-1}R_1)^{-\frac12}+2^{j-1}R_1 
 \Vert u \Vert_{L^2((2^{j-1}-1)x_1, (2^{j+1}-1)x_1)}  \\ & \quad  + \Vert v \Vert_{L^2((2^ {j-1}-1)x_1,(2^{j+1}-1)x_1)}\big). \end{split} \]

Let $J$ be the maximal natural number so that $(2^{J+1}-1)x_1 \le x_1 + R_2$, and, to simplify the notation assume that we have equality.  
Then 
\[ \begin{split} \Vert w_+ \Vert^2_{L^2(x_1,x_2)}
\, & = \sum_{j=1}^J \Vert w_+ \Vert^2_{L^2((2^j-1)x_1, (2^{j+1}-1)x_1)} + \|w_+\|_{L^2(x_1+R_2,x_2)}^2
 \\ & \lesssim R_1^{-1} + R_2^2 \Vert u \Vert^2_{L^2(0,x_2)} + \Vert v  \Vert^2_{L^2(0,x_2)}
 \end{split} 
 \] 
 We can repeat the process to get the general estimate
 \begin{align}
     \|w_+\|_{L^2(x_j, x_{j+1})}^2
     \lesssim  \sum_{k=j}^{j+1}\left( R_k^{-1}+ R_k^2\|u\|_{L^2(x_k, x_{k+1})}^2+\|v\|_{L^2(x_k,x_{k+1})}^2\right).
 \end{align}
 Notice that at most one $R_j=\frac12$, and while $R_j<\frac12$,  $(2R_j)^{-1}= \|v\|_{L^2(x_{j-1}, x_{j})}^2$. We arrive at \eqref{eq:elementary} using $R_k < 1$ and can sum up the estimate to obtain 
 \[ \Vert w_+ \Vert_{L^2(0,1)}\le c ( 1+ \Vert u+ v_x \Vert_{H^{-1}(-1,2)}). \]
 The reason we get the estimate on $(-1,2)$ is that we could have $x_{N+1} > 1$.

\noindent{\bf Proof of  Claim 1: Small data.} 
Again we replace $u$ by $u+ v_x$ with 
\[ \Vert u \Vert_{L^2(a-1,b+1)} + \Vert v \Vert_{L^2(a-1,b+1)} 
\le 2 \Vert u+ v_x \Vert_{H^{-1}(a-1,b+1)} \]
and 
\[ \sup_{I \subset (a-1,b+1), |I|=2} 
\Vert u \Vert_{L^2(I)} + \Vert v \Vert_{L^2(I} 
\le 2\sup_{I \subset (a-1,b+1), |I|=2}  \Vert u+ v_x \Vert_{H^{-1}(I )} \]
It is useful to consider the more symmetric formulation with $\omega= w+ 1-v$ which satisfies 
\begin{equation}\label{eq:omegax}   \omega_x + \omega^2+2v\omega - 1 =  u-v^2, \quad \text{ on } I=(a-1,b+1)  \end{equation}  
and we recall that we may restrict to $ \tau = 1$.
If $v_x+u=0$ it can be solved explicitly and the set of all  global solutions  are given by \[ \omega =  \tanh(  x-y) \]
or $ \omega = \pm 1$.  Suppose that $ \omega(y)=0$ for a point $ y \in I $ and let $ \dot \omega(x) = \omega(x) - \tanh(x-y)$.  It satisfies 
\[
    \partial_x \dot \omega + (\dot \omega + 2\tanh+ 2v )\dot \omega  = u - v^2 - 2v \tanh
\]
hence, if $ x >y$ (the argument for $ x<y$ being similar)  
\[ \dot \omega(x) = \int_y^x  \exp\Big( \int_s^x -2\tanh(\sigma-y) -2v(\sigma) -\dot \omega(\sigma)) d\sigma \Big) ( u(s)-v^2(s) -2 v(s) \tanh(s-y)) ds \]   
 Suppose that $ \Vert \dot \omega \Vert_{L^\infty(y,x)} < \frac12$
and 
\[ \Vert v \Vert_{L^2(I)}+ \Vert u \Vert_{L^2(I)}  < \frac1{1000} \]  
on unit sized intervals $I \subset (a,b)$. 
Then
\[ \int_s^x -2\tanh(\sigma-y) -2v(\sigma) -\dot \omega(\sigma)) d\sigma 
\le 3- |x-s|, 
\] 
hence, by \eqref{eq:uniformsmall} 
\[ | \dot \omega(x) |\le e^3
\sum_{k=0}^\infty  
 e^{-k} \Big(\Vert u \Vert_{L^2((x-k-1,x-k)\cap (a,b) )}+ 3\Vert v \Vert_{L^2((x-k-1,x-k)\cap (a,b) }\Big)\le \frac{8 e^3}{1000} \frac14. \] 
By a continuity argument 
\[ \Vert\dot \omega \Vert_{L^\infty(a,b)} \le \frac14 \] 
and 
\[ | \dot \omega(x) |\le  c \sup_{I\subset (a,b), |I|=1}  ( \Vert u \Vert_{L^2(I)}+ \Vert v \Vert_{L^2(I)}) \]
and again by using the $L^\infty$ bound in the exponential
\[ \Vert \dot \omega \Vert_{L^\infty(x,x+1)}
\le e^4 \sum_{k=0}^\infty e^{-k}  \Big(\Vert u \Vert_{L^2((x-k-1,x-k)\cap (a,b) )}+ 3\Vert v \Vert_{L^2((x-k-1,x-k)\cap (a,b) }\Big),  \] 
we estimate the $L^2$ norm of $ \dot \omega$ on the unit size interval by the $L^\infty$ norm and apply Schur's lemma 
to arrive at 
\[ \Vert   w-\tanh(x-y) \Vert_{L^2(a,b)}
\le \Vert v \Vert_{L^2(a,b)}+ \Vert \dot \omega \Vert_{L^2(a,b)} \le 8e^4  (\Vert v \Vert_{L^2(a,b)}+\Vert u \Vert_{L^2(a,b) }).\]  
The two estimates \eqref{eq:connectsol2} and \eqref{eq:connectsol}
are an immediate consequence. 
We can easily adapt the argument to the case when $ |\omega(y)| \le \frac12 $ at one point. Suppose that $ \omega\ge \frac12 $ on $(a-1,b+1)$, the case $ \omega \le -\frac12 $ being similar.
By Claim 2 $\Vert \omega \Vert_{L^2(I)} \le \Vert w \Vert_{L^2(I)}+1 \le C$ on unit sized intervals in $(a,b)$. Since $\omega$ satisfies equation \eqref{eq:omegax}  
\[ \frac12 \le \omega \le C \quad 
\text{ on } (a,b) \] 
for some universal constant $C$ and hence 
\[ - \frac12 \le w-v \le C \qquad \text{ on } (a,b).\]
Let $ \dot \omega = w-v$. Then 
\[ \partial_x \dot \omega + \dot \omega^2 +2 \dot \omega +2 v \dot \omega    = u - v^2      \]
and 
\[ \dot \omega(x) =  \exp\Big(-\int_a^x 2+  \dot \omega +2v  ds \Big)   \dot \omega(a)
+ \int_a^x \exp\Big( -\int_s^x 2+ \dot \omega+2 v   d\sigma \Big) (u+ v^2) ds \]  
Since $ \sup_I \Vert v \Vert_{L^2(I)} \le \frac1{1000}$ 
\[\int_a^x 2+  \dot \omega +2v ds \ge    (x-a)-\frac1{500}   \]
and, with a small  modification of the previous argument for $a\le x \le b$
\[  \dot \omega(x) \le  2 C e^{-(x-a)} +  500 \sup_I \Big( \Vert u \Vert_{L^2(I)}+ \Vert v \Vert_{L^2(I)} \Big)    \]
and 
\[ \Vert \dot \omega \Vert_{L^2(a,b)} \le C + 500 \Vert u + v_x \Vert_{H^{-1}(a-1,b+1)}. \]
Again \eqref{eq:connectsol2} and \eqref{eq:connectsol}
are an immediate consequence.
\end{proof}

\bigskip 

\subsection{The diffeomorphism \texorpdfstring{$w\to u $}{w to u }}
\label{subsec:diffeo} 
Properties of the modified Miura map are collected in the next proposition. All constants will depend on $ \tau^{-1} \Vert w \Vert^2_{L^2} $ which by Lemma \ref{lem:uw} is equivalent to  having them depend on $ \tau^{-1/2} \Vert u \Vert_{H^{-1}_\tau}$ and the norm of 
$ (\tau^{-2}(-\partial^2 + u) + 1)^{-1} $ as operator on  $L^2$. For simplicity we write $ c( \tau^{-1/2} \Vert w \Vert_{L^2})$.
Given $ \tau$ we call a subset $Q_U \in H^{-1}_\tau$ bounded if there exists $C$ and $ \delta$ so that
\[ \tau^{-1/2} \Vert u \Vert_{H^{-1}_\tau}\le C  \]
and 
\[  \Vert -\phi_{xx} + (u+\tau^2) \phi \Vert_{L^2} \ge (\delta \tau)^2 \Vert \phi \Vert_{L^2}.   \]
By Lemma \ref{lem:uw} a set $Q_W \subset L^2$ is bounded if and only if 
\[ Q_U= \{ w_x + 2\tau w + w^2 : w \in Q_W \}\]
is bounded.
We will use this notation and the notions below. To cover later needs we formulate the next result in larger generality than needed as this point. Let $ W^{n,p}_\tau$ be the standard Sobolev space if $ N \ge 0$ equipped with the norm
\[ \Vert f \Vert_{W^{N,p}_\tau}^p = \sum_{n=0}^N \tau^{p(N-n)} \Vert f^{(n)} \Vert_{L^p}^p \]  
with obvious modifications if $p=\infty$. If $N=-1$ we define 
\[ \Vert f \Vert^p_{W^{-1,p}_\tau(I)} =
\inf \Big\{ \Vert g \Vert^p_{L^p}+ \Vert h \Vert^p_{L^p}: f = g_x + \tau h \Big\}. \]

\begin{definition}\label{def:slowlyvarying} We say that a nonnegative function $ \gamma \in C^\infty( \R) $ is slowly varying of rate $ \alpha$ if 
\[ |\gamma'| \le \alpha \gamma  \] 
and if for every $k$ there exists $c_k$ so that 
\[   |\gamma^{(k)} (x) | \le c_k \alpha^k \gamma(x). \] 
\end{definition} 
Typically examples are $\gamma = e^{\alpha (x-x_0)} $, $ \cosh( \alpha (x-x_0)) $  and $ \sech( \alpha(x-x_0) ) $.

The following lemma provides an alternative description of a weighted version of $W^{-1,p}$.

\begin{lemma} \label{lem:w1infty} 
Let $ \gamma $ be slowly  varying of rate $ \tau/2$. 
Let $ \gamma g \in W^{-1,p}_\tau$. There exists a unique solution $f$ with $ \gamma f  \in L^p$ to 
\[ f' + \tau f = g \] 
which satisfies
\[ \frac25\Vert \gamma  g \Vert_{W^{-1,p}_\tau} \le \Vert \gamma f \Vert_{L^p}  \le 4 \Vert \gamma g \Vert_{W^{-1,p}_\tau  }. \]  
\end{lemma} 
\begin{proof} 
First 
\begin{equation} \label{eq:weightedestimate}  \gamma g = \tau \gamma f- \frac{\gamma'}{\gamma} (\gamma f)  + (\gamma f)'  \end{equation}  
and 
\[ \Vert \gamma g \Vert_{W^{-1,p}_\tau} \le \frac52\Vert \gamma f \Vert_{L^p}. \] 
We decompose according to the definition
$ \gamma g= \tau g_0+g'_1$
and rewrite the equation above as 
\[ (\gamma f)' + (\tau- \frac{\gamma'}{\gamma}) (\gamma f)= \tau g_0+ g_1'. \] 
Then 
\[\begin{split} \gamma  f(x)\, &  = \gamma(x) \int_{-\infty}^x \exp(-\tau( x-y))\gamma^{-1}(y)  ( \tau  g_0 + g_1' ) dy
\\ & =  \tau \gamma(x) \int_{-\infty}^x \exp(-\tau( x-y))\gamma^{-1}(y)  g_0 dy 
\\ & \quad - \gamma(x) \int_{-\infty}^x \exp(-\tau( x-y))\gamma^{-1}(y)\big( \tau + \frac{\gamma'}\gamma\big)   g_1  dy + g_1.
\end{split} 
\]
and by Young's inequality 
\[ \Vert \gamma f \Vert_{L^p} \le 2\Vert g_0 \Vert_{L^p} + 4 \Vert g_1 \Vert_{L^p}. \]
Indeed, by using Lipschitz continuity of $\ln(\gamma(x))$,
\[
    \int_{-\infty}^x e^{-\tau(x-y)}\frac{\gamma(x)}{\gamma(y)} |h(y)|\, dy \leq \int_{-\infty}^x e^{-\frac{\tau}2(x-y)}|h(y)|\, dy,
\]
which is estimated in $L^p$ using Young's inequality.
\end{proof} 

\begin{definition}
	A subset $Q \subset H^N(\R)$ is equicontinuous if and only if
	\begin{equation} \label{eq:defequicontinuity} 
		\lim_{h\to 0} \sup_{w\in Q} \|w(\cdot{} + h)-w\|_{H^N} \to 0.
	\end{equation}
	A subset $Q \subset H^N(\R)$ is called tight, if  for every $ \varepsilon >0$ there exists $R$ so that
	\begin{equation} \label{eq:deftight} 
	\sup_{w\in Q} \Vert  w \Vert_{H^N(\R\backslash [-R,R])} < \varepsilon. 
	\end{equation} 
	 A subset $Q$ is precompact  if it is tight and equicontinuous. 
\end{definition}

We collect the analytic properties of the modified Miura map in the following proposition. 

\begin{proposition}\label{prop:equivalencewu}  Let $ \tau \ge 1$ and $ \gamma$ slowly  varying of rate $ \tau/2$. 
The implicit constants in the sequel are independent of $\tau$ but they depend on constants in Definition \ref{def:slowlyvarying}  \begin{enumerate} 
\item \label{equi:norm}  Let $u=w_x+2\tau w+w^2$ and  $ N \ge 0 $. Then following estimates hold: 
\begin{equation}\label{eq:norm}  \begin{split}  \Vert \gamma u \Vert_{H^{N-1}_\tau}\, &  \le c(\tau^{-1/2}  \Vert w \Vert_{L^2} )  \Vert \gamma w \Vert_{H^N_\tau} \\ \Vert \gamma w \Vert_{H^N_\tau}\, &  \le c\big( \tau^{-1/2} \Vert w\Vert_{L^2}  \big) \Vert \gamma u \Vert_{H^{N-1}_\tau}.  \end{split}   \end{equation} 
\item \label{equi:diffeo}   The map 
\[ \Theta: H^N_\tau(\R) \ni w \to w_x+ 2\tau w + w^2 \in \{ u \in H_\tau^{N-1} : -\partial^2+ u +\tau^2 \text{ p.d.}\}  \] 
is a diffeomorphism with all (including higher) Fr\'echet derivatives of $ \Theta$ and $ \Theta^{-1}$
 bounded by a constant depending only on  $ \tau^{-N-\frac12} \Vert w \Vert_{H^N_\tau}$.
 \item \label{equi:wcontinuity} 
 Let $Q_W \subset L^2$ be a bounded subset and $Q_U$ its image under the modified Miura map. 
 Then 
$Q_U \subset H^{-1}_\tau $ is equicontinuous if and only 
$Q_W \subset L^2$ is equicontinuous. A set $ Q_U \subset H^{-1}$
is equicontinuous  if and only for every $ \varepsilon>0$  there exists $\tau_1$ so that for $ w \in L^2$ with 
\[ w_x + 2 \tau_1 w + w^2 \in Q_U \] 
we have $ \Vert w \Vert_{L^2} < \varepsilon$.
\item \label{equi:wtight}   Suppose that $Q_W \subset L^2$ is bounded.  Then 
$Q_U \subset H^{-1}_\tau$ is tight if and only if $Q_W \subset L^2$
is tight. 
 \item \label{equi:wprecompact} Suppose that $Q_W \subset L^2 $ is bounded. Then 
$Q_U \subset H^{-1}_\tau$ is precompact if and only if $Q_W \subset L^2$
is precompact. 
 \end{enumerate} 
\end{proposition}
\begin{proof} 
We consider the linear equation  
\[ L_w \psi := \psi_x + 2\tau \psi + 2w \psi = f \]
in considerable  detail.  
\begin{lemma} \label{lem:linearesteq} 
Let $ \gamma$, $ \gamma_1$ and $ \gamma_2$ be slowly  varying of rate $ \tau/2$, let $ n \ge 0 $, $ w \in L^2$
and $ \gamma w \in H^n$. Then 
\begin{equation}\label{eq:linforward}   \Vert \gamma L_w  \psi \Vert_{H^{n-1}_\tau} \le c (1+  \tau^{-1/2} \Vert w \Vert_{L^2} )  \Vert \gamma \psi \Vert_{H^n_\tau}
+ c \tau^{-1/2} \Vert \psi \Vert_{L^2} \Vert \gamma w \Vert_{H^n_\tau},
 \end{equation}  
\begin{equation}\label{eq:linwinverse}  \Vert \gamma_1\gamma_2 \psi \Vert_{H^n_\tau}  \le   c_n \exp ( 2 \tau^{-1} \Vert w \Vert^2_{L^2}) \Big( \Vert \gamma_1 \gamma_2 L_w \psi \Vert_{H^{n-1}_\tau} + \tau^{-1/2}  \Vert \gamma_1 w \Vert_{H^{n}_\tau}  \Vert \gamma_2 L_w \psi \Vert_{H^{-1}_\tau}\Big). \end{equation} 
\end{lemma} 

\begin{proof} 
We begin with the case $n=0$ and estimate 
\[
\begin{split} 
\Vert \gamma L_w  \psi \Vert_{H^{-1}_\tau} \, & \le 2\Vert \gamma \psi  \Vert_{L^2} + \tau^{-1/2} \Vert \gamma w\psi \Vert_{L^1}+ \tau^{-1}  \Vert \gamma' \psi \Vert_{L^2} 
\\ & \le  (3+ \tau^{-1/2} \Vert w \Vert_{L^2}) \Vert \gamma \psi \Vert_{L^2} 
\end{split} 
\] 
where we used \eqref{eq:l1h} and  Lemma \ref{lem:w1infty}. 
To bound the inverse we consider the representation \eqref{eq:linpresentation}
\[  \psi(x) = \int_{-\infty}^x  \exp\Big( -2\tau (x-y) - 2\int_y^x w dt \Big) f(y) dy, \] 
and  write the integral kernel as $ \exp(g(x,y)) $ where as in \eqref{eq:kernelbound} $g(x,y) \le -\tau (x-y) + \frac1\tau \Vert w \Vert_{L^2}^2$. We again bound 
\begin{equation} \label{eq:weightedl2}  \Vert \gamma(L_w)^{-1} \gamma^{-1} \Vert_{L^2\to L^2} \le \frac1\tau \exp\Big( \tau^{-1}\Vert w  \Vert_{L^2}^2\Big) \end{equation}  
using Schur's lemma and the obvious estimates
\[\begin{split} \hspace{2cm} & \hspace{-2cm}  \max\Big\{  \sup_x \gamma(x) \int_{-\infty}^x  \gamma(y)^{-1} \exp\Big(-2 \tau (x-y) - 2\int_y^x w dt \Big)\, dy , \\ & \quad 
\sup_y  \gamma^{-1}(y) \int_{y}^\infty  \gamma(x) \exp\Big(-2 \tau (x-y) - 2\int_y^x w dt \Big)\, dx  \Big\} 
\\ & \le \frac1\tau \exp\Big( \frac1\tau \Vert w \Vert_{L^2}^2 \Big) 
\end{split} 
\] 
 We next bound  $ \Vert \gamma L_w^{-1}\gamma^{-1} \Vert_{H^{-1}_\tau \to L^2}$ 
 and write using Lemma \ref{lem:w1infty} with $ 2\tau$ instead of $\tau$
\[   \psi' + 2\tau \psi + 2w \psi = 2\tau  f + f' \] 
so that $ \phi =\psi-f $ satisfies %\friedrich{why does one write the RHS that way? if one chooses $2\tau f + f'$, the RHS becomes only $-2wf$.}
\[ \phi' + 2\tau \phi + 2w \phi =  -2w f \]
and using \eqref{eq:weightedl2},
\[ \Vert \gamma \phi \Vert_{L^2} \le  c \exp\Big( \tau^{-1} \Vert w \Vert_{L^2}^2\Big) 
\Big(1+ \tau^{-1/2} \Vert  w \Vert_{L^2} \Big)  \Vert \gamma f \Vert_{L^2}.  \]

 Before we turn to $ n >1$ we collect calculus type estimates in Sobolev spaces.

\begin{lemma} \label{lem:calculus1} 
Let $ \gamma$, $ \gamma_1$ and $\gamma_2$ be slowly $ \tau/2$ varying. 
The following estimates hold:
\begin{equation}\label{eq:productn1}  \Vert \gamma_1\gamma_2 fg  \Vert_{H^{n}_\tau} \le c_n \Big(\Vert \gamma_1 \gamma_2 f \Vert_{W^{n,\infty}_\tau} \Vert g \Vert_{L^2}+ 
\Vert \gamma_1 f \Vert_{L^2} \Vert \gamma_2 g \Vert_{W^{n,\infty}_\tau} \Big) 
\end{equation}  
%and the consequence of $ \Vert f \Vert_{L^\infty} \le \Vert f \Vert_{L^2}^{\frac12} \Vert f' \Vert_{L^2}^{\frac12} $
\begin{equation}\label{eq:sobolevn1}   \Vert \gamma f \Vert_{W^{n,\infty}_\tau} \le \Vert \gamma  f \Vert^{\frac12}_{H^n_\tau} \Vert \gamma f \Vert^{\frac12}_{H^{n+1}_\tau} \le c \tau^{-1/2} \Vert \gamma f \Vert_{H^{n+1}_\tau}. \end{equation} 
\end{lemma} 
\begin{proof} \eqref{eq:productn1}  simply follows from the Leibniz rule if $ \gamma=1$:
\begin{align*}
    \|fg\|_{H^n_\tau}\lesssim &\sum_{j=0}^n \tau^{n-j}\|(fg)^{(j)}\|_{L^2}
    \lesssim \sum_{j=0}^n\tau^{n-j} \left(\|f^{(j)}g\|_{L^2}+\|fg^{(j)}\|_{L^2}\right)\\
    \lesssim & \sum_{j=0}^n\tau^{n-j} \left(\|f^{(j)}\|_{L^\infty}\|g\|_{L^2}+\|f\|_{L^2}\|g^{(j)}\|_{L^\infty}\right)\\
    \lesssim & \Vert f \Vert_{W^{n,\infty}_\tau} \Vert g \Vert_{L^2}+ 
\Vert f \Vert_{L^2} \Vert g \Vert_{W^{n,\infty}_\tau}.
\end{align*}
\eqref{eq:sobolevn1} is similarly a consequence of the following estimate for $ \gamma=1$
\[\|f\|_{L^\infty}\le \|f\|^{\frac12}_{L^2}\|f'\|^{\frac12}_{L^2}.\]
The case of general $\gamma$ follows from 
\[ \Vert f \Vert^p_{W^{n,p}_\gamma} 
\sim \sum_{k \in \mathbb{Z}} 
\Vert f \Vert^p_{W^{n,p}(\tau(k-1),\tau(k+1))} \] 
and 
\[ \Vert \gamma f \Vert_{W^{n,p}_\tau(\tau(k-1), \tau(k+1))} \sim \gamma(\tau k) 
\Vert f \Vert_{W^{n,p}_\tau(\tau(k-1),\tau(k+1))}. \] 
\end{proof}
Let $n \ge 1$. Inequalities \eqref{eq:productn1} and \eqref{eq:sobolevn1} immediately give
the forward estimate \eqref{eq:linforward}:
\[ \Vert \gamma \partial \psi  \Vert_{H^{n-1}_\tau} \le C\Vert \gamma w \Vert_{H^n_\tau}\] 
\[ 2\tau \Vert  \gamma  \psi \Vert_{H^{n-1}_\tau} \le 2 \Vert \gamma \psi  \Vert_{H^n_\tau}   \] where the first estimate follows by commuting the derivative with $\gamma$ and estimating lower order terms, and 
\[\begin{split}  \Vert \gamma w \psi \Vert_{H^{n-1}_\tau}
\, & \le c_n  \Big( \Vert \gamma w \Vert_{W^{n-1,\infty}}\Vert \psi \Vert_{L^2} + \Vert w \Vert_{L^2}
\Vert \gamma \psi \Vert_{W^{n-1},\infty} \Big)
\\ & \lesssim \tau^{-1/2} \Vert \gamma w \Vert_{H^n} \Vert \psi \Vert_{L^2} + \tau^{-1/2} \Vert w \Vert_{L^2} \Vert \gamma \psi \Vert_{H^n_\tau}
\end{split} 
\] 
To prove \eqref{eq:linwinverse} consider 
\[ \psi' + 2\tau \psi + 2w \psi = f. \] The case $n=0$ has been proven above. Then, taking $n-1$ derivatives of the equation, and using \eqref{eq:productn1}
\[  \Vert \gamma_1\gamma_2  \psi^{(n)} \Vert_{L^2}  \le 2 \tau \Vert \gamma_1\gamma_2 \psi^{(n-1)} \Vert_{L^2} 
+ \Vert \gamma_1\gamma_2 f^{(n-1)} \Vert_{L^2} 
+ c\Vert \gamma_1\gamma_2 w \psi \Vert_{H^{n-1}_\tau}. \] 
By Lemma \ref{lem:calculus1} and Young's inequality
\[\begin{split}  \Vert \gamma_1\gamma_2 w \psi \Vert_{H^{n-1}_\tau}\, &  \le \tau^{-1/2} \Vert \gamma_1 \psi \Vert_{L^\infty} \Vert \gamma_2 w \Vert_{H^{n-1}_\tau} + \Vert w \Vert_{L^2} \big(\Vert \gamma_1\gamma_2 \psi \Vert_{H^{n-1}_\tau}\Vert \gamma_1\gamma_2 \psi \Vert_{H^n_\tau} \big)^{1/2}  \\  &  \le \tau^{-1/2} \Vert \gamma_1 \psi \Vert_{L^2} \Vert \gamma_2 w \Vert_{H^n_\tau} + 2c\Vert w \Vert^2_{L^2} \Vert \gamma_1\gamma_2 \psi \Vert_{H^{n-1}_\tau}  
+ \frac1{2c} \Vert \gamma_1\gamma_2 \psi \Vert_{H^n_\tau}.
\end{split} 
\] 
We subtract $\frac1{2} \|\gamma_1\gamma_2 \psi\|_{H^n_\tau}$ from the combined estimate and iterate. Altogether,
\[ \Vert \gamma_1 \gamma_2 \psi \Vert_{H^n_\tau}  \lesssim \Vert \gamma_1\gamma_2 f \Vert_{H^{n-1}_\tau} + \tau^n\big(1+\tau^{-1/2} \Vert  w \Vert_{L^2}\big)^{2n-1} \Vert \gamma_1 w \Vert_{H^{n-1}}  \Vert \gamma_2 \psi \Vert_{L^2}. 
\] 
Together with the $L^2$ bounds  above  this implies \eqref{eq:linwinverse}. 
\end{proof} 

We return to the proof of  \eqref{eq:norm} with $N\ge 1$. We recall 
\[ L_{\frac{w}{2}} w = w_x+ 2\tau w + w^2 = u \] 
and  \eqref{eq:linforward} implies the first inequality of \eqref{eq:norm}.
For the second inequality   \eqref{eq:linwinverse} we observe
\[\begin{split}  \Vert \gamma w \Vert_{H^N_\tau} 
\, & \le \Vert \partial_x ( \gamma w) + 2\tau \gamma w \Vert_{H^{N-1}_\tau} 
\\ \, & \le  \Vert \gamma ( \partial_x w+ 2\tau w ) \Vert_{H^{n-1}_\tau} + \Vert \gamma' w \Vert_{H^{n-1}_\tau}
\\ & \le 
\Vert \gamma u \Vert_{H^{N-1}_\tau} + \Vert \gamma w^2 \Vert_{H^{N-1}_\tau}
+ c_n \tau \Vert \gamma w \Vert_{H^{n-1}_\tau} 
\\ & \le  \Vert \gamma u \Vert_{H^{N-1}_\tau} + c_n \Vert w \Vert_{L^2} \Big( \Vert \gamma w \Vert_{H^{N-1}_\tau}\Vert \gamma w \Vert_{H^{N}_\tau}\Big)^{1/2}+ c_n \tau \Vert \gamma w \Vert_{H^{n-1}_\tau}
\\ & \le  \Vert \gamma u \Vert_{H^{N-1}_\tau} + c^2_n \Vert w \Vert^2_{L^2}  \Vert \gamma w \Vert_{H^{N-1}_\tau}
+ \frac14 \Vert \gamma w \Vert_{H^{N}_\tau}+ c_n  \tau \Vert \gamma w \Vert_{H^{n-1}_\tau}
\end{split} 
\] 
and we complete the argument in the same fashion as for the linear estimate.

To see injectivity we consider for $ j=1,2$
\[ \partial w_j + 2\tau w_j + w_j^2 = u_j \] 
hence with $ w = w_2-w_1$
\[ \partial w + 2 \tau w + (w_1+w_2) w = u_2-u_1 \] 
and \eqref{eq:linwinverse} provides a bound 
\[ \Vert w^2- w^1 \Vert_{L^2} \le c \exp(\tau^{-1} ( \Vert w_1 \Vert_{L^2}^2 + \Vert w_2 \Vert_{L^2}^2) \Vert u_2-u_1 \Vert_{H^{-1}_\tau}. \] 
Moreover \eqref{eq:linwinverse} provides an estimate for the Lipschitz norm of the inverse restricted to  bounded  sets corresponding to $w$ in a ball in $L^2$.
Fr\'echet differentiability of the inverse is immediate and its differential is given by the inverse of $L_w$. The inverse function theorem gives Fr\'echet smoothness of the inverse and bounds can be obtained from differentiating  $(L_w)^{-1}$
as in multivariate calculus. 

The range of
$\Theta$ 
is $H^{N-1}_\tau \cap \{-\partial^2 + u + \tau^2 > 0\}$. Indeed, one inclusion follows from \eqref{eq:equivalent2}. For the other inclusion note that if $u \in H^{-1}$ with $- \partial^2 + u + \tau^2 > 0$, we can define $w$ via its Jost function as $w = \partial_x \ln \psi - \tau$. By positivity of $\-\partial^2 + u + \tau^2$ and with the notation as in \eqref{eq:equivalent3}, $\tau^2 - \tau_1^2 > 0$, and hence $w \in L^2$. $w \in H^N_\tau$ then follows from \eqref{eq:norm}.

We turn to the proof of Proposition \ref{prop:equivalencewu} \eqref{equi:wcontinuity}.
The maps $ \Xi$ and $ (\Xi)^{-1} $ are uniformly Lipschitz on bounded sets satisfying the equivalent conditions in Lemma \ref{lem:uw}, hence $Q_u$ in \eqref{equi:wcontinuity} is equicontinuous iff $Q_v$ is equicontinuous.   
 
It is easily seen that condition \eqref{eq:defequicontinuity} for $N=0$ implies tightness in the Fourier space, more precisely $Q \subset H^{-1}$ is equicontinuous
if and only if for every $ \varepsilon>0$ there exists $R$ so that 
\[ \sup_{f\in Q} \Vert |\xi|^{-1}  \hat f \Vert_{L^2( \R \backslash (-R,R) ) } < \varepsilon \] 
which is equivalent to the existence of $ \tau_0$  so that 
\[  \sup_{f\in Q, \tau \ge \tau_0} \Vert f \Vert_{H^{-1}_\tau} < \varepsilon. \]
For details we refer to \cite{MR3990604}.

Together with Lemma \ref{lem:uw}
and Proposition \ref{prop:equivalencewu} \eqref{eq:norm} with $N=0$
we see that equicontinuity is equivalent to the $L^2$ norm of preimages of the modified Miura map with parameter $\tau_1$ being small for $ \tau_1$ large, which is  the second part of the equicontinuity claim of Proposition \ref{prop:equivalencewu} \eqref{equi:wcontinuity}. 

To prove equivalence of  tightness we use the  following equivalent characterization.

\begin{lemma} 
A set $ Q \subset H^s$ is bounded and tight if and only if for all $ \alpha>0$ there  exists an $ \alpha$-slowly varying function $ \eta\ge 1$ with $\lim_{x\to \pm \infty} \eta(x)= \infty$
so that 
\[ \sup_{w\in Q} \Vert \eta w \Vert_{H^s} < \infty. \] 
\end{lemma} 

\begin{proof} 
Let 
\[ \alpha_n = \sup_{w\in Q} \Vert w \Vert_{H^s(\R\backslash [- 2^n ,2^n] ) }. \] 
A bounded set $Q$ is tight if and only if $ \alpha_n \to 0$. Choose $n_k\ge k $ so that 
$ \alpha_{n_k} \le \exp^{-2\alpha k} $. It is not hard to define an $ \alpha $- slowly varying function $\eta\ge 1 $, monotone in radial direction, so that 
\[ \eta( \pm 2^{n_k}) = \exp( \alpha k). \] 
Then 
\[ \sup_{w\in Q} \Vert \eta w \Vert_{L^2} < \infty.\]
\end{proof} 

We turn to the proof of Proposition \eqref{prop:equivalencewu} \eqref{equi:wtight}. Suppose that 
$Q_w \in L^2$ is bounded and tight. Then, if $ \gamma $ is slowly varying with rate $ \tau/2$ than by Proposition \ref{prop:equivalencewu} \eqref{equi:norm}   
\[ \Vert\gamma  (\partial_x w + 2\tau w + w^2)\Vert_{H^{-1}_\tau} \sim  \Vert \gamma w \Vert_{L^2}\]  with implicit constants depending on $ \tau^{-1/2} \Vert w \Vert_{L^2}$.  
Finally Proposition \ref{prop:equivalencewu} \eqref{equi:wprecompact} follows from Proposition \ref{prop:equivalencewu} \eqref{equi:wcontinuity} and \eqref{equi:wtight}. \end{proof}

Let $ 1 \le \tau < \im z$. The map 
\[   L^2\ni w \to     w(z) \in L^2\]
defined by inverting the Miura map resp. solving the equation on the left via the left Jost function
\[    w'(z) -2iz w(z) + w^2(z)= w' +2 \tau w - w^2 \]
which will be useful at several points. 

\begin{theorem} \label{thm:locallipschtz} 
Let $1 \le \tau < \tau_1$, $ N \in N$  and  $ \gamma$ be slowly $\tau/2$ varying. Then 
\[      \Big\Vert \gamma \frac{\delta}{\delta w} \int w^2(i\tau_1,x) dx \Big\Vert_{H^{N}} \le c(\Vert w \Vert_{L^2})  \Vert \gamma w \Vert_{H^N} \] 
and the Lipschitz constant of the variational derivative is bounded by
\[
\begin{split}\Big\Vert &\gamma\Big(\frac{\delta}{\delta w}\int w^2(i\tau_1,x;w) dx(w_1) -  \frac{\delta}{\delta w}\int w^2(i\tau_1,x;w) dx(w_2)\Big)\Big\Vert_{H^{N}} \\
&\le c( \Vert w_1 \Vert_{L^2}, \Vert w_2 \Vert_{L^2})\Big(  \Vert \gamma (w_2-w_1) \Vert_{H^N} +\Vert w_2-w_1 \Vert_{L^2} ( \Vert \gamma w_2 \Vert_{H^N}+ \Vert \gamma w_1 \Vert_{H^N})\Big). \end{split}\]    
\end{theorem} 

\begin{proof} 
We first compute the variational derivative. Let $ \phi$ be a test function. Then 
\[  \int \frac{\delta}{\delta w} \Vert  w^2(i\tau_1) \Vert_{L^2}^2 \phi dy = 2  \int (-\partial+ 2\tau +2 w) (-\partial + 2\tau_1 + 2 w(i\tau_1))^{-1} w(i\tau) \phi dx.  \]
The first estimate is a consequence of \eqref{eq:norm}. Moreover 
\[ \begin{split} u_2-u_1\, & =      w'_2 +2\tau w_2 + w^2_2 - (w'_1+2\tau w_1 + w_1^2 )
\\ & = \partial (w^2-w^1)+ 2 \tau (w^2-w^1) + (w^2+w^1) (w^2-w^1)
\\ & = \partial (w^2(i\tau_1)-w^1(i\tau_1)+ 2 \tau_1 (w^2(i\tau_1)-w^1(i\tau_1) \\ & \qquad + (w^2(i\tau_1)+w^1(i\tau_1)) (w^2(i\tau_1)-w^1(i\tau_1) )
\end{split} 
\]
and estimate follows by the linear estimates  by Lemma 
 \ref{lem:linearesteq} , \eqref{eq:linforward} and  \eqref{eq:linwinverse}.
\end{proof}

\subsection{The good variable hierarchy}\label{subsec:goodvariable}

The map
\[ v \to \tau v -\frac{v_x}{2(1+v)} \] 
is a diffeomorphism between 
subsets of Banach spaces. It relates the good variable hierarchy and the Gardner hierarchy in a very similar fashion as the modified Miura map related the Gardner hierarchy and the KdV hierarchy. The results in this section complement the previous results. Proofs are similar, technical but to a large extend standard, which was different for Lemma \ref{lem:uw}. We only provide part of the  proofs here which we consider nonstandard.  The following theorem describes properties of this map. 

\begin{theorem}\label{thm:equivalence}   
\begin{enumerate} 
\item \label{equi:wv} If $ s>-\frac12$, $Q_v \subset \{ v \in H^{s+1}: v>-1\}$, $Q_w = \{ \tau v - \frac12 \partial \ln(1+v)\} \subset H^s$ the following is equivalent
\begin{itemize} 
\item There exists $r$ so that $ \Vert w \Vert_{H^s_\tau } \le r $ for all  $ w \in Q_w$.
\item There exists $R$ and $ \delta >0 $ so that $ 1+v \ge \delta $
and $ \Vert v \Vert_{H^{s+1}_\tau} \le R $
\end{itemize} 
\item \label{equiB:norm} Let $ s \ge -1 $. Then following estimates hold for $ 0<\varepsilon \le 1$ 
\[ \Vert \gamma w \Vert_{H^{s}_\tau} \le c\Big( \Vert v \Vert_{L^\infty}, \Vert (1+v)^{-1} \Vert_{L^\infty}  \Big)  \Vert \gamma w \Vert_{H^{s+1}_\tau} \] 
\[ \Vert \gamma v \Vert_{H^{s+1}} \le c(  \Vert w  \Vert_{H^{-\frac12+\varepsilon}_\tau})  \Vert \gamma w \Vert_{H^{s}_\tau}.    \] 
\item \label{equi:thetav}  Let $s > -\frac12$. The map 
\[ \Theta_v: \{ v \in H^{s+1}_\tau(\R): v >-1\}  \ni v \to \tau v - \frac12 \partial \ln(1+v) \in  H^s_\tau(\R) \] 
is a diffeomorphism with all (including higher) Fr\'echet derivatives of $ \Theta$ and $ \Theta^{-1}$
 bounded by a constant depending on $ \tau^{-s} \Vert v \Vert_{H^s_\tau}$ and $ \inf v+1$.
\item \label{equi:vuniform} Suppose that $ s>-\frac12$ and that the equivalent conditions of \ref{equi:wv} hold. Then $ Q_w\subset H^{s}_\tau  $ is equicontinuous if and only if $Q_v \subset H^{s+1}_\tau$ is equicontinuous. 
\item \label{equi:vtight} Suppose that the equivalent condition of \eqref{equi:wv} hold, then $Q_w$
is tight in $H^s_\tau$ if and only if $Q_v$ is tight in $H^{s+1}_\tau$. 
\item \label{equi:vprecompact} Suppose that the equivalent condition of \eqref{equi:wv} hold. Then $Q_w \subset H^s_\tau$ is precompact if and only if $Q_v \in H^{s+1}_\tau$ is precompact. 
\end{enumerate} 
\end{theorem}

\subsubsection{Proof of Theorem \ref{thm:equivalence}} 

\begin{proof} 
We  study the map 
\[  v  \to  \tau v   - \frac12 \partial_x  \ln( 1+ v ) = : w \] 
with inverse 
\[ v = \frac{1}{\tau (L_w)^{-1} (1)} -1. \]  
The lower bound $1+v \ge \delta$ of \eqref{equi:wv} is equivalent to the upper bound
\[ (L_w)^{-1}1 \le \delta^{-1}. \] 
Let  $-\frac12<s< \frac12 $.  We combine the   embedding $H^{s+1} \subset C^{s+\frac12}$ and $x\to  \int_{y}^x w \in H^{s+1}_{loc}$ if $ w \in H^s$ with Young's inequality 
\[ \left| \int_{y}^x w d\sigma \right| \le c  |x-y|^{s+\frac12}    \Vert w \Vert_{H^s} \le \tau  |x-y| +    c \tau^{-\frac{1+2s}{1-2s}}   \Vert w \Vert_{H^s}^{\frac2{1-2s}}.   \]

Again by Schur's lemma 
\[ \Vert L_w^{-1} f \Vert_{L^\infty} \le   \frac1\tau \exp\Big( c \tau^{-\frac{1+2s}{1-2s}}   \Vert w \Vert_{H^s}^{\frac2{1-2s}}\Big) \Vert f \Vert_{L^\infty} \]
which implies the lower bound
\begin{equation} 
\label{eq:deltavw} 
1+ v \ge \tau \exp(-c  \tau^{-\frac{1+2s}{1-2s}}   \Vert w \Vert_{H^s}^{\frac2{1-2s}}\Big).
\end{equation} 
Trivially  we have 
\begin{equation} \label{eq:omegawinfty} 
\Vert w \Vert_{W^{-1,\infty}_\tau} \le  \Vert v \Vert_{L^\infty} + \Big\Vert \frac1{1+v} \Big\Vert_{L^\infty},   \end{equation} 
We write 
\[ \ln(1+v)= \int_0^1 \frac1{1+tv} dt v = \psi^{-1} v \]
with $ \psi \ge \delta$.
Let $ \phi= \ln(1+v)$. Then 
\[ \phi_x + 2\tau \psi  \phi      =  2 w \]
hence, for $1\le p \le \infty $ 
\begin{equation}\label{eq:phiest}   \Vert \phi \Vert_{L^p} \le \frac{2}{\delta \tau}  \Vert w \Vert_{W^{-1,p}_\tau} \end{equation}  
hence 
\begin{equation} \label{eq:vtow} 
\Vert v \Vert_{L^\infty} \le \exp\Big( \frac{1}{\delta \tau}\Vert w \Vert_{W^{-1,\infty}}\Big).
\end{equation} 
The bounds \eqref{eq:deltavw}, \eqref{eq:omegawinfty} and \eqref{eq:vtow} are the basis for the remaining estimates. 
\end{proof}

We simplify the estimate a bit by substituting 
$ w = \tau  v -  \frac{ v_x}{1+ v} $.
We obtain linear combinations of 
\[   ( \tau + \partial \frac1{2(1+\tilde v)} )^{-1}\partial^{1+j_0}
(1+v)^{-M} \tau^{2l} \prod_{k=1}^K  v^{(j_k)}
\] 
with (denoting by $D= 1+ \sum_{k=0}^K j_k$ the total number of derivatives) $2L+D = 2N+2$.

\section{Weak solutions}
\label{sec:weak} 

In this section we study weak solutions to equations of the $N$th KdV equation and the $N$th Gardner equation. Under weak regularity conditions $w$ is a weak solution to the $N$th Gardner equation if and only if $u= w_x + 2\tau w + w^2$ is a weak solution to the $N$th KdV equation. 
This reduces the proof of the main theorem to a study of the weak solutions to the Gardner hierarchy.

\subsection{Calculus estimates in Sobolev spaces} 
In almost all sections below we will need estimates of differential monomials in $L^p$ or weighted $L^p$.
\begin{definition} Let $1\le p\le \infty$, $N \ge 0$, $ \tau >0$ and $I=(a,b)$. We define 
\[   W^{N,p}_\tau(I) = \Big\{ f= \sum_{j=0}^N \tau^{N-j} \partial^j  f_j :  f_j \in L^p(I) \Big\} \] 
with 
\[ \Vert f \Vert_{W^{-N,p}_\tau(I)} = \inf \Big\{\Big(\sum_{j=0}^N \Vert f_j  \Vert^p_{L^p(I)} \Big)^{1/p} : f= \sum_{j=0}^N \tau^{N-j} \partial^j  f_j \Big\}. \]
We  define for $ N \ge 0 $ and  $ \tau > 0$
\[ \Vert g \Vert_{W^{N,p}_\tau(I)} = \sum_{j=0}^N \tau^{N-j} \Vert g^{(j)} \Vert_{L^p(I)}\]
and for $ \tau >0 $ and $ s \in \R$
\[ \Vert f \Vert_{H^s_\tau(\R)} = \Vert  ( \tau^2 +\xi^2)^{s/2}  \hat f \Vert_{L^2(\R) }  \]
and for an interval $I$ 
\[ \Vert f \Vert_{H^s_\tau(I)} =  \inf\{  \Vert \tilde f \Vert_{H^s_\tau}, f = \tilde f  \text{ on } I \} \]
\end{definition} 

It is obvious that the restriction to smaller intervals
is a
bounded linear operator of norm $1$.
There exists a bounded extension operator to functions supported on twice the interval. 

\begin{lemma} 
The following norms are equivalent: Suppose that 
$n \in \Z$, $ \tau >0$ $\tau r \ge 1$ and $ 1 \le p \le \infty$. Then 
\[ \Vert f \Vert_{W^{n,p}_\tau(\R)} \sim \big\Vert \Vert f \Vert_{W^{n,p}_\tau((k-1)r,(k+1)r))} \big\Vert_{l^p}. \] 
\end{lemma} 
\begin{proof}
It suffices to verify the claim for $ \tau = 1$ and $r=1$. The case $n\ge 0$
is obvious. Let $n>0$, $f= \sum_{j=0}^n \partial^j f_j $ with $\sum_{j=0}^n \Vert f_j  \Vert^p _{L^p(\R)} \sim \Vert f \Vert^p_{W^{-n,p}(\R)}$. Then 
\[ \big\Vert \Vert f \Vert_{W^{-n,p}(k-1,k+1)} \big\Vert^p_{l^p}
\le \sum_k \sum_{j=0}^{n}  \Vert f_j  \Vert^p_{L^p(k-1,k+1)}  = 2 \sum_{j=0}^n \Vert f_j \Vert_{L^p(\R)}^p \] 
if $1\le p < \infty $ with obvious modifications if $p=\infty$. 

For the opposite direction choose $f_{k,j} \in L^p(k-1,k+1) $ so that $f= \sum_{j=0}^n \partial^j f_j$ on $[k-1,k+1]$,  $ \sum_{j=0}^n \Vert f_{k,j} \Vert^p_{L^p(k-1,k+1)}  \le 2 \Vert f \Vert^p_{W^{-n,p}(k-1,k+1)} $, choose a partition of unity $ \sum \eta(x-k) = 1$ with $ \supp \eta \subset (-1,1)$, $ \eta=1$ on $(-1/4,1/4)$. Then 
\[  f = \sum_{k\in \Z} \eta(x-k) \sum_{j=0}^n \partial^j f_{k,j} 
= \sum_{k\in \Z} \sum_{j=0}^n \partial^j \eta f_j    -\sum_{l=1}^{j-1} \binom{j}{l} \partial^{j-l}  ( \eta^{(l)}(x-k) f_j). \] 
We set 
\[ f_j= \sum_{k\in \Z}  \eta(x-k) f_j - \sum_{l=j+1}^n \binom{l}{j}  \eta^{(l-j)}(x-k) f_l \] 
and obtain 
\[ \sum_{j=0}^n  \Vert f_j \Vert_{L^p}^p \le    c  \sum_{k} \Vert f \Vert^p_{L^p(k-1,k+1)}\]
again with obvious modifications if $p=\infty$.
\end{proof} 

We turn to an  interpolation inequality, for which we  provide a proof for completeness.

\begin{lemma} \label{lem:elementary} 
Let $ 0 \le j < s$, $2 \le q \le  \infty$, $ 2\le r < \infty $ and  
\[  \frac{s}r = \frac{s-j}q + \frac{j}2.    \]
Then 
\[ \Vert f^{(j)} \Vert^s_{L^r} \le c \Vert f \Vert^{s-j}_{L^{q}} \Vert f^{(s)} \Vert^j_{L^2}. \]
\end{lemma} 

\begin{proof} 
The lemma relies on three elementary estimates. 
Suppose that 
\[ \frac1p +\frac1q = \frac2r, \,2\le r\le  \infty, \,1\le p,q \le \infty \] 
Then 
\begin{equation}\label{elementaryinterpolation}   \Vert f' \Vert^2_{L^r} \le (r-1) \Vert f \Vert_{L^p} \Vert f'' \Vert_{L^q}. \end{equation} 
This follows from 
\[
\begin{split} 
\int |f'|^r dx \,  & = \int \partial_x ( |f'|^{r-2}f'  f) dx - (r-1) \int f |f'|^{r-2} f'' dx 
\\ & \le  (r-1)\Vert f' \Vert_{L^r}^{r-2} \Vert f \Vert_{L^p} \Vert f'' \Vert_{L^q} 
\end{split} 
\] 
which implies \eqref{elementaryinterpolation}.  There is a version for fractional derivatives.  Let $0< s < 1$ and $1\le p, q \le \infty$. We define the homogeneous Besov norm 
\[ \Vert f \Vert_{\dot B^s_{pq}}  =   \left( \int_0^\infty  \big( h^{-s} \Vert f(.+h)- f \Vert_{L^p}\big)^q  \frac{dh}{h} \right)^{1/q} \]
with the obvious modification for $q=\infty$. Then $\dot B^s_{pq} \subset \dot B^2_{p\tilde q} $ whenever $ \tilde q \le q$ and 
\[ \Vert f \Vert_{\dot B^s_{22}} = c \Vert |\xi|^s \hat f \Vert_{L^2}=\Vert f \Vert_{\dot H^s}.\]
Let $ \frac1r= \frac{1-s}q + \frac{s}p$ then  
\[  \Vert h^{-s} f(.+h)- f \Vert_{L^r} =  \Vert  |f(.+h)-f|^{1-s}   (h^{-1} |f(.+h) -f(.)|^s     \Vert_{L^r}    \le  ( 2 \Vert f \Vert_{L^q} )^{1-s} \Vert f' \Vert_{L^p}^ s \] 
hence
\begin{equation} \label{eq:elementary2}  \Vert f \Vert_{\dot B^s_{r,\infty}}  \le 2 \Vert f \Vert_{L^p} \Vert f' \Vert_{L^q}.  \end{equation} 
Here we used  
\[  \Vert f' \Vert_{L^r} = \sup_{h >0 }   h^{-1} \Vert f(.+h) - f \Vert_{L^r}.\]
Similarly 
\[ \begin{split}\hspace{1cm} & \hspace{-1cm} h^{-r} \int |f(x+h)-f(x)|^{r-2} (f(x+h)-f(x) (f(x+h)-f(x))  dx
 \\ &  = h^{-r} \int f(x) \Big( |f(x+h)-f(x) |^{r-2}(f(x+h)-f(x)) \\ & \qquad- |f(x)-f(x-h)|^{r-2}(f(x) -f(x-h) \Big) dx  
\\ & \le  \Vert f \Vert_{\dot H^{s-1} } \Big( h^{-1} \Vert f(.+h)-f \Vert_{L^r}\Big)^{r-2} \Vert |h|^{-1-s} f(.+h) -2 f(x) + f(.-h) \Vert_{L^q} 
\end{split} 
\]
hence, if 
\[ \frac1p +\frac1q = \frac2r, \quad  r\ge 2 \]
\begin{equation}\label{eq:interpolationfrac}   \Vert f' \Vert_{L^r}^2\le \Vert f \Vert _{\dot B^{1-s}_{q,\infty}}  \Vert f' \Vert_{\dot B^s_{p,\infty}}. \end{equation}

Recursively we obtain for 
$ j \leq n$, $ 2 \le p, q , r$ satisfying 
\[ \frac{n}{r} = \frac{n-j}p + \frac{j}{q} \]
\begin{equation}\label{interpolationinequality}  \Vert f^{(j)} \Vert^n_{L^r}  \le c \Vert f \Vert^{n-j}_{L^p} \Vert f^{(n)} \Vert^j_{L^q}. \end{equation} 
Indeed, suppose this estimate holds for $n$.  
Let $ 1\le j \le n+1$ and 
\[  \frac{1}{r_j} = \frac{n+1-j}{n+1} \frac1p  + \frac{j}n \frac1q \] 
Then 
\[ 
\begin{split} 
\Vert f^{(j)} \Vert^{nj}_{L^{r_j}} \, & \lesssim \left( \Vert f' \Vert^{n+1-j}_{L^{r_1}} \Vert f^{(n+1)} \Vert^{j-1}_{L^q} \right)^j 
\\ & \le c \Vert f \Vert^{(n+1-j)(j-1)}_{L^p}  \Vert f^{(j)} \Vert^{n+1-j}_{L^{r_j}} \Vert f^{(n+1)} \Vert^{(j-1)j}_{L^q} 
\end{split} 
\] 
where 
\[ nj  -( n+1-j)  =(j-1)(n+1)   \]   
which implies estimate \eqref{interpolationinequality} for $s \in \N$. The general case follows by using \eqref{eq:elementary2} and \eqref{eq:interpolationfrac} in addition.  
\end{proof} 

A particular instance is the following. 
Let $ 0\le d\le N  $ and $ p \ge 2$. 
 Then 
\begin{equation} \label{eq:interpolationa} 
\Vert u^{(d)} \Vert_{L^{\frac{N+2}{1+d}}} 
\le c \Vert u \Vert^{\frac{N-2d}{N}}_{L^{N+2}} \Vert u^{(N/2)} \Vert^{\frac{2d}{N}}_{L^2}. 
\end{equation} 
A weighted variant is
\begin{equation} \label{eq:interpolationb} 
\Vert \sech( x) u^{(d)} \Vert_{\frac{N+2}{1+d}}
\le c \Vert u \Vert^{\frac{N-2d}{N}}_{L^{N+2}} \Vert  \sech(x)  u \Vert^{\frac{2d}{N}}_{H^{N/2}}. 
\end{equation}

To see this we apply a standard extension argument to deduce from \eqref{eq:interpolationa} on intervals $I$
of length $2$ 
\[   \Vert u^{(d)} \Vert_{L^{\frac{n*2}{1+d}}(I)} 
\le c \Vert u \Vert^{\frac{N-2d}{N}}_{L^{N+2}(I)} \Vert u \Vert^{\frac{2d}{N}}_{H^{N/2}(I)}.  \]
We multiply by the weight, take the power $\frac{n+2}{1+d}$, add the intervals to arrive at \eqref{eq:interpolationb}.

To proceed, we show the following multilinear  estimates.

\begin{lemma}\label{lem:estimatedifferentialmonomial2} 
If  $h =  \prod_{j=1}^{N+2-d} u^{(\alpha_j)}$ 
is a product  with a total number of derivatives $d= \sum \alpha_j\le N$
 then, if no term carries more than $N/2$ derivatives,
\begin{equation} \label{eq:esthamon} \prod_{j=1}^{N+2-d}  \Vert u^{(\alpha_j)}\Vert_{L^{\frac{N+2}{1+\alpha_j}}} 
\le c\left\{ \begin{array} {l} \Vert u \Vert_{L^{N+2}}^{N+2-d-\frac{d}{N}} \Vert u^{(N/2)} \Vert^{\frac{d}{N}}_{L^2}  \\[2mm] \Vert u \Vert_{L^{N+2}}^{N+2} + \Vert u^{(N/2)} \Vert_{L^2}^2
\\[2mm] \Vert u \Vert_{L^2}^N \Vert u^{(\frac{N+d}4)} \Vert^2_{L^2}
\end{array}\right.
\end{equation} 
and 
\begin{equation} \label{eq:esthbmon}  \prod_{j=1}^{N+2-d}  \Vert \sech^{\frac{2(1+\alpha_j)}{N+2} }(x) u^{(\alpha_j)}\Vert_{L^{\frac{N+2}{1+\alpha_j}}}
\le c\left\{ \begin{array} {l} \Vert u \Vert_{L^2}^{N-d} \Vert \sech(x) u \Vert_{L^2}^{1-\frac{d}{N}} \Vert \sech(x) u \Vert^{1+\frac{d}{N}}_{H^{N/2}}  \\ \Vert u \Vert_{L^2}^N \Vert \sech(x) u^{(\frac{N+d}4)} \Vert^2_{H^{\frac{N+d}4}}. 
\end{array}\right.
\end{equation} 
If  $h =  \prod_{j=1}^{N+2-d} u_j^{(\alpha_j)}$ 
is a product  with a total number of derivatives $d= \sum \alpha_j\le N$
 then, if no term carries more than $N/2$ derivatives
 \begin{equation} \label{eq:estha} 
 \begin{split} 
\int  |h| dx\, &  \le \prod_{j=1}^{N+2-d} \Vert u_j^{(\alpha_j)}\Vert_{L^{\frac{N+2}{1+\alpha_j}}}  \\ &  \le c  \prod_{j=1}^{N+2-d} \left( \Vert u_j \Vert_{L^{N+2}}^{N+2-d-\frac{2\alpha_j}{N}} \Vert  u_j^{(N/2)}  \Vert_{L^2}^{\frac{2\alpha_j}{N}}  \right)^{\frac1{N+2-d}}
\\    & \le c \prod_{j=1}^{N+2-d} \Big( \Vert u_j \Vert_{L^{N+2}}^{N+2} + \Vert u_j^{(N/2)} \Vert_{L^2}^2  \Big)^{\frac{1}{N+2-d} }
\\   \Vert u_j^{(\alpha_j)}\Vert_{L^{\frac{N+2}{1+\alpha_j}}}  & \le c \Vert u_j \Vert_{L^2}^{1- \frac{4\alpha_j(N+1) +2N}{(N+2)(N+d)}}  \Vert u_j^{(\frac{N+d}4)} \Vert_{L^2}^{\frac{4\alpha_j(N+1) +2N}{(N+2)(N+d)}} 
\end{split} 
\end{equation}
and 
\begin{equation} \label{eq:esthb} 
\begin{split} 
\int \sech^2(x)  |h| dx \, & \le \prod_{j=1}^{N+2-d} \Vert \sech^{\frac{2(1+\alpha_j)}{d+m}  } (x)  u_j^{(\alpha_j)} \Vert_{L^{\frac{N+2}{1+\alpha_j}} }
\\ & \le c  \prod_{j=1}^{N+2-d} \left( \Vert   u_j \Vert_{L^{2}}^{N-d}\Vert \sech(x) u_j \Vert^{2-\frac{2\alpha_j}{N}}_{L^2} \Vert \sech(x)  u_j^{(N/2)}  \Vert_{H^{N/2}}^{\frac{2\alpha_j}{N}}  \right)^{\frac1{N+2-d}}
\\   \Vert u_j^{(\alpha_j)}\Vert_{L^{\frac{N+2}{1+\alpha_j}}}  & \le c \Vert u_j \Vert_{L^2}^{1- \frac{4\alpha_j(N+1) +2N}{(N+2)(N+d)}}  \Vert \sech (x) u_j^{(\frac{N+d}4)} \Vert_{L^2}^{\frac{4\alpha_j(N+1) +2N}{(N+2)(N+d)}}. 
\end{split} 
\end{equation}
\end{lemma}

\begin{proof} 
It suffices to prove \eqref{eq:estha} and \eqref{eq:esthb}.  The inequalities \eqref{eq:esthamon} and \eqref{eq:esthbmon} are immediate consequences.  The first inequality in \eqref{eq:estha}  is a consequence of H\"older's inequality,  for the second estimate we apply lemma \ref{lem:elementary}.  The third estimate follows by Young's inequality. 

To prove the fourth estimate we recall that for $ 0< s<1$
\[ f =  c_s |x|^{-1+s} * f^{(s)} \]
from which we obtain the fractional 
 Sobolev embedding 
\[ \Vert f \Vert_{L^p} \le c \Vert u^{(\frac12-\frac1p)} \Vert_{L^2} \]
by the Hardy Littlewood Sobolev inequality.
Let $d \ge 0 $ and $ p\ge 2$. Then by the fractional Sobolev embedding 
\[ \Vert u^{(d)} \Vert_{L^p} \le   c \Vert u^{(d+\frac12-\frac1p)}\Vert_{L^2} \]
and the interpolation inequality for $ 0 \le \sigma \le s$
\[ \Vert f \Vert_{\dot H^\sigma} \le \Vert f \Vert^{1-\frac{\sigma}s}_{L^2} \Vert f \Vert^{\frac{\sigma}s} _{\dot H^s}  \]
implies the fourth estimate. 

We turn to \eqref{eq:esthb}. Again the first inequality is H\"older's inequality. We continue with the Hardy-Littlewood Sobolev inequality. Again w interpolate on intervals for fixed length, multiply by the weight, square and sum over the intervals to arrive at the second inequality of \eqref{eq:esthb}.  
\end{proof}

\subsection{The nonlinearities of the KdV and the Gardner equations}\label{subsec:nonlinearity}

We recall the structure of the $N$th equation in each of the two hierarchies involved. Each of them can be written as
\[
    \psi_t - (-1)^{N}\psi^{(2N+1)} = \partial F_N(\psi),
\]
where the structure of $F_N$ is described by Theorem \ref{thm:formofkdv}. For later purposes we want to pull out as many derivatives as possible:
\begin{lemma}\label{lem:gardnerstructure} 
For KdV, we have,
\[
\begin{split}
    &F_N^{\KdV}(u) = \sum_{K,(j_k)_{0\le k \le K}}c_{K,(j_k)}\partial^{j_0}  \prod_{k=1}^K u^{(j_k)},\\
    2 &\leq K \leq N+1, \qquad K + \frac{1}{2} \sum_{k=0}^K j_k = N+1,
    %\\
%    0 &\leq i \leq N+1-K-\Big[\frac{N+1-K}{2}\Big], \qquad 
\quad j_k \leq N+1-K-\frac{j_0}{2}, \text{ if } j \geq 1
\end{split}
\]
and for Gardner,
\[
\begin{split}
    &F_N^{\Gardner}(w) = \sum_{K,(j_k)_k,l}c_{(j_k)_k,K,l}\partial^{j_0}  \tau^l\prod_{k=1}^K w^{(j_k)},\\
    2 &\leq K \leq 2N+1, \qquad l+K + \sum_{k=0}^K j_k  = 2N+1,\\
    %\\
    %0 &\leq i \leq %2N+1-K-l-\Big[\frac{2N+1-K-l}{2}\Big], 
    j_k &\leq \frac{2N+1-K-l-j_0}{2} \text{ if } j \geq 1
\end{split}
\]
where $K + l$ is always odd.
\end{lemma}
\begin{proof}
Theorem \ref{thm:formofkdv} describes the structure of the Hamiltonians. 
For KdV, $F_N(u)$ is a sum over differential monomials $\prod_{k=1}^K u^{(j_k)}$ with
\[
    2 \leq K \leq N+1, \qquad K + \frac{1}{2} \sum_{k=1}^K j_k = N+1,
\]
and for Gardner, $F_N(w)$ is a sum over $\tau^l \prod_{k=1}^K w^{(j_k)}$ with 
\[
    2 \leq K \leq 2N+1, \qquad l + K + \sum_{k=1}^K j_k = 2N+1.
\]
For KdV and Gardner we reduce the highest number of derivatives falling on one factor until there are at least two factors with the highest number of derivatives. This can be done as follows: Consider a differential monomial $ \prod_{k=1}^K u^{(j_k)}$. We order the factors 
so that $ k \to j_k $ decreases monotonically. If $ j_1 > j_2= j_l > j_{l+1}$ we write 
\[  \prod_{k=1}^K  u^{(j_k)} = \frac1{l} \partial  (u^{(j_2)})^l \prod_{k=l+1}^K u^{(j_k)} 
- \frac1l  (u^{(j_2)})^l \partial \prod_{k=l+1}^K  u^{(j_k)}. \]
$K+l$ being odd follows from the formula
\[
    \partial_x \frac{\delta H^{\KdV}_N}{\delta u}(w_x + 2\tau w + w^2) = (\partial_x + 2\tau + 2w)\frac{\delta H_N^{\Gardner}}{\delta w}(w).
\]
\end{proof}

The Gardner Hamiltonian $H_N^{\Gardner}$ is an integral over $ \frac12(w^{(N)})^2$ plus 
 a sum of differential monomials 
\[   \tau^l  \prod_{k=1}^K w^{(j_k)}  \] 
with $K \geq 3$, $ l+\sum_k j_k+ K = 2N+2$, $l+K $ even, and no term carries more then $2N-l-K$ derivatives. We apply \eqref{eq:estha} of Lemma \ref{lem:estimatedifferentialmonomial2}: 
\[  \Big\Vert \tau^l  \prod_{k=1}^K w^{(j_k)} \Big\Vert_{L^1} \le c \tau^l \Vert w \Vert_{L^2}^{K-2}
\Vert w^{(N-\frac{K-2}4 -\frac{l}2 )}   \Vert^2_{L^2}
\le c  (\tau^{-1/2} \Vert w \Vert_{L^2})^{K-2}  \Vert w \Vert^2_{H^N_\tau}.
\] 
hence 
\begin{equation}    \left| H_N^{\Gardner}(w)- \Vert w^{(N)} \Vert_{L^2}^2 \right| \le c (1+ \tau^{-1/2}\Vert w \Vert_{L^2})^{2N-1} (\tau^{-1/2} \Vert w \Vert_{L^2})  \Vert w \Vert_{H^N_\tau}^2.\end{equation}
The same argument shows with $ u = w_x + \tau w$
\begin{equation} 
\begin{split} 
\left| H_N^{\KdV} - \Vert u^{(N)} \Vert_{L^2}^2\right| 
\,  & \le c (1+ \tau^{-1/2}\Vert w \Vert_{L^2})^{N-1} (\tau^{-1/2} \Vert w \Vert_{L^2})  \Vert w \Vert_{H^{N+1}_\tau}^2
\\ & \le  c (1+ \tau^{-1/2}\Vert u \Vert_{H^{-1}_\tau})^{N-1} (\tau^{-1/2} \Vert u \Vert_{H^{-1}_\tau})  \Vert u \Vert_{H^{N}_\tau}^2.
\end{split} 
\end{equation} 
In the same fashion
\[ \left| \int \frac{\delta H^{\Gardner}_N}{\delta w} \phi dx \right| \le c (1+ \tau^{-1/2} \Vert w \Vert_{L^2} )^{2N-1} \Big( \tau^{-1/2} \Vert \phi \Vert_{L^2} \Vert w \Vert_{H^N_\tau}^2 
+ \tau^{-1/2} \Vert w \Vert_{L^2} \Vert w \Vert_{H^N_\tau}  \Vert \phi \Vert_{H^N_\tau}\Big).  \] 
Let $ \gamma $ be slowly varying of rate $ \tau $. We can localize this estimate and add up the intervals of length $ \tau^{-1}$ to obtain 
\begin{equation}\label{eq:gardnerboundN}   \Big\Vert \gamma^2 \big(\frac{\delta}{\delta w} H^{\Gardner}_N- w^{(2N)} \big) \Big\Vert_{H^{-N}_\tau} 
\le c\tau^{-N-\frac12} \Big(1+ \tau^{-1/2} \Vert w \Vert_{L^2}\Big)^{2N-1}    \Vert \gamma w \Vert_{H^N_\tau}^2. 
\end{equation} 
We claim the similar estimate 
\begin{equation}\label{eq:KdVboundN}   \Big\Vert \gamma^2 \big(\frac{\delta}{\delta u} H^{\KdV}_N- u^{(2N)} \big) \Big\Vert_{H^{-N}_\tau} 
\le c\tau^{-N+1/2} \Big(1+ \tau^{-1/2} \Vert u \Vert_{H^{-1}_\tau}\Big)^{N-1}     \Vert \gamma u \Vert_{H^{N-1}_\tau}^2. 
\end{equation} 
By Lemma \ref{lem:gardnerstructure} we have to bound a sum of monomials 
\[  \partial^{j_0} \prod_{k=1}^K u^{(j_k)}  \]
where $2\le  K \le N+1$, the total number of derivatives being 
$2(N+1-K)$, with at most half the derivatives on a single factor. 
We set $ u = \partial_x v + \tau v$ so that the total number of derivatives becomes $d=2(N+1)-K-j_0$ (we count factors $ \tau$ like derivatives).  By Lemma \ref{lem:estimatedifferentialmonomial2}, 
\[
\begin{split}
\tau^{1/2} \Big\Vert  \prod_{k=1}^K v^{(1+\alpha_k)} \Big\Vert_{H^{-1}_\tau}\, &  \le \Big \Vert \prod_{k=1}^K v^{(1+\alpha_k)} \Big\Vert_{H^{-1}_\tau} 
\\ & \le c \Vert v \Vert_{L^2}^{K-2}   \Vert v^{  (N+\frac12-\frac{K}4 -\frac{j_0}2 )} \Vert_{L^2}^2
\\ & \le c \tau^{-j_0-\frac{K-2}2} \Vert u \Vert^{K-2}_{H^{-1}_\tau} \Vert u \Vert^2_{H^{N-1}_\tau}  
\end{split}
\] 
which implies \eqref{eq:KdVboundN} is the same fashion as we proved \eqref{eq:gardnerboundN}. 

The estimates \eqref{eq:KdVboundN} and \eqref{eq:gardnerboundN} allow to define weak solutions to equations of the KdV and the Gardner hierarchy. 

\begin{proposition}  The following estimates hold for $ R \ge  1$ and an interval $I$ of length $1$
\begin{equation} \Big\Vert \sech^2(x/R)  \partial \frac{\delta}{\delta u} H^\KdV_N (u) \Big\Vert_{H^{-N-2}}  \le  c( \Vert u \Vert_{H^{-1}} \Big( \Vert \sech^2 u \Vert_{H^{N-1} } + \Vert \sech u \Vert^2_{H^{N-1}} \Big) \end{equation} 
\begin{equation} \Big\Vert \sech^2(x/R)  \partial \frac{\delta}{\delta w} H^\Gardner_N (w) \Big\Vert_{H^{-N-2}}  \le  c(\Vert w \Vert_{L^2} ) \Big( \Vert \sech^2 w \Vert_{H^{N} } + \Vert \sech w \Vert^2_{H^{N}} \Big)  \end{equation} 
The following estimates hold for $ R \ge  1$ and an interval $I$ of length $1$
\begin{equation} \Big\Vert \sech^2(x/R)  \partial \frac{\delta}{\delta u} H^\KdV_N (u) \Big\Vert_{H^{-N-3}}  \le  c( \Vert u \Vert_{H^{-1}}) \Big( \Vert \sech^2 u \Vert_{H^{N-2} } + \Vert \sech u \Vert^2_{H^{N-2}} \Big) \end{equation} 
\begin{equation} \Big\Vert \sech^2(x/R)  \partial \frac{\delta}{\delta w} H^\Gardner_N (w) \Big\Vert_{H^{-N-2}}  \le  c(\Vert w \Vert_{L^2} )( 1+ \Vert w \Vert_{L^\infty})  \Big( \Vert \sech^2 w \Vert_{H^{N-1} } + \Vert \sech w \Vert^2_{H^{N-1}} \Big)  \end{equation} 
\end{proposition} 

The proposition is a consequence of Theorem \ref{thm:formofkdv} and \eqref{eq:estha} resp \eqref{eq:esthb} of Lemma \ref{lem:estimatedifferentialmonomial2}.

The proposition allows to define weak solutions in natural regularity classes.

\begin{definition} 
Let $I=(a,b)$ be an open interval. We call 
\[ u \in L^\infty(I,H^{-1}) \qquad \text{ with }   u^{(N-1)}\in L^2_{loc}(I \times \R) \] 
a weak solution to the $N$th KdV equation if it satisfies the equation in the distributional sense. We call 
\[w \in L^\infty(I,L^2) \qquad \text{ with }   w^{(N)}\in L^2_{loc}(I \times \R) \] 
a weak solution to the $N$th Gardner equation if it satisfies the equation in a distributional sense. 
We call $ u \in L^\infty(H^{N-1}(\X)) $  a weak solution  to the $N$th KdV equation and  $ w \in L^\infty H^{N-1}$  a weak solution to the $N$th Gardner equation if they are weak solutions to th corresponding equation. 
\end{definition} 

We define  $L^2_u(I\times \R )\subset L^2_{loc}(I\times \R)$
by 
\[ \Vert f  \Vert_{L^2_u(I\times \R^n)} = \sup_k \Vert f \Vert_{L^2(I \times (k,k+1))}.\]
We will relate weak solutions to different Gardner equations to another, and to weak solutions to the KdV equation via the modified Miura map. Since 
\[
\begin{split} 
\Vert \partial_x \psi \Vert_{L^2}^2 + \int ( u\psi   - v \partial_x \psi)\psi dx 
\, & \ge \Vert \psi_x \Vert_{L^2}^2 - \Vert u \Vert_{L^2} \Vert \psi \Vert_{L^2}^{3/2} \Vert \psi_x \Vert^{\frac12}_{L^2} 
- \Vert v \Vert_{L^2} \Vert \psi \Vert_{L^2} \Vert \psi_x \Vert_{L^2} 
\\ & \ge  -\big( \frac12 \Vert v \Vert^2  + \frac3{2^{\frac{10}{3} }}  \Vert u \Vert^{\frac43}_{L^2}\big)
\Vert \psi \Vert_{L^2}^2 
\end{split} 
\]
we see that $u+ v_x $ lies in the range of the Miura map if 
\begin{equation}\label{eq:condtau}     \tau^2 >  \frac12 \sup_{t\in I} \Vert u(t) \Vert_{H^{-1}}^2  + \frac3{2^{\frac{10}{3} }} \sup_{t\in I}  \Vert u(t)  \Vert^{\frac43}_{H^{-1}} \end{equation}
Suppose that $ 0< \tau_1 < \tau_2$. Then $ w+ 2\tau_1 w + w^2$ is in the range of the $ \tau_2$ Miura map.

\begin{theorem} \label{thm:equivalenceweaksolutionsuw}  
Let $ 1\le \tau_1 < \tau_2$, assume  $ w_{1}\in L^\infty(I; L^2)$ with $w^{(N)}_1\in L^2_u(I \times \R)$  is a weak solutions to the $N$th $ \tau_1$ Gardner equation. Then $u= \partial_x w_{1} + 2\tau_1 w_{1} + w^2_{1}$ satisfies $u \in L^\infty(I; H^{-1}) $ and  $u^{(N-1)} \in L^2_u(I \times \R)$. Moreover 
 it is  a weak solution to the $N$th KdV equation. Define $w_2$ by 
\[  \partial_x w_2 + 2\tau_2 w_2 + w_2^2 = \partial_x w_1 + 2\tau_1 w_1 + w_1^2. \]
It satisfies $w_2 \in L^\infty(I; L^2)$ and  $ w_2^{(N)} \in L^2_u( I \times \R) $ if  this holds for $w_1$. Moreover it  is a weak solution to the $N$th $\tau_2 $ Gardner equation if and only if 
$w_1$ is a weak solution to the $N$th $\tau_1$ Gardner equation. 

Suppose that $u\in L^\infty(I; L^2)$ with $ u^{(N-1)} \in L^2_u(I \times \R) $ is a weak solution to the $N$th KdV equation, 
\[ -\partial^2 + u(t) + \tau^2 \] 
is positive definite uniformly in $t$ (which holds if $ \tau >0$ satisfies \eqref{eq:condtau}) and $w$ is defined by 
\[ w_x + 2\tau w + w^2 = u.\]
Then $w\in L^\infty(\R; L^2) $, $ w^{(N)} \in L^2_u(I \times \R) $.
If moreover 
$ w \in L^\infty$ or $ w \in L^\infty (I; H^{N})$ then $w$ is a weak solution to the $N$th Gardner equation. \end{theorem}

\begin{proof}[Proof of Theorem \ref{thm:equivalenceweaksolutionsuw}] 
\noindent{\bf Step 0: The spaces.}
Let $ w_1 \in L^\infty(I; L^2)$. Lemma \ref{lem:uw}  implies 
\[  u:= \partial_x w_1 + 2\tau_1 w_1 + w_1^2  \in L^\infty(I; H^{-1}(\R)) \]   and 
\[ -\partial^2+ u(t) + \tau_1^2 \]
is uniformly positive definite and hence in the range of the $ \tau_2 $ modified Miura map and $ w_2\in L^\infty(I; L^2) $ if $ w_1 \in L^\infty(I; L^2)$.
Since 
\[ \Vert f \Vert_{L^2_u(I \times \R) } \sim \sup_{x_0} \Vert \sech(\kappa( x-x_0))f \Vert_{L^2(I \times \R) }
\] 
we can apply  Proposition \ref{prop:equivalencewu} with $N$
and $ \gamma =\sech( \kappa (x-x_0))  $ and $ \kappa < \tau $
to  see that in addition 
\[ w^{(N)}_1 \in L^2_u(I \times \R) 
\Longleftrightarrow 
u^{(N-1)} \in L^2_u \Longleftrightarrow w_2^{(N)} \in L^2_u. \]

\noindent{\bf Step 1: Weak solution to Gardner define weak solutions to KdV.} 
We first prove that the modified Miura map maps weak solutions to $N$th Gardner equation to weak solutions to $N$th KdV.
Suppose that $w\in L^\infty(I; L^2(\R))$ with $w^{(N)} \in L^2_u(I \times \R) $ is a weak solution to the $N$th $\tau$ Gardner equation and $\tau_1 \ge \tau$.
We define $w_\varepsilon = J_\varepsilon w := j_\varepsilon * w$, for a mollifier $j_\varepsilon$, and 
\begin{equation}\label{eq:defuepsilon}
    u_\varepsilon = w_{\varepsilon,x} + 2\tau w_{\varepsilon} + w_\varepsilon^2
\end{equation}
We assumed that $w$ is a weak solution hence  $w_\varepsilon$ satisfies 
\[
    \partial_t w_\varepsilon = \partial_x J_\varepsilon \frac{\delta H_N^{\Gardner}}{\delta w}(w)
\]
Let $ \phi\in C^\infty_c(\R)$. Then 
\[ t \to \int \frac{\delta H_N^{\Gardner}}{\delta w} \phi dx \in L^1(I), \]
hence, since  $w$ is assumed to be a weak solution to the $N$th Gardner equation 
\[ (t \to \int w \phi dx ) \in W^{1,1}(I) \] 
and  for all $n >0$ 
\[ \sup_x \Vert \partial^n w_\varepsilon(t,x)\Vert_{W^{1,1}(I)} < \infty.  \]  We calculate using the chain rule for functions in $W^{1,1}$
\[
\begin{split}
    \partial_t u_\varepsilon &= \partial_t ( \partial_x w_\varepsilon + 2\tau w_\varepsilon + w_\varepsilon^2) 
    \\ & = ( \partial + 2 \tau +2 w_\varepsilon ) J_\varepsilon \partial_t w
\\ & =     
        (\partial_x + 2\tau + 2w_\varepsilon)\partial_x J_\varepsilon \frac{\delta H_N^{\Gardner}}{\delta w}(w).
\end{split}
\]
Since $w_\varepsilon$ and $u_\varepsilon$ are smooth for almost all $t \in I$, \eqref{eq:defuepsilon} implies
\[
    (\partial_x + 2\tau + 2w_\varepsilon)\partial_x \frac{\delta H_N^{\Gardner}}{\delta w} (w_\varepsilon) = \partial_x \frac{\delta H_N^{\KdV}}{\delta u}(u_\varepsilon).
\]
Using this identity and the fact that $J_\varepsilon$ commutes with the linear part of the equation, we obtain
\begin{equation} \label{eq:weakconvergenceeq} 
    \partial_t u_\varepsilon - \partial \frac{\delta H^{\KdV}_N}{\delta u}(u_\varepsilon) = (\partial + 2\tau + 2w_\varepsilon)\partial_x\big(J_\varepsilon F_N^{\Gardner}(w) - F_N^{\Gardner}(w_\varepsilon)\big).
\end{equation} 
The modified Miura map is continuous as a map from 
\[ L^\infty(I; L^2) \cap L^2 H^{N}_u \to L^\infty(I;H^{-1}) \cap L^2 H^{N-1}_u. \] 
By \eqref{eq:KdVboundN}  
\[ L^\infty(I;H^{-1}) \cap L^2 H^{N-1}_u \ni u \to  \frac{\delta H^{\KdV}_N}{\delta u } (u) \in  L^2  H^{-N}_u \]
hence in $L^\infty (H^{-N-1}_{loc})$ 
\[ \partial \frac{\delta H^{\KdV}_N}{\delta u}(u_\varepsilon) 
\to \partial \frac{\delta H^{\KdV}_N}{\delta u}(u).    \]
By construction $J_\varepsilon F_N^{\Gardner}\to F_N^{\Gardner} $
and $F^{\Gardner}_N(w_\varepsilon) \to F^{\Gardner}_N(w) $ in $L^2H^{-N}_u$ by the bound \eqref{eq:gardnerboundN}. This gives in the sense of distributions 
\[ (\partial+2\tau) \partial ( J_\varepsilon F_N^{\Gardner}(w) - F_N^{\Gardner}(w_\varepsilon)  ) \to 0 \]  
It remains to verify 
\[  w_\varepsilon \partial_x (   \big(J_\varepsilon F_N^{\Gardner}(w) - F_N^{\Gardner}(w_\varepsilon)\big)     \to 0       \]
in a distributional sense. We pull out the derivative and the claim follows from the bound 
\begin{equation}   \Vert \gamma^2 (\partial_x w_\varepsilon) (F_N^{\Gardner}(w)-w^{(2N)} )  \Vert_{L^1}
\le c  \tau (1+ \tau^{-\frac12}\Vert w \Vert_{L^2})^{2N} \Vert \gamma w \Vert^2_{H^N_\tau}.     \end{equation} 
Again we need do bound differential monomials of $K+1\ge 3$ factors with $d= 2N+3-K-1-l-\alpha_0$ derivatives with at most $1$ or half the derivatives on one factor, 
\begin{equation}  \begin{split} \tau^{l} \Big\Vert w_x \partial_x \prod_{k=1}^K w^{(\alpha_k)} \Big\Vert_{L^1}
\, & \le c\tau^{l} \Vert w \Vert_{L^2}^{K-1} 
\Vert w^{(N+\frac12 -\frac{K-1}4- \frac{l+\alpha_0}{2}  ) }\Vert^2_{L^2} 
\\ & \le   \tau^{1-\frac{K-1}2-\alpha_0}  \Vert w \Vert_{L^2}^{K-1} \Vert w \Vert^2_{H^N_\tau}. 
\end{split} 
\end{equation}  
Thus $u$ is a weak solution to the $N$th KdV hierarchy. 

\bigskip 

\noindent{\bf Step 2: Weak solutions to the $N$th Gardner equation with different $\tau $} 
Changing the notation slightly  we assume $ 2\le \tau  \le \tau_1$, that $w$ lies in the $N$th Gardner Kato smoothing space  and that it is a weak solution of the $N$ th $\tau$ Gardner equation,  and 
 \[ w^{\tau_1}_x + 2\tau_1 w^{\tau_1} + (w^{\tau_1})^2= w_x + 2\tau w  + w^2. \] 
 We want to prove  that $ w^{\tau_1}$ is a weak solution to the $N$th $\tau_1$ Gardner equation.
We regularize the solution as above $w_\varepsilon = j_\varepsilon * w =: J_\varepsilon w$.  Again 
 \begin{equation} \label{eq:wepsilon}
 \partial_t w_\varepsilon - \partial \frac{\delta H^{\Gardner}_N}{\delta w}( w_\varepsilon)  = \partial  \Big( J_\varepsilon  F^{\Gardner}_N(w) -   F^{\Gardner}_N(w_\varepsilon) \Big)  
  \end{equation} 
  Clearly $ w_\varepsilon $ is smooth in space and by \eqref{eq:gardnerboundN}
  \[ \Vert J_\varepsilon \partial (F^{\Gardner}_N - \partial^{2N} w) \Vert_{L^1_u(I\times \R)} \le c (\tau^{-1/2}  \sup_t \Vert w(t) \Vert_{L^2})  \Vert w \Vert^2_{L^2(H^N_{\tau,u})}. \] 
    Define $ w_\varepsilon^{\tau_1}$ by 
    \begin{equation} \partial  w^{\tau_1}_\varepsilon + 2 \tau_1 w^{\tau^1}_\varepsilon + ( w^{\tau_1}_\varepsilon)^2 = \partial_x w_\varepsilon + 2\tau w_\varepsilon + w_\varepsilon^2 \end{equation} 
We differentiate both sides of the  equation with respect to $t$, use the chain rule for $W^{1,1}$ functions for fixed $x$, and invert one operator to arrive at 
    \begin{align*}
    \partial_t w^{\tau_1}_\varepsilon &  = ( \partial + 2\tau_1 + 2 w^{\tau_1}_\varepsilon )^{-1} (\partial + 2\tau + 2 w_\varepsilon) \partial_t w_\varepsilon \\  &  = \partial_t w_\varepsilon + (\partial+ 2 \tau_1 + 2 w^{\tau_1}_\varepsilon )^{-1} (w_\varepsilon - w^{\tau_1}_\varepsilon +\tau-\tau_1) \partial_t w_\varepsilon. 
    \end{align*}
    We use the identity
     \[ \partial \frac{\delta H^{\KdV}_N}{\delta u} = (\partial +2 \tau + 2 w) \partial \frac{\delta H^{\Gardner}}{\delta w}  \]
     twice, once for $ \tau$ and then for $ \tau_1$ to see that 
\[ \partial \frac{\delta H^{\Gardner}_{\tau_1, N}}{\delta w}(w^{\tau_1}_\varepsilon) = 
    \Big\{ 1+     (\partial+ 2 \tau_1 + 2 w^{\tau_1}_\varepsilon )^{-1} (w_\varepsilon - w^{\tau_1}_\varepsilon + \tau - \tau_1)\Big\}  \partial \frac{\delta H^{\Gardner}_N}{\partial w} ( w_\varepsilon). 
\]
Altogether 
     \[ \partial_t w_\varepsilon^{\tau_1}- \partial \frac{\delta}{\delta w} H^{\Gardner}_{\tau_1,N} ( w^{\tau_1}_\varepsilon) = \Big\{ 1+     (\partial+ 2 \tau_1 + 2 w^{\tau_1}_\varepsilon )^{-1} (w_\varepsilon - w^{\tau_1}_\varepsilon + \tau - \tau_1)\Big\} 
 \partial  \Big( J_\varepsilon  F^{\Gardner}_N(w) -   F^{\Gardner}_N(w_\varepsilon) \Big)  \] 
By the continuity of $w \to \frac{\delta H^{\Gardner}}{\delta w}(w)$ the left hand side converges to 
\[ w_t - \partial \frac{\delta H^{\Gardner}_N}{\delta w} (w) \] 
in the sense of distributions. 
We will show that the right hand side converges to $0$ in the sense of distributions which implies that $ w^{\tau_1} $ is a weak solution to the $N$th $\tau_1 $ Gardner equation. 
We immediately turn to the most difficult term 
\[      (\partial+ 2 \tau_1 + 2 w^{\tau_1}_\varepsilon )^{-1} (w_\varepsilon - w^{\tau_1}_\varepsilon) 
\partial \Big\{ \Big( J_\varepsilon  F^{\Gardner}_N(w) -   F^{\Gardner}_N(w_\varepsilon) \Big)   \] 
which contains no linear term. 
By the same continuity argument as before it suffices that the map from $w$ to 
\[  \tau^l (\partial+ 2 \tau_1 + 2 w_1 )^{-1} w_2 \partial^{j_0+1}   \prod_{k=1}^K  w^{(j_k)}  
\] 
is continuous in $w_1$, $w_2$ and $w$. 
We want to pull the derivatives $ \partial^{(j_0+1)} $ recursively to the left.
In the first step we obtain a term were the derivative falls on $w_2$, and one with the derivative between the integral operator and $w_2$. We recall 
\[ ( \partial+ 2\tau + 2 w_1)^{-1} \partial = 1 - (\partial+ 2\tau + 2 w_1)^{-1}( 2\tau + 2 w_1). \]
By a repeated application of this argument we obtain a linear combination  of terms 
\[ \tau_1^{l}  \partial^{l_0} \prod_{k=1}^{K_1}  w_1^{(j_k)}
w_2^{(l_{K_1+1})} \prod_{k=K_1+2}^{K_1+K+1} w^{(j_k)} \] where 
\[ l +    \sum_{k=0}^{K+K_1+1} j_k  +K+K_1+1  = 2N+1  \] and 
\[ \tau_1^{l} (\partial+ 2\tau + 2 w_1)^{-1}  \prod_{k=1}^{K_1}  w_1^{(j_k)}
w_2^{(l_{K_1+1})} \prod_{k=K_1+2}^{K_1+K+1} w^{(j_k)} \]
where 
\[    \sum_{k=1}^{K+K_1+1} j_k  +K+K_1+1  = 2N+2. \]
 It suffices to bound the differential monomials in $L^1$, with $D= \sum_{k=0}^{K_1}l_k+ \sum_{k=1}^K j_k$
 the total number of derivatives
 and $M= K_1+1+K_0$ the number of factors
 \[ \Big\Vert \sech^2(x) \prod_{j=1}^{K_1}  w_1^{(l_k)}
w_2^{(l_0)} \prod_{j=1}^K w^{(j_k)}
 \Big\Vert_{L^1} \le c \prod_{j=1}^3  \left(  \Vert w_j \Vert^{M-2}_{L^2} \Vert  \sech(x) w_j  \Vert^2_{H^{D+\frac{M-2}4}} \right)^{\alpha_j}
\] 
with $ \alpha_1=\frac{K_1}{M}$, $\alpha_2 =\frac1{M} $, $ \alpha_3 =\frac{K}{M}$.
Altogether we arrive at
\[ \begin{split} \hspace{1cm} & \hspace{-1cm} 
\Big\Vert \gamma  (\partial + 2 \tau_1 + 2w_1)^{-1} w_2 \partial F^{\Gardner}_N \Vert_{L^1} \\ &  \le c(\Vert w_1 \Vert_{L^2}, \Vert w_2 \Vert_{L^2}, \Vert w \Vert_{L^2}) \Big( \Vert  \gamma w \Vert^2_{H^N_\tau}+ \Vert \gamma w_1 \Vert^2_{H^N_{\tau_1}}+ \Vert \gamma w_2 \Vert^2_{H^N_{\tau_1}}\Big).
\end{split} 
\]         
This completes the proof that if $w$ is a weak solution to the $N$th $\tau $ Gardner equation then the modified Miura maps define a weak solution to the $N$th $\tau_1$ Gardner equation, and vice versa. 

\bigskip 

\noindent{\bf Step 3: Weak solutions to KdV define weak solutions to Gardner under the additional assumption 
$ \Vert w \Vert_{L^\infty} < \infty$.} 
Now suppose that $u$ is a weak solution to the $N$th KdV equation. We want to prove that $w$ is a weak solution to the Gardner equation.  Let $ u_\varepsilon = J_\varepsilon u $ and define $w_\varepsilon$ by the Miura map, 
\[ \partial_x  w_\varepsilon + 2 \tau w_\varepsilon + w_\varepsilon^2 = u_\varepsilon. \]  
We apply the chain rule and argue as above to obtain 
\[
        \partial_t w_\varepsilon - \partial \frac{\delta H_N^{\Gardner}}{\delta w}(w_\varepsilon) = (\partial + 2\tau + 2w_\varepsilon)^{-1}\partial\big(J_\varepsilon F_N^{\KdV}(u)  - F_N^{\KdV}(u_\varepsilon)\big).
\]

The Gardner terms on the left hand side are covered by our previous considerations. 
Only the right-hand side needs a new consideration.  Again it suffices to provide bounds for 
\[ (\partial +2 \tau +2 \tilde w)^{-1} \partial^{1+j_0} \prod_{k=1}^K u^{(j_k)}. \] 
As above we pull derivatives to the left, making use of 
\[ (\partial + 2 \tau + 2 \tilde w)^{-1} \partial  = 1 - ( \partial+ 2\tau +2 \tilde w)^{-1} ( 2\tau + 2 \tilde w). \]  
We substitute $u=v_x+ 2\tau v$ and expand so that we obtain a linear combination of terms, and, by an abuse of notation we write $ v$ for $ \tilde w$,
\begin{equation} \label{eq:direct}   \tau^l \partial^{j_0} \prod_{k=1}^K v^{(j_k)} \end{equation} 
where 
\[ 2\le K, \quad \sum_{k=0}^K j_k =  2N+2 - K -l, \quad j_k \le N  \] 
and 
\begin{equation} \label{eq:inverse}   \tau^l(\partial +2 \tau +2v)^{-1}\prod_{k=1}^K v^{(j_k)} \end{equation}  
where 
\[ 2\le K, \quad   \sum_{k=1}^K j_k = 2N+3-K-l, \quad  j_k \le  N. \]  
The estimate for \eqref{eq:direct} is 
\[ \Vert  \tau^l \partial^{j_0} \prod_{k=1}^K v^{(j_k)} \Vert_{L^1_{u}(\R) }
\le c(\tau^{-1/2} \Vert v \Vert_{L^2})    \Vert v \Vert_{H^N_{\tau,u}}^2 \] 
and for \eqref{eq:inverse} 
\[ 
\begin{split} 
\Vert \tau^l(\partial +2 \tau +2 v)^{-1}  \prod_{k=1}^K v^{(j_k)}\Vert_{L^1_u} 
\, & \lesssim  \tau^{l-1} \Big\Vert \prod_{k=1}^K v^{(j_k)}\Big\Vert_{L^1_u} 
\\ & \lesssim       c( \tau^{-1/2} \Vert v \Vert_{L^2}) (1+ \tau^{-1}\Vert v \Vert_{L^\infty}) 
\Vert v \Vert_{H^N_{\tau,u}}^{2}. 
\end{split} 
\] 
Indeed the first inequality follows from writing $(\partial +2 \tau +2 v)^{-1}$ as an integral operator which is bounded as an operator on $L^1_u$. For the second inequality, note that by \eqref{eq:estha} from Lemma \ref{lem:estimatedifferentialmonomial2} we have
\[
    \tau^{l-1} \Big\Vert \prod_{k=1}^K v^{(j_k)}\Big\Vert_{L^1} \lesssim \tau^{l-1}\|v\|_{L^2}^{K-2}\|v^{(N-\frac{K-2}{4}-\frac{l-1}{2})}\|_{L^2}^2 \lesssim (\tau^{-1/2}\|v\|_{L^2})^{K-2}\|v\|_{H^N_\tau}^2,
\]
in all situations but when $l=0, K =3$ because of $N-\frac{K-2}{4}-\frac{l-1}{2} < N$ and since when $l = 0$, then $K \geq 3$ in \eqref{eq:inverse}. Now assume $j_1 + j_2 + j_3 = 2N, j_i \leq N$, and $j_1 \geq j_2 \geq j_3$. When $j_3 = 0$ then we can use Hölder's inequality directly, otherwise
\[
    \big\|v^{(j_1)}v^{(j_2)}v^{(j_3)}\big\|_{L^1} \leq \big\|v^{(j_1)}\big\|_{L^2}\big\|v^{(j_2)}\big\|_{L^2}\big\|v^{(j_3)}\big\|_{L^\infty} \lesssim \|v\|_{L^\infty}  \|v^{(N)}\|_{L^2}^2,
\]
from the Gagliardo-Nirenberg inequality. Performing a weighted summation implies
\[
    \tau^{-1}\big\|v^{(j_1)}v^{(j_2)}v^{(j_3)}\big\|_{L^1_u} \lesssim \tau^{-1}\|v\|_{L^\infty}\|v\|_{H^N_{\tau,u}}^2,
\]
which proves the second inequality.
\end{proof}

\subsection{The good variable hierarchy}

The fundamental Sobolev estimate in this setting is 
\begin{equation} \label{eq:fundamental}  \Big\Vert \prod_{k=1}^K  v^{(j_k)} \Big\Vert_{L^1_u} 
\le c \Vert v \Vert_{L^\infty}^{K-2}  \Vert v \Vert^2_{H^{\sigma}_u} \quad \text{ 
where } \quad 
 2\sigma = \sum_{j=1}^K j_k. \end{equation}  
The index $u$ refers to uniform Lebesgue or Sobolev spaces: We the supremum of the norm on intervals of length $1$.
Let $\frac12 < s \le 1$, $J$ an interal of length $1$ $ w \in L^\infty(J, H^{s-1}) \cap L^2(J, H^{s+N-1}_u) $.
In the last term there an abuse of notation: 
\[ \Vert w \Vert_{L^2 H^{s+N-1}_u} = 
\sup_{|I|\le 1 }  \Vert w^{(s+N-1)}\Vert_{L^2( J \times I)}. \] 
We claim that this suffices to define weak solutions to the $N$th Gardner equation on the time interval $J$.
We write $ w = \partial_x v - \tau v$ and we have to bound 
\[    \tau^l \partial^{j_0} \prod_{k=1}^K  \partial^{j_k} ( \partial_x v - \tau v) \] 
which can be written as a linear combination of terms (changing the meaning of $l,j_0,K$)  
\[ \tau^l \partial^{j_0} \prod_{k=1}^K v^{(j_k)} \] 
where
\[ l + \sum_{k_0}^{K} j_k = 2N+1 \] 
and no term carries more than $N$ derivatives. We obtain  a bound
\[  \Vert \prod_{k=1}^K v^{(j_k)} \Vert_{L^1_u} \le C \Vert v \Vert_{L^\infty}^{K-2} \Vert v \Vert_{H^{D}}^2 \]
with $ 2D = \sum_{j=1}^K j_k $.
The claims follow now from embeddings and estimates for the operator $ \partial-2\tau$.

We recall the structure of the $N$th equation good variables hierarchy. It can be written as
\[
    \psi_t - (-1)^{N}\psi^{(2N+1)} = \partial F_N^{GV}(\psi),
\]
where the structure of $F_N$ is the following:
\begin{lemma}\label{lem:gardnerstructuregood} 
\[
\begin{split}
    &F_N^{\operatorname{GV}}(w) = \sum_{K,(j_k)_k,l,M}c_{K,(j_k)_k,l,M}\partial^{j_0}  (1+v)^{-M}\tau^{2l}\prod_{k=1}^K v^{(j_k)},\\
    0 &\leq M \leq 2N+1, \qquad 2 \leq K \leq 3N+1  
    \qquad  \sum_{k=0}^K j_k  = 2N-2l,
    %\\
    %1 &\leq i \leq 2N-l-1 - \Big[\frac{2N-l}{2}\Big], 
    \qquad j_k \leq N-l.
\end{split}
\]
\end{lemma}
\begin{proof}
 $F_N^{GV}(v)$ is a sum over $(1+v)^{-M}\tau^l \prod_{k=1}^K v^{(j_k)}$ with (cf. Theorem \ref{thm:formofgoodvariableequation})
\[
\begin{split}
    0 \leq M &\leq 2N-1, \qquad 2 \leq K \leq 2N+1, \\
    \qquad l + \sum_{k=1}^K j_k &= 2N, \qquad   \sum_{k = 1}^K j_k \geq M+1 \text{ if } M \geq 1.
\end{split}
\]
For the good variable equation, we proceed as in Lemma \ref{lem:formofgoodvariableequation}.\end{proof}

Let $ s > \frac12$. Then we can define the notion of a weak solution to the good variable equation under the regularity assumption 
\[  v \in L^\infty(I, H^s) \cap L^2 (I,H^{s+N})_u. \] 

\begin{theorem}\label{thm:vwweak}
 Let $s>\frac12$,$ \tau\ge 2$, $I$ an open interval. Suppose $ v \in L^\infty(H^s) \cap L^2 H^{s+N}_u$ satisfies $ v > -1$. Let  $ w=  v -\frac12\partial \ln v $. Then $ w \in L^\infty(H^{s-1}) \cap L^2 H^{s-1+N}_u $. Vice versa: Suppose that 
 $ w \in L^\infty H^{s-1} \cap L^2 H^{s-1+N}_u$. Then there exists a unique $v \in L^\infty$ with $v>-1$ and $ (v+1)^{-1} $ which satisfies 
 \[  \tau v - \partial_x (\ln( 1+ v)) = w. \] 
 Moreover $ v \in L^\infty H^s \cap L^2 H^{s+N}_u $. 
 Under these assumptions $w$ is a solution to the $N$th $ \tau $ Gardner equation iff $v$ is a weak solution to the $N$th $ \tau$ Good variable equation. 
 \end{theorem}
\begin{proof}
We argue in the same fashion for the second part. We recall the equation for $v$ (see Theorem \ref{thm:formofgoodvariableequation} 
\[  v_t = (-1)^N v^{(2N+1)} + \partial_x F^{gv}_N \] 
with    
\[ \partial \frac{\delta H^{\Gardner}_N}{\partial w} \Big|_{w= \tau v - \frac12 \partial \ln(1+v) } = 
    \Big( \tau - \partial \frac1{2(1+v)}\Big) \partial \left( F^{GV}_N(v)+ \partial^{2N} v \right) . \] 
    Suppose that $v$ is a weak solution the $N$th equation and   set $v_\varepsilon = J_\varepsilon v$, $ w_\varepsilon= \tau v_\varepsilon -\frac12 \partial \ln(1+v_\varepsilon)$.  Then 
\[ 
 \partial_t   w_\varepsilon - \partial \frac{\delta H^{\Gardner}_N}{\delta w}(w_\varepsilon) 
= \Big( \tau -  \partial \frac1{2(1+v_\varepsilon)}\Big)   \partial ( J_\varepsilon F^{GV}_N(v)- F_N^{GV}(v_\varepsilon)). 
\] 
The left hand side converges to 
\[ \partial_t w - \partial \frac{\delta H^{\Gardner}}{\delta w} \] 
in a distributional sense. We claim that the right hand side converges to zero. 
To prove that is suffices to consider convergence for summands as in Lemma \ref{lem:gardnerstructuregood}: 
\[ \tau^{2l} \Big( \tau -  \partial \frac1{2(1+v_\varepsilon)}\Big)   \partial^{1+j_0}  \Big( (1+v)^{-M}\prod_{k=1}^K v^{(j_k)}      \Big). \]  
Again we pull the derivatives in front. We obtain a sum of terms 
\[  \tau^L \partial^{j_0} (1+v)^{-M} \prod_{k=1}^K v^{(j_k)}\] 
where 
\[  0 \leq M \leq 2N+3, \qquad 2 \leq K \leq 3N+2  
    \qquad  \sum_{k=0}^K j_k  = 2N+1-L
    \qquad j_k \leq N-L/2.  \] 
We may ignore the powers of $\tau$ and of $(1+v)$ for the question of distributional convergence, which now follows from \eqref{eq:fundamental}

Suppose that $w$ is a weak solution to the $N$ Gardner equation, $ w_\varepsilon = J_\varepsilon w$ and $v_\epsilon$ satisfies $ \tau v - \frac12 \partial \ln(1+v)$. Then 
\[ 
 \partial_t  v_\varepsilon - (-1)^N v^{(n+1)}_\varepsilon -\partial F^{GV}_N(v_\varepsilon)) =   \Big(\tau-\partial \frac1{2(1+v_\varepsilon)}\Big)^{-1} \partial 
\Big\{ J_\varepsilon  F^{\Gardner}_N(w ) - F^{\Gardner}_N (w_\varepsilon) \Big\}.  
\]
The convergence of the left hand side to the good variable equation follows from the bounds which ensure that the notion of a weak solutions is well defined. 
  We again have to provide bounds for 
\[ \Big(\tau-\partial \frac1{2(1+v_\varepsilon)}\Big)^{-1} \partial 
 F^{\Gardner}_N(w ),  \] 
more precisely for its summands 
\[ \Big(\tau-\partial \frac1{2(1+ v_\varepsilon)}\Big)^{-1} \partial^{1+m}   \prod_{k=1}^K w^{(j_k)}. \] 
We write $ w = \tau \tilde v- \tilde v_x$ and ignore the difference between $ \tilde v$, $v$ and $v_\varepsilon$ in the sequel. 
We pull derivatives to the left:
\[
\begin{split} 
\Big(\tau-\partial \frac1{2(1+v)}\Big)^{-1} \partial \, & = \Big(\tau-\partial \frac1{2(1+ v)}\Big)^{-1} \partial (2(1+ v))^{-1} (2(1+ v)) 
\\ & = - 2(1+ v) +  \tau   \Big(\tau-\partial \frac1{2(1+v)}\Big)^{-1}  (2(1+\tilde v)). 
\end{split} 
\] 
Consider 
\[  \tau \phi + \partial \frac{\phi}{(1+v)} = f. \] 
Let $ \psi = \frac{\phi}{1+v}$ and rewrite it as 
\[ \partial \psi + \tau ( 1+ v) \phi = f \] 
and deduce 
\[ \Vert \phi \Vert_{L^\infty} \le \Vert \psi \Vert_{L^\infty} \Vert 1 + v \Vert_{L^\infty} \le c \Vert 1+ v \Vert_{L^\infty} \Vert (1+ v)^{-1} \Vert_{L^\infty}
\Vert f \Vert_{L^1_u}. \] 
The final estimate now follows from \eqref{eq:fundamental}. 
\end{proof}

\section{Regularity of weak solutions and Kato smoothing} 
\label{sec:tightness} 

We analyze the structure of the conservation laws for smooth solutions to the $N$th KdV  Gardner equation. Moreover, we use the conservation laws in their weak form to show smoothing estimates and precompactness of orbits for weak solutions.

\subsection{The energy-flux identity}

We obtain the energy flux identity for smooth solutions by multiplying the equation by%  
$ \eta w$ where $\eta$ is a test function in $\R^2$:
\[ 0 = \int_{\R^2}  \eta w ( w_t - \partial_x \frac{\delta}{\delta w} H^{\Gardner}_N ) dx dt = \frac12 \int_{\R^2} -w^2 \eta_t  + \frac{\delta }{\delta w} H^{\Gardner}_N \partial_x( w \eta)  dx dt. \] 
We rewrite the last term for fixed $t$
\[ \int \frac{\delta}{\delta w} H^{\Gardner}_N
\partial (\eta w) dx = \frac{d}{ds} H^{\Gardner}_N(w+ s\eta  \partial_x  w) )\Big|_{s=0}
+ \int  \frac{\delta}{\delta w} H^{\Gardner}_N
w \partial_x \eta dx.  \] 
We study the first term on the right hand side for differential monomials $I= \int \prod w^{(j_k)} dx $ 
\[\begin{split} \frac{d}{ds}  I(w+ s \eta w_x) |_{s=0}
\, & = \int \eta \sum_{k=1}^k  w_x^{(j_k)} \prod_{i\ne k } w^{(j_i)} dx 
= \int \eta \partial_x  \prod_{k=1}^K  w^{(j_k)} dx 
\\ & = 
- \int \eta_x \prod_{k=1}^K w^{(j_k)}  dx 
.  \end{split}   \] 
We obtain the flux 
\begin{equation}\label{eq:flux}   Fl_N (w) =       \sum_{K=2}^{2N+2} \sum_{l=0}^{2N+2-K} \partial^{j_0}  \sum_{\sum j_k = 2N+2 -K-l } c_{K,l,j_k} \tau^l    \prod_{j=1}^K w^{(j_k)}  \end{equation}  
which has almost the same structure as the energy density so that the energy flux equation  
\begin{equation}\label{eq:energyflux}  \partial_t \frac12 w^2 = \partial_x Fl_N(w) \end{equation} 
holds. In particular the quadratic part of the flux is 
\[
   \partial Fl_{N,2} = 2w w^{(2N+1)} =\partial\big((2N+1)|w^{(N)}|^2 + \sum_{k=1}^n a_k \partial^{2k} |w^{(N-k)}|^2\big),
\]
for some combinatorical constants $a_k \in \R$.

\begin{lemma}\label{lem:energyfluxweak} Suppose that $ w \in L^\infty(\R, L^2 ) $ with $  w^{(N)} \in L^2_u( \R \times  \R) $ is a weak solution to the $N$th Gardner equation. Then
  \[ \partial_t w^2 = \partial_x Fl_N(w)  \]
  in the sense of distributions.
\end{lemma}

\begin{proof}
  Let $ \rho \in C^\infty_c ( (-1,1)) $ be nonnegative and even  with $ \int \rho dx =1$. We define $\rho_\varepsilon= \varepsilon \rho( x/ \varepsilon) $ and $
  J^x_\varepsilon f  = \rho_{\varepsilon} *_x f$ and similarly $J^t_\varepsilon$ 
  is the convolution with respect to time.
 We recall that 
 weak solutions $w$ satisfy 
 $ J_\varepsilon^x w \in W^{1,1}(I, L^2(J) )$ 
  whenever $I,J$ are compact intervals. As a consequence 
  \[  t \to \int w(t) \phi dx \] 
  is continuous for every test function $ \phi$, or, equivalently, after modifications on the set of times of measure zero $ t \to w(t) \in L^2$ 
  is weakly continuous. 
  
  Let $w$ be a weak solution satisfying the regularity assumptions. Then with $w_\varepsilon= J^t_{\varepsilon_1} J^x_{\varepsilon_2} w $, similar to above  
  \begin{equation} \label{eq:fluxreg} \begin{split} 0\, & = \int \Big(w_t- \partial \frac{\delta H^{\Gardner}_N}{\delta w}( w) \Big)  J_{\varepsilon_1}^t J_{\varepsilon_2}^x (\phi  w_\varepsilon) dx dt  \\ &  = \int  - \frac12 w_\varepsilon^2 \partial_t \phi  +  \Big(J^t_{\varepsilon_1}J^x_{\varepsilon_2} \frac{\delta H^{\Gardner}_N}{\delta w} (w)\Big) \Big(w_\varepsilon \partial_x \phi+ (\partial_x w_\varepsilon) \phi \Big)  \, dx dt.
  \end{split} 
    \end{equation}  
  All terms converge as the time regularization $ \varepsilon_1$ tends to $0$. We set $\varepsilon_1=0 $ and write $ \varepsilon= \varepsilon_2$ in the sequel.
  By the continuity in $ \varepsilon$ we have to justify 
  \[  \int_{\R^2}  \Big(J^t_{\varepsilon_1}J^x_{\varepsilon_2} \frac{\delta H^{\Gardner}_N}{\delta w} (w)\Big) \Big(w_\varepsilon \partial_x \phi+ (\partial_x w_\varepsilon) \phi \Big)  \, dx dt   \to \int FL_N \partial_x \phi dx dt.     \] 
  We obtain for the linear resprectively quadratic term 
  \[ \int -\frac12 w_\varepsilon^2 \partial_t \eta + FL_{N,2}(w_\varepsilon) \partial_x \eta\,  dx\, dt  \to \int - \frac12 w^2+ Fl_{N,2}(w) \partial_x \eta\, dx \, dt  \] 
  since   $ w^{(j)}$ is in $L^2_{loc}$ of space time for $j \le N$. 
   
  Let $ h(w)=\tau^l \prod_{k=1}^K w^{(j_k)} $ be monomial in the energy density up to a coefficient  with $K \ge 3$ and let $FL$ be the corresponding part of the flux. Let $D$ be the total number of derivatives. Then 
  \[ K+D +l= 2N+2 \]
Then 
  \[  \lim_{\varepsilon \to 0 } \int_{\R^2}  J_\varepsilon h(w) ( w_\varepsilon  \partial_x \phi + \partial_x w_\varepsilon \phi) \,  dx\, dt 
  = \int_{\R^2}   h(w) ( w  \partial_x \phi + \partial_x w \phi dx) \, dt. 
  \] 
  Summation and comparison with \eqref{eq:fluxreg} 
  implies  the energy flux identity. 
  \end{proof}

\subsection{Consequences of the energy-flux identity: Kato smoothing for weak solutions} \label{subsec:Katosmoothing}

Weak solutions satisfy the energy flux identity
\eqref{eq:energyflux} by Lemma \ref{lem:energyflux} which can be written as 
\begin{equation} \label{eq:integrated}   \int w^2 \partial_t \eta - \flux_N(w) \partial_x \eta dx dt = 0 \end{equation}  
for smooth compactly supported functions $ \eta$.  
We say that $u\in L^\infty(I; H^{-1}) $
with $ u^{(N-1)} \in L^2_{loc}( I \times \R) $
lies in the $N$ th KdV Kato smoothing space if  for all compact subintervals $J\subset I$
\begin{equation}\label{eq:KdVKato}
\limsup_{x_0 \to \pm \infty} \Vert u^{(N-1)} \Vert_{L^2(J\times (x_0-1,x_0+1))} = 0. 
\end{equation}
We say that $w\in L^\infty(I; L^2) $
with 
$ w^{(N)} \in L^2_{loc}( I \times \R) $ 
lies in the $N$ the Nth Gardner  Kato smoothing space if  for all compact subintervals $J\subset I$
\begin{equation}\label{eq:GardnerKato} 
\limsup_{x_0 \to \pm \infty} \Vert w^{(N)} \Vert_{L^2(J\times (x_0-1,x_0+1))} = 0. 
\end{equation}

The condition \eqref{eq:KdVKato} is equivalent to 
\[
\limsup_{x_0 \to \pm \infty} \Vert \sech(\sigma( x-x_0)) u \Vert_{L^2(J;H^{N-1}(\R))} = 0. 
\]
and \eqref{eq:GardnerKato} is equivalent to 
\[
\limsup_{x_0 \to \pm \infty} \Vert \sech(\sigma(x-x_0)) w \Vert_{L^2(J;H^{N}(\R))} = 0. 
\]
If $ u = w_x + 2\tau w + w^2 $ (for fixed time)  and $ \sigma< \tau $
then 
\[ \Vert \sech(\sigma(x-x_0)) w \Vert_{H^N}  \sim \Vert \sech(\sigma(x-x_0)) u \Vert_{H^{N-1}}    \]
with constants depending on $\tau$ and $ \Vert w \Vert_{L^2} $.

\begin{lemma}\label{lem:L2cont}  Let $ N \ge 1$, $ \tau\ge 1$.  Suppose  that $w$  in the $N$th Gardner Kato smoothing space (\eqref{eq:KdVKato}, \eqref{eq:GardnerKato}) on $\R$ 
is a weak solution to the $N$th Gardner equation.
Then  $ w \in C(\R,L^2(\R))$ and there exists $\rho=\rho(\tau^{-1/2} \Vert w \Vert_{L^2})$ so that  it satisfies with $ w_0 = w(0)$ the Kato smoothing estimate for $0 < \kappa \le \tau$
\begin{equation}
\begin{split} 
\hspace{2cm} & \hspace{-2cm} 
\sup_t \int (1+ \tanh(\kappa  (x- x_0- \rho  \tau^{2n} t)) |w(t,x)|^2 dx     
 \\ &  + \kappa \int_0^\infty   \int \sech^2(\kappa (x-x_0-\rho \tau^{2n} t))   
(|w^{(n)}|^2+ \tau^{2n} |w|^2)   dx dt \\  \lesssim  \, &  \int \big(1+\tanh(\kappa( x-x_0))\big)w^2_0 \, dx.
\end{split} \label{localsmoothinghigherorder} 
	\end{equation} 
	The implicit constant depends only on $N$.
\end{lemma}
\begin{proof} 
Let $ \eta_j\in C^\infty_c$, $j=0,1$ and assume
that $\eta_0$ is nonnegative and satisfies  $ \eta_0(0)=1$. We  define for $R>0$
\[ \eta(t,x) = \eta_1(t) \eta_0(x/R). \]
Then 
\[ \lim_{R \to \infty } \int_{\R^2}  w^2 \partial_t \eta - FL_N \partial_x \eta dx = \int_{\R} \Vert w(t) \Vert_{L^2}^2 \partial_t \eta_0(t) dt \] 
 hence $ \Vert w \Vert_{L^2}$ is independent of time (up to sets of measure $0$).
 We have seen that $t\to w_x^\varepsilon\in L^2$ is weakly continuous which together with norm continuity implies that $ w \in C(\R; L^2)$. 
 
 We want plug in 
 \[ \eta = \chi_{[0,T]} (1+\tanh(\kappa (x-x_0-\rho \tau^{2n} t))) 
 \] 
 into the integrated energy-flux identity \eqref{eq:integrated}. There are two obstacles: The characteristic function is not smooth, and $1+\tanh$ is not compactly supported. We deal with the second problem by multiplying in addition by $\eta_0(x/R)$ as above, and with the first obstacle by regularization. 
  We obtain 
\[ 
\begin{split} \hspace{2cm} & \hspace{-2cm} 
\frac12  \int (1+ \tanh( \kappa (x-x_0-\rho \tau^{2n} t))) w^2 dx \Big|_0^T 
\\ &  =   -\kappa \int_0^T   \int   \sech^2(\kappa(x-x_0-\rho \tau^{2n} t))) \big( \frac{\rho \tau^{2n}}2 w^2    -   Fl_N(w) \big)  dx.
\\ & = : \int_0^T A(t) dt 
\end{split} 
\] 
The quadratic term on the RHS is - omitting the argument of $\sech^2$ as well as the time integration -
\[ - \kappa \int \sech^2 \big( \frac{\rho \tau^{2n}}2 w^2+ (2N+1) |w^{(N)}|^2\big) dx 
+ a_k \sum_{k=1}^N (\sech^2)^{2k+1}) |w^{(N-k)}|^2 dx \]  
and 
\[\begin{split}  \left|  \int \sech^2  Fl_{N\ge 3} dx \right|\, &  \le c (1+ \Vert w \Vert_{L^2}^2)^N 
\Vert \sech w \Vert_{H^N_\tau} \Vert \sech w \Vert_{H^{N-1}_\tau} 
\\ & \le             \Vert \sech(x) w \Vert_{H^N_\tau}^2 
+ c^2 (1+ \Vert w \Vert_{L^2})^{2N^2}   \Vert \sech^2 w \Vert^2_{L^2}) 
\end{split} 
\]
Together the integrand can be estimated
\[ A(t) \le  -\kappa  \int \sech^2( \rho \tau {2n}  w^2 + 2N     |w^{(N)}|^2 dx + c^2(1+\Vert w \Vert_{L^2}^{N^2}) \Vert \sech w \Vert^2 + \Vert \sech w \Vert_{H^{N-1}_\tau}^2.    \] 
We choose 
\[ \rho  \ge (C + c^2 (1+ \Vert w \Vert^2_{L^2})^{N^2+1} \tau^{-2N}  ) \]  and subtract  the second and the third term on the RHS from both sides.
This gives  \eqref{localsmoothinghigherorder}.
\end{proof}

\subsection{Precompactness of orbits of weak solutions}\label{subsec:precompactness}

In this section we show that the orbit of weak solutions  of the Gardner hierarchy with initial data in a  precompact subset of $L^2(\R)$ are precompact over compact time intervals.

\begin{theorem}\label{thm:gardnerprecompactness}
    Let $Q \subset \Sc(\R)$ be a precompact subset of $L^2(\R)$, let $A$ be a set of weak solutions with initial data in $Q$ and let $I$ by a bounded interval. Then 
    \[ \{  w(t): w\in A, t\in I \} \subset L^2\] 
    is precompact. 
    \end{theorem}

This will be essential in the proof of our main theorem, as we want to upgrade weak to strong convergence using a compactness argument.
To check precompactness, one can check boundedness, equicontinuity, and tightness.

\begin{proof} We begin with boundedness. Since $Q$ is precompact it is bounded and there exists $R>0$ so that 
\[ \Vert u_0 \Vert_{L^2} < R \] 
for $u_0\in Q$. Weak solution conserve the $L^2$ norm by Lemma \ref{lem:L2cont}.   

We turn to tightness on the  the right. Let $ \varepsilon>0$.
Since $Q$ is precompact and hence tight there exists $R$ so that for $x_0>R$ 
\[ \sup_{w_0 \in Q} \int (1+ \tanh(\kappa(x-x_0) ) |w_0|^2dx <\varepsilon.   \]
By Lemma \ref{lem:L2cont} again 
\[ \sup_{w_0\in Q} \sup_{t} 
\int ( 1+ \tanh( \kappa ( x-x_0-\gamma \tau^{2N} t) )) |w(t)|^2 dx \le c \varepsilon. \] 
Given $ \varepsilon>0$ and a bounded time interval $I$ we find $x_0$ so that \begin{equation}   \sup_{t\in I}  \sup_{w_0 \in Q}  \Vert w  \Vert_{L^2(x_0,\infty)} <\varepsilon. 
\end{equation}  

We turn to the proof of equicontinuity.
The set 
\[ Q_u = \{ w_x+2\tau w + w^2 : w \in Q \} \] 
is precompact since $Q$ is. If $ \tau_1 \ge \tau$ it lies in the range of the $ \tau_1$ Miura map. $Q_v$  is equicontinuous hence there there exists $ \tau_1$ so that 
\[ \Vert u \Vert_{H^{-1}_{\tau_1}} < \varepsilon \] 
Let $ Q_{\tau_1} = M_{\tau_1}^{-1}(Q_u) $, Then 
\[ \Vert w_0 \Vert_{L^2} \le 2\varepsilon \] 
for all $ w_0 \in Q_{\tau_1}$.
Under these mapping every weak solution to the $\tau $ Gardner equation is mapped to a weak solution to the $ \tau_1 $ Gardner equation. This equation preserves the $L^2$ norm, hence
\[ \Vert w_x(t) + 2 \tau w^{\tau_1}(t)+ (w(t))^2 \Vert_{H^{-1}_{\tau_1}}
=\Vert w^{\tau_1}_x(t) + 2 \tau_1 w^{\tau_1}(t)+ (w^{\tau_1}(t))^2 \Vert_{H^{-1}_\tau} < \varepsilon \] 
and the orbits are equicontinuous on compact time intervals. 

The Kato smoothing estimate for $ w^{\tau_1}$ implies 
\[ \sup_{x_0} \kappa\int_0^\infty  \Vert  \sech(\kappa (  x-x_0 - \gamma \tau_1^{2N} t) )  w \Vert^2_{H^N_{\tau_1}} dt \le \varepsilon,  \]   
hence 
\[ \sup_{x_0} \kappa\int_0^\infty  \Vert  \sech(\kappa (  x-x_0 - \gamma \tau_1^{2N} t) )  w_x +\tau w \Vert^2_{H^{N-1}_{\tau_1}} dt \le c\varepsilon\] 
for all weak solutions with initial datum in $Q$.
This implies the high frequency bound
\begin{equation} \label{eq:highfrequency}  \sup_{x_1} \int_0^\infty \int \sech^2 (\kappa_1( x-x_1-\rho \tau_1^{2N}  t) ) 
 (| w^{(N)}|^2 + \tau_1^{2(N-1)}(|w_x|^2+w^2) )   dx dt < \frac{\varepsilon}{\kappa_1} \end{equation}
 for all weak solutions $w$ to the $N$-Gardner equation  with $w(0) \in Q$
 and $ \kappa_1 $ much smaller than $\tau$.

  With the high frequency estimate in place we turn to tightness on the left. We claim in the setting of the theorem: 
 Given $ I=[0,T]$ and  $ \varepsilon>0$ there exists $x_0$ so that 
 \begin{equation} \label{eq:tightnessleft} 
 \sup_{0\le t \le T } \sup_{w_0 \in Q} 
 \Vert w(t) \Vert_{L^2(-\infty, x_0) } \le \varepsilon \end{equation} 
To simplify the notation we set $T=1$.
We use 
\[ \eta(t,x) = \chi_{[0,T]} \Big( 1- \tanh( \kappa (x-x_0+\rho \tau^{2n} t))\Big) \]  with $\kappa = \delta  \tau^{-2n}$
in the integrated energy-flux identity. 
Again this function is not admissible since it is not a test function, neither in $t$ nor in $x$. In the same way as for the $L^2$ conservation we obtain for $ T \le 1$
and $x_0 \in \R$.
\begin{equation}\label{eq:wrongside}   
\int \eta  w^2 dx \Big|_0^T
 \le c \kappa_0    \int_0^1   \Vert  \sech(\tau^{-2N} (x-x_0+\rho\tau^{2n} t))) w^{(N)}  \Vert_{L^2}^2.
\end{equation} 
We choose $ \tau_1$ so that in the notation above 
\[ \Vert w^{\tau_1}(0) \Vert_{L^2} \le \varepsilon.\] 
By tightness of $Q$ we can choose $R$ so that the initial term  satisfies 
\[ \int (1+ \tanh( \tau_1^{-2N}  ( x-x_0+ \tau_1^{2N} )) ) |w|^2 dx < \varepsilon \] 
We set $x_1 = x_0 - 3 \tau_1^{2N}$ and 
$ \kappa_1 = \frac19 \tau_1^{-2N} $ in \eqref{eq:highfrequency} to estimate the  
right hand side in \eqref{eq:wrongside} 
by $ c \varepsilon$. 
\end{proof}

\section{Smoothing and convergence for the difference flow}
\label{sec:differenceflow}

We recall the formula \eqref{eq:GardnerWadati} 
\[ \T_{-1}^{\Gardner}(z,w,\tau)= \frac{-iz}{4z^2+4\tau^2}
\ln{\det}_2(1+ K).
\] 
where to shorten the notation we define 
\begin{equation} \label{eq:Kdef} 
\begin{split} 
K(w) \, & = \Big( \begin{matrix}0 & (-iz-\partial)^{-1} w \\ -(-iz +\partial)^{-1} w & (-iz+\partial)^{-1} 2\tau (-iz-\partial)^{-1}\end{matrix} \Big)
\\ & =  \Big( \begin{matrix}0 & (-iz-\partial)^{-1} w \\ -(-iz +\partial)^{-1} w & -\frac{\tau}{iz} \big[(-iz+\partial)^{-1}+  (-iz-\partial)^{-1}\big]w \end{matrix} \Big)
\end{split} 
.\end{equation} 
Then 
\[  -\ln {\det}_2 \big( 1+ K(w)\big)   = \sum_{n=2}^\infty \frac{(-1)^n}n  \trace{K(w)}^n
\]
and 
\[ \T_{N}^{\Gardner}(z,w,\tau) = \frac{-iz(2z)^{2N+2}}{4z^2+4\tau^2}
\ln{\det}_2\big( 1+ K(w) \big)- \sum_{n=0}^N (2z)^{2(N-n)}H^{\Gardner}_n(w,\tau) .
\]
In this section we prove Proposition \ref{prop:diffham}, in particular the bound on $ \T^{\Gardner}_N$ \eqref{eq:TnGardner}, the characterization of equicontinuity of subsets $Q$ in  $H^N$ \eqref{eq:propequi} and the weak convergence statement \eqref{eq:propuni}. As a consequence we also obtain the corresponding bound \eqref{eq:tauNdec} and \eqref{eq:approximate}. 
A central step consists in the  study of $ (2z)^{N+1}\trace K(w)^n $.

\subsection{Schatten classes and the case of large \texorpdfstring{$n\ge N+2$}{n>N+1}}
We first deal with the sum over large $n\ge 2N+4$  and collect simple estimates for the operator $K$. The central estimates are 
\begin{lemma} \label{lem:largen}
Let $N\ge 3$,  $N +2 \le n $ and $ \im z \ge \tau$. Then 
 \begin{equation} \label{eq:largen1} 
\begin{split} (\im z)^{N+1} |\trace K^n(w) | \,& \le c  ((\im z)^{-1/2} \Vert w \Vert_{L^2})^{n-N-2} \Vert w \Vert_{L^{N+2}}^{N+2} \\ 
& \le  c  ((\im z)^{-1/2} \Vert w \Vert_{L^2})^{n-N-2} \Vert w \Vert_{L^2}^{N}\Vert w^{(N/4)} \Vert^2_{L^2}
\\
&\le  c  ((\im z)^{-1/2} \Vert w \Vert_{L^2})^{n-N-2} \Vert w \Vert_{L^2}^{N+1}\Vert w^{(N/2)} \Vert_{L^2} 
\end{split}
\end{equation}
where $ \Vert w^{(s)} \Vert_{L^2}$ is defined using the Fourier transform.
\end{lemma} 

Since for $ N \ge 4$
\[ \int |w|^{N+2} dx \le \Vert w \Vert_{L^2}^2 \Vert w \Vert_{L^\infty}^N 
\le  \Vert w \Vert_{L^2}^{2+\frac{N}2} \Vert w_x \Vert_{L^2}^{\frac{N}2} 
\le \Vert w \Vert_{L^2}^{2+\frac{N}2} \Big( \Vert w \Vert^{1-\frac4N}_{L^2} \Vert w^{(N/4)} \Vert^{\frac{4}N}_{L^2} \Big)^{\frac{N}2},  
\]
\[   \Vert w^{(N/4)} \Vert_{L^2}^2 \le \Vert w \Vert_{L^2} \Vert w^{(N/2)} \Vert_{L^2} \]
and 
\[  \int |w|^3 dx \le \Vert w \Vert \Vert  w^{(1/2)} \Vert_{L^2}^2\]
it suffices to prove the first inequality.

To this end recall the Schatten classes $\I_p$ of compact operators $A$ with $l^p$-summable singular values,
\begin{equation}\label{eq:Schattenclass}
    \|A\|_{\I_p}^p = \sum_{n=1}^\infty \mu_n(A)^p.
\end{equation}
Here we always consider $L^2$ and its powers as the underlying Hilbert space. Special cases are the trace class operators $\I_1$, the Hilbert-Schmidt operators $\I_2$, and the class of compact operators $\I_\infty$ which are defined with the obvious adaption in \eqref{eq:Schattenclass}. We will make use of the fact that the $\I_p$ are $*$-ideals in the sense that for all $A \in \I_p, B: L^2 \to L^2$
\[
    AB, BA \in \I_p, \qquad \|AB\|_{\I_p}, \|BA\|_{\I_p} \leq \|A\|_{\I_p}\|B\|_{L^2 \to L^2},
\]
and of the Hölder-like inequality
\[
    \|AB\|_{\I_r} \leq \|A\|_{\I_p}\|B\|_{\I_q}, \qquad \frac1r=\frac1p+\frac1q.
\]
In particular for all $\sum_i \frac1{p_i} = 1$,
\[
    \big|\trace\big(\prod_i A_i\big)\big| \leq \big\|\prod_i A_i\big\|_{\I_1} \leq \prod_i \|A_i\|_{\I_{p_i}}.
\]
We refer to \cite{MR3364494} for a thorough introduction to these spaces. Finally we note that the classes $\I_p$ admit the interpolation property, see \cite[Proposition 2.1]{Benyamini}.

\begin{lemma}\label{lem:opest} 
We have with $ \sigma = \im z $
\begin{equation}\label{eq:KWest}
\begin{split} 
    \Vert   K(w) \Vert_{\I_2} \, & =  \Big(2+ \frac{\tau^2}{|z|^2}\Big)^{\frac12} \sigma^{-\frac12}\Vert  w \Vert_{L^2}, \\ 
    \Vert K(w) \Vert_{\I_p} \, & \le  \Big( 2+ \frac{\tau}{|z|} \Big)    \sigma^{-1+\frac1 p} \Vert w \Vert_{L^p}, \qquad p \geq 2
    \end{split} 
\end{equation}
\end{lemma}
\begin{proof} 
For the first estimate we calculate $K^* K$ and obtain
\[\begin{split}
    \|K\|_{\I_2}^2 &= \trace(K^*K) \\
    &= \|(-iz-\partial)^{-1}w\|_{\I_2}^2 + \|(-iz+\partial)^{-1}w\|_{\I_2}^2 + 4\tau^2\|(-z^2-\partial^2)^{-1}w\|_{\I_2}^2
    \\ &= ( 1+ \tau^2/z^2)     (2iz)  \int ((-\partial + 2iz)w)^2 dx . 
    \end{split}
\]
By using the integral kernel of $(-iz-\partial)^{-1}$ and Fubini's theorem we obtain
\[
    \|(-iz-\partial)^{-1}w\|_{\I_2}^2 = \int_{x>y} |e^{iz(x-y)}w(y)|^2\,dydx = \sigma^{-1} \|w\|_{L^2}^2.
\]
The same holds for the second summand. For the third summand we calculate, using Fubini and then Plancherel,
\[
    \|(-z^2-\partial^2)^{-1}w\|_{\I_2}^2 = (2\pi)^{-1}\|(-z^2+\xi^2)^{-1}\|_{L^2_\xi}^2 \|w\|_{L^2}^2 = \frac{1}{4\sigma |z|^2}\|w\|_{L^2}^2.
\]
This proves first identity of \eqref{eq:KWest} for the Hilbert-Schmidt norm. Next 
\[
    \|(-iz\pm\partial)^{-1}w\|_{L^2 \to L^2} \leq  \|(-iz\pm\partial)^{-1}\|_{L^2 \to L^2}\|w\|_{L^\infty} =  \sigma^{-1}\|w\|_{L^\infty}
\]
which implies 
\[ \| K(w) \|_{L^2 \to L^2} \le  (1+ \frac{\tau}{\sigma} )      \sigma^{-1}   \]

and  for the second inequality we interpolate by viewing the operator $K(w)$ as a map $L^p \to \I_p$, where for $p = 2$ we use our first inequality and for $p=\infty$ the operator bound. 
\end{proof} 
 Let $N \in \N $, $ N \le n$. Then 
\begin{equation} \label{eq:KnN} 
\begin{split} |\trace K^n| \le \Vert K^n \Vert_{\I_1}\, & \le   \Vert K^{N} \Vert_{\I_1} \Vert K \Vert_{\I_2}^{n-N} 
\\ & \le  c_N(\im z)^{1-N} \Vert w \Vert_{L^{N}}^{N} \Big(\frac{3}{(\im z)^{1/2}} \Vert w \Vert_{L^2}\Big)^{n-N}. 
\end{split} 
\end{equation} 
There is one more structural observation: 
The integral kernel of the lower right entry of $K$ 
is 
\[   2\tau  \int^\infty_{\max\{x,y\}} e^{-2iz (2t-x-y)} dt =   \frac{\tau}{iz} e^{-iz|x-y|}    \]
hence 
\[ 2 (iz-\partial)^{-1} \tau (iz+\partial)^{-1} 
= -\frac\tau{iz} \big( (-iz-\partial)^{-1} + (-iz +\partial)^{-1}    \big) 
\]
If in an expansion the lower right entry of 
\[ \left( \begin{matrix} 0 & (-iz -\partial)^{-1} w \\ (-iz +\partial)^{-1}w & -\frac{\tau}{iz} \big( (- iz-\partial)^{-1} + ( -iz+\partial)^{-1}\big) w \end{matrix} \right)  \]
is involved (which is always the case if $n $ is odd) then we gain a factor $ \frac{\tau}{\im z}$ in all the estimates above.

\subsection{The  case of \texorpdfstring{$n < N+2$}{n<N+2}: The structure of the terms.}\label{sec:diffhamiltonian} 

We will encounter a special class of multidimensional integrals frequently in this section. 
\begin{definition}[Primitive integrals]\label{def:primitive}  We call integrals of the type 
\begin{equation}
    \Big(\frac{2\tau}{2iz}\Big)^{2m-n}\int_{x_1 < \dots < x_n} \prod_{j=1}^{n} e^{2iz x\cdot y} w_j(x_j) dx_j, \quad n/2 \leq m \leq n
\end{equation} 
where 
\[ \sum_{j} y_j=0 , \quad y\cdot x \ge x_n-x_1 \qquad \text{ for all } x \in A= \{ x_1<x_2 < \dots < x_n\} \]
 primitive integrals. 
\end{definition} 

We decompose 
\[ \T_{-1}^\Gardner(z,w,\tau) = \sum_{n=2}^\infty  \T_n(z,w,\tau) \]
where
\[ \T_n(z,w,\tau) =  \frac{(-1)^n}{n}  \frac{iz}{4z^2 + 4 \tau^2} \trace K^n(w) \]
is the $n$-homogeneous part with respect to $w$. In the previous section we have seen that if $ n \ge M+2 $ and $ \im z \ge \tau $
\begin{equation} \label{eq:highTN}  |\T_n(z,w,\tau)| \le (c (\im z)^{-1/2} \Vert w \Vert^2_{L^2})^{M+2-n}\frac{|z|}{\im z}   (\im z)^{-M}  \Vert w \Vert^{M+2}_{L^{M+2}} . \end{equation}

We begin with an algebraic decomposition of $\T_n $.

\begin{lemma} \label{lem:tracestructure} 
Suppose that $n < M+2$. Then we can write
\begin{equation}  (2z)^{M+1} \T_{n} = \sum_{j \le M-2}  (2z)^{M-j} \int h_{n,j} dx +   (2z)^{-1}    \T_{M+1,n} \end{equation} 
(so that the  leading term is $ (2z)^M \trace K^n(w)$)
where $h_{n,j}$ are differential polynomials independent of $M$ and 
\begin{equation} 
\T_{M+1,n} = \sum_{j=1}^J  \Big( \frac{z}{\im z} \Big)^{D_j-1} \tau^{k_j}   B_{M+1,n,j}
+  \frac{4z^2 \tau^{M+2-n}}{4z^2+ 4 \tau^2 }  \Big( \frac{z}{\im z} \Big)^{n-1}  B_{M+1,n} 
\end{equation} 
where the $h_{n,j}$ are differential polynomials, 
the $B_{M+1,n,j}$ are  $D_j$ dimensional primitive integrals, $D_j\ge 2$, 
\[  (2\im z)^{D-1} \int_{x_1<x_2 \dots < x_D} e^{2iz y \cdot x}  P_1(x_1) w( x_2) \dots w(x_{D-1}) P_D(x_D) dx_1 \dots dx_D \]
where $P_1$ and $P_D$ are differential monomials. In total there are $n$ $w$ factors (including those with derivatives) and there are at most 
\[ d \le M+2 - n -k_j \]
derivatives which are evenly distributed over $P_1$ and $P_D$ (equal if $2m-k$ is even, with a difference $1$ otherwise). 
If $n$ is odd then $k_j \ge 1 $.

 The $B_{M+1,n}$ are sums of $n$ dimensional primitive integrals 
\[   (2\im z)^{n-1} \int_{A_n} \exp( iz y_j\cdot x) \prod_{k=1}^n w(x_k) dx_k. \] 
\end{lemma}

\begin{proof} 
To shorten the notation we write $ R_{\pm} = (-iz\pm \partial)^{-1}$ so that $K$ (see \eqref{eq:Kdef}) can be written as 
\[  K(w) = \left( \begin{matrix} 0 & R_-w \\ -R_+ w & -\frac{\tau}{iz} ( R_-+R_+) w  \end{matrix} \right). \]
Then, since $\trace R_+ w R_+ w = \trace R_-w R_-w = 0 $ 
\[ 
\begin{split} 
\trace K^2\, &  = \trace \Big( -R_+w R_- w - R_-wR_+ w - \frac{\tau^2}{z^2} ( R_+ w R_+ w + R_-w R_- w+  R_+ w R_-w + R_-w R_+ w)  \Big)
\\ & = - \frac{ 4z^2+ 4\tau^2}{4z^2} \int e^{2iz|x-y|}   w(x) w(y)    dx dy, 
\end{split} 
\]
 At $z=i\tilde \tau$  we obtain 
\begin{equation} \label{eq:trK2} \begin{split}  
\trace K^2\, & = \frac{4\tau^2 - 4\tilde \tau^2 }{4\tilde \tau^2}
\int \int_{\R^2} e^{-2\tilde\tau |x_1-x_2|} w(x_1) w(x_2) dx_1 \, dx_2 =   \frac{4\tau^2 - 4\tilde \tau^2 }{\tilde \tau}     \Vert w \Vert_{H^{-1}_{2\tilde \tau}}^2 
\\ & =   \frac{4\tau^2 - 4\tilde \tau^2 }{\tilde \tau}     \Big( \sum_{n=0}^N  \frac{(-1)^n}{(2\tilde \tau)^{2n+2}} \Vert w^{(n)} \Vert_{L^2}^2 + \frac{(-1)^{N+1}}{ (2\tau)^{2N+1}}
\Vert w^{(n+1)} \Vert_{H^{-1}_{2\tau}}^2 \Big) 
\end{split} 
\end{equation} 
where we used 
\[ 
\frac{1}{(2\tau)^2+|\xi|^2} =  \sum_{j=0}^n (-1)^j  \frac{\xi^{2j}}{(2\tau)^{2j+2}} + (-1)^{n+1} (2\tau)^{-2(n+1)} \frac{|\xi|^{2(n+1)}}{(2\tau)^2 + |\xi|^2}. 
\]   
Similarly
\[\begin{split}  \trace K^3\, &  = \Big( \frac{3\tau^3}{i z^3}+ \frac{6\tau}{iz} \Big)  \trace \Big((R_-w)^2 R_+w +  (R_+ w)^2 R_-w) \Big) 
\\ &=  \Big( \frac{3\tau^3}{i z^3}+ \frac{6\tau}{iz} \Big) \int_{x_1<x_2<x_3} e^{ 2iz (x_3-x_1)} w(x_1) w(x_2) w(x_3) dx_1 dx_2 dx_3, 
\end{split} 
\]
\[ \begin{split} 
\trace K^4 \, & = \trace\Big( \big(2+ \frac{2\tau^2}{z^2} + \frac{\tau^4}{z^4} \big) R_-w R_+ w R_-w R_+ w
+ \big( \frac{2\tau^2}{z^2}+\frac{4\tau^4}{z^4}  \big) R_-w R_-w R_+w R_+ w 
\\ & \qquad + \Big( \frac{2\tau^2}{z^2} + \frac{4\tau^4}{z^4}\Big)   ((R_-w)^3R_+w + R_-w(R_+w)^3 )\Big)  
\\ & =\!\! \int\limits_{x_1<x_2<x_3<x_4} \!\!\Big(\big(8+ \frac{8\tau^2}{z^2} + \frac{4\tau^4}{z^4} \big)e^{\gamma_1 }  + \big( \frac{6\tau^2}{z^2}+  \frac{12\tau^4}{z^4} \big) e^{\gamma_2} \Big)   w(x_1) w(x_2) w(x_3) w(x_4) dx_1 dx_2 dx_3 dx_4 
\end{split} 
\]
with 
\[ 
\gamma_1 = 2iz ( x_4+x_3-x_2-x_1), \qquad  
\gamma_2 = 2iz ( x_4-x_1).
\]
Expanding the $2\times 2$ matrices in general we see that $ \trace K^n$ is a linear combination
of products of $n$ factors $R_\pm w$.
 Let $\Sigma$ be the  permutations of $n$ elements. Since 
 \[ \R^m \backslash \{ x_j = x_k \text{ for some } j\ne k\} =  \bigcup_{\sigma \in \Sigma}
 \{ ( x_{\sigma_j})_{1\le j \le m}: x_1 < x_2 < \dots < x_m \} \]  
 $\trace K^n$ can be written as a linear combination of expressions of type
\begin{equation}\label{eq:sumkn} 
    \Big(\frac{2\tau }{2iz}\Big)^{2m-n}\int_{x_1 < \dots < x_n} \prod_{j=1}^{n} e^{2iz x\cdot y} w(x_j) dx_j, \quad n/2 \leq m \leq n
\end{equation} 
where 
\[ \sum_{j} y_j=0 , \quad y\cdot x \ge x_n-x_1 \qquad \text{ for all } x \in \Omega_n= \{ x_1<x_2 < \dots < x_n\}. \]

In Lemma \ref{genhoelder0} it will be shown that $n$ dimensional primitive integrals are bounded by 
\[ C (2\im z)^{1-n} \prod_j \Vert w_j \Vert_{L^n} \]  and $(2\im z)^{n-1}  $ times a primitive integral with $z = i\tau$ converges to 
\[        \int_{\Omega_n, x_1=0} e^{- y\cdot x} dx_2 \dots dx_n \int \prod w_j(x) dx_j. \]

\bigskip

We illustrate the cancellations with an example before stating a general algorithm. Consider 
\begin{equation*}
    T = (2iz)^{3} \int_{\Omega_4} e^{-2iz(x_1 + x_2 - x_3 - x_4)} w(x_1) \dots w(x_4) \, dx.
\end{equation*}
Define $\phi = -2iz(x_1 + x_2 - x_3 - x_4)$. Thus $T$ does not decay as $\im z  \to \infty$.  We do a partial integration in $x_4$ to obtain
\begin{equation*}
\begin{split}\hspace{1cm} & \hspace{-1cm} 
   - (2iz)^2 \int_{\Omega_4} e^{\phi} w(x_1) w(x_2) w(x_3) w'(x_4) \, dx \\ & - (2iz)^2 \int_{\Omega_3} e^{-2iz(x_1 + x_2 - 2x_3)}w(x_1) w(x_2) w^2(x_3)\, dx.
   \end{split} 
\end{equation*}
The first term has increased decay compared to $T$, for the price of one derivative. For the second term, we partially integrate again from the right in $x_3$. We iterate this procedure and arrive at
\begin{equation*}
    \begin{split}
        %\frac{1}{\friedrich{c}}
        T = & \, -(2iz)^2 \int_{\Omega_4} e^{\phi} w(x_1) w(x_2) w(x_3) w'(x_4) \, dx\\
        &\quad  +   iz  \int_{\Omega_3} e^{-2iz(x_1 + x_2 - 2x_3)} w(x_1) w(x_2) (w^2)'(x_3) \, dx\\
        &\quad -\frac{1}{2} \int_{\Omega_2} e^{-2iz(x_1 - x_2)} w(x_1) (w^3)'(x_2) \, dx + \frac{1}{2}\int_\R w^4 \, dx.
    \end{split}
\end{equation*}
The one-dimensional integral has no decay, all the higher-dimensional integrals have decay $(\im z)^{-1}$ if $w$ is sufficiently regular and integrable. We can further iterate, but now from the left, and apply the same procedure in $x_1$ to all of the remaining multidimensional integrals. We obtain 
\begin{equation*}
    \begin{split}
        T = & \, 2iz \int_{\Omega_4} e^{\phi} w'(x_1) w(x_2) w(x_3) w'(x_4) \, dx -\frac{1}{2} \int_{\Omega_3} e^{-2iz(2x_1 - x_2 - x_3)} (w^2)'(x_1) w(x_2) w'(x_3) \, dx\\
        &\quad + \frac{1}{4iz} \int_{\Omega_2} e^{2iz(x_1 - x_2)} (w^3)'(x_1) w'(x_2) \, dx - \frac{1}2 \int_{\Omega_3} e^{-2iz(x_1 + x_2 - 2x_3)} w'(x_1) w(x_2) (w^2)'(x_3) \, dx \\
        &\quad +\frac{1}{8iz} \int_{\Omega_2} e^{-4iz(x_1-x_2)} (w^2)'(x_1) (w^2)'(x_2) \, dx +\frac{1}{4iz} \int_{\Omega_2} e^{-2iz(x_1 - x_2)} w'(x_1) (w^3)'(x_2) \, dx\\
        &\quad + \frac{1}{2}\int_\R w^4 \, dx  - \frac{1}{4iz} \int_{\R} (w^3)'w\, dx- \frac{1}{8iz} \int_\R w^2 (w^2)' \, dx -\frac{1}{4iz}\int_\R w^3 w' \, dx.
    \end{split}
\end{equation*}
Now all the multidimensional integrals here have decay $(\im z)^{-2}$. Note that the one-dimensional integrals in this step all vanish, which is expected since we know that $\T_{-1}(i\tau) \sim \sum (2i\tau)^{-2j} H_{j-1}^{\Gardner}$ is only nonzero for even powers of $\tau$. Iterating even further will give the quartilinear term from $H_2^{\Gardner}$, and, of $H_N^{\Gardner}$ in general.

The following algorithm for connected integrals describes how to treat the cancellations in general:
Let $  n<M+2$. The term
\[   T :=   (2iz)^{M}  \int_{\Omega_n}  e^{2i  z y \cdot x} \prod w(x_j) dx_j \]
has homogeneity $n$ in $w$ and is written as a $n$ dimensional integral.  We integrate by parts in the $x_n$  integral which gives 
\begin{equation}\label{eq:firststep} 
\begin{split} \hspace{1cm}& \hspace{-1cm}  
T=(2iz)^{M-1} y_n^{-1} \int_{\Omega_n}     e^{2iz y \cdot x}  w'(x_n)  \prod_{j=1}^{n-1}  w(x_j) dx_j dx_1 
\\ & + (2iz)^{M-1} \int_{\Omega_{n-1}}     e^{2iz \Big( \sum\limits_{j=1}^{n-2}y_j x_j+ (y_{n-1}+y_{n}) x_{n-1}\Big) } w^2(x_{n-1}) \prod_{j=1}^{n-2} w(x_j) dx_j. 
\end{split} 
\end{equation} 
Let $ \hat y = (y_1, \dots , y_{m-1}+y_m)$. Then 
\[ \hat y \cdot x \ge x_{m-1}-x_1 \quad \text{ for } x \in \Omega_{m-1}. \]   
At this stage we keep the first term on the right hand side of \eqref{eq:firststep} and repeat the argument with the second term iteratively if $n<M+1$, ending up with a sum of integrals of dimension 
$n, n-1, n-2, \dots 1$. For all multidimensional (all besides the one dimensional) integrals we repeat the iterative procedure from $x_1$. We obtain two one-dimensional integrals plus multidimensional integrals where we gained a decay of $z^{-2}$, 

\[
\begin{split} \hspace{1cm} & \hspace{-1cm} 
 T =  c_1    (2iz)^{M+1-n}     \int w^n dx  + c_2 (2iz)^{M-n}  \int w^{n-1} \partial_x w dx  
\\ & + (2iz)^{M-2} \sum_{D=2}^n  \sum_{k_1+k_2=n-D} \sum_l c_{M,l,k_1} (2iz)^{D-1}  \int_{\Omega_D} e^{2iz \langle y_{D,l}, x\rangle} 
\partial_{x_1}  w^{k_1}(x_1) \partial_{x_D} w^{k_2} (x_D)
\\ & \qquad \times \prod_{j=2}^{D-1} w(x_j)  dx_j dx_1 dx_D  
\end{split} 
\] 
where  the homogeneity in $w$ is always $n$. 
The second one dimensional integral clearly vanishes, but we keep it at this point since the vanishing of the corresponding terms in the next steps is true but not immediately obvious. 

We repeat this iterative procedure $
M+2-n$ times and arrive at
\[ T =  \sum_{k=0}^{M+1-n}(2iz)^{k}\int P_{M,k}  dx +  (2iz)^{-1} A_{M+1,n} \] 
where $P_{M,k}$ are differential polynomials of homogeneity $n$ with a total of $M+1-n-k$ derivatives. 
The terms $A_{M+1,n}$ are sums over $D$ dimensional integrals, 
\[ A_{M+1,n} = \sum_{D=2}^n  z^{D-1}    \sum_j 
 \int_{\Omega_D} e^{iz \langle  y_{D,j} \cdot x \rangle} P_{1,D,j}(x_1) w(x_2) \dots P_{D,D,j}(x_D) dx_1 \dots dx_D.  \]
Here $P_{D,j}$ are homogeneous differential monomials of homogeneity $n_j$, $n_1+n_D+D-2= n$,    with a total of $M+2-n$ derivatives. 

We complete the proof of Lemma \ref{lem:tracestructure}: We expand
$\trace K_n $ into a sum  $(\tau /z)^k $,
$0 \le k \le n$
times the trace of the product of $n$ operators $R_\pm$, which we expand into a sum of $n$ dimensional primitive integrals. 

We multiply by the prefactor which we expand as far as we need
\[ \frac{(-1)^n}{n} \frac{iz(2z)^{2N+2}}{4z^2+4\tau^2}= \frac{(-1)^ni}{2n} \sum^L_{l=0}  
 (-1)^{l} (2\tau)^{2l}  (2z)^{2N-2l}
 + \frac{(-1)^ni}{2n}  \frac{(2z)^{2N-2L}}{4z^2+4\tau^2}.
\] 
Finally  we do the  integrations by parts.
\end{proof}

\subsection{The case \texorpdfstring{$n < M+2$}{n<M+2}: Estimates} 

In the previous section we have seen that we can decompose $(2z)^M \trace K^n(w)$ into integrals over differential polynomials and a sum of primitive integrals. In the next lemma we provide estimates for the primitive integrals.

\begin{lemma}\label{lem:traceestimates}  
Let $N \ge 0$, $M+2= d+n$, $1\le D \le n$,
\[I= (2\im z)^{D-1} \int_{\Omega_D} e^{2izy \cdot x} 
P_1(x_1) w(x_2) \dots P_D(x_D) dx_1 \dots dx_D  \]
be a primitive integral (see Definition \ref{def:primitive}) 
with $ n$ terms $w$, a total number of derivatives $d$, evenly distributed among $P_1$ and $P_D$ as above. 
Then following estimates hold: 
\[ |I| \le  c \left\{  \begin{array}{l} \Vert w \Vert_{L^{M+2}}^{M+2} + \Vert w^{(\frac{M}2)} \Vert_{L^2}^2    \\[2mm]     \Vert w \Vert^{n+2-\frac{d}{M}}_{L^2}  \Vert w^{(\frac{M}2)} \Vert^{1+\frac{d}{M}}_{L^2}  
\\[2mm] \Vert w \Vert^{n-2}_{L^2} \Vert u^{(\frac{M+1}2 -\frac{n}4)} \Vert^2_{L^2}.
\end{array}\right. \] 
\end{lemma} 
We apply these estimates to $ \T_{M+1,n} $ of Lemma \ref{lem:tracestructure}.
\begin{corollary}\label{cor:smalln} 
If $n < M+2$ then 
\begin{equation} \label{eq:even}  |\T_{M+1,n}| \le c \Big(\frac{|z|}{\im z} \Big)^{n-1}   \left\{ \begin{array}{l}  \Vert w \Vert^{M+2}_{L^{M+2}}  +   \Vert w \Vert^2_{H^{\frac{M}2}_\tau} 
\\[2mm] \Vert w \Vert_{L^2}^{n-2} \Vert w \Vert^2_{H^{\frac{M+1}2-\frac{n}4}_\tau}.
\end{array} 
\right.  \end{equation} 
and if $n$ is in addition odd 
\begin{equation} \label{eq:odd}  |\T_{M+1,n} | \le c \tau  \Big(\frac{|z|}{\im z} \Big)^{n-1}
\times \left\{ \begin{array}{l}  \Vert w \Vert^{M+1}_{L^{M+1}}  +   \Vert w \Vert^2_{H^{\frac{M-1}2}_{\tau}} 
\\ \Vert w \Vert^{n-2}_{L^2} \Vert w \Vert^2_{H^{\frac{M}2-\frac{n}4}_\tau} 
.\end{array} \right. 
\end{equation} 
\end{corollary} 

\begin{proof} 
We apply Lemma \ref{lem:traceestimates} to the decomposition 
of \ref{lem:tracestructure}. We have to estimate $I$ of homogeneity $n$ with a total number  
$ d \le M+2 -n-k$ derivatives where $k$ is the power of $\tau$.
The corollary follows by H\"older's inequality 
\[   \int \tau^j |w|^{M+2-j} dx \le \Big( \int w^{M+2} dx \Big)^{1-\frac{j}{M}  } \Big(\int \tau^{M} w^2 dx\Big)^{\frac{j}{M}}     \]
and interpolation 
\[ \tau^{j/2} \Vert  w^{((M-j)/2)} \Vert_{L^2} \le  \Vert \tau^{M/2} w \Vert^{\frac{j}{M}}_{L^2} \Vert w^{(M/2)} \Vert^{1-\frac{j}{M}}. \]
The sum starts at $j=1$ if $n$ is odd. Tracing the change in the estimates implies the result in that case. 
\end{proof}

We will use H\"older type estimates for simpler primitive integrals first. 
To shorten notation, we recall the notation \begin{equation} \label{def:Am} \Omega_D = \{x\in \R^D: x_1 < \dots < x_D\}. \end{equation} 

\begin{lemma}  \label{genhoelder0} 
Let $D \ge 2$ and $ y \in \R^D$ satisfying 
\begin{equation}\label{condy}   \sum_{j=1}^D y_j = 0 , \qquad y \cdot x \ge x_D-x_1 \quad \text{ if }  x\in \Omega_D. \end{equation}
Define the primitive integral 
\[ I(z) =   \int_{ \Omega_D}  e^{2i z y \cdot x} \prod_{j=1}^D f_j(x_j) dx_j.    \]
Let $1\le p_j \le \infty$ satisfy $\sum\limits_{j=1}^D \frac1{p_j} \ge  1$. Then 
\begin{equation}  |I(z) | \le   (\im z)^{\sum\limits_{j=1}^D \frac1{p_j}-D} \prod \Vert f_j \Vert_{L^{p_j}}. \label{ineq:I0} \end{equation} 
If $ \sum \frac1{p_j}=1$ and $Q_j \subset L^{p_j}$  are equicontinuous sets then 
\[   \lim_{\tau \to \infty} (2\tau)^{D-1}  I(i\tau)  =   \int_{ \Omega_D: x_1=0}  e^{-y \cdot x} dx_2 \dots dx_D    \int_{\R} \prod_{j=1}^D  f_j(x) \,  dx  \]
uniformly for  $f_j \in Q_j$.
\end{lemma} 
\begin{proof} 
Fubini's theorem shows the trivial estimate 
\[ |I| \le \prod \Vert f_j \Vert_{L^1}.  \] 
We consider the endpoint case $\sum_{j=1}^M \frac1{p_j} = 1$. We observe that \eqref{condy} implies $y_M\ge 1$. Let $\frac1{\tilde p_{M-1}} = \frac1{p_{M-1}} + \frac1{p_M} \in [1,\infty]$ and 
\[ \tilde f_{M-1} = f_{M-1} \Big( (e^{2iz y_M x}\chi_{x>0})*f_M\Big) \] 
We estimate 
\[ \Vert  \tilde f_{M-1} \Vert_{L^{\tilde p_{M-1}}} \le 
\Vert f_{M-1} \Vert_{L^{p_{M-1}}} \Vert (e^{2iz y_Mx} \chi_{x>0})* f_M \Vert_{L^{p_M}}   
\]
and 
\[ \begin{split} \Vert (e^{2iz y_Mx} \chi_{x>0}) *f_M \Vert_{L^{p_M}}\, &  \le (2 \im z )^{-1} \Vert (e^{-y_Mx} \chi_{x>0})\Vert_{L^1} \Vert  f_M \Vert_{L^{p_M}} \\ & \le  (2\im z)^{-1} \Vert  f_M \Vert_{L^{p_M}}. \end{split}  \]
Let $ (\tilde y)_j = y_j $ if $ j < M-1$ and $ \tilde y_{M-1} = y_{M-1}+ y_M$. Then \eqref{condy} is satisfied: 
\[ \tilde y\cdot  \tilde x = y \cdot (x_1, \dots x_{M-1}, x_{M-1}) \ge x_{M-1} - x_1 \qquad \text{ in } \Omega_{M-1}. \]
We complete this case by induction on $M$. The general case follows now either by a simple variant or an iterative application of complex interpolation by Riesz-Thorin. 
To see the uniform convergence on equicontinuous sets we observe 
\begin{equation} \label{eq:convdist} 
(2\tau)^{M-1}   \chi_{\Omega_M}  e^{-2\tau   y \cdot x} \to   \int_{ \Omega_M: x_1=0}  e^{-y \cdot x} dx_2  \dots dx_M     \prod_{j=2}^M \delta_{x_j-x_1} 
\end{equation} 
as $ \tau \to \infty$ in the sense of functionals on $C_c(\R^M)$ which is seen by writing 
\[ (2\tau)^{M-1} I(i\tau)    =  \int_{0 < s_2 < s_3 \dots s_M} e^{-2\tau \sum\limits_{j=2}^M y_j s_j}      \int  f_1(x) \prod_{j=2}^M f_j(x+ s_j )ds_j \,   dx. 
\]
hence 
  \[ \begin{split} \hspace{1cm} & \hspace{-1cm} \left| (2\tau)^{M-1}I(i\tau) - \int \prod f_j(x) dx \right| \\ &   =  \left| \int_{0 < s_2 < s_3 \dots s_M} e^{-2\tau \sum\limits_{j=2}^M y_j s_j}      \int  f_1(x) \left( \prod_{j=2}^M f_j(x+ s_j )-\prod_{j=2}^M f(x) \right)ds_j \,   dx \right| 
  \\ & \le c ( \supp f_j) \tau^{-1} \prod_{j=1}^M \Vert f_j \Vert_{C^1} . 
\end{split} 
\]   

Let $Q_j \subset L^{p_j}(\R)$ be bounded equicontinuous sets. Given $\varepsilon >0$ there exists $R$ and $C>0$ so that for $ f_j \in Q_j $
there exists $F_j \in C^1((-R,R)$  with $ \Vert F_j \Vert_{C^1}\le C$ with $ \Vert F_j - f_j\Vert_{L^{p_j}} < \varepsilon$. Then 
\[  I(i\tau, f_j )  - I ( i\tau, F_j)  = \sum_{k=1}^M I(i\tau, f_1, \dots, f_{j-1},  f_j-F_j, F_{j+1},\dots , F_M)   
\]
hence 
\[ \left|  (2\tau )^{M-1}  I(i\tau, f_j ) -\int_{ \Omega_M: x_1=0}  e^{-y \cdot x} dx_2  \dots dx_M    \int \prod f_j dx \right| \le 
C \tau^{-1} + cM \varepsilon  
\]
uniformly for $ f_j \in Q_j $.
\end{proof}

\begin{proof}[Proof of Lemma \ref{lem:traceestimates}]
Let $n_j$  be the homogeneity of $P_j $ in $w$ and its derivatives. Then $n_1+n_D + M-2 = n$
We estimate using Lemma \ref{genhoelder0} to estimate 
\[ |I| \le  \Vert w \Vert^{m-2}_{L^{n+d}} \Vert P_1 \Vert_{L^{\frac{n+d}{n_1+d_1}}} \Vert P_D \Vert_{L^{\frac{n+d}{n_2+d_2}}}  \]
and, with $ P_1 = \prod w^{(\alpha_j)} $ by H\"older's inequality 
 \[ 
\Vert P_1 \Vert_{L^{\frac{n+d}{n_1+d_1}}}\le \prod \Vert w^{(\alpha_j)} \Vert_{L^{\frac{n+d}{1+\alpha_j} }}\]
with a similar estimate for $p_M$.
The two estimates of \ref{lem:traceestimates} follow by Lemma \ref{lem:estimatedifferentialmonomial2}. 
\end{proof}

\subsection{Bounding the difference Hamiltonian \texorpdfstring{ $ \T^\Gardner_N$}{TNGardner}}  
\label{subsec:diffHamiltonian}

The central estimate in this subsection is \eqref{eq:TnGardner}. Here we prove the claims on $\T^\Gardner_N$ and $\T^{\KdV}_N$ stated in Section \ref{sec:proof} which we repeat here.  

\begin{proposition}\label{lem:TNest}  The following estimate holds if $ \Vert w \Vert_{L^2} < \frac14\sqrt{\im z }$: 
\[ |\T^{\Gardner}_N (iz,w,\tau)|
\le c \left( \frac{|z|}{\im z}\right)^{2N+1}  (\im z)^{-2}  
\Big( \Vert w \Vert^{2N+4}_{L^{2N+4}}
+ \Vert w \Vert^2_{H^{N+1}_\tau}\Big). 
\]
\end{proposition}

\begin{proof}
We write
\begin{align} \notag
    2z\T_N^{\Gardner} &= (2z)^{2N+3}\T_{-1}^{\Gardner} - \sum_{k=0}^N(2z)^{2(N-k)+1}H_k^{\Gardner}\\ \notag
    &= (2z)^{-1}\sum_{n=2}^{2N+3} \T_{2N+3,n}  + \sum_{n > 2N+3}\frac{(-1)^n}{n} 
    \frac{iz(2z)^{2N+3}}{4z^2+4\tau^2}
    \trace K^n\\ & \quad - \sum_{k=0}^N(2z)^{2(N-k)+1}H_k^{\Gardner} + \sum_{n=2}^{2N+3} \sum^{2N+2}_{j=0} (2z)^{2N+2-j}\int 
  h_{n,j}  \label{eq:third}
\end{align}
and apply \eqref{eq:KnN} to the traces and estimate
\[ \Big| \frac{iz(2z)^{2N+3}}{4z^2+4\tau^2}
    \trace K^n\Big| \le  c_N \Big( |z|/\im z\Big)^{2N+2}(\im z)^{-1} \Big( \frac{3}{(\im z)^{1/2}} \Vert w \Vert_{L^2} \Big)^{n-2N-4} \Vert w \Vert^{2N+4}_{L^{2N+4}}.  \]
Summing over $n$ we obtain the desired estimate for the traces. By Lemma \ref{lem:tracestructure}
and Corollary \ref{cor:smalln}, if $n \le 2N+3$
\[ |\T_{2N+3,n}| \le c \Big( \frac{|z|}{\im z} \Big)^{n-1} \Big( \Vert w \Vert_{L^{2N+4}}^{2N+4} + \Vert w \Vert^2_{H^{N+1}}\Big) \]

Comparison of the expansion of Lemma \ref{lem:tracestructure} with the expansion \eqref{eq:Gardnerexpansion}
allows to identify the coefficients:
The $h_j$  vanish if $j$ is odd and 
$\int \sum h_{2j,n} dx= H^{\Gardner}_j$. In particular \eqref{eq:third} vanishes and
\begin{equation} 
 \T_{2N+2,n} = 2z \T_{2N+1,n} = (2z)^2 \T_{2N,n} + H^\Gardner_{N,n}.
\end{equation}
\end{proof} 

\begin{lemma}
The quadratic part of $\T^\Gardner_{-1}$ is 
\[   \frac12   \frac{-iz(2z)^{2N+2}}{4z^2+4\tau^2_0} \trace K^2(w) =  \frac12 \int \Big(  (-\partial + 2iz)^{-1} w \Big)^2\, dx.  \]
The quadratic part of $ \T^\Gardner_N(z,w,\tau_0)$ is 
\[  \frac12 \int ((\partial -2iz)^{-1} w^{(N+1)})^2 dx  . \]
\end{lemma}
\begin{proof} 
We give two different arguments. First the quadratic part of 
\[ \int w^2(z) - w^2 dx =  \int [ (\partial -2iz + w(iz))^{-1}(\partial +2\tau + w)w ]^2 -w^2 dx 
\]
is 
\[  \int    w (-\partial^2 -4z^2)^{-1}(-\partial^2  +4\tau^2)w  -w^2 dx = - (4 z^2+4 \tau^2 )  \int \big[ ( -\partial + 2iz) w \big]^2 .  \]
Alternatively we computed $\trace K^2$ explicitly in 
\eqref{eq:trK2}.
\end{proof}

We recall that we denote the higher order part for $z = i\tau$ by 
\[ \T^{NL}_{N} = \T^\Gardner_{N} - \frac12 \Vert w^{(N+1)} \Vert^2_{H^{-1}_{2\tau}}  . \]

\begin{lemma}\label{lem:decomposition}
We have 
\[ |\T^{NL}_{N}(i \tau,w,\tau_0) | \le c(\tau_0, \Vert w  \Vert_{L^2}) \tau^{-1}  \Vert w \Vert^2_{H^{N}_{\tau_0}}.
 \]
 \end{lemma} 
\begin{proof}
By the proof of Proposition \ref{lem:TNest} we have to estimate 
$\T_{2N+2,n}$.
 The claim follows by Corollary \ref{cor:smalln}. More precisely, with $M=2N+1$ and maximal a total number of $d=2N+3-n $ derivatives 
\[  |\T_{2N+2,n}| \le c    \left( \frac{|z|}{\im z}\right)^{n-1} 
\Vert w \Vert_{L^2}^{n-2} \Vert w \Vert^2_{H^{N+1- \frac{n}4 }_\tau}    
\]     
Then $ N+1- \frac{n}4 \le N $ if $n \ge 4$. Finally 
\[ |\T_{2N+2,3}| \le  c \tau  \left( \frac{|z|}{\im z}\right)^{2} \Vert w \Vert_{L^2}
 \Vert w \Vert^2_{H^{N-\frac14}_\tau}  \] 
 since $3$ is odd. The claim follows by summation with respect to $n$.
\end{proof}

 We obtain the claim on equicontinuity \eqref{eq:propequi} of Proposition \ref{prop:diffham} as corollary. 
\begin{corollary}\label{cor:Gardnerequi} 
A subset $Q \subset H^N$ is bounded and equicontinuous  iff 
\[ \lim_{\tau \to \infty} \sup_{w\in Q} |\T^\Gardner_N(i\tau, w,\tau_0) |=0. \]
\end{corollary} 
\begin{proof}
By Lemma \ref{lem:decomposition} is suffices to consider the quadratic part.  Then 
\[ \Vert  w^{(N+1)} \Vert^2_{H^{-1}_\tau} 
= \int  \frac{\xi^2}{\tau^2+ \xi^2}   |\xi^N \hat w|^2 dx \to 0  \]
uniformly for $w^{(N) } \in \tilde Q$ iff $ \tilde Q \subset L^2$ is bounded and equicontinuous.   
\end{proof}

We turn to bounding variational derivatives and prove \eqref{eq:propuni}, which we formulate as a lemma. 

\begin{lemma}\label{lem:diffN1}
Let $N \ge  2$, $Q \subset H^{N-1} $ be bounded and equicontinuous. Then 
\[ \Big\Vert \frac{\delta}{\delta w} \T^\Gardner_N(i\tau , w , \tau_0) \Big\Vert_{H^{-N-1}} \to 0 \quad \text{ as } \tau \to \infty  \]
uniformly in $Q$.
\end{lemma}
\begin{remark} 
The proof applies also to $\T^\Gardner_{app,N}$ essentially without change, taking into account Lemma \ref{lem:approximate}.
\end{remark} 

\begin{proof}
We observe that 
\[  \Vert (4\tau^2- \partial^2)^{-1}\partial^{2N+2} w \Vert_{H^{-N-1}}  \]
converges to $0$  uniformly in bounded and equicontinuous sets. 
From the definitions and the properties of the trace 
\[ \int \phi \frac{\delta}{\delta w} \frac1n \trace K^n = \trace K(\phi) K^{n-1}(w)\]
and hence all terms with $n \ge 2N+3$ converge uniformly to $0$  on bounded sets in $H^1$.
If $ n \ge 2N+3$ we hence obtain $\tau = \im z$ in this estimate
\[ \Big| \frac1n \int \phi \frac{\delta}{\delta w} \trace K^n dx\Big| \le (c \tau^{-1/2} \Vert w \Vert_{L^2})^{2N+3-n}(\im z)^{-2N-2} \Vert \phi \Vert_{L^\infty} \Vert w \Vert^{2N+2}_{L^{2N+2}}.   \]
For $ 3 \le n \le 2N+2$
we use the decomposition of Lemma \ref{lem:tracestructure} which reduces the claim to the corresponding claim for primitive integrals as in Lemma \ref{lem:tracestructure}. Let $M=2N+1$, $I$
be an $D$ dimensional  primitive integrals with a prefactor   $(\im z)^{D-1}$
with homogeneity $n$ in $w$ and a total of $d\le  2N+3-n$ derivatives which are evenly distributed. The variational derivative plus duality amounts to replacing $w$'s by $\phi$, resp. $n$ by $n-1$ and $M$ by $M-1$. Of course the situation improves further if derivatives fall on $\phi$.
Hence, if $n>2$ by Lemma \ref{genhoelder0} 
\[  \Big| \int \phi  \frac{\delta}{\delta w } I  dx  \Big| \le c   \sum_{m=0}^{d/2}
\Vert \phi^{(m)}  \Vert_{L^\infty} 
  \Vert w \Vert^{n-3}_{L^2} \Vert w \Vert^2_{H^{N+\frac14-\frac{n}4}_\tau}.
\]
 Thus all terms of homogeneity  $n\ge 5$ are bounded by $\tau^{-1} \Vert w \Vert^2_{H^{N-1}}$ with a constant depending on $\Vert w \Vert_{L^2}$.

It remains to consider the quadratic ($n=3$) and cubic contributions ($n=4$) with no derivative falling on $\phi$. We consider 
\[ \T^\Gardner_{N,n} := (2z)\T_{2N-1,n} + H^\Gardner_{N,n}   \]
with $n=3$ or $n=4$ with suggestive notation. If at least one derivative falls on $\phi$ or the power of $\tau$ is at least one in the cubic case ($n=4$) or at least 2 in the quadratic case ($n=3$) we argue as above and obtain a decay like $ (\im z)^{-1} $. For the remain case we omit the final integration by parts and write 
$ (2\tau)^{-1} \T_{2N+2,3} $ as a linear combination of terms of the type   
\[\begin{split} \hspace{1cm} & \hspace{-1cm} \tau_0 (2\tau)^{-2} \int_{\Omega_3} e^{-2\tau( x_3-x_1) }  w^{(N-1)}(x_1) \phi(x_2) w^{(N-1)}(x_3) dx_1dx_2dx_3 \\ & - \int_{0<x_2<x-2} e^{-x_3} dx_2dx_3    \int     (w^{(N-1)}(x))^2 \phi(x) dx  
\end{split}
\]
which converge to $0$ uniformly on bounded equicontinuous sets in $H^{N-1}$ by Lemma \ref{genhoelder0}.
For $n=4$ we obtain similarly a linear combination of terms of the type 
\[ 
\begin{split}\hspace{1cm} & \hspace{-1cm} 
\tau_0 (2\tau)^{-2} \int_{\Omega_3} e^{-2\tau( x_3-x_1) }  w^{(d_1)}(x_1)w^{(d_2)}(x_1)  \phi(x_2) w^{(N-1)}(x_3) dx_1dx_2dx_3 \\ & - \int_{0<x_2<x-2} e^{-x_3} dx_2dx_3    \int  w^{(d_1)} w^{(d_2)}     w^{(N-1)} \phi\, dx  \end{split}   \]
which converge to $0$ uniformly on bounded equicontinuous sets in $H^{N-1} \cap L^\infty$, and hence on bounded equicontinuous sets in $H^{N-1}$ if $ N \ge 2$.

\end{proof}

\subsubsection{Proof of Proposition \ref{eq:kdv:Lipschitz}}
\label{subsub:propeq:kdv}
We translate these bounds to the KdV side by the Miura map and prove the estimates of Propositions \ref{eq:kdv:Lipschitz}. We recall \eqref{eq:Gardner-KdV}
\[ 4z^2 \T^{\KdV}_N(z, w_x +2\tau w + w^2)  = 4 \tau^2 \T^\Gardner_{N}(z, w, \tau) + 4 z^2 \T_{N+1}^\Gardner(z,w,\tau). \] 
We estimate
\[ \begin{split} |\T^{\KdV}_N(z, u) | 
\, &\le |4 (\tau/z)^2 \T^\Gardner_{N}(z, w, \tau)| + | \T_{N+1}^\Gardner(z,w,\tau)| 
\\ & \le c((\im z)^{-1/2} \Vert w \Vert_{L^2})  (\im z)^{-2}\Big( \frac{|z|}{\im z}\Big)^{2N+3} \Big(\Vert w \Vert^2_{H^{N+2}_\tau} + \Vert w \Vert^{2N+4}_{L^{2N+4}} \Big)
\end{split}
\]
Together with the fact that $\|w\|_{H^{N+2}} \sim \|u\|_{H^{N+1}}$ with a constant depending on $\tau_0$ (see Proposition \ref{prop:equivalencewu}), and 
\[ \Vert w \Vert_{L^{2N+4}}^{2N+4} \le   c(\Vert w \Vert_{L^2}) \Vert u \Vert^{N+2}_{L^{N+2}},\]
 shows \eqref{eq:tauNdec}. %The bound \eqref{eq:tauNdecder} follows by similar arguments. 

\subsection{Weighted estimates} 

In this subsection we prove the central estimates used to prove local smoothing of the difference flow \eqref{eq:diffflowKato} and convergence of the difference flow \eqref{eq:diffflow}. 

\begin{proposition}\label{prop:weighted} 
A)  The following estimates hold for $ n \ge M+2$.
\begin{equation}\Big| \int \sech^2 \phi \frac{\delta}{\delta w} \trace K^n(w)dx \Big|  \le \big(c (\im z)^{-1/2} \Vert w \Vert_{L^2} \big)^{n-M-2}(\im z)^{-M-1}
\Vert \phi \Vert_{L^\infty} \Vert w \Vert_{L^2}^{M-1} 
\Vert \sech w \Vert^2_{H^{\frac{M-1}4}_\tau}
\end{equation} 
\begin{equation}\Big| \int \tanh w \partial  \frac{\delta}{\delta w} \trace K^n(w)dx \Big|  \le \big(c (\im z)^{-1/2} \Vert w \Vert_{L^2} \big)^{n-M-2}(\im z)^{-M-2}|z|
 \Vert w \Vert_{L^2}^{M} 
\Vert \sech w \Vert^2_{H^{\frac{M}4}_\tau}.
\end{equation} 
B)  Let $2\le  n < M+3$ and let 
\begin{equation} \label{eq:primitive}   I = (\im z)^{D-1} \int_{\Omega_D} e^{2iz y \cdot x} P_1(x_1) w(x_2) \dots P_D(x_D) dx_1 \dots dx_D  \end{equation} 
be a $D\ge 2$ dimensional primitive integral with homogeneity $n$, a total number of  derivatives $d=M+3-n $ which are evenly distributed. 
Then
\begin{equation} 
\Big| \int \sech^2 \phi \frac{\delta}{\delta w}I  dx \Big|  
\le \Vert \phi \Vert_{H^M} \Vert w \Vert_{L^2}^{n-3} 
\Vert \sech w \Vert^2_{H^{\frac{2M+3-n}4}_\tau}\qquad  \text{ if }n \ge 3
\end{equation}
and
\begin{equation} 
\Big| \int \tanh w \partial_x  \frac{\delta}{\delta w } I dx \Big| \le c  
 \Vert w \Vert_{L^2}^{n-2} \frac{|z|}{\im z}  
\Vert \sech w \Vert^2_{H^{\frac{2M+4-n}4}_\tau}. 
\end{equation}
\end{proposition}

We postpone the proof and deduce \eqref{eq:diffflow} and \eqref{eq:diffflowKato}.
\begin{proof}[Proof of \eqref{eq:diffflow} and \eqref{eq:diffflowKato}]
In this section constants are allowed to depend on $\tau_0$ which simplifies the notation.
Let $M=2N+2$ and $n \ge 2N+4$. By duality 
\[ \Big\Vert \sech^2 \frac{\delta}{\delta w}\T_{2N+3,n} \Big\Vert_{L^1}
\le (c(\im z)^{-1/2} \Vert w \Vert_{L^2})^{n-2N-3} 
\Big( \frac{|z|}{\im z} \Big)^{2N+3} 
\Vert w \Vert_{L^2}^{2N} \Vert \sech w \Vert^2_{H^{\frac{N}2+\frac14}}.\]
If $2<n \le M+2=2N+3$,
\[ 
\Big\Vert \sech^2 \frac{\delta}{\delta w} \T_{2N+3,n}(i\tau, w , \tau_0)  \Big\Vert_{H^{-N}} \le c
\Vert w \Vert_{L^2}^{n-3} \Vert \sech w \Vert^2_{H^{N+ \frac{5-n}4}} 
\]
and, since  $n=3,5$ is 
\[
\Big\Vert \sech^2 \frac{\delta}{\delta w} \T_{2N+3,5}(i\tau, w , \tau_0)  \Big\Vert_{H^{-N}} \le c
  \Vert w \Vert_{L^2}^{n-3} \Vert \sech w \Vert^2_{H^{N- \frac{1}2}} 
\]
\[
\Big\Vert \sech^2 \frac{\delta}{\delta w} \T_{2N+3,3}(i\tau, w , \tau_0)  \Big\Vert_{H^{-N}} \le c
\Vert w \Vert_{L^2}^{n-3} \Vert \sech w \Vert^2_{H^N} 
\]
Let  $ \T^{>5}_N$ be the contributions of homogeneity $n\ge 5$ to $\T_N$. Then, if $ \Vert w \Vert_{L^2} \le \frac14 \tau$ 
\begin{equation} \Vert \sech^2 \frac{\delta}{\delta w} \T_N^{>5}(i\tau, w, \tau_0) \Vert_{H^{-N}}\le c \tau^{-2}  (1+ \Vert w \Vert^{2N+1}_{L^2}) \Vert \sech w \Vert_{H^{N-\frac14}}^2  . \end{equation}

Let $ \chi \in \C^\infty_c([-2,2])$ be an even function, identically $1$ on $[-1,1]$. 
We decompose 
\[ w = w_< + w_> \qquad   w_<=    \mathcal{F}^{-1}    \chi(|\xi|/\tau)   \hat w. \]
If $ 0 \le s \le N$
\[ \Vert \sech w \Vert_{H^s}  \le \Vert \sech w_< \Vert_{H^s} + \Vert \sech w_> \Vert_{H^s}, \]
\[ \begin{split} \Vert \sech w_> \Vert_{H^s}\, &  \le \Vert \sech w_> \Vert^{1-\frac{s}{N}}_{L^2}
\Big( \Vert \sech w_> \Vert_{L^2} + \Vert \sech w_>^{(N)} \Vert_{L^2}\Big)^{\frac{s}{N} }
\\ & \le c (\Vert \sech w \Vert_{L^2} + \Vert \sech  w^{(N+1)} \Vert_{H^{-1}_{2\tau}}) 
\end{split} 
\] 
and 
\begin{equation}\label{eq:lowfreq}\begin{split}    \Vert \sech  w_< \Vert_{H^s_\tau}\, &   \le  \Vert \sech w \Vert^{1-\frac{s}{N+1}}_{L^2} \Vert \sech w_< \Vert^{\frac{s}{N+1}}_{H^{N+1}} 
\\ & \le c \tau^{ \frac{s}{N+1} } ( \Vert \sech w \Vert_{L^2} + \Vert \sech w^{(N+1)}_< \Vert_{H^{-1}_{2\tau}}).
\end{split} 
\end{equation} 
These weighted interpolation estimates are a consequence of the interpolation estimates without weight by first proving them on bounded intervals by restriction and extension, and summation over the intervals. 
If $n> 4 $ or $n=3$ we arrive at the desired bound by 
\[ \begin{split} \tau^{-2} \Vert  \sech w \Vert^2_{H^{N}}
\, & \le \tau^{-2} \Big( \Vert \sech w_< \Vert_{H^{N}}
+ \Vert \sech w_> \Vert_{H^{N}}\Big)^2
\\ & \le c \tau^{-\frac{2}{N+1}} \Big( \Vert  w  \Vert_{L^2}
+ \Vert \sech w^{(N+1)} \Vert_{H^{-1}_\tau} \Big)^2
\end{split} 
\]
The case $n=4$ is more delicate. We have the two bounds 
for the quartic terms in $\T^\Gardner_N$, which, using a different notation to above, 
 \[ \Big\Vert \sech^2\frac{\delta}{\delta w} \T^\Gardner_{N,4}(\tau, w, \tau_0)  \Big\Vert_{H^{-N}} \le c \tau^{-1} \Vert w \Vert_{L^2} \Vert \sech w \Vert^2_{H^{N-\frac14} }  .\]
\[ \Big\Vert \sech^2\frac{\delta}{\delta w} \T^\Gardner_{N,4}(\tau, w, \tau_0)   \Big\Vert_{H^{-N}} \le c \tau^{-2} \Vert w \Vert_{L^2} \Vert \sech w \Vert^2_{H^{N+\frac14}}.\]
We split $w = w_< + w_>$, polarize  the cubic  term ($n=4$)   and expand. There can be 0, 1, 2 or 3 high frequency factors $w_>$. Then, with a suggestive notation 
\[ \begin{split} \hspace{1cm} & \hspace{-1cm} \Big\Vert \sech^2 \frac{\delta}{\delta w} \T^\Gardner_{N,4,>>>}(\tau, w, \tau_0)  \Big\Vert_{H^{-N}} + \Big\Vert \sech^2\frac{\delta}{\delta w} \T^\Gardner_{N,4,<>>}(\tau, w, \tau_0)  \Big\Vert_{H^{-N}} \\ & \le c \tau^{-1} \Vert w \Vert_{L^2} \Vert \sech w_> \Vert^2_{H^{N-\frac14} } 
\end{split}
\]
since $ \Vert w_< \Vert_{H^s} \le 2\Vert w_> \Vert_{H^s}$ (and commuting the weight with the Fourier multiplier and putting more derivatives on the low frequency part), 
\[ \Big\Vert \frac{\delta}{\delta w} \T^\Gardner_{N,4,<<>}(\tau, w, \tau_0)  \Big\Vert_{H^{-N}} \le c \tau^{-2} \Vert w \Vert_{L^2} \Vert \sech w_< \Vert_{H^{N+\frac12}} \Vert \sech w_> \Vert_{H^{N}} \] 
and 
\[ \Big\Vert \frac{\delta}{\delta w} \T^\Gardner_{N,4,<<<}(\tau, w, \tau_0)  \Big\Vert_{H^{-N}} \le c \tau^{-2} \Vert w \Vert_{L^2} \Vert \sech w_< \Vert^2_{H^{N+\frac14}}.\]  
We observe that 
\[ \Vert \sech w_> \Vert_{H^N} \le c \Big( \Vert  \sech w \Vert_{L^2} + \Vert \sech  w^{(N+1)} \Vert_{H^{-1}_{2\tau}} ). \]
Together with \eqref{eq:lowfreq} and the easy case $n=2$ this implies \eqref{eq:diffflow}. 

\bigskip

The analogous estimates for Kato smoothing are 
\[ \Big| \int \tanh(x) w\partial \frac{\delta}{\delta w} \trace K^n dx \Big|  \le \big(c (\im z)^{-1/2} \Vert w \Vert_{L^2} \big)^{n-2N-4} (\im z)^{-2N-3}    
 \Vert w \Vert_{L^2}^{2N+2} 
\Vert \sech w \Vert^2_{H^{\frac{N+1}2}_\tau}
\]  
for $n \ge  2N+4$ and
\[ \Big| \int \tanh(x) w \partial  \frac{\delta}{\delta w} \T_{2N+3,n} dx \Big|  \le c  
 \Vert w \Vert_{L^2}^{n-2} 
\Vert \sech w \Vert^2_{H^{N + \frac{6-n}4}_\tau}
\]  
We argue as above if $n \ge 6$. The case $n=5$ is simple due to the gain if $n$ is odd. The case $n=3$
follows by the very same arguments as for $n=4$ above. There are only trivial changes if there are two or more high frequency terms. In total we get a bound with $c$
depending on $N$,$\tau_0$ and $ \Vert w \Vert_{L^2}$, assuming $ \Vert w \Vert_{L^2} \le c \tau^{1/2} $,

\[ \Big| \int \tanh(x) w \partial  \frac{\delta}{\delta w} \T^{>2}_{N} dx \Big|  \le c \tau^{-\frac1{2(N+2)} }  
( \Vert \sech w \Vert^2_{L^2}  +  \Vert \sech w^{(N+1)} \Vert_{H^{-1}_{2\tau} })
\]  
which implies \eqref{eq:diffflowKato}.
\end{proof}

\begin{proof} [Proof of Proposition \ref{prop:weighted}]
We have
\[ \int \phi \sech^2(x/R) \frac{\delta}{\delta w} \trace K^n(w) dx 
=  n\trace \Big( K(\sech^2(x/R) \phi) K^{n-1}(w)\Big). 
\]
Again we estimate the traces for $n \ge 2N+2$. Suppose that  $ 2\le p \le \infty$,  $ \eta$ be slowly $\tau$ varying  and $ 1 \le \tau \le \im z$. We claim 
\begin{equation}\label{eq:weighted}
    \Vert \eta K(w) \eta^{-1} \Vert_{\I_p} \le (1+ \frac{\tau}{\im z} ) \Vert w \Vert_{L^p}. 
\end{equation}
To verify the claim we bound  the Hilbert-Schmidt norm 
\[  \Vert \eta K(w) \eta^{-1} \Vert^2_{\I_2}
\le c (\im z)^{-1/2} \Vert w \Vert_{L^2} 
\]
and 
\[  \Vert \eta K(w) \eta^{-1} \Vert^2_{L^2\to L^2}
\le c (\im z)^{-1} \Vert w \Vert_{L^\infty}. \]
Estimate \eqref{eq:weighted} follows by interpolation.
An algebraic manipulation gives 
\[ \begin{split}
\hspace{1cm}& \hspace{-1cm} K( \sech^2(x/R)\phi)  K^{M}(w) = K(\phi) \big(\sech^2 K(\sech^{\frac{2}M}w)\sech^{-2}\big) \\ & \big( \sech^{\frac{2(M-1)}M}K(\sech^{\frac2M}w) \sech^{-\frac{2(M-1)}M}\big)   \dots \big(\sech^{\frac2M} K(\sech^{\frac2M} w) \sech^{-\frac2M}\big)
\end{split}
\]
and
\begin{equation}
\begin{split} \hspace{1cm} & \hspace{-1cm} 
(\im z)^M\Big|\trace K(\sech^2 \phi) K^n(w) \Big| \le  
   ( 2 (\im z)^{-1/2} \Vert w \Vert_{L^2})^{n-M} 
\im z   \Vert K(\phi)\Vert_{L^2\to L^2}
\\  & \times \prod_{j=1}^M  (\im z)^{\frac{M-1}M} \Vert  \sech^{\frac{2j}M} K^M(\sech^{\frac2M} w) \sech^{-\frac{2j}M}  \Vert_{\I_M}
\\ \le &  \Vert \phi \Vert_{L^\infty} ( 2 (\im z)^{-1/2} \Vert w \Vert_{L^2})^{n-M} 
\Vert \sech^{\frac2M} w \Vert_{L^M}^M.
\end{split}
\end{equation}
The argument is clearly more flexible: We may distribute the weight among the factors as we wish. Primitive integrals can be understood as traces, and we may estimate them in the same fashion. We consider a primitive integral defined in \eqref{eq:primitive} 
 Then, by H\"older's inequality  
\[ \Big| \int \sech^2 \phi \frac{\delta}{\delta w} I dx\Big| 
\le c  \sum_{k=0}^{d/2} \Vert \phi^{(k)} \Vert_{L^\infty} \sum_{\sum \alpha_j = d-k}   \prod_{j=1}^{n-1}  \Vert \sech^{\frac{2(1+\alpha_j)}{n+d-k-1}}   w^{(\alpha_j)} \Vert_{L^{\frac{n+d-k-1}{1+\alpha_j}}}.  \] 
We apply  \eqref{eq:esthbmon} of Lemma \ref{lem:estimatedifferentialmonomial2}
and  arrive at 
\[ \Big| \int \sech^2(x/R) \phi \frac{\delta}{\delta w} I dx\Big| \le c \Vert \phi \Vert_{C^{d/2}} \Vert w \Vert_{L^2}^{n-3} \Vert \sech(x/R) w \Vert^2_{H^{\frac{n+2d-3}4}}.  \]

\bigskip

We turn to  the proof of \eqref{eq:diffflowKato} and begin with the calculation, setting $R=1$ for simplicity, 
\[ \begin{split} 
\hspace{.5cm} & \hspace{-.5cm} 
\int  w \tanh^2 \partial_x  \frac{\delta}{\delta w} \trace K^n(w) dx 
=  -\frac1n \trace \Big( K((\tanh  w)_x ) K^{n-1}(w) \Big)   
\\ & =  -  \frac{n-1}{n} \trace \Big(K ( \sech^2 w ) K^{n-1}(w)\Big) 
\\ & \quad + 
\frac{n-1}n \sum_{j=0}^{n-2}  \trace   \Big( K( \tanh w) K^j(w) K(w_x) K^{N-j-2}(w) - K( \tanh w_x )   K^{n-1}(w) \Big).  
\end{split} 
\] 
In the fashion as above with obvious modifications we obtain for $n \ge M+2$
\[    \Big|   \trace K ( \sech^2 w ) K^{n-1}(w)\Big| \le c ((\im z)^{-1/2} \Vert w \Vert_{L^2})^{n-M-2}
(\im z )^{1-M} \Vert w \Vert_{L^2}^{M} \Vert \sech w \Vert^2_{H^{\frac{M}4}}
\]
and for $I$ as above 
\[   \Big| \int \sech^2 w \frac{\delta}{\delta w} I dx\Big| \le c \Vert \phi \Vert_{C^{d/2}} \Vert w \Vert_{L^2}^{M-2} \Vert \sech w \Vert^2_{H^{\frac{n+2d-2}4}}.  \]

The last line above requires  more effort. To remove the derivative (and pay by a factor $ z $) 
we observe 
\[  K(w_x) = [\partial , K(w) ]  \]
and, with $ \tilde K(w)$ of the same structure as $K(w)$, but with some sign changes
\[ \partial K (w) = -iz \tilde K(w)  + \left( \begin{matrix} 0 & w \\  w & 0   \end{matrix} \right).  \]
By a small abuse of notation the first summand corresponds to replacing $w_x$ by $izw$. We estimate  (again $\eta$ is assumed to be $\tau$ slowly varying)  
\[  \Big\Vert  \eta \left( \begin{matrix} 0 & w_1 \\ w_1 & 0 \end{matrix} \right) K(w_2) \eta^{-1}  \Big\Vert_{\I_2} \le   c \left\{ \begin{array}{l}  \Vert w_1 \Vert_{L^2}\Vert w_2 \Vert_{L^2} \\ 
  (\im z)^{-1/2}\Vert w_1 \Vert_{L^\infty }\Vert w_2 \Vert_{L^2}     \end{array} \right. \]
and 
 \[  \Big\Vert  \eta \left( \begin{matrix} 0 & w_1 \\ w_1 & 0 \end{matrix} \right) K(w_2) \eta^{-1}  \Big\Vert_{L^2\to L^2} \le   c    \tau^{-1}  \Vert w_1 \Vert_{L^\infty}\Vert w_2 \Vert_{L^\infty} . \] 
We  interpolate the estimates for $ p \ge 2$ 
\begin{equation} \Big\Vert  \eta \left( \begin{matrix} 0 & w_1 \\ w_1 & 0 \end{matrix} \right) K(w_2) \eta^{-1}  \Big\Vert_{\I_p} \le   c \tau^{-1+ \frac{2}p}   \Vert w_1 \Vert_{L^{2p}}\Vert w_2 \Vert_{L^{2p}} .       \end{equation}

\bigskip 
First $ K(\tanh w) = K(w) \tanh$. 
In the last line above we want to commute $\tanh $ and $K(w)$. 
 The integral kernel  of the commutator of $ \tanh(x/R)$ and $R_{\pm}$ 
\[ g(x,y) = (\tanh(x) -\tanh(y) ) \chi_{\pm (x-y) > 0 }   e^{2iz |x-y|} \min  \{ \sech^2(x), \sech^2(y)\}  \]   
is bounded by 
\[  \frac1{\im z} e^{-\im z |x-y|}( \cosh(x) +\cosh(y) )^{-1}  \]
hence for $ p \ge 2$ 
\[    \Vert  \eta   [\tanh , R(w)]  \eta^{-1}  \Vert_{\I_p} \le  \frac{c}{\im z} \Vert w \Vert_{L^p}. \]
 Together, if $ n \ge M+2$,
\[ \begin{split} \hspace{1cm} & \hspace{-1cm} | \trace (K(w) \tanh K(w)\dots K(w_x) \dots        ) -\trace{K(w)\dots K(w_x)\tanh   } | 
\\ & \le c ((\im z)^{-1/2} \Vert w \Vert_{L^2} )^{n-M-2} |z| (\im z)^{-M-2}   \Vert \sech^{\frac2{M+2}} w \Vert_{L^{M+2}}^{M+2}.   
\end{split} 
\]
Let $I $ be a primitive integral as in \eqref{eq:primitive}. We apply the very same arguments and obtain 
\[ \left| \int \tanh \frac{\delta}{\delta w} I dx \right| \le c \frac{|z|}{\im z} \Vert w \Vert_{L^2}^{M-2} \Vert \sech  w\Vert^2_{H^{\frac{n+2d-2}4}}.    \]
This completes the proof. 
\end{proof}

\appendix

\section{Calculations}

\subsection{Proof of Theorem \ref{thm:formofgoodvariableequation}}\label{sec:proofofgveq}

In this section we prove Theorem \ref{thm:formofgoodvariableequation}. To do so we introduce some notation. Let $\omega(\tau) = 2\tau v$ be a rescaling of the good variable. Then if $u$ solves the $N$th KdV equation, $\omega$ solves
\begin{equation}\label{eq:Nthomega}
    \omega_t = 2\partial\Big[(\omega+2\tau)\sum_{n=-1}^{N-1}(2i\tau)^{2(N-1-n)}\frac{\delta H_n^{\KdV}}{\delta u}\Big].
\end{equation}
The relation between $u$ and $\omega$ now reads
\[
    u = -\frac12\frac{\omega_{xx}}{\omega+2\tau} + \frac34\frac{\omega_{x}^2}{(\omega+2\tau)^2}  + \frac14\omega^2 + \tau \omega.
\]
And in the case of the KdV equation $N=1$ we find e.g.
\begin{equation} 
\omega_t =  \partial_x\Big[-\omega_{xx} + 3\tau \omega^2 + \frac{1}{2}\omega^3 + \frac{3}{2} \frac{\omega_x^2}{\omega + 2\tau}\Big].
\end{equation}
Thus in order to prove Theorem \ref{thm:formofgoodvariableequation} we are lead to prove that the $N$th equation for $\omega$ can be written in the form $\omega_t = \partial_x \tilde F_N$, where
    \begin{equation}
        \tilde F_N = \sum_{n,l,K,d}(\omega+2\tau)^{-n} \tau^l \tilde f_{N,n,k,d}(\omega),
    \end{equation}
    where $\tilde f_{N,n,k,d}$ has homogeneity $k$ in $\omega$ and a total number of derivatives $d$, and the sum is restricted by
    \[
        \begin{split}
            0 &\leq n \leq 2N-1, \qquad l + k + d = 2N+1+n, \qquad l + k \geq n+1,\\
            \#&\{\text{factors of } \omega \text{ with at least 1 derivative}\} \geq n+1 \quad \text{if}\quad n\geq 1.
        \end{split}
    \]
Moreover, the linear part of the equation is $(-1)^N w^{(2N+1)}$, and $\tau^l\tilde f_{N,n,k,d}$ contains no term of the form $\tau^l \omega^k \omega^{(2N)}$ with $l + k = n+1$. The number of derivatives $d$ is always even.

Define a generalized differential monomial in $\omega$ to be an expression of the form
\begin{equation}
    \tau^l (\omega+2\tau)^{-n} \prod_{i=1}^k \omega^{(\alpha_i)}.
\end{equation}
We call $k$ its homogeneity, $n$ its negative homogeneity, and denote by $d = \sum_{i=0}^k \alpha_i$ its total number of derivatives. Moreover, we define the degree
\[
    D = l + k - n + d.
\]
A generalized differential polyonomial of degree $D$ is then a sum of generalized differential monomials of degree D. Note that expanding the fraction by powers of $\omega + 2\tau$ leaves its degree invariant. What the first part of Theorem \ref{thm:formofgoodvariableequation} now says is that the $N$th equation for $\omega$ has the form
\[
    \omega_t = \partial_x F_N,
\]
where $F_N$ is a generalized differential polynomial in $\omega$ of degree $2N+1$ with some further properties.

\begin{proof}[Proof of Theorem \ref{thm:formofgoodvariableequation}]
    First we note that
    \[
        \partial_x \frac{\prod_{i=1}^k \omega^{(\alpha_i)}}{(\omega+2\tau)^n} = \sum_{j=1}\frac{\omega^{(\alpha_j+1)}\prod_{i\neq j} \omega^{(\alpha_i)}}{(\omega+2\tau)^n} - n \frac{w_x\prod_{i=1}^k \omega^{(\alpha_i)}}{(\omega+2\tau)^{n+1}}.
    \]
    Each derivative falling onto $(\omega+2\tau)^{-1}$ produces a new factor of $\omega$ with at least one derivative. The degree of the corresponding generalized monomial gets increased by one. Moreover it is clear that the degree of the product of generalized monomials is just the sum of the two degrees.
    
    For KdV, we know
    \[
        \frac{\delta H_m^{\KdV}}{\delta u} = (-1)^m u^{(2m)} + \sum_{k=2}^{m+1}\sum_{\alpha_1+\dots+\alpha_k = 2(m+1-k)} c_{\alpha_i} \partial_x^{\alpha_1}u\dots\partial_x^{\alpha_k} u.
    \]
    We plug in
    \[
    u = -\frac12\frac{\omega_{xx}}{\omega+2\tau} + \frac34\frac{\omega_{x}^2}{(\omega+2\tau)^2}  + \frac14\omega^2 + \tau \omega,
    \]
    which is a generalized differential polynomial in $\omega$ of degree $2$. Thus for each $k$ and $\alpha_1+\dots+\alpha_k = 2(m+1-k)$, the degree of
    \[
        \partial_x^{\alpha_1}u\dots\partial_x^{\alpha_k} u
    \]
    in $\omega$ is $2k + 2(m+1-k) = 2(m+1)$. The same holds for the linear term $u^{(2m)}$. This shows that
    \[
    \sum_{n=-1}^{N-1}(2i\tau)^{2(N-1-n)}\frac{\delta H_n^{\KdV}}{\delta u}
    \]
    has degree $2N$ in $\omega$, and the factor $(\omega+2\tau)$ in the equation for $\omega_t$ increases the degree again by one. It remains to prove that $0 \leq n \leq 2N-1$, $l+k \geq n+1$, and to analyze the term for $n = 0$. 
    
    To prove $0 \leq n \leq 2N-1$ we see that in order to create higher $n$, there has to be either a high power or a big amount of derivatives falling onto $(\omega+2\tau)^{-2}\omega_{x}^2$. The highest possible power comes from $u^{N}$ (for $m = N-1$) and leads to $n = 2N$, which after multiplication with $(\omega+2\tau)$ becomes $2N-1$. The term with most derivatives is $u^{(2N-2)}$ which again could leads to a term with $n = 2N$, or $2N-1$ after multiplication. The terms in between the extremes are handled likewise.
    
    The inequality $l+k \geq n+1$ is equivalent to $d \leq 2N$ by the degree condition. This in fact follows again by looking at the worst terms: if the full powers of $u^{N}$ hit either $(\omega+2\tau)^{-1}\omega_{xx}$ or $(\omega+2\tau)^{-2}\omega_{x}^2$ we arrive at $d = 2N$, and so do we if all derivatives from $u^{(2N-2)}$ hit $(\omega+2\tau)^{-1}\omega_{xx}$ or $(\omega+2\tau)^{-2}\omega_{x}^2$.
    
    To see that the number of factors with at least one derivative is $\geq n+1$ for $n \geq 1$, note that the only way to produce powers of $(\omega + 2\tau)^{-1}$ is if $u$ is one of the factors $(\omega + 2\tau)^{-1}\omega_{xx}$ and $(\omega + 2\tau)^{-2}\omega_{x}^2$. In both cases, the number is $\geq n$. Multiplication with other factors will leave these properties invariant, as does taking derivatives. Multiplication by $(\omega + 2\tau)$ in the end leads to the statement. This also gives $d \geq n+1$ if $n \geq 1$, and $l + k \leq 2N$ if $n \geq 1$, respectively $l + k \leq 2N+1$ if $n = 0$, because in this case $d = n = 0$ may happen. It also shows $k \geq n+1$.
    
    The number of derivatives $d$ is always even, because this statement holds for the KdV monomials (which in turn follows from the Lenard recursion), and because writing $u$ in terms of $\omega$ only contains terms with an even number of derivatives.
    
    We analyze the term with $n = 0$. Consider $\delta H_m^{\KdV}/\delta u $ with $u$ as above. We can only reach $n = 0$ if all copies of $u$ are $\frac14\omega^2 +\tau \omega$, or if exactly one $u$ is $(\omega+2\tau)^{-1}\omega_{xx}$ (because in this case the factor of $(\omega+2\tau)$ decreases $n$ by one to zero). In the latter case, all possible derivatives have to fall onto $\omega_{xx}$. A linear term in $\omega$ can only be reached by
    \[
        2(\omega + 2\tau)(\tau \omega)^{(2m)} = 4\tau^2 \omega^{(2m)} + \text{bilinear},
    \]
    and by
    \[
        2(\omega + 2\tau)\Big(-\frac12\frac{\omega_{xx}}{\omega+2\tau}\Big)^{(2m)} = -\omega^{(2m+2)} + \text{higher homogeneities}.
    \]
    Using the summation over $-1 \leq m \leq N-1$ from \eqref{eq:Nthomega} cancels out all linear terms but the one of order $2N$.
    
    We consider the term with $d = 2N$ derivatives and where $k + l = n+1, n \geq 1$. $2N$ derivatives can only be reached if we are in the case
    \[
        2(\omega + 2\tau)u^{(2N-2)}
    \]
    where $u$ is either $\frac{\omega_{xx}}{\omega+2\tau}$ or $\frac{\omega_{x}^2}{(\omega+2\tau)^2}$. In both cases, the only way to create more homogeneity is by factors of $\omega_x$, which shows the impossibiliy of a term $(\omega+2\tau)^{-n}(\omega^n \omega^{(2N)})$.
    
    We turn to the form of the equation for $v = \omega/(2\tau)$. In general if we have a $k$-linear form $A_k$ and an equation
    \[
        \omega_t = \frac{A_k(\omega)}{(\omega + 2\tau)^n},
    \]
    then the corresponding equation for $v$ will be
    \[
        v_t = (2\tau)^{k-1-n}\frac{A_k(v)}{(v + 1)^n}.
    \]
    Thus the new power of $\tau$ in the transition from $\omega$ to $v$ will be $l + k - n - 1$, and the degree condition becomes
    \[
        l + d = 2N.
    \]
    The upper bound on the homogeneity $k$ for $\omega$, $k \leq 2N+1$ does not change when going to $v$, because it does not depend on $l$. Similarly, the number of derivatives $d$ stays even. Due to the condition $d + l = 2N$, $l$ has to be even as well.
\end{proof}

\subsection{Good Variables}

For $N = 1$ the good variable equation is
\begin{equation*} 
    v_t =  \partial_x\Big[-v_{xx} + 6\tau^2 v^2 + 2\tau^2 v^3 + \frac{3}{2} \frac{v_x^2}{v+1}\Big],
\end{equation*}
and for $N = 2$
\begin{equation*}
\begin{split}
    v_t =&  \partial_x\Big[v_{xxxx} - 7\tau^2v_{xx}v^2 - 4\tau^2v v_x^2 -14\tau^2 v_{xx}v\\
    &-4\tau^2v_x^2 + 6\tau^4v^5 +30\tau^4v^4 +40\tau^4v^3\\
    &+ (v+1)^{-1}\big(-\frac52v_{xx}^2 - 5v_{xxx}v_x+18\tau^2v_x^2v + \frac92\tau^2v_x^2v^2 -6\tau^2v_x^2\big)\\
    &+(v+1)^{-2}\Big(\frac{25}{2}v_{xx}v_x^2\Big) + (v+1)^{-3}\Big(-\frac{45}{8}v_x^4\Big)\Big].
\end{split}
\end{equation*}

\section{The AKNS hierarchy}
\subsection{The hierarchy}

The nonlinear Schrödinger equation (NLS) is the equation
\begin{equation}\label{eq:NLS}
    iq_t = -q_{xx} \pm 2|q|^2 q,
\end{equation}
where the sign $+$ corresponds to the defocusing NLS and $-$ to the focusing NLS. The complex modified Korteweg-de Vries equation (we will abbreviate both the complex and the real mKdV just by mKdV) is the equation
\begin{equation}\label{eq:cmKdV}
    q_t = -q_{xxx} \pm 6|q|^2 q_x,
\end{equation}
and again we talk about defocusing and focusing mKdV. Both equations are Hamiltonian equations of the form
\begin{equation*}
    iq_t = \left.\frac{\delta H}{\delta r}\right|_{r = \pm \bar q}
\end{equation*}
with Hamiltonian $H = \int q'r' + q^2r^2 dx$ for NLS and Hamiltonian $H = -i \int q'r'' + 3q^2rr'$ for mKdV. They are special cases of the AKNS system of equations
\begin{equation*}
    \begin{split}
        iq_t &= -q_{xx} + 2q^2 r,\\
        ir_t &= r_{xx} - 2r^2 q,
    \end{split}
\end{equation*}
for NLS, and
\begin{equation*}
    \begin{split}
        iq_t &= -iq_{xxx} + 6iqq'r,\\
        ir_t &= -ir_{xxx} - 6irr' q,
    \end{split}
\end{equation*}
for mKdV, with the additional restriction $r = \pm \bar q$. The AKNS equations can be written as the Hamiltonian system of equations 
\begin{equation}\label{eq:AKNS}
        iq_t = \frac{\delta H}{\delta r}, \quad ir_t = -\frac{\delta H}{\delta q},
\end{equation}
where one takes the Hamiltonians as above with the Poisson structure
\begin{definition}\label{defpoissonstructure}
The AKNS symplectic form on $L^2(\R;  \C^2) $ is  given by 
\[  \omega( (q_1,r_1) , (q_2, r_2)) = \int q_1 r_2 - q_2 r_1 dx. \]  
It defines the  AKNS Poisson structure 
    \[\{F,G\}=- i \int \frac{\delta F}{\delta q}\frac{\delta G}{\delta r}  - \frac{\delta F}{\delta r}  \frac{\delta G}{\delta q} dx .\]
\end{definition}
A Lax operator for the AKNS equations is given by
\begin{equation}\label{eq:laxode}
     L(q,r) \phi := i\left( \begin{matrix}  - \partial & q \\ -r &  \partial \end{matrix} \right)\phi 
\end{equation}
with  a Lax equation 
\[ L(q,r) \phi = z \phi. \]
The  Wadati Laplace operator is a special case (see Subsection \ref{subsec:determinants} and it seems  worthwhile to explore the various and strong relations between all the hierarchies.

We assume at first that $ q$ and $r$ decay fast and $ z = \xi \in \R$. There exist two fundamental systems, $ \psi_{-+}, \psi_{--}$  normalized at $-\infty$
\[  \lim_{x\to -\infty}  e^{i\xi x}\psi_{-+}(x) = \left( \begin{matrix} 1 \\ 0 \end{matrix} \right), \qquad   \lim_{x\to -\infty}  e^{-i\xi x}\psi_{--}(x) = \left( \begin{matrix} 0 \\ 1 \end{matrix} \right) \] 
and a fundamental system normalized at $ \infty$
\[  \lim_{x\to \infty}  e^{-i\xi x}\psi_{+-}(x) = \left( \begin{matrix} 0 \\ 1 \end{matrix} \right), \qquad  \lim_{x\to \infty}  e^{i\xi x}\psi_{++}(x) = \left( \begin{matrix} 1 \\  0 \end{matrix} \right). \] 
The solution space of the problem \eqref{eq:laxode} is a two-dimensional vector space. As a consequence
these solutions are connected on the real line by 
\[ \left( \begin{matrix} \psi_{+-}\\  \psi_{-+} \end{matrix} \right)  = \left( \begin{matrix} a_+(\xi) & b_+(\xi) \\      b_-(\xi) & a_-(\xi)         \end{matrix}  \right)          \left( \begin{matrix} \psi_{--}\\ \psi_{++}\end{matrix} \right). \] 
There  are  simple alternative expression 
\[ a_+(\xi) = W( \psi_{-+}, \psi_{+-})= \det ( \psi_{-+}, \psi_{+-} ) , \quad a_-(\xi) = W(\psi_{++}, \psi_{--}) .\] 
Here $ W$ is the Wronskian, which is independent of $x$.
The solutions $ \psi_{-+}$ and $\psi_{+-}$  have a holomorphic extensions to the upper half plane called left and right Jost functions and we define the transmission coefficient 
\[ T(z)^{-1} = W( \psi_{-+}, \psi_{+-}) = \lim_{x \to \infty}e^{izx}\psi^1_{-+}(x) = \lim_{x \to -\infty} e^{-izx}\psi^2_{+-}(x).\] 
for $ z$ in the upper half plane. The solutions $ \psi_{--}$ and $ \psi_{++} $ have a holomorphic extension to the lower half plane and we define 
\[ T(z) = W(\psi_{++} , \psi_{--} ) = \lim_{x\to -\infty}e^{izx}\psi_{++}^1(x) = \lim_{x \to  \infty}e^{-izx}\psi_{--}^2(x) \] 
on  the lower half plane.  The reason for this choice are the same as for KdV: It is a choice which gives simultaneous asymptotic series in the lower and the upper half plane. We will see later that
\begin{equation}\label{eq:AKNSasymptotic}
    - \log T^{\AKNS}(z,q,r) \sim \frac{i}{2z}\int qr\, dx + \frac{1}{(2z)^2} \int qr_x\, dx + \frac{i}{(2z)^3}\int q_x r_x + q^2 r^2\, dx + \dots.
\end{equation}
\begin{definition} \label{defT} 
The functions $\alpha(x,z,q,r),\beta(x,z,q,r)$ and $\gamma(x,z,q,r)$ are defined by
\[
\left( \begin{matrix} \frac{\gamma}2&\! \alpha \\ \beta &\! \frac{\gamma}2 \end{matrix} \right) 
= \left\{ \begin{array}{ll} \displaystyle \frac{T}2\left( \begin{matrix} \psi_{-+}^1\psi^2_{+-} + \psi_{-+}^2\psi_{+-}^1 & 2 \psi_{-+}^1 \psi_{+-}^1 \\ 2\psi^2_{-+} \psi_{+-}^2 &\psi_{-+}^1\psi^2_{+-} + \psi_{-+}^2\psi_{+-}^1 \end{matrix} \right) & \text{ if } \im z >0 \\[5mm] \displaystyle \frac{1}{2T} \left( \begin{matrix} \psi_{--}^1\psi^2_{++} + \psi_{--}^2\psi_{++}^1 & 2 \psi_{--}^1 \psi_{++}^1 \\ 2\psi^2_{--} \psi_{++}^2 &\psi_{--}^1\psi^2_{++} + \psi_{--}^2\psi_{++}^1 \end{matrix} \right)& \text{ if } \im z < 0. 
\end{array} 
\right. 
\]
\end{definition}
Observe that for each $z$, we have $\gamma\to 1, $ and $ \alpha \to 0$ and $ \beta \to 0 $ when $x\rightarrow \pm\infty.$

\begin{lemma}\label{lemma:GFAKNS}
For $z$ in the upper half-plane, the Green's function for the operator $L -z 1$ is 
\begin{equation*}  G(x,y,z) = -iT(z)\left\{ \begin{array}{rl} \left( \begin{matrix} \psi^1_{-+}(x,z) \psi^2_{+-}(y,z)& \psi^1_{-+}(x,z) \psi^1_{+-}(y,z)  \\ \psi^2_{-+}(x,z) \psi^2_{+-}(y,z) & \psi^2_{-+}(x,z) \psi^1_{+-}(y,z)    \end{matrix} \right) & \text{ if } x<y, \\[5mm] 
\left(  \begin{matrix} 
\psi^1_{+-}(x,z)  \psi^2_{-+}(y,z) & \psi^1_{+-}(x,z) \psi^1_{-+}(y,z) \\ \psi^2_{+-}(x,z) \psi^2_{-+}(y,z) & \psi^2_{+-}(x,z) \psi^1_{-+}(y,z)    \end{matrix} \right) & \text{ if } y<x.\end{array}\right. \end{equation*}
\end{lemma}
\begin{proof}
 We observe that the columns considered as functions of $x$ satisfy 
 \[ L G =zG \] 
 whenever $ x \ne y$. It is the Green's function since, for $x^+$ being the limit from above and $x^-$ being the limit from below,
 \[ G(x^+,x)-G(x^-,x) = -i\left( \begin{matrix} -1 & 0 \\ 0 & 1 \end{matrix} \right), \]
 which shows that $(L-z1)G(x,y) = \delta(x-y)$.
\end{proof}

The function $ \alpha$, $ \beta $ and $ \gamma$ are closely related to the Green's function. 

\begin{lemma} Let 
\[ g(z,x) = \lim_{h\to 0+} \frac{i}2 (G(x+h,x,z) + G(x, x+h,z) ) \] 
be the diagonal Green's function. Then 
\[ g(z,x) = \left\{ 
\begin{array}{cl} 
\left( \begin{matrix}  \gamma/2 & \alpha \\ \beta &  \gamma/2 \end{matrix} \right) & \text{ if } \im z>0 \\[4mm]  
-\left( \begin{matrix} \gamma/2 & \alpha \\ \beta & \gamma/2 \end{matrix} \right) & \text{ if } \im z <0. \end{array} \right. \] 
\end{lemma}

The importance of these objects is that they characterize $\log T$ and its functional derivatives. We introduce the resolvents $R_{\pm} = (-iz \pm  \partial)^{-1}$ for $ \im z >0$,
\begin{equation} \label{eq:integralkernelsresolvent}
R_+ f(x) = \int_{-\infty}^x   e^{iz(x-y)} f(y) dy, \quad R_- f(x) = 
\int_x^\infty e^{-iz(x-y)} f(y) dy. 
\end{equation} 
Then 
\[ L(q,r)-z1 = (L(0,0)-z1) \left( 1+ \left( \begin{matrix} 0 & R_- q \\ - R_+ r & 0 \end{matrix} \right) \right).\]
Unfortunately the operator in the bracket is only Hilbert-Schmidt for $ q,r \in L^2$, but not trace class, even for Schwartz functions. For trace class operators $K$  one has the expansion 
\begin{equation}\label{eq:logdeta}   \ln \det (1 - K ) =  \sum_{n=1}^\infty \frac1n \trace K^n \end{equation} 
where $ \trace K^n $ is defined for $K$ in the $L^n$ Schatten class. In particular only the first term is problematic for the bracket above. 
 On the other hand, formally at least this trace should be zero due to the off-diagonal block matrix form of the operator. This motivates the use of the renormalized determinant
\[ {\det}_2 ( 1+ K ) = \det( I +K) \exp( - \trace K) \] 
for trace class functions, which has a unique extension to Hilbert Schmidt operators $K$, see . We refer to Subsection \ref{subsec:determinants} for details.

The first set of statements below are elementary results from the theory of ODEs. For the reader's convenience we collect the proofs in the Appendix \ref{sct:functionalderivatives}.
\begin{lemma}\label{lemma:NLSfunctionalderivative} 
    $\alpha,\beta,\gamma$ are connected to the transmission coefficient via 
    \begin{equation} 
   \frac{d}{dz} \log T^{\AKNS} = i \int \gamma-1 \, dx,\quad 
        \alpha = \frac{\delta}{\delta r}\log T,\quad
        \beta = -\frac{\delta}{\delta q}\log T. \label{vardeqlnT}
    \end{equation}
Let $ \im z >0$ and $ r,q \in L^2$. Then the Fredholm determinant below is well defined (with the interpretation described above) and, if $ \im z >0$,
\begin{equation} \label{functionaldeterminant}  {\det}_2\left( 1+ i(L(0,0)-z1)^{-1}  \left(\begin{matrix} 0 & q \\ -r & 0 \end{matrix} \right) \right) = T(z,q,r)^{-1}. \end{equation}
\end{lemma}

\begin{proof} We provide short conceptional proof for  \eqref{functionaldeterminant} (which is well known, see \cite{MR1723386} )     by calculating the derivative of the functional determinant with respect to the potentials. 
This requires a bit of care. 
 We observe that both sides are identically $1$ if $q=r=0$. By an abuse of notation, $G(z,q,r) := -i(L(q,r)-z1)^{-1}$ whenever it is defined. 
We approximate $r$ and $q$ by Schwartz functions and replace $R_\pm$ by the trace class operators $R_{\pm}^\tau = [ (1- \frac1{\tau} \partial)(-iz \pm \partial) ]^{-1}$ which have integral kernels 
\[ k^\tau_{\pm}(x,y) = (e^{\pm iz(\cdot)}\chi_{\{\pm(\cdot)>0\}}) * (\tau e^{\tau(\cdot)}\chi_{\{-(\cdot)>0\}})(x-y).\]
Then, with the obvious notation 
\[ \ln {\det}_2 \left( 1+ G(z,0,0) \left( \begin{matrix} 0 & tq \\ - tr & 0 \end{matrix} \right) \right) 
= \lim_{\tau\to \infty}\ln \det  \left( 1+ G^\tau(z,0,0) \left( \begin{matrix} 0 & tq \\ - tr & 0 \end{matrix} \right) \right) 
\] 
since the trace of the second summand on the right hand side vanishes, and the operator converges in the Hilbert-Schmidt norm. Moreover,
\[ \int \beta^\tau q - \alpha^\tau r\, dx \to \int \beta q - \alpha r \,dx \] 
since the Green's function converges.  We use  the operator identity $1+(t+s)A = (1+tA)(1 + (1+tA)^{-1}sA)$ below, and calculate
\[\begin{split} \hspace{0.5cm} & \hspace{-0.5cm}  \frac{d}{dt} \log {\det}\left(1+ G^\tau(z,0,0) \left( \begin{matrix} 0 & tq \\ -tr & 0 \end{matrix} \right) \right)  
\\ & = \frac{d}{ds} \log {\det}\left(1+ G^\tau(z,0,0) \left( \begin{matrix} 0 & (t+s)q \\ -(t+s)r & 0 \end{matrix} \right) \right)_{s=0}  
\\ & = \frac{d}{ds} \log \det \left[ 1 + s\left( 1+ G^\tau(z,0,0)\left( \begin{matrix} 0 & tq \\ -tr & 0 \end{matrix} \right)\right)^{-1} G^\tau(z,0,0) \left( \begin{matrix} 0 & q \\ -r & 0 \end{matrix} \right) \right]_{s=0} 
\\ & = \frac{d}{ds} \log \det \left[ 1 + s G^\tau(z,tq,tr) \left( \begin{matrix} 0 & q \\ -r & 0 \end{matrix} \right) \right]_{s=0} 
=\trace \left[ G^\tau(z,tq,tr) \left( \begin{matrix} 0 & q \\ -r & 0 \end{matrix} \right) \right] 
\\ & = \int \beta^\tau(z,tq,tr) q - \alpha^\tau(z,tq,tr) r dx \qquad   \to \quad  -\frac{d}{dt} \log T(z, tq,tr) \quad \text{ as  } \tau\to \infty .
\end{split} 
\]

\end{proof}

\begin{lemma}\label{lemma:NLSODE}
    $\alpha,\beta,\gamma$ satisfy the ODE
    \begin{equation}\label{eq:NLSODE} \begin{split}  \gamma'\, & = 2(q \beta  + r \alpha ) \\  
    \alpha'\, &  = -2iz \alpha    + q \gamma \\  
    \beta'\, &  = 2iz  \beta   + r\gamma. 
    \end{split} 
    \end{equation}
\end{lemma}
\begin{proof}

The statements follow from  differentiating products of components of solutions to the ODE \eqref{eq:laxode}. We carry this out for one term, the others follow likewise.
 \[ \frac{d}{dx} ( \psi^1_{-+}\psi^2_{+-} + \psi^2_{-+} \psi^1_{+-} ) = 2 q \psi^2_{-+} \psi^2_{+-} +2r \psi^1_{-+} \psi^1_{+-}. \]
\end{proof}
Combining the equations \eqref{eq:NLSODE} gives
\begin{equation*}
    (\alpha(z_1)\beta(z_2))' = (2iz_2 - 2iz_1)\alpha(z_1)\beta(z_2) + q\gamma(z_1)\beta(z_2) + r\gamma(z_2)\alpha(z_1).
\end{equation*}
This has two important consequences. The first follows from setting $z_1 = z_2$ and yields an alternative equation for $\gamma$,
\begin{equation} \label{alternativegamma} 
    \gamma^2 = 1+4\alpha \beta,
\end{equation}
The second consequence is that
\begin{equation}\label{eq:ODEAKNScommutation}
\begin{split} 
\\ \hspace{3cm} & \hspace{-3cm}   \left(\alpha(z_1)\beta(z_2) + \alpha(z_2)\beta(z_1)-\frac{1}{2}\gamma(z_1)\gamma(z_2)\right)' 
\\ &   = 2i(z_2 - z_1)(\alpha(z_1)\beta(z_2) - \alpha(z_2)\beta(z_1))
\end{split}
\end{equation}
which can be used to show that the transmission coefficients are Poisson commuting, see Theorem~\ref{thm:poissonbrackets}.

\subsubsection{Symplectic forms and Poisson structures}

%It is easy to see that the Jacobi identity is satisfied 
%\[ \{ \{ F,G\},H\} + \{ \{ G,H\}, F \} +\{ \{ H,F\}, G\} = 0. \] 

%Note that while the above Poisson bracket is only defined for functionals, we can extend the definition to functions $f(x)$ by identifying them as the point evaluation functionals $f \mapsto f(x)$ having $x$ as a parameter. In particular we find that $\frac{\delta q(x)}{\delta q}(y) = \delta(x-y)$, $\frac{\delta r(x)}{\delta q}(y) = 0$ and the $n$th AKNS equations take the form
%\begin{equation*}
%    q_t = \{q,H_n\}, \quad r_t = \{r,H_n\}.
%\end{equation*}
It is convenient to consider Poisson brackets of operators $A(p,q): X \to Y $ with functions as follows: Let $\phi\in X $ and $ \psi \in Y^*$. Then 
$(p,q) \to \psi(A(p,q)(\phi)) $ is a function, and we define the Poisson product of $A(p,q)$ with a function $H$ as the operator defined by 
\[  \psi(   \{ A, H \} (\phi)) = \{ \psi( A(\phi)) , H\} \] 
whenever this is defined. In particular, if $A$ is the multiplication by a differential polynomial then the Poisson product of the multiplication operator is the multiplication  by the Poisson products. 

We compute some Poisson brackets. The proof of the next theorem is strongly inspired by \cite{harropgriffiths2020sharp}.\footnote{Note that our notation slightly differs from theirs. In particular, their $\gamma+1$ corresponds to our $\gamma$.}

\begin{theorem}\label{thm:poissonbrackets} The transmission coefficients $T(z_1)$ and $T(z_2)$ and its logarithms Poisson commute: 
\begin{equation} \label{PoissoncommutelnT}   \{  \log T(z_1), \log T(z_2) \}  = 0. \end{equation}  
Moreover, we have
\begin{equation}\label{Poissonqr} \{  q(x), \log T(z)\}  = -i\alpha(x,z), \quad \{  r(x), \log T(z)\} = -i  \beta(x,z). \end{equation} 
\begin{equation}  \label{laxAKNS} 
\begin{split} 
\hspace{.5cm} & \hspace{-.5cm} \{ -i (L(q,r)-z_1 1), \log T(z_2)\}  = -i\left( \begin{matrix} 0 & \alpha(z_2) \\ -\beta(z_2) & 0 \end{matrix} \right)  \\ & =  -\frac{1}{2(z_2-z_1)} \left( (L-z_11) \Big( \begin{matrix} 0 & \alpha(z_2) \\ -\beta(z_2) & 0 \end{matrix} \Big) +  \Big( \begin{matrix} 0 & \alpha(z_2) \\ -\beta(z_2) & 0 \end{matrix} \Big) (L-z_11) \right) 
\\ & \quad + \frac1{4(z_2-z_1)}\left(  (L-z_11) \Big( \begin{matrix} \gamma(z_2) & 0 \\ 0 & -\gamma(z_2) \end{matrix} \Big) -  \Big( \begin{matrix} \gamma(z_2) & 0 \\ 0 & - \gamma(z_2) \end{matrix} \Big) (L-z_11) \right). 
\end{split} 
\end{equation}  
Also 
\begin{equation}  \label{poissongamma} 
\begin{split}
    \{  \alpha(z_1), \log T(z_2) \} &= \frac{1}{2(z_1-z_2)}(\alpha(z_1)\gamma(z_2)-\alpha(z_2)\gamma(z_1)),\\
    \{  \beta(z_1), \log T(z_2) \} &= \frac{1}{2(z_1-z_2)}(-\beta(z_1)\gamma(z_2)+\beta(z_2)\gamma(z_1)),\\
    \{  \gamma(z_1), \log T(z_2) \} &= \frac{1}{2i (z_1-z_2)^2} \partial_x \Big\{ \alpha(z_1) \beta(z_2)+ \beta(z_1) \alpha(z_2) -\frac12 \gamma(z_1) \gamma(z_2) \Big\}. 
\end{split}    \end{equation} 
\end{theorem}

\begin{proof}
$\log T(z_1)$ and $\log T(z_2)$ Poisson commute as a consequence of \eqref{eq:ODEAKNScommutation}. Since the functional derivatives of the logarithms are proportional to the functional derivatives of the functions, the same argument works for the transmission coefficients itself.

Equations \eqref{Poissonqr} follow from $\frac{\delta q(x)}{\delta q}(y) = \delta(x-y)$, the Dirac measure, and $ \frac{\delta }{\delta r} q(x) = 0 $. Then 
\[ \{ q(x), \log T(z) \}= -i \int  \delta(x-y)\alpha(y) dy = -i \alpha(x), \]
and similarly for $\{ r(x), \log T(z)\} $. Also the  first line of \eqref{laxAKNS} follows from the definition of the Poisson structure (we drop the evaluation in the notation). Using \eqref{eq:NLSODE} we see that
\[ \begin{split} \hspace{1.5cm} & \hspace{-1.5cm}  L(z_1) \Big( \begin{matrix} 0 & \alpha(z_2) \\ -\beta(z_2) & 0 \end{matrix} \Big) +  \Big( \begin{matrix} 0 & \alpha(z_2) \\ -\beta(z_2) & 0 \end{matrix} \Big) L(z_1)\\ 
& = \left( \begin{matrix} - q\alpha(z_2)  - r \beta(z_2) & -2iz_1 \alpha(z_2) -\alpha'(z_2) \\ 2iz_1 \beta(z_2) -\beta'(z_2) & -q \alpha(z_2) - r \beta(z_2) \end{matrix} \right)  
\\ & =  
 \left( \begin{matrix} -\frac12 \gamma'(z_2) & - q \gamma(z_2) \\  - r \gamma(z_2)  & -\frac12 \gamma'(z_2) \end{matrix}\right)  + 2i (z_2-z_1) \left( \begin{matrix} 0 & \alpha(z_2) \\ -\beta(z_2) & 0\end{matrix} \right) 
 \end{split} 
 \] 
and 
\[  L(z_1) \Big( \begin{matrix} \gamma(z_2) & 0 \\ 0 & -\gamma(z_2) \end{matrix} \Big) -  \Big( \begin{matrix} \gamma(z_2) & 0 \\ 0 & - \gamma(z_2) \end{matrix} \Big) L(z_1)= 
 \left( \begin{matrix} - \gamma'(z_2)  & -2q \gamma(z_2)\\ -2r \gamma(z_2)  &  -\gamma'(z_2) \end{matrix} \right). 
\] 
We sum the terms to arrive at \eqref{laxAKNS}.

We begin with the differentiation of the resolvent $L^{-1} = (L_0-iz)^{-1}$ to prove the last Poisson brackets. Then using \eqref{laxAKNS} and the resolvent identity $(A-B)^{-1}- A^{-1} =    (A-B)^{-1}  B         A^{-1}$,
we see that
\[\begin{split} \hspace{.5cm} & \hspace{-.5cm}  \{L^{-1}(z_1) , \log T(z_2) \}  = -i L^{-1}(z_1)  \left( \begin{matrix} 0 & \alpha(z_2)  \\ - \beta(z_2) & 0 \end{matrix} \right) L^{-1}(z_1) 
\\ & = \frac{1}{2(z_2-z_1)}\left\{  \left( \begin{matrix} 0 & \alpha(z_2)  \\ - \beta(z_2) & 0 \end{matrix} \right) L^{-1}(z_1) + L^{-1}(z_1) \left( \begin{matrix} 0 & \alpha(z_2)  \\ - \beta(z_2) & 0 \end{matrix} \right) \right\} 
\\ & \qquad - \frac1{4(z_2-z_1)} \left\{ \left( \begin{matrix} \gamma(z_2) & 0 \\ 0 & - \gamma(z_2)  \end{matrix} \right) L^{-1}(z_1)- L^{-1}(z_1) \left( \begin{matrix} \gamma(z_2) & 0 \\ 0 &- \gamma(z_2) \end{matrix} \right)   \right\}
\end{split} 
\] 
We evaluate  the integral kernel at $ x=y$, which is possible since the kernel of RHS is continuous. 

\[\begin{split}  \{ g(z_1) , \log T(z_2)\}\, &  = \frac{1}{2(z_2-z_1)}\left( \begin{matrix} \alpha(z_2) \beta(z_1) - \alpha(z_1) \beta(z_2)\hspace{-12pt}   & \alpha(z_2) \gamma(z_1) \\ - \beta(z_2) \gamma(z_1) & \alpha(z_2) \beta(z_1) - \alpha(z_1) \beta(z_2) \end{matrix} \right) 
\\ & \qquad - \frac{\gamma(z_2)}{2(z_2-z_1)} \left( \begin{matrix} 0 & \alpha(z_1) \\  -\beta(z_1)&0\end{matrix}  \right).  
\end{split} 
\] 
The claimed equality follows by \eqref{eq:ODEAKNScommutation}. 
\end{proof}

From now on we assume that at least one of the two functions $q,r$ is decaying. Otherwise the coefficients in the asymptotic expansion below might be undefined. In the situation where both functions are nondecaying, one may still be able to do an asymptotic expansion in a different spectral parameter, see \cite{MR4186010}. The AKNS Hamiltonians are defined as the coefficients in the asymptotic series
\begin{equation}\label{eq:defNLShierarchy}
	\log  T(z) \sim -i \sum_{n=1}^{\infty} (2 z)^{-n} H_{n}^{\AKNS}.
\end{equation}
As above we define
\[  \T^{\AKNS}_N =  (2z)^N \log  T(z)+ i \sum_{n=1}^{N} (2 z)^{N-n} H_{n}^{\AKNS}.  \]

Similarly, we use $\alpha$, $\beta$ and $\gamma$ as generating functions, 
\[ \gamma \sim 1 + \sum_{n=1}^\infty (2z)^{-n}\gamma_n, \quad \alpha \sim \sum_{n=1}^\infty (2z)^{-n}\alpha_n, \quad \beta_n \sim \sum_{n=1}^\infty (2z)^{-n}\beta_n.  \] 
Note that these expansions are not absolutely convergent but only asymptotic, and they hold for $z \in \C \backslash \R$.
Their precise meaning is contained in the following theorem.

\begin{theorem}\label{thm:asymptoticexpansion}
The transmission coefficient is a meromorphic function in $z$ for $ z \in \C \backslash \R$. The poles coincide with the eigenvalues in the upper half plane. 
The following estimates hold for $N \in \N$: 
\begin{equation} \label{eq:hn} 
|H_N(q,r) |   \le c_N  (\Vert q^{(\frac{N-1}2)} \Vert^2_{L^2}  + \Vert r^{( \frac{N-1}2) } \Vert_{L^2}^{2} + \Vert q \Vert^{2N}_{L^2} + \Vert r \Vert^{2N}_{L^2}  ) 
\end{equation} 
Let $z \in \C$ and $q,r \in L^2$ be such that 
\[   |\im z|^{-\frac12} ( \Vert q \Vert_{L^2} + \Vert r \Vert_{L^2} ) \le \frac1{100}. \] 
\begin{equation} \label{eq:TNAKNS}   \left| \T_N^{\AKNS} \right|   \le c_N\left( \frac{|z|}{\im z} \right)^{2(N-1)}  \Big(\Vert q^{(\frac{N-1}2)}  \Vert_{L^2}^2 + \Vert r^{(\frac{N-1}2)}  \Vert_{L^2}^2 + \Vert q \Vert^{2N}_{L^2} + \Vert r \Vert^{2N)}_{L^2} \Big)  
\end{equation} 
\begin{equation} \label{eq:difftnakns}  \begin{split} 
  \left\| \frac{\delta}{\delta (q,r)} \T^ {\AKNS}_{N}   \right\|_{H^{-\frac{N}2}}
 \le \, & c_N\left( \frac{|z|}{\im z} \right)^{2(N-1)}   \Big(1+ \Vert q^{(\frac{N-1}2)} \Vert_{L^2} + \Vert r^{(\frac{N-1}2)}  \Vert_{L^2}^2 + \Vert q \Vert_{L^2}^{2N-1} + \Vert r \Vert_{L^2}^{2N-1} \Big)
 \\ & \times \Big(\Vert q \Vert_{H^{\frac{N-1}2}} + \Vert r \Vert_{H^{\frac{N-1}2} } \Big). 
 \end{split} 
\end{equation} 
All expression here are holomorphic in $z$, $r$ and $q$. Derivatives can be controlled by the Cauchy integral formula.  
\end{theorem} 
\begin{proof} 
The Lax equation  has the following symmetries. Suppose that $\left(z,q,r, \left( \begin{matrix} \phi^1 \\ \phi^2 \end{matrix} \right)\right) $ satisfy \eqref{eq:laxode}. Then 
\begin{enumerate} 
\item (Translation symmetry) $\big(z, q(.+h), r(.+h), \phi(.+h)\big)$, $ h \in \R $ 
\item (scaling symmetry) $(\lambda z , \lambda q( \lambda . ) , \lambda r( \lambda .), \phi( \lambda .) ) $, $ \lambda >0 $ and 
 \item  (Galilean symmetry) $\left(z + \xi,  e^{-2i\xi x} q, e^{2i\xi x} r, \left(\begin{matrix} e^{-i\xi x} \phi^1 \\ e^{i\xi x} \phi^2 \end{matrix} \right) \right)$, $ \xi \in \R$ 
\end{enumerate} 
all satisfy \eqref{eq:laxode}. 

As a consequence 
\begin{equation} T^{\AKNS} ( \lambda z, \lambda q(\lambda x) , \lambda r(\lambda x)) = T( z,q,r)   \end{equation} 
\begin{equation} 
T^{\AKNS}( z  + \xi, e^{-2i\xi x} q, e^{2i\xi x} r) = T( z, q,r). 
\end{equation} 
We compare the asymptotic series: 
\[ \log T^{\AKNS}(\lambda z, \lambda q(\lambda x), \lambda r(\lambda x)) \sim -i \sum_{n=1}^\infty (2\lambda z)^{-n} H^{\AKNS}_n ( \lambda q (\lambda x) , \lambda r(\lambda x)) \sim    -i \sum_{n=1}^\infty (2z)^{-n} H^{\AKNS}_n(q,r)   \]
and 
\[       H^{\AKNS}_N( \lambda q(\lambda .), \lambda r(\lambda .) ) = \lambda^{N} H_N^{\AKNS}(q,r).  \]
Similarly 
\[\begin{split}  \log T^{\AKNS}(i\tau+\xi, e^{-2i\xi x} q,  e^{2ixi x}  r)\, &  \sim -i \sum_{n=1}^\infty (2(i\tau +\xi))^{-n} H^{\AKNS}_n ( \lambda q (\lambda x) , \lambda r(\lambda x))
\\ & \sim    -i \sum_{n=1}^\infty (2i\tau)^{-n} H^{\AKNS}_n(q,r)  
\end{split} 
\]
and,
\[\begin{split}  \sum_{n=1}^\infty (2i\tau)^{-n} H^{\AKNS}_n(q,r)\, &  \sim  \sum_{n=1}^\infty   (2i\tau)^{-n} ( 1+\frac{2\xi}{2i\tau})^{-n} H^{\AKNS}_n(e^{-2i\xi x} q, e^{2i\xi x} r ) 
\\ & \sim \sum_{n=1}^\infty \sum_{m=0}^\infty (-2\xi)^m  (2i\tau)^{-n-m} H^{\AKNS}_n(e^{-2i\xi x} q,e^{2i\xi x} r)
\end{split} 
\] 
which implies 
\[   H^{\AKNS}_N ( e^{2i\xi x} q , e^{-2i\xi x } r ) 
=   \sum_{n=0}^N   (-2\xi)^{N-n}      H^{\AKNS}_n (q,r) 
\]
Together 
\[  \T_N^{\AKNS}( i\tau + \xi,q,r)  
= (\tau-i \xi)^N  \T_N^{\AKNS}( i, e^{2\xi x/\tau} q(x/\tau) /\tau, e^{-2\xi x/ \tau} r(x/\tau )/\tau)     \]
which  reduces estimate \eqref{eq:TNAKNS} to the case $z=i$ under  the assumptions $\Vert q \Vert_{L^2} + \Vert r \Vert_{L^2} \ll 1$,
\[ |\T^{\AKNS}_N(i) | \le   c_N  (\Vert q^{(\frac{N-1}2)} \Vert^2_{L^2} + \Vert r^{(\frac{N-1}2}  \Vert_{L^2}^{2} + \Vert q \Vert^{2N}_{L^2} + \Vert r \Vert^{2N}_{L^2}  ).     \]
 The  proof if this estimate is analogous 
 to Section \ref{sec:differenceflow} 
 but simpler and we omit it.  
It also implies estimates \eqref{eq:hn} and \eqref{eq:difftnakns}.

\end{proof}

With these definitions and Theorem \ref{thm:asymptoticexpansion}, \eqref{eq:NLSODE} becomes the equivalent of the Lenard recursion 
\begin{equation}\label{eq:AKNSrecursion}
    \begin{split}
        \gamma_n'\, &  = 2 ( q \beta_n + r \alpha_n),\\
        \alpha_{n+1}\, &  = i\alpha_{n}' - iq\gamma_{n}, \\
        \beta_{n+1} \, & = - i\beta_{n}' + ir \gamma_n
    \end{split} 
\end{equation} 
with $\gamma_0 = 1, \alpha_0 = \beta_0 = 0$. Using the alternative equation for $\gamma$ \eqref{alternativegamma}, we also find
\begin{equation}
    2\gamma_n = \sum_{k=1}^{n-1} 4\alpha_{k}\beta_{n-k} - \gamma_{k}\gamma_{n-k}.
\end{equation}
Note that since there are no anti-derivatives involved in the recursion, we can directly conclude that $H^{\AKNS}_n$ is an integral over a differential polynomial in $q$ and $r$. The first Hamiltonians and iterates of $\alpha, \beta$ and $\gamma$ can be found in Appendix \ref{iterates}.

The next theorem shows that the $N$th equation in the AKNS hierarchy $q_t = \{q,H_N\}, r_t = \{r,H_N\}$ takes the simple form
\begin{equation}\label{eq:AKNSalphabeta}
    q_t = \alpha_N, \qquad r_t = \beta_N.
\end{equation}

\begin{theorem}\label{thm:formofakns}
    The Hamiltonians $H_N^{\AKNS}$ Poisson commute with $ \log T(z)$. Any two Hamiltonians Poisson commute. They are given by
    \begin{equation}\label{eq:hamiltonian}   H_N^{\AKNS} =  \frac1{2N}  \int \gamma_{N+1}  dx, \quad \frac{\delta}{\delta q} H_N = -i\beta_N,  \quad \frac{\delta}{\delta r }H_N = i\alpha_N. \end{equation} 
    Moreover,
    \begin{equation}
        \label{poissonHn} \{ q, H_N^{\AKNS} \} = \alpha_N  ,  \quad \{ r, H_N^{\AKNS} \} = \beta_N,
    \end{equation}
    
    \begin{equation}
        \label{poissonalpha} 
        \begin{split} \{ \alpha(z), H_N^{\AKNS} \}\, &  =  -i\sum_{l+j = N-1} (2z)^{l}\Big[ \alpha(z)\gamma_j - \gamma(z) \alpha_j \Big],
    \\ \{ \beta(z), H_N^{\AKNS} \}\, &  =  -i\sum_{l+j=N-1}  (2z)^l \Big[  -\beta(z) \gamma_j + \gamma(z) \beta_j \Big],
    \\  \{\gamma(z), H_N^{\AKNS} \}\, &  =  2\sum_{l+j=N-2} (l+1) (2z)^{l}\partial_x  \left[  \alpha(z) \beta_{j} + \beta(z) \alpha_{j} - \frac12 \gamma(z) \gamma_{j} \right],
    \\ \{\gamma_k, H_N^{\AKNS}\}\, &  = 2\sum_{m+j = N+k-2} (m-k+1)\partial_x ( \alpha_{m} \beta_j + \alpha_j  \beta_{m} - \frac12 \gamma_m \gamma_{j} ).
    \end{split} 
    \end{equation}
\end{theorem} 

\begin{proof}The statements follow from Theorem \ref{thm:poissonbrackets} and by making the asymptotic expansions using Theorem \eqref{thm:asymptoticexpansion}.
\end{proof}

\subsection{Nonvanishing limits}\label{sec:nonvanishing}
In the previous part we assume that $ q,r \in L^2$. Surprisingly many results carry over for $L(z,a + q,b+r)$ 
with $a,b \in \C$ and $q,r \in L^2$ instead of $(q,r)$ as arguments.  The reference operator is now 
\begin{equation}\label{eq:L0} L_0 =\left( \begin{matrix} -iz -\partial & a \\ -b & -iz+\partial \end{matrix} \right).  \end{equation}
We write the equation $L\psi=0$ as 
\begin{equation}\label{eq:diff}   \psi' = \left( \begin{matrix} -iz & q + a \\ r+b & iz \end{matrix} \right) \psi \end{equation}  
The characteristic exponents are the roots of $\lambda^2 +z^2 - ab = 0$,  $ \lambda = \pm iz \sqrt{1-ab/z^2}$ which motivates the definition  $  \zeta= -z \sqrt{1-ab/z^2}$ . 
It cannot be purely imaginary  if $|\im z|\ge c(a,b) $.
We choose a basis of the eigenspaces as columns of
\[ U= \left( \begin{matrix} z+\zeta  &-ia  \\ ib &    z+\zeta  \end{matrix} \right), \quad \det U = 2z (z+\zeta) , \quad U^{-1} = \frac{1}{2z(z+\zeta)} \left( \begin{matrix} z+\zeta & ia \\ -ib & z+\zeta \end{matrix} \right). \] 
The Ansatz $\psi = U \phi $ gives
\begin{equation}\label{eq:diff2}   \phi' = \left( \begin{matrix} -i \zeta & 0 \\ 0 & i \zeta \end{matrix} \right) \phi 
+ \frac1{2z} \left( \begin{matrix} i(bq+ar) & (z+\zeta) q + \frac{a^2}{z+\zeta}r \\ (z+\zeta) r+ \frac{b^2}{z+\zeta}q & -i(bq+ar) \end{matrix} \right) \phi. 
\end{equation} 
Similar to the construction in Section \ref{sec:mmiura}, there exist unique Jost solutions to \eqref{eq:diff} provided $ |\im z| $ is sufficiently large and $ r, q \in L^1$, normalized by 
\[ \lim_{x\to -\infty} e^{i\zeta x} \psi_l(x) = \left( \begin{matrix} z+\zeta \\ ib \end{matrix} \right), \quad \lim_{x\to \infty} e^{-i\zeta x} \psi_r(x) = \left( \begin{matrix} -ia \\ z+\zeta \end{matrix} \right) \]
or, equivalently, if $ \phi_l$ satisfies \eqref{eq:diff2} with the 'initial' condition$ \lim\limits_{x\to -\infty} e^{izx} \phi_l = \left( \begin{matrix} 1 \\ 0 \end{matrix} \right)$ resp.  $ \lim_{x\to \infty} e^{-izx} \phi_r= \left( \begin{matrix} 0 \\ 1 \end{matrix} \right) $ this allows to define the transmission coefficient as in the fast decaying case,
\[ T^{-1}(\zeta) := \lim_{x\to \infty} e^{i\zeta x} \phi_l^1(x) = \lim_{x\to -\infty} e^{-i\zeta x}\phi_r^2(x) =  W( \phi_l, \phi_r)= \frac1{2z(z+\zeta)} W( \psi_l, \psi_r) . \] 
We renormalize the transmission coefficient in order to be able to define it for $r,q \in L^2$, without the integrability condition and observe that 
\[  \rho_l = \exp i\Big(\zeta x - \frac1{2z}\int_0^x bq+ar dy \Big)  \phi_l 
\]
satisfies
\begin{equation}\label{eq:diff3}   \rho' = \left( \begin{matrix} 0 & 0 \\ 0 & 2i \zeta \end{matrix} \right) \rho
+ \frac1{2z} \left( \begin{matrix} 0 & (z+\zeta) q + \frac{a^2}{z+\zeta}r \\ (z+\zeta) r+ \frac{b^2}{z+\zeta}q & -2i(bq+ar) \end{matrix} \right) \rho. 
\end{equation} 
We again normalize 
\[ \lim_{x\to -\infty} \rho_l = \left( \begin{matrix} 1 \\ 0 \end{matrix} \right) \] 
and define  the renormalized transmission coefficient $ T_r(z)$
\[ (T_r )^{-1} = T^{-1} \exp\Big( - \frac1{2z}\int bq + ar dy \Big) 
=\lim_{x\to \infty} \rho_l^1. \]
The quantity on the right hand side is now defined for $r,q \in L^2$. 

The resolvent is defined for $ |\im z | $ large and we obtain the same relation between resolvent and the logarithm of the transmission coefficient. The effect of the renormalization is transparent: 
\[ \tilde \alpha = \alpha - \frac{b}{2z},\quad \tilde \beta = \beta - \frac{a}{2z}.\] 
Exactly as in the decaying case (we diagonalize $L_0$) 
\[ {\det}_2 \left( 1+ L_0^{-1} \left( \begin{matrix} 0 & q \\ -r & 0 \end{matrix} \right) \right) =  T_r^{-1}, \]
where $L_0$ is given by \eqref{eq:L0}. Indeed, a close inspection of the argument used in Lemma \ref{lemma:NLSfunctionalderivative} shows that it relies on the decomposition $L_0^{-1} L = 1 + L_0^{-1}(L-L_0)$ rather than on the form of $L_0$.

In particular we obtain the recursion formulas and  the calculation of Poisson brackets carries over to this situation.

\subsection{Functional Derivatives of the Transmission Coefficient}\label{sct:functionalderivatives}
This section contains the proof of Lemma \ref{lemma:NLSfunctionalderivative}.

\begin{proof}[Proof of Lemma \ref{lemma:NLSfunctionalderivative}.]

The equation 
\[ L(z) \psi = f \] 
has a forward fundamental solution  $G(x,y;z)$  given by  $0$ if $ x<y$ and otherwise (observe that $ W( \psi_{++}, \psi_{+-}) = 1 $)
\[   
-  \left( \begin{matrix}  \psi^1_{++}(x)\psi^2_{+-}(y) - \psi^1_{+-}(x)\psi^2_{++}(y)    & \psi^1_{++}(x)\psi_{+-}^1(y) - \psi^1_{+-}(x)\psi^1_{++}(y)   \\
\psi^2_{++}(x)\psi^2_{+-}(y) - \psi^2_{+-}(x)\psi^2_{++}(y)  & \psi^2_{++} (x) \psi^1_{+-}(y) - \psi^2_{+-}(x) \psi^1_{++}(y)   \end{matrix} \right).
\]
This can be seen by checking that $L_xG(x,y;z) = 0$ away from the diagonal and by the jump condition
\begin{equation*}
    G(x^+,x) - G(x^-,x) = \left(\begin{matrix} -1 & 0\\0 & 1\end{matrix}\right),
\end{equation*}
as this implies $L_x G(x,y) = \delta(x-y)$. To determine $ \frac{\delta T}{\delta q} $ and $ \frac{\delta T}{\delta r} $ resp. for $ \im z >0  $, recall

\[ \frac{d}{dt} T(z; q+ t \dot q ,r+ t\dot r)|_{t=0}  = : \int \frac{\delta T}{\delta q} \dot q + \frac{\delta T }{\delta r } \dot r dy \] 
We differentiate the equation with respect to $t$ (dots are $t$-derivatives) and consider $\dot \psi = \dot \psi_{-+}$
\[ L(z) \dot \psi = \left( \begin{matrix}- \dot q \psi^2_{-+} \\  \dot r \psi^1_{-+}\end{matrix} \right).  \] 
Hence 
\[
\begin{split} 
\dot \psi (x)\, &  =  \psi_{++}(x)  \int_{-\infty}^x \psi_{+-}^2(y) \psi_{-+}^2(y) \dot q(y) 
- \psi_{+-}^1(y) \psi_{-+}^1(y) \dot r(y) dy \\ & \qquad 
- 
\psi_{+-}(x)   \int_{-\infty}^x
 \psi^2_{++}(y)  \psi^2_{-+}(y) \dot q(y)
 - \psi^1_{+-}(y) \psi^1_{-+} \dot r dy
\end{split} 
\]
and 
\[ \frac{d}{dt} T^{-1} ( q+ t \dot q , r+\dot r) = \lim_{x\to \infty} e^{izx} \dot \psi^1(x) =  \int  \psi_{+-}^2 \psi_{-+}^2  \dot q dy - \int \psi_{+-}^1 \psi_{-+}^1 \dot r dy \] 
Here, the second summand vanishes due to the assumption $\im z > 0$. Thus 
\[ \frac{\delta T^{-1}}{\delta q} =   
 \psi_{+-}^2 \psi_{-+}^2, \quad 
\frac{\delta T^{-1}}{\delta r} = -
 \psi_{+-}^1 \psi_{-+}^1.    
\] 
We turn to the derivative in the spectral parameter. Define $\tilde \psi = e^{izx}\psi$. We are interested in $T^{-1} = \lim_{x \to \infty} \tilde \psi_1(x)$. We calculate
\[\left(\begin{matrix} -\partial & q\\ -r & -2iz + \partial\end{matrix}\right)\tilde \psi = 0.\]
Thus
\[\left(\begin{matrix} -\partial & q\\ -r & -2iz + \partial\end{matrix}\right)\dot{\tilde \psi} = \left(\begin{matrix}0 \\ 2i \tilde \psi_2 \end{matrix}\right).\]
Then $\Psi = e^{izx}\dot{ \tilde \psi}$ solves
\[\left(\begin{matrix} -iz-\partial & q\\ -r & -iz + \partial\end{matrix}\right)\Psi = \left(\begin{matrix}0 \\ 2i \psi^2_{-+}  \end{matrix}\right).\]
Hence
\begin{equation} \label{dzDT1} 
\begin{split} \frac{d T^{-1}}{dz}\, &  = \lim_{x\to \infty}
e^{izx}  \Psi^1(x) = -2i \int \psi_{+-}^1 \psi_{-+}^2  dy 
\\ & =  -  i \int \psi_{+-}^2 \psi_{-+}^1 + \psi_{+-}^1 \psi_{-+}^2 - T^{-1} dy 
\end{split} 
\end{equation} 
If $\im z <0 $  we use the backward fundamental solution $G(x,y;z)$ which is $ 0 $ for $ x>y$ and otherwise 
\[   
\left( \begin{matrix}  \psi^1_{-+}(x)\psi^2_{--}(y) - \psi^1_{--}(x)\psi^2_{-+}(y)    & \psi^1_{-+}(x)\psi_{--}^1(y) - \psi^1_{--}(x)\psi^1_{-+}(y)   \\
\psi^2_{-+}(x)\psi^2_{--}(y) - \psi^2_{--}(x)\psi^2_{-+}(y)  & \psi^2_{-+} (x) \psi^1_{--}(y) - \psi^2_{--}(x) \psi^1_{-+}(y)   \end{matrix} \right) 
\] 
We differentiate again the equation and consider 
\[ L(z) \dot \psi = \left( \begin{matrix} -\dot q \psi^2_{++} \\  \dot r \psi^1_{++}\end{matrix} \right).  \]
Hence 
\[
\begin{split} 
\dot \psi (x)\, &  =  -\psi_{-+}(x)  \int_x^\infty  \psi_{--}^2(y) \psi_{++}^2(y) \dot q(y) 
- \psi_{--}^1(y) \psi_{++}^1(y) \dot r(y) dy \\ & \qquad 
+ 
\psi_{--}(x)   \int^\infty_x
 \psi^2_{-+}(y)  \psi^2_{++}(y) \dot q(y)
 + \psi^1_{-+}(y) \psi^1_{++}(y) \dot r(y) dy
\end{split} 
\] 
We obtain 
\[ \frac{d}{dt} T(z,q+t\dot q, r+t\dot r) |_{t=0} = \lim_{x\to -\infty} e^{iz x} \dot \psi(x)  = \int -\psi_{--}^2 \psi^2_{++} \dot q + \psi_{--}^1 \psi_{++}^1 \dot r dy \] 
and 
\[ \frac{\delta}{\delta q}T(z,q,r)  =
-\psi_{--}^2 \psi_{++}^2, \qquad \frac{\delta}{\delta r} T(z,q,r) = \psi_{--}^1 \psi_{++}^1. \]
We observe that 
\[ G(z,q,r) L(z,\tilde q, \tilde r) = 1 +   G(z,q,r) \left( \begin{matrix} 0 & \tilde q-q \\ -(\tilde r-r) & 0 \end{matrix} \right) \] 
where the second summand on the right hand side is a trace class operator if $p,q \in L^1$. We calculate for $\im z>0$
\[ 
\begin{split} 
\frac{d}{dt} \det\Big( G(z,q,r) L(z, q+ t\dot q, r+t\dot r\Big)\Big|_{t=0} & =  \trace \left( G \left( \begin{matrix} 0 & \dot q \\ -\dot r & 0 \end{matrix}  \right)\right) \\ & = \int \trace_{\R^2}\left[   \left(\begin{matrix}  \gamma/2  & \alpha \\ \beta & \gamma/2 \end{matrix} \right) \left( \begin{matrix} 0 & \dot q \\ -\dot r & 0 \end{matrix} \right) \right] dx 
\\ &= \int  \beta \dot q - \alpha \dot r  dx 
. 
\end{split} 
\]
The left hand side can be rewritten as 
\[ \det ( G(z ,q,r) L(z,0,0)  ) \frac{d}{dt} \det  G(z,0,0) L(z,q+t\dot q, r+ t \dot r)  \Big|_{t=0}
= \int \beta \dot q - \alpha \dot r 
\]  
and we arrive at 
\[ \frac{d}{dt} \ln \det(G(z,0,0) L(z,q+t \dot  q, r+ t \dot  r) = 
\int \beta \dot q - \alpha \dot r dx.  \] 
Since $ \det( G(z,0,0) L(z,0,0) ) =1$ we see that 
\begin{equation}  \det\left( 1+ G(z,0,0) \left( \begin{matrix} 0 & q \\ -r & 0 \end{matrix} \right)\right)  = T(z,q,r) \end{equation} 
on the upper half plane. 

In the same fashion as above 
\begin{equation} \label{dzDT2} 
\begin{split} \frac{d T}{dz}\, &   = 2i \int \psi_{--}^1 \psi_{++}^2 dy 
\\ & =    i \int \psi_{++}^2 \psi_{--}^1 + \psi_{++}^1 \psi_{--}^2  -  T  dy.
\end{split} 
\end{equation} 
We arrive at 
\[ \frac{\delta \ln T}{\delta q} =-  \left\{ \begin{array}{rl} T \psi_{-+}^2 \psi_{+-}^2 & \quad \text{ if } \im z > 0  \\ T^{-1} \psi_{--}^2 \psi_{++}^2 & \quad \text{ if } \im z < 0 \end{array} \right.  \] 
\[ \frac{\delta \ln T}{\delta r} = \left\{ \begin{array}{rl} T\psi_{-+}^1 \psi_{+-}^1 & \quad \text{ if } \im z > 0  \\ T^{-1} \psi_{--}^1 \psi_{++}^1 & \quad \text{ if } \im z < 0 \end{array} \right.  \] 
\[ \frac{d}{dz} \ln T = \left\{ \begin{array} {rl}\displaystyle  i  \int T  (\psi_{-+}^1 \psi_{+-}^2 + \psi_{-+}^2 \psi_{+-}^1)-1 dx     & \qquad \text{ if } \im z > 0 \\[4mm]
\displaystyle i  \int T^{-1}  (\psi_{--}^1 \psi_{++}^2 + \psi_{--}^2 \psi_{++}^1)-1 dx 
& \qquad \text{ if } \im z < 0 ,
\end{array} \right. \]
which finishes the proof.
\end{proof}

\subsection{Embedding other hierarchies into the AKNS hierarchy}
\label{sec:embedding}

In this additional section we will show that the AKNS hierarchy contains the NLS $(r= \pm \bar q)$, and as a part of it, the real mKdV hierarchy. There is also the complex KdV hierarchy which is obtained by setting $r=1$, and the Gardner hierarchy related to  the Wadati Lax operator (see Subsection \ref{subsec:determinants}. In all cases we specialize the transmission coefficient and their variational derivatives, study structural  properties like symmetry and relations between transmission coefficients for different Lax operators - the most important being the connection between the Wadati operator and KdV via the modified Miura map. Then we deduce real symplectic and Poisson structures from the complex Poisson structure for AKNS and relate that in many cases to Gardner Poisson structures. In particular we find three Hamiltonian structures for KdV: The Gardner structure, the Magri structure, and the interpretation as AKNS Hamiltonian vector field restricted to a subset of functions.
 
In fact, we could define the Gardner hierarchy resp. the generating function $\T_{-1}^{\Gardner}$ of the Gardner Hamiltonians and a study of their structure merely from our knowledge on AKNS. In the end we took a shortcut to find the Gardner generating functions by studying the good variables in more detail (see Lemma \ref{lem:poissonbrackets}), but we decided to keep the AKNS approach since it helps in understanding the various connections between classical integrable hierarchies.

%Interestingly the fundamental conserved quantity $\Vert w \Vert_{L^2}^2$ is independent of the coefficients of the asymptotic series for the generating function which is a surprising fact. This connection can be extended to all of the functions $\T_N^{\Gardner}$ (see part E) in Theorem \ref{thm:gardnermiura2}.

\subsection{Complex KdV hierarchy} \label{section:kdvembedded}
In Sections \ref{sec:proof} and \ref{sec:AKNS}, we constructed the KdV hierarchy by means of its transmission coefficient $T^{\KdV}$ in the upper half-plane. Now we consider the AKNS Lax operator 
\[ L(z,u,1) = \left( \begin{matrix} -iz-\partial & u \\ -1 & - iz +\partial \end{matrix} \right), \] 
the Hamiltonians and the diagonal Green's function evaluated at $q=u$ and $r=1$. Note that the Lax operator is holomorphic in $u$. Taking  a fresh start we now define 
\[ T^{\KdV}(z,u) = T^{\AKNS}(z, u,1). \]  
An intriguing consequence of the proposition below is that equations of the KdV hierarchy can be understood as three different Hamiltonian evolutions
\begin{enumerate} 
\item As Hamiltonian equations with the Hamiltonian  $H^{\KdV}_n$ with respect to the Gardner Poisson bracket. 
\item As Hamiltonian equations with the Hamiltonian $H^{\KdV}_n$ with respect to the Magri Poisson bracket.
\item As restriction of the Hamiltonian equations with the Hamiltonian $H^{\AKNS}_{2n+2}$ with respect to the AKNS symplectic structure, restricted to the set $(u,1)$.
\end{enumerate} 

\begin{proposition} \label{prop:complexkdv} 
\noindent {\bf A) Recursion relations.}
The functions $\log T^{\KdV} $ and  $\beta(z,u,1)$ are odd in $z$.  We define $H^{\KdV}_n(u) = \frac12 H^{\AKNS}_{2n+3}(u,1)$.
The recursion relations can be written as  
\begin{equation} 
\label{eq:lenardodeakns}\begin{split}   \gamma\, &  = \beta' - 2iz\beta, \qquad  \alpha  = \frac12 ( \beta'' -2iz \beta' - 2u \beta), \\
    0\, & = \beta''' + 4 (z^2-u)  \beta'  - 2u' \beta, \quad \beta_{2n-1}''' - 4 u \beta_{2n-1}' -2u'\beta_{2n-1} =-\beta_{2n+1}.
    \end{split} 
\end{equation} 
where $\beta$,$\beta_n$ and their derivatives are evaluated at $(z,u,1)$.\hfill\ \linebreak
\noindent{\bf B) Poisson brackets.}  Let $z_1,z_2 \in \C\backslash \R $ and $n,m \in \N$. Then 
 $ \log T^{\KdV}(z_1)$ and  $\log T^{\KdV}(z_2)$ all  Poisson commute with respect to the Gardner structure, whenever $T^{\KdV}$ is defined (i.e. when $z_{1,2}$ is outside the spectrum). As a consequence 
$H_n^{\KdV}$, $H_m^{\KdV}$ Poisson commute on sufficiently regular functions.
The Gardner Poisson brackets of $\beta$  with $ \log T^{\KdV} $ and $H^{\KdV}_{n}$ satisfy 
\begin{equation}\label{eq:kdvpoissonbrackets} \begin{split} 
   \Big\{ \frac{1}{\beta(z_1)} ,\log T^{\KdV}(z_2) \Big\}_{\Gardner}& = \frac{\displaystyle \partial\Big[ \frac{1}{\beta(z_1)}\frac{\delta}{\delta u}\log T(z_2)\Big]}{(2z_1)^2-(2z_2)^2}
    \  ,\\
    \Big\{\frac{1}{\beta(z)}, H^{\KdV}_n\Big\}_{\Gardner} &= -\partial\Big[\frac1{\beta(z)} \frac{\delta}{\delta u}\sum_{m=-1}^{n-1}  (2z)^{2(n-1-m)}H_{m}^{\KdV}\Big]
\end{split}
\end{equation}
and with 
\[ 
    \mathcal{T}_N^{\KdV} = \sum_{n=-1}^N H_n^{\KdV}  (2z)^{2N-2n}  +\frac{(2z)^{2N+3}}{2i}  \log T^{\KdV},
\]
we have (again with the Gardner bracket) for $N \geq 0$,
\begin{equation}\label{eq:kdvpoissonbrackets2}
\Big\{ \frac{1}{\beta(z_1)}, \mathcal{T}^{\KdV}_N(z_2)\Big\}   = \frac{\displaystyle \partial \Big[ \frac1{\beta} \frac{\delta}{\delta u} \Big((2z_2)^{2}   \mathcal{T}^{\KdV}_{N-1}(z_2) + \sum_{m=-1}^{N-1}(2z_1)^{2(N-m)}  H^{\KdV}_m \Big)
\Big]}{(2z_1)^2-(2z_2)^2}.
\end{equation}

\noindent{\bf C) Hamiltonian structures.} 
The three Hamiltonian structures (KdV with Gardner Poisson structure, KdV with Magri Poisson structure and the restriction of the even AKNS Hamiltonian vector field) are 
%relation between the KdV flows and the even AKNS flows is 
expressed by
\begin{equation} \label{eq:relkdvakns} 
\begin{split} 
 2iz \partial_x \frac{\delta}{\delta u} \log T^{\KdV} 
\, & = 2iz \partial_x \frac{\delta}{\delta q} \log T^{\AKNS}(z)\Big|_{(q,r)=(u,1)} 
\\ & =\frac{\delta}{\delta r} \Big( \log T^{\AKNS}(z) + \log T^{\AKNS}(-z)\Big)\Big|_{(q,r)=(u,1)}
\\ & = (2iz)^{-1} \Big(\partial^{(3)} - 2(u \partial +\partial u) \Big)         \frac{\delta}{\delta u} \log T^{\KdV} 
\end{split} 
\end{equation} 
which implies the identification of the Gardner, Magri and AKNS Hamiltonian vector fields  
\begin{equation}\label{eq:threehamiltonian} 
\partial \frac{\delta}{\delta u} H^{\KdV}_{n} = -i\frac{\delta}{\delta r} H^{\AKNS}_{2(n+1)}\Big|_{(q,r)= (u,1)}
\\  = -\Big( \partial^{(3)}-2(u \partial+\partial u) \Big) \frac{\delta}{\delta u} H^{\KdV}_{n-1}.
\end{equation} 
Finally
\begin{equation} \label{eq:stationarity} \frac{\delta}{\delta q} ( \ln T^{\AKNS}(z) + \ln T^{\AKNS}(-z) )\Big|_{(u,1)}=0,  \qquad  \frac{\delta}{\delta q} H^{\AKNS}_{2n}\Big|_{(u,1)}=0  . \end{equation}

%\alpha(z,u,1) + \alpha(-z,u,1) = -2iz \beta'(z, u,1). \end{equation}
%Equivalently, we have with $G(u) = F(u,1) $
%\begin{equation}\label{eq:kdvpoissonbracket}
%    \{ F, \ln  T^{\AKNS}(z) + \ln T^{\AKNS}(-z) \} \big|_{(u,1)} = 2z %\{ G, \ln T^{\KdV} \}_{\Gardner}, 
%\end{equation}
%hence  set $\{ (u,1)\} $ is invariant under the flow of $H_{2(n+1)}^{\AKNS}$ and the flow on it coincides with Hamiltonian flow of $H_{n}^{\KdV}$ with respect to the Gardner Poisson bracket. 
\end{proposition}

\begin{proof} 
A) Let $ \phi_l$ be a left Jost function for $L(z,u,1)$. Then $\left( \begin{matrix} \phi_{l,1}+2iz \phi_{l,2} \\ \phi_{l,2} \end{matrix} \right)$ is a left Jost function for $L(-z,u,1)$. This implies that $ \log T^{\KdV}$ is an odd function of $z$, and the same is true for 
$\beta = - \frac{\delta}{\delta u} \log T^{\KdV}$.
Similarly $ \phi_2$ satisfies $ -\phi_2''+ u \phi = z^2 \phi  $ which shows that the definition of $T^{\KdV}$ is consistent with the first subsection. We solve the third equation in \eqref{eq:NLSODE} for $\gamma$, substitute $\gamma$ in the first equation and solve for $ \alpha$ to obtain the first two identities of \eqref{eq:lenardodeakns}. We substitute $\gamma$ and $ \alpha$ in the second line in \eqref{eq:NLSODE} to obtain the recursion equation for $ \beta$. The  asymptotic expansion $ \beta \sim \sum \beta_n (2z)^{-n}$ implies the last identity in \eqref{eq:lenardodeakns}.

B) From AKNS Poisson commutativity and \eqref{eq:lenardodeakns},
\begin{align*}
    0 &= \int \frac{\delta \log T^{\AKNS} (z_1)}{\delta r}\frac{\delta \log T^{\AKNS} (z_2)}{\delta q} - \frac{\delta \log T^{\AKNS}(z_1)}{\delta q}\frac{\delta \log T^{\AKNS}(z_2)}{\delta r} \, dx\Big|_{(u,1)}  \\
    & = - \int  \alpha(z_1) \beta(z_2) - \alpha(z_2) \beta(z_1) \, dx \\
    & = \frac{1}{2}\int 2i z_1 \beta'(z_1) \beta(z_2) - 2i z_2 \beta'(z_2) \beta(z_1)\, dx \\ 
    & = i \int (z_1+z_2) \beta'(z_1) \beta(z_2) \, dx \\
    &= -i(z_1 + z_2) \int \frac{\delta}{\delta u}\log T^{\AKNS}(z_1,u,1) \partial_x \frac{\delta}{\delta u}\log T^{\AKNS}(z_2,u,1) \, dx\Big|_{(u,1)} \\ 
    & = -i (z_1+z_2) \{ \log T^{\KdV}(z_1), \log T^{\KdV}(z_2) \}_{\Gardner}. 
\end{align*}
Poisson commutation of $\log T^{\KdV}(z_1)$ and $ \log T^{\KdV}(z_2)$ implies Poisson commutation of $\log T^{\KdV}(z)$, $H_n^{\KdV}$ and $H_m^{\KdV}$.

We turn to \eqref{eq:kdvpoissonbrackets}. We specialize \eqref{poissongamma} to deduce for the AKNS Poisson brackets,
\[
\Big\{ \beta(z_1),\log T^{\AKNS}(z_2)\Big\}\Big|_{(u,1)}\, =  \frac{\beta'(z_1) \beta(z_2) - \beta'(z_2) \beta(z_1)}{2(z_1-z_2)} - i\beta(z_1)\beta(z_2),
\]
and a similar formula for $\log T^{\AKNS}(-z_2)$. Summing these two, the term $i\beta(z_1)\beta(z_2)$ drops out since $\beta$ is odd in $z_2$. 
Now we claim that for $G(u) = F(u,1)$ where $F = F(q,r)$, we have
\[
    \{F, \log T(z) + \log T(-z)\}|_{(u,1)} = 2z \{G, \log T^{\KdV}(z)\}_{\Gardner}
\]
Indeed, let $A(z) = \log T(z) + \log T(-z)$. We will later see that by \eqref{eq:relkdvakns},
\[2iz \partial \frac{\delta}{\delta u} \log T^{\KdV} = \frac{\delta}{\delta r} A(z)|_{(u,1)},\]
and by \eqref{eq:stationarity} $\frac{\delta}{\delta q} A(z)|_{(u,1)} = 0$. Thus
\[ 
\begin{split}
    \{F, \log T(z) + \log T(-z)\}|_{(u,1)} &= \frac{1}{i} \int \frac{\delta F}{\delta q} \frac{\delta A}{\delta r} - \frac{\delta F}{\delta r} \frac{\delta A}{\delta q} \, dx |_{(u,1)}\\
    &= \frac{1}{i} \int \frac{\delta G}{\delta u} 2iz \partial \frac{\delta}{\delta u} \log T^{\KdV}(z) \, dx\\
    &= 2z \{G, \log T^{\KdV}\}_{\Gardner}.
\end{split}
\]
Applying this to $\beta(z_1) = F$ shows
\[
    \{\beta(z_1),\log T^{\KdV}(z_2)\}_{\Gardner} = \frac{2}{4z_1^2-4z_2^2}(\beta'(z_1) \beta(z_2) - \beta'(z_2) \beta(z_1)),
\]
and we arrive at the first part of \eqref{eq:kdvpoissonbrackets} by using
\[
    \big\{\frac{1}{\beta},\log T\big\}_{\Gardner} = \frac{-1}{\beta^2}\{\beta,\log T\}_{\Gardner}.
\]
The second part of \eqref{eq:kdvpoissonbrackets} follows by using the asymptotic expansion. We turn to the proof of \eqref{eq:kdvpoissonbrackets2} which we prove inductively. For $N=0$, the claim holds as can be seen using \eqref{eq:kdvpoissonbrackets} directly. Now since
\[
    \T_{N+1}(z_2) = -(2z_2)^2\T_N(z_2) + H_{N+1},
\]
we have
\begin{align*}
    \{\beta^{-1},\T_{N+1}(z_2)\} &= \frac{-(2z_2)^2\partial \Big[ \frac1{\beta} \frac{\delta}{\delta u} \Big((2z_2)^{2}\T^{\KdV}_{N-1} + \sum_{m=-1}^{N-1}(2z_1)^{2(N-m)}  H^{\KdV}_m \Big)
    \Big]}{(2z_1)^2-(2z_2)^2} \\
    &\quad - \partial \big(\beta^{-1} \sum_{m=-1}^N \frac{\delta}{\delta u} (2z_1)^{2(N-m)}H_m\big)\\
    &= \frac{\partial \Big[ \frac1{\beta} \frac{\delta}{\delta u} \Big((2z_2)^2\T^{\KdV}_{N} - (2z_2)^2\sum_{m=-1}^{N}(2z_1)^{2(N-m)}  H^{\KdV}_m \Big)
    \Big]}{(2z_1)^2-(2z_2)^2} \\
    &\quad + \frac{\partial \Big[\beta^{-1} (-(2z_1)^2+(2z_2)^2)\sum_{m=-1}^N \frac{\delta}{\delta u} (2z_1)^{2(N-m)}H_m\Big]}{(2z_1)^2-(2z_2)^2}\\
    &= \frac{\partial \Big[\beta^{-1}\big((2z_2)^2\T^{\KdV}_{N} -(2z_1)^2\sum_{m=-1}^N \frac{\delta}{\delta u} (2z_1)^{2(N-m)}H_m\big)\Big]}{(2z_1)^2-(2z_2)^2},
\end{align*}
which is the right-hand side of \eqref{eq:kdvpoissonbrackets2} for $N+1$.

C) The first identity of \eqref{eq:relkdvakns} is the definition. The second identity can be equivalently written as 
\[ -2iz \beta' =  \alpha(z) + \alpha(-z) \]  
which follows from the second equation of \eqref{eq:lenardodeakns} and the observation that the first and the last term on the right hand side are odd. The identity of the left hand side and the right hand side in \eqref{eq:relkdvakns} is equivalent to the third line in \eqref{eq:lenardodeakns}. The asymptotic series give \eqref{eq:threehamiltonian}. The last identity \eqref{eq:stationarity} is a direct consequence of the fact that $\beta$ is odd. 
\end{proof}

\subsection{Defocusing NLS hierarchy} 

The defocusing NLS hierarchy contains the (complex) mKdV hierarchy as the even part, and hence also the real mKdV hierarchy. The relations between the hierarchies are interesting in themselves, but we will not use them outside of this subsection. 

The defocusing NLS hierarchy is  the case $ r = \bar q$. We choose the standard real symplectic form 
\[ \omega(f,g) = \im \int f \bar g - g \bar f dx \] 
and define
\[ T^{\NLS} (z,q) = T^{\AKNS}(z, q, \bar q).  \] 

\begin{proposition} 
\noindent{\bf A) The transmission coefficient.} 
The transmission coefficient has the following  properties
\begin{equation} \label{eq:symmetriesdefocNLS}
\begin{split} 
1=\, &  T^{\NLS} (z,q) \overline{T^{\NLS} ( \bar z, q)}, \\
\log T^{\NLS}(z,q) = \, & - \log \overline{T^{\NLS}( \bar z,q)},  \\ 
\gamma(z, q, \bar q)  =\, &  \overline{\gamma( \bar z, q, \bar q )} \\
\beta(z,q,\bar q ) = \, & \overline{\alpha( \bar z, q , \bar q) ) }. 
\end{split} 
\end{equation} 
\noindent{\bf B) The recursion formula.}  The following identities hold
\[ \gamma'(z) = 2( q \overline{ \alpha( \bar z)} + \bar q \alpha(z)), \quad \alpha' = - 2i z \alpha + q \gamma, \]
hence $H_n \in \R$, $ \gamma_n(x)  \in \R$, $ \alpha_n = \bar \beta_n$
and  $ \gamma_n' = 4 \real (q \bar \alpha_n)$ and 
\[ \alpha_{n+1} = i \alpha_n' - i q \gamma_n. \] 
\end{proposition}
As a consequence the Hamiltonian flow of $H_n$ preserves the structure $ r= \bar q $.

\begin{proof} 
Let $ \phi_l(z,q,\bar q) $ be a Jost solution. Then also 
$ \left( \begin{matrix} \bar \phi^2_{l} \\ \bar \phi_l^1   \end{matrix} \right) $ is a Jost solution to the spectral value 
$ \bar z$. This implies the first line and the second line in \eqref{eq:symmetriesdefocNLS}. 
Since 
\[ (L(z,q,\bar q))^* = - L(\bar z,q,\bar q)   \]   
we obtain $ \alpha(z) = \overline{\beta( \bar z)}$ which implies the remaining assertions. 
\end{proof} 
It is remarkable that the fourth equation is the complex mKdV equation 
\[  q_t + q_{xxx}- 6 |q|^2 q_x = 0. \]
We note that unlike in the case of real potentials, we do not expect to be able to write the even flows with respect to the Gardner Poisson structure. This can already be seen from the complex defocusing mKdV flow, since $ |q|^2 q_x$ can in general not be written as a total derivative, and hence not as Hamiltonian equations with respect to the Gardner structure.

\subsection{Defocusing real mKdV} 

The real defocusing mKdV hierarchy is a special degenerate version of the Gardner hierarchy. This is the case $r=q=v$ with real valued functions $v$. 
It is a special case of the NLS hierarchy, note however that real functions are contained in a Lagrangian subspace of the symplectic form, i.e. it vanishes identically on it. The relevant Poisson structure is the Gardner Poisson structure. The recursion relations allow to relate the symplectic structure of the AKNS hierarchy to the Gardner Poisson structure of mKdV (in particular to deduce Poisson commutation of $\log T$ with respect to the Gardner Poisson structure from the Poisson commutation of $\log T$ for the AKNS structure. 
The mKdV hierarchy  is connected to the KdV hierarchy via the Miura map. This is useful, however it is difficult to work with it, since  it is not even a local diffeomorphism from $L^2$ to $H^{-1}$, see \cite{MR2189502}. 
 We define
\[T^{\mKdV}(z,v) = T^{\AKNS}(z,v,v)\qquad ( = T^{\NLS}(z,v)). \]

\begin{proposition} 
\noindent{\bf A) Properties of the transmission coefficient.}
The generating function  $\ln  T^{\mKdV}$ is odd and real on the imaginary axis. Moreover 
\begin{equation}\label{eq:mkdvmiurat}
    T^{\mKdV}(z,v) = T^{\KdV} (z,v_x+ v^2 ).
\end{equation} 
We define $H_n^{\mKdV}(v)= H_{n-1}^{\KdV}(v_x+v^2)$ so that 
\[ \ln T^{\mKdV} \sim -2i\sum_{n=0}^\infty H_n^{\mKdV}  (2z)^{-1-2n}.  \]

\noindent{\bf B) Recursion relations. }
The functions $\alpha, \beta$ and $ \gamma$ are real on the imaginary axis.  $\alpha-\beta$ is odd and $\gamma$ and $ \alpha+ \beta$ are even, and $\alpha(z) = \beta(-z)$. The recursion  formula can be written as 
\begin{equation}\label{eq:mKdVrecursion}  \begin{split} ( \alpha + \beta)''\, &  =  -4z^2 (\alpha+\beta) + 2 (v \gamma)', \\ 
  ( \alpha-\beta)'\, &  = -2iz  (\alpha+ \beta),
  \\ \gamma' \, & = 2v ( \alpha+ \beta) ,
  \\ \alpha_{2n+2}+\beta_{2n+2} \, & = (\alpha_{2n}+\beta_{2n})'' -2(v \gamma_{2n})', 
  \\   \gamma_{2n}'\, &  = 2v (\alpha_{2n}+\beta_{2n}). 
 \end{split} 
\end{equation}
\noindent{\bf C) Poisson brackets.}  Let $ z_1,z_2 \in \C \backslash \R$, $n,m\in \N$. The transmission coefficient 
$T^{\mKdV}(z_1)$, $T^{\mKdV}(z_2)$, $ H^{\mKdV}_n$ and $H^{\mKdV}_m$ all Poisson commute with respect to the Gardner Poisson structure. The Magri structure and Gardner structure are related by 
\begin{equation}\label{eq:mkdvpoisson}   \big\{ f(v_x+v^2), g(v_x+v^2)\big\}_{\Gardner}   =  \big\{f,g\big\}_{\Magri}\big|_{v_x+v^2}.   \end{equation}

\noindent{\bf D) Hamiltonian structures.}
The relation between the mKdV flows, KdV flows  and the even AKNS flow is expressed by
\begin{equation}\label{eq:mKdVrelations}
\begin{split} 
-2iz\Big(\frac{\delta}{\delta r} - \frac{\delta}{\delta q}\Big)  \log T^{\AKNS}\Big|_{(z,v,v)} \, & = \partial \frac{\delta}{\delta v} \log T^{\mKdV}= \partial(-\partial+2v) \frac{\delta}{\delta u} \log T^{\KdV}\Big|_{v_x+v^2}  ,\\
\frac{\delta}{\delta q} H^{\AKNS}_{2(n+1)}|_{(z,v,v)}=      \frac{\delta}{\delta r}  H^{\AKNS}_{2(n+1)} |_{(z,v,v)}   &= \partial \frac{\delta}{\delta v} H_{n}^{\mKdV}(v)
     = \partial (-\partial+2v) \frac{\delta}{\delta u} H_{n-1}^{\KdV}\Big|_{v_x+v^2} .
\end{split} 
\end{equation}
Moreover, we have the following chain rule for the time derivative
\begin{equation}
    (\partial + 2v)\partial \frac{\delta}{\delta v} \log T^{\mKdV} = 4z^2 \partial \frac{\delta}{\delta u} \log T^{\KdV}\big|_{v_x + v^2}
\end{equation}
In particular, $u=v_x+v^2$ solves the $n$th KdV equation whenever $v$ solves the $n$th mKdV equation, and $v$ solves the $n$th mKdV equation if and only if $(v,v)$ satisfies the $2(n+1)$th AKNS Hamiltonian equation.
\end{proposition} 

\begin{proof} 
A) Let $ \phi_l$ be a left Jost function for $L(z,v,v)$. Then $\psi = \phi_{l}^1 + \phi_l^2$ satisfies
\[ (-z^2- \partial^2 + v^2 + v_x ) \psi = 0 \] 
and hence by inverting the reasoning from the KdV case,
\[\left( \begin{matrix} (\phi_{l}^1+\phi_{l}^2)' - iz(\phi_{l}^1+\phi_{l}^2) \\ \phi_l^1 + \phi_l^2 \end{matrix}\right)
\]
is a left Jost function for $L(z,v_x+v^2,1)$.
Thus
\[ T^{\mKdV}(z,v) = T^{\KdV}(z, v_x+v^2). \] 
In particular $\ln T^{\mKdV}$ and $ \alpha-\beta = \frac{\delta}{\delta w}T^{\mKdV} $ are odd in $z$. Specializing the symmetries of defocusing NLS \eqref{eq:symmetriesdefocNLS} to $q = \bar q = v$, we see that $\alpha(i\tau,v,v) = \beta(-i\tau,v,v)$ and $\gamma(i\tau,v,v) = \gamma(-i\tau,v,v)$. Combining these facts we find that $\alpha(i\tau) + \beta(i\tau)$ and $\gamma(i\tau)$ are even. Since $\alpha$ and $\beta$ are holomorphic also  $\alpha(z) + \beta(z)$ and $\gamma(z)$ are even.

B) The recursion relations are an immediate consequence of A) and the recursion relations for the NLS hierarchy. 

C) The Poisson commutativity holds, because
    \begin{align*}
        0 &= \int \frac{\delta \log T^{\AKNS} (z_1)}{\delta q} \frac{\delta \log T^{\AKNS}(z_2)}{\delta r} - \frac{\delta \log T^{\AKNS}(z_1)}{\delta r} \frac{\delta \log T^{\AKNS}(z_2)}{\delta q} \\
        &= -\int \beta(z_1)\alpha(z_2) - \alpha(z_1)\beta(z_2)\\
        &= \frac12 \int -(\alpha+\beta)(z_1)(\alpha-\beta)(z_2)+(\alpha-\beta)(z_1)(\alpha+\beta)(z_2)\\
        &= \Big(\frac1{2iz_1} +\frac{1}{2iz_2}\Big) \int (\alpha-\beta)'(z_1)(\alpha-\beta)(z_2)\\
        &= \frac{i(z_1+z_2)}{2z_1z_2} \int \frac{\delta}{\delta v} T^{\mKdV}(z_1,v) \partial_x \frac{\delta} {\delta v}T^{\mKdV}(z_2,v,v)
        \\ &= \frac{i(z_1+z_2)}{2z_1z_2} \Big\{\log T^{\mKdV}(z_1), \log T^{\mKdV}(z_2) \Big\}_{\Gardner}.
    \end{align*}
Moreover, from the operator identity
\begin{equation}\label{eq:operatorequivalencemkdv}
    (-\partial^3+4(v_x+v^2)\partial + 2(v_x+v^2)_x) = (\partial + 2v)\partial (-\partial + 2v),
\end{equation}
we see
\[  \big\{ f , g \big\}_{\Magri}  = \int (-\partial+2v )  \frac{\delta}{\delta u} f \partial (-\partial + 2 v) \frac{\delta}{\delta u } g dx 
= \big\{ f(v_x+v^2),g(v_x+v^2) \big\}_{\Gardner} .
\]

%The recursion formulas follow by setting $z = i\tau$ in \eqref{eq:NLSODE} and an explicit calculation. The second recursion equation together with the definition of $ \alpha$ and $\beta$ give the equation \eqref{eq:mKdVrelations}. 

%The second equation of \eqref{eq:mKdVGardnerstructure} follows from the definition of $\alpha$ and $\beta$, $\alpha(z) = \beta(-z)$, the  first equation of \eqref{eq:mKdVrelations} and the asymptotic series. 

%The statement for the Hamiltonians \eqref{eq:mkdvmiurah} is seen by going the the asymptotic expansion in \eqref{eq:mkdvmiurat}. 
D) The first part of the first identity in \eqref{eq:mKdVrelations} is just a consequence of the second recursion relation \eqref{eq:mKdVrecursion}. The identity 
\[(\partial + 2v)(\alpha(z,v,v)+\alpha(-z,v,v)) = 2iz\partial \beta(z,u,1),
\]
is an equivalent formulation of the second part of the first identity in \eqref{eq:mKdVrelations}. It is a consequence of the second recursion relation \eqref{eq:mKdVrecursion} to which we 
apply the operator $(\partial+2v)$, use the chain rule and \eqref{eq:operatorequivalencemkdv}.  
The second identity of the first line follows by the chain rule, and the second line spells the conseqences out for the aymptotic series. 
Using the second equality in \eqref{eq:mKdVrelations}, the operator identity \eqref{eq:operatorequivalencemkdv} and the Lenard recursion, the last equality follows.
\end{proof} 

We specialise Theorem \ref{thm:formofakns} to $q = r \in \R$ and obtain the following

\begin{lemma} The mKdV Hamiltonian $H^{\mKdV}_n(v) = \frac{1}{2}H_{2n+1}^{\AKNS}(v,v)$ can be written as an integral over a sum of homogeneous differential polynomials 
\[ H^{\mKdV}_n(v) =  \frac12 \int |v^{(n)}|^2 dx   + \sum_{k=2}^{2n+2}  \int e^{\mKdV}_{n,k} dx \] 
where $e^{\mKdV}_{n,j}$ is a sum of monomials of  degree $2n+2$,  homogeneity $2k$ and weight $2(n+1-k)$, each factor in the monomials having order at most $[(n+1-k)/2] + 1$. Moreover 
\[ e^{\mKdV}_{n,2n+2} =     \frac1{2(2n+1)}   \binom{2n+2}{n+1}   v^{2n+2}. \]
\end{lemma}

\section{Hamiltonians, Equations and Recursions}\label{iterates}

For convenience, we list here some of the equations and Hamiltonians of the hierarchies contained in the AKNS hierarchy. Recall, \eqref{eq:AKNSrecursion},
\begin{equation*}
    \begin{split}
        \gamma_n'\, &  = 2 ( q \beta_n + r \alpha_n),\\
        \alpha_{n+1}\, &  = i\alpha_{n}' - iq\gamma_{n}, \\
        \beta_{n+1} \, & = - i\beta_{n}' + ir \gamma_n
    \end{split} 
\end{equation*} 
with $\gamma_0 = 1, \alpha_0 = \beta_0 = 0$, and the alternative equation for $\gamma$
\begin{equation*}
    2\gamma_n = \sum_{k=1}^{n-1} 4\alpha_{k}\beta_{n-k} - \gamma_{k}\gamma_{n-k}.
\end{equation*}
The first few iterates are
\begin{align*}
    \alpha_0 &= 0, \quad \beta_0 = 0, \quad \gamma_0 = 1,\\
    \alpha_1 &= -iq, \quad \beta_1 = ir, \quad \gamma_1 = 0,\\
    \alpha_2 &= q',\quad \beta_2 = r',\quad \gamma_2 = 2qr,\\
    \alpha_3 &= iq'' -2iq^2r, \quad \beta_3 = -ir'' + 2ir^2 q, \quad \gamma_3 = -2i(qr' - q'r),\\
    \alpha_4 &= -q''' + 6qq'r, \quad \beta_4 = -r''' + 6rr'q, \quad \gamma_4 = -2(qr''+rq''-q'r')+ 6q^2r^2,\\
    \alpha_5 &= -i(q^{(4)} - 8qq''r-6(q')^2r-4qq'r'-2q^2r''+6q^3r^2),\\
    \beta_5 &= i(r^{(4)} - 8rr''q-6(r')^2q-4rr'q'-2r^2q''+6r^3q^2),\\
    \gamma_5 &= -2i(q'''r-r'''q -q''r'+r''q' + 6(-qq'r^2+rr'q^2))\\
    \gamma_6 &= 2[qr^{(4)}+rq^{(4)}-(q'r'''+q'''r')+q''r'']\\ &-10((q')^2r^2+q^2 (r')^2)-20(q^2rr''+r^2qq'')+20q^3r^3
\end{align*}
For the first Hamiltonians we construct the functional antiderivatives by hand and find
\begin{align*}
    H_{1}^{\AKNS} &= \int qr \,dx,\\
    H_2^{\AKNS} &= -\frac{i}{2}\int qr'-q'r \,dx = -i\int qr'\,dx,\\
    H_3^{\AKNS} &= \int q'r'+q^2r^2\,dx,\\
    H_4^{\AKNS} &= -\frac{i}{2}\int q'r''-q''r' + 3(q^2rr' - r^2qq')\,dx = -i\int q'r'' + 3q^2rr' \,dx,\\
    H_5^{\AKNS} &= \int q''r'' + \frac{3}{2}(q^2)'(r^2)' + ((qr)')^2 + 2q^3 r^3.
\end{align*}
\subsubsection{Complex KdV}
We set $r = 1$. Then,
\begin{align*}
    \alpha_0 &= 0, \quad \beta_0 = 0, \quad \gamma_0 = 1,\\
    \alpha_1 &= -iq, \quad \beta_1 = i, \quad \gamma_1 = 0,\\
    \alpha_2 &= q',\quad \beta_2 = 0,\quad \gamma_2 = 2q,\\
    \alpha_3 &= iq'' -2iq^2, \quad \beta_3 = 2i q, \quad \gamma_3 = 2iq',\\
    \alpha_4 &= -q''' + 6qq', \quad \beta_4 = 0, \quad \gamma_4 = -2q''+ 6q^2,\\
    \alpha_5 &= -i(q^{(4)}-6(q')^2 -8qq''+6q^3), \quad \beta_5 = -i(2q''+6q^2), \quad \gamma_5 = 2i(-q'''+6qq')
\end{align*}

\subsubsection{Defocusing NLS, complex mKdV}
We set $r = \bar q$. Then,
\begin{align*}
    \alpha_0 &= 0, \quad \beta_0 = 0, \quad \gamma_0 = 1,\\
    \alpha_1 &= -iq, \quad \beta_1 = i\bar q, \quad \gamma_1 = 0,\\
    \alpha_2 &= q',\quad \beta_2 = \bar q',\quad \gamma_2 = 2|q|^2,\\
    \alpha_3 &= iq'' -2i|q|^2 q, \quad \beta_3 = -i\bar q'' + 2i|q|^2 \bar q, \quad \gamma_3 = 4\im(q\bar q'),\\
    \alpha_4 &= -q''' + 6|q|^2q', \quad \beta_4 = -\bar q''' + 6|q|^2 \bar q', \quad \gamma_4 = -2(2\real(q\bar q'')-|q'|^2)+ 6|q|^4,\\
    \alpha_5 &= -i(q^{(4)} - 8|q|^2q''-6(q')^2\bar q-4q|q'|^2-2q^2\bar q''+6|q|^4q),\\
    \beta_5 &= i(\bar q^{(4)} - 8|q|^2\bar q''-6(\bar q')^2q-4\bar q|q'|^2-2\bar q^2q''+6|q|^4 \bar q),\\
    \gamma_5 &= 4\im(q'''\bar q -q''\bar q') + 12\im(|q|^2 q\bar q')
\end{align*}
and
\begin{align*}
    H_{1} &= \int |q|^2 \,dx,\\
    H_2 &= \im \int q\bar q'\,dx,\\
    H_3 &= \int |q'|^2+|q|^4\,dx,\\
    H_4 &= \im\int q'\bar q'' + 3|q|^2q\bar q' \,dx,\\
    H_5 &= \int |q''|^2 + \frac{3}{2}|(q^2)'|^2 + ((|q|^2)')^2 + 2|q|^6.
\end{align*}

\subsubsection{Defocusing real mKdV}
We set $r = q \in \R$. Then,
\begin{align*}
    \alpha_0 &= 0, \quad \beta_0 = 0, \quad \gamma_0 = 1,\\
    \alpha_1 &= -iq, \quad \beta_1 = iq, \quad \gamma_1 = 0,\\
    \alpha_2 &= q',\quad \beta_2 = q',\quad \gamma_2 = 2q^2,\\
    \alpha_3 &= iq'' -2iq^3, \quad \beta_3 = -iq'' + 2iq^3, \quad \gamma_3 = 0,\\
    \alpha_4 &= -q''' + 6q^2q', \quad \beta_4 = -q''' + 6q^2 q', \quad \gamma_4 = -2(2qq''-(q')^2)+ 6q^4,\\
    \alpha_5 &= -i(q^{(4)} - 10q^2q''-10(q')^2 q + 6q^5),\\
    \beta_5 &= i(q^{(4)} - 10q^2q''-10(q')^2 q + 6q^5), \quad \gamma_5 = 0
\end{align*}
and
\begin{align*}
    H_{1} &= \int q^2 \,dx, \quad H_2 = 0,\\
    H_3 &= \int (q')^2+q^4\,dx, \quad H_4 = 0,\\
    H_5 &= \int (q'')^2 + 10q^2(q')^2 + 2q^6.
\end{align*}

\subsubsection{Gardner}
We set $q = w, r = w+2\tau_0 \in \R$. Then,
\begin{align*}
    \alpha_0 &= 0, \quad \beta_0 = 0, \quad \gamma_0 = 1,\\
    \alpha_1 &= -iw, \quad \beta_1 = i(w+2\tau_0), \quad \gamma_1 = 0,\\
    \alpha_2 &= w',\quad \beta_2 = w',\quad \gamma_2 = 2w(w+2\tau_0),\\
    \alpha_3 &= iw'' -2iw^2(w+2\tau_0), \quad \beta_3 = -iw'' + 2i(w+2\tau_0)^2 w, \quad \gamma_3 = 4i\tau_0w',\\
    \alpha_4 &= -w''' + 6w^2w' + 12\tau_0 ww', \quad \beta_4 = -w''' + 6w^2w' + 12\tau_0 ww',\\
    \gamma_4 &= -2(2\tau_0w''-(w')^2)+ 6w^2(w+2\tau_0)^2,\\
    \alpha_5 &= -i(w^{(4)} - 8ww''(w+2\tau_0)-6(w')^2(w+2\tau_0) \\
    &\qquad -4w(w')^2-2w^2w''+6w^3(w+2\tau_0)^2),\\
    \beta_5 &= i(w^{(4)} - 8(w+2\tau_0)w''w-6(w')^2w \\
    &\qquad -4(w+2\tau_0)(w')^2-2(w+2\tau_0)^2w''+6(w+2\tau_0)^3w^2),\\
    \gamma_5 &= -2i(2\tau_0w'''+ 6(-ww'(w+2\tau_0)^2+(w+2\tau_0)w'w^2))
\end{align*}
and the Hamiltonians $H_n^{\Wadati}(w,\tau_0) = H_n^{\AKNS}(w,w+2\tau_0)$ become
\begin{align*}
    H_{1}^{\Wadati} &= \int w^2 + 2\tau_0 w \,dx, \\
    H_2^{\Wadati} &= 0,\\
    H_3^{\Wadati} &= \int (w')^2 +w^4 + w^2(w+2\tau_0)^2\,dx, \\
    H_4^{\Wadati} &= 0,\\
    H_5^{\Wadati} &= \int (w'')^2 + \frac{3}{2}(w^2)'((w+2\tau_0)^2)' + ((w(w+2\tau_0))')^2 + 2w^3(w+2\tau_0)^3.
\end{align*}
The Gardner Hamiltonians $H_{n-1}^{\Gardner}(w,\tau_0) = \frac{1}{2}H_{2n+1}^{\Wadati}(w,\tau_0) - 4\tau_0^2 H_{n-1}^{\Gardner}(w,\tau_0)$ if $n \geq 1$ are
\begin{align*}
    H_{0}^{\Gardner} &= \int w^2 \,dx, \\
    H_1^{\Gardner} &= \int w_x^2 +w^4 + w^4 +4\tau_0 w^3,\\
    H_2^{\Gardner} &= \int w_{xx}^2 +10w^2w_x^2 + 2w^6 + 4\tau_0(5ww_x^2 + 3w^5) + 24\tau_0^2 w^4\,dx.
\end{align*}

\newpage
\printbibliography

@book {MR1747916,
    AUTHOR = {Cannas da Silva, Ana and Weinstein, Alan},
     TITLE = {Geometric models for noncommutative algebras},
    SERIES = {Berkeley Mathematics Lecture Notes},
    VOLUME = {10},
 PUBLISHER = {American Mathematical Society, Providence, RI; Berkeley Center
              for Pure and Applied Mathematics, Berkeley, CA},
      YEAR = {1999},
     PAGES = {xiv+184},
      ISBN = {0-8218-0952-0},
   MRCLASS = {58B34 (16S80 17B63 46L89 53D17 53D55 58B32)},
  MRNUMBER = {1747916},
MRREVIEWER = {Johannes Huebschmann},
}

@article {MR4448993,
    AUTHOR = {Bahouri, Hajer and Perelman, Galina},
     TITLE = {Global well-posedness for the derivative nonlinear
              {S}chr\"{o}dinger equation},
   JOURNAL = {Invent. Math.},
  FJOURNAL = {Inventiones Mathematicae},
    VOLUME = {229},
      YEAR = {2022},
    NUMBER = {2},
     PAGES = {639--688},
   ISSN = {0020-9910,1432-1297},
   MRCLASS = {35Q55},
  MRNUMBER = {4448993},
MRREVIEWER = {Marco\ Gipo\ Ghimenti},
       DOI = {10.1007/s00222-022-01113-0},
       URL = {https://doi.org/10.1007/s00222-022-01113-0},
}

@article {MR4565673,
    AUTHOR = {Harrop-Griffiths, Benjamin and Killip, Rowan and Vi\c{s}an,
              Monica},
     TITLE = {Large-data equicontinuity for the derivative {NLS}},
   JOURNAL = {Int. Math. Res. Not. IMRN},
  FJOURNAL = {International Mathematics Research Notices. IMRN},
      YEAR = {2023},
    NUMBER = {6},
     PAGES = {4601--4642},
      ISSN = {1073-7928,1687-0247},
   MRCLASS = {35Q55},
  MRNUMBER = {4565673},
       DOI = {10.1093/imrn/rnab374},
       URL = {https://doi.org/10.1093/imrn/rnab374},
}

@misc{killip2023sharp,
      title={Sharp well-posedness for the Benjamin--Ono equation}, 
      author={Rowan Killip and Thierry Laurens and Monica Visan},
      year={2023},
      eprint={2304.00124},
      archivePrefix={arXiv},
      primaryClass={math.AP}
}

@article {gerard2022explicit,
    AUTHOR = {G\'{e}rard, Patrick},
     TITLE = {An explicit formula for the {B}enjamin-{O}no equation},
   JOURNAL = {Tunis. J. Math.},
  FJOURNAL = {Tunisian Journal of Mathematics},
    VOLUME = {5},
      YEAR = {2023},
    NUMBER = {3},
     PAGES = {593--603},
      ISSN = {2576-7658,2576-7666},
   MRCLASS = {35Q53 (35C05 37K15 47B35)},
  MRNUMBER = {4662323},
       DOI = {10.2140/tunis.2023.5.593},
       URL = {https://doi.org/10.2140/tunis.2023.5.593},
}

@article {MR3652066,
    AUTHOR = {G\'{e}rard, Patrick and Grellier, Sandrine},
     TITLE = {The cubic {S}zeg\H{o} equation and {H}ankel operators},
   JOURNAL = {Ast\'{e}risque},
  FJOURNAL = {Ast\'{e}risque},
    NUMBER = {389},
      YEAR = {2017},
     PAGES = {vi+112},
      ISSN = {0303-1179,2492-5926},
      ISBN = {978-2-85629-854-1},
   MRCLASS = {37K15 (30J10 35F20 35Q55 35R11 47B35)},
  MRNUMBER = {3652066},
MRREVIEWER = {Jens\ Wirth},
}

@article {MR2173592,
    AUTHOR = {Praught, Jeffery and Smirnov, Roman G.},
     TITLE = {Andrew {L}enard: a mystery unraveled},
   JOURNAL = {SIGMA Symmetry Integrability Geom. Methods Appl.},
  FJOURNAL = {SIGMA. Symmetry, Integrability and Geometry. Methods and
              Applications},
    VOLUME = {1},
      YEAR = {2005},
     PAGES = {Paper 005, 7},
   MRCLASS = {37K10 (01A70 35Q51 35Q53 37-03 37K05 37K15)},
  MRNUMBER = {2173592},
MRREVIEWER = {Marco Pedroni},
       DOI = {10.3842/SIGMA.2005.005},
       URL = {https://doi.org/10.3842/SIGMA.2005.005},
}

@article {BOActa,
    AUTHOR = {G\'{e}rard, Patrick and Kappeler, Thomas and Topalov, Peter},
     TITLE = {Sharp well-posedness results of the {B}enjamin-{O}no equation
              in {$H^s(\mathbb T,\mathbb R)$} and qualitative properties of its
              solutions},
   JOURNAL = {Acta Math.},
  FJOURNAL = {Acta Mathematica},
    VOLUME = {231},
      YEAR = {2023},
    NUMBER = {1},
     PAGES = {31--88},
      ISSN = {0001-5962,1871-2509},
   MRCLASS = {58J40 (35Q53)},
  MRNUMBER = {4652410},
       DOI = {10.4310/acta.2023.v231.n1.a2},
       URL = {https://doi.org/10.4310/acta.2023.v231.n1.a2},
}

@article {BObirkhoff,
    AUTHOR = {G\'{e}rard, Patrick and Kappeler, Thomas},
     TITLE = {On the integrability of the {B}enjamin-{O}no equation on the
              torus},
   JOURNAL = {Comm. Pure Appl. Math.},
  FJOURNAL = {Communications on Pure and Applied Mathematics},
    VOLUME = {74},
      YEAR = {2021},
    NUMBER = {8},
     PAGES = {1685--1747},
      ISSN = {0010-3640,1097-0312},
   MRCLASS = {35Q53 (58J40)},
  MRNUMBER = {4275336},
MRREVIEWER = {Mircea\ Crasmareanu},
       DOI = {10.1002/cpa.21896},
       URL = {https://doi.org/10.1002/cpa.21896},
}

@article {MR1215780,
    AUTHOR = {Bourgain, J.},
     TITLE = {Fourier transform restriction phenomena for certain lattice
              subsets and applications to nonlinear evolution equations.
              {II}. {T}he {K}d{V}-equation},
   JOURNAL = {Geom. Funct. Anal.},
  FJOURNAL = {Geometric and Functional Analysis},
    VOLUME = {3},
      YEAR = {1993},
    NUMBER = {3},
     PAGES = {209--262},
      ISSN = {1016-443X},
   MRCLASS = {35Q55},
  MRNUMBER = {1215780},
MRREVIEWER = {Yun Mei Chen},
       DOI = {10.1007/BF01895688},
       URL = {https://doi.org/10.1007/BF01895688},
}

@article {MR1780703,
    AUTHOR = {Christodoulou, Demetrios and Lindblad, Hans},
     TITLE = {On the motion of the free surface of a liquid},
   JOURNAL = {Comm. Pure Appl. Math.},
  FJOURNAL = {Communications on Pure and Applied Mathematics},
    VOLUME = {53},
      YEAR = {2000},
    NUMBER = {12},
     PAGES = {1536--1602},
      ISSN = {0010-3640},
   MRCLASS = {76B15 (35Q53 35Q55 35Q58)},
  MRNUMBER = {1780703},
MRREVIEWER = {Walter Craig},
       DOI = {10.1002/1097-0312(200012)53:12<1536::AID-CPA2>3.3.CO;2-H},
       URL =
              {https://doi.org/10.1002/1097-0312(200012)53:12<1536::AID-CPA2>3.3.CO;2-H},
}

@article {MR1969209,
    AUTHOR = {Colliander, J. and Keel, M. and Staffilani, G. and Takaoka, H.
              and Tao, T.},
     TITLE = {Sharp global well-posedness for {K}d{V} and modified {K}d{V}
              on {$\mathbb R$} and {$\mathbb T$}},
   JOURNAL = {J. Amer. Math. Soc.},
  FJOURNAL = {Journal of the American Mathematical Society},
    VOLUME = {16},
      YEAR = {2003},
    NUMBER = {3},
     PAGES = {705--749},
      ISSN = {0894-0347},
   MRCLASS = {35Q53 (35B30 37K10)},
  MRNUMBER = {1969209},
MRREVIEWER = {Vladislav G. Dubrovsky},
       DOI = {10.1090/S0894-0347-03-00421-1},
       URL = {https://doi.org/10.1090/S0894-0347-03-00421-1},
}

@article {MR1780702,
    AUTHOR = {Schneider, Guido and Wayne, C. Eugene},
     TITLE = {The long-wave limit for the water wave problem. {I}. {T}he
              case of zero surface tension},
   JOURNAL = {Comm. Pure Appl. Math.},
  FJOURNAL = {Communications on Pure and Applied Mathematics},
    VOLUME = {53},
      YEAR = {2000},
    NUMBER = {12},
     PAGES = {1475--1535},
      ISSN = {0010-3640},
   MRCLASS = {76B15 (35Q53 35Q55 35Q58)},
  MRNUMBER = {1780702},
MRREVIEWER = {Walter Craig},
       DOI = {10.1002/1097-0312(200012)53:12<1475::AID-CPA1>3.0.CO;2-V},
       URL =
              {https://doi.org/10.1002/1097-0312(200012)53:12<1475::AID-CPA1>3.0.CO;2-V},
}

@article {MR3400442,
    AUTHOR = {Buckmaster, Tristan and Koch, Herbert},
     TITLE = {The {K}orteweg--de {V}ries equation at {$H^{-1}$} regularity},
   JOURNAL = {Ann. Inst. H. Poincar\'{e} Anal. Non Lin\'{e}aire},
  FJOURNAL = {Annales de l'Institut Henri Poincar\'{e}. Analyse Non Lin\'{e}aire},
    VOLUME = {32},
      YEAR = {2015},
    NUMBER = {5},
     PAGES = {1071--1098},
      ISSN = {0294-1449},
   MRCLASS = {35Q53 (35B30 35B35 35C08 35D30)},
  MRNUMBER = {3400442},
MRREVIEWER = {John Albert},
       DOI = {10.1016/j.anihpc.2014.05.004},
       URL = {https://doi.org/10.1016/j.anihpc.2014.05.004},
}

@article {MR512420,
    AUTHOR = {Deift, P. and Trubowitz, E.},
     TITLE = {Inverse scattering on the line},
   JOURNAL = {Comm. Pure Appl. Math.},
  FJOURNAL = {Communications on Pure and Applied Mathematics},
    VOLUME = {32},
      YEAR = {1979},
    NUMBER = {2},
     PAGES = {121--251},
      ISSN = {0010-3640},
   MRCLASS = {34B25 (35P25 58F07)},
  MRNUMBER = {512420},
MRREVIEWER = {R. C. Gilbert},
       DOI = {10.1002/cpa.3160320202},
       URL = {https://doi.org/10.1002/cpa.3160320202},
}

@book {MR3364494,
    AUTHOR = {Simon, Barry},
     TITLE = {Operator theory},
    SERIES = {A Comprehensive Course in Analysis, Part 4},
 PUBLISHER = {American Mathematical Society, Providence, RI},
      YEAR = {2015},
     PAGES = {xviii+749},
      ISBN = {978-1-4704-1103-9},
   MRCLASS = {47-01 (34-01 35-01 42B35 42B37 43-01 46-01 81-01)},
  MRNUMBER = {3364494},
MRREVIEWER = {Fritz Gesztesy},
       DOI = {10.1090/simon/004},
       URL = {https://doi.org/10.1090/simon/004},
}

@article {MR2189502,
    AUTHOR = {Kappeler, Thomas and Perry, Peter and Shubin, Mikhail and
              Topalov, Peter},
     TITLE = {The {M}iura map on the line},
   JOURNAL = {Int. Math. Res. Not.},
  FJOURNAL = {International Mathematics Research Notices},
      YEAR = {2005},
    NUMBER = {50},
     PAGES = {3091--3133},
 %     ISSN = {1073-7928},
   MRCLASS = {37K15 (35Q53 37K10)},
  MRNUMBER = {2189502},
MRREVIEWER = {Dmitry G. Shepelsky},
     %  DOI = {10.1155/IMRN.2005.3091},
    %   URL = {https://doi.org/10.1155/IMRN.2005.3091},
}

@book {MR1992536,
    AUTHOR = {Gesztesy, Fritz and Holden, Helge},
     TITLE = {Soliton equations and their algebro-geometric solutions.
              {V}ol. {I}},
    SERIES = {Cambridge Studies in Advanced Mathematics},
    VOLUME = {79},
      NOTE = {$(1+1)$-dimensional continuous models},
 PUBLISHER = {Cambridge University Press, Cambridge},
      YEAR = {2003},
     PAGES = {xii+505},
 %     ISBN = {0-521-75307-4},
   MRCLASS = {37K40 (14H70 34L05 35Q53 37K10 37K15 37K20)},
  MRNUMBER = {1992536},
MRREVIEWER = {Emma Previato},
     %  DOI = {10.1017/CBO9780511546723},
    %   URL = {https://doi.org/10.1017/CBO9780511546723},
}

@article {bringmann2019global,
    AUTHOR = {Bringmann, Bjoern and Killip, Rowan and Visan, Monica},
     TITLE = {Global well-posedness for the fifth-order {K}d{V} equation in
              {$H^{-1}(\mathbb R)$}},
   JOURNAL = {Ann. PDE},
  FJOURNAL = {Annals of PDE. Journal Dedicated to the Analysis of Problems
              from Physical Sciences},
    VOLUME = {7},
      YEAR = {2021},
    NUMBER = {2},
     PAGES = {Paper No. 21, 46},
      ISSN = {2524-5317,2199-2576},
   MRCLASS = {35Q53},
  MRNUMBER = {4304314},
       DOI = {10.1007/s40818-021-00111-4},
       URL = {https://doi.org/10.1007/s40818-021-00111-4},
}

@article {MR2018661,
    AUTHOR = {Christ, Michael and Colliander, James and Tao, Terrence},
     TITLE = {Asymptotics, frequency modulation, and low regularity
              ill-posedness for canonical defocusing equations},
   JOURNAL = {Amer. J. Math.},
  FJOURNAL = {American Journal of Mathematics},
    VOLUME = {125},
      YEAR = {2003},
    NUMBER = {6},
     PAGES = {1235--1293},
      ISSN = {0002-9327},
   MRCLASS = {35Q53 (35Q55 35R25 37K10 37K15)},
  MRNUMBER = {2018661},
MRREVIEWER = {John Albert},
       URL =
              {http://muse.jhu.edu/journals/american_journal_of_mathematics/v125/125.6christ.pdf},
}

@article {MR783348,
    AUTHOR = {Segal, Graeme and Wilson, George},
     TITLE = {Loop groups and equations of {K}d{V} type},
   JOURNAL = {Inst. Hautes \'{E}tudes Sci. Publ. Math.},
  FJOURNAL = {Institut des Hautes \'{E}tudes Scientifiques. Publications
              Math\'{e}matiques},
    NUMBER = {61},
      YEAR = {1985},
     PAGES = {5--65},
      ISSN = {0073-8301},
   MRCLASS = {58F07 (14K99 35Q20 58G35)},
  MRNUMBER = {783348},
MRREVIEWER = {A. M. Vinogradov},
       URL = {http://www.numdam.org/item?id=PMIHES_1985__61__5_0},
}

@book {MR900587,
    AUTHOR = {Pressley, Andrew and Segal, Graeme},
     TITLE = {Loop groups},
    SERIES = {Oxford Mathematical Monographs},
      NOTE = {Oxford Science Publications},
 PUBLISHER = {The Clarendon Press, Oxford University Press, New York},
      YEAR = {1986},
     PAGES = {viii+318},
      ISBN = {0-19-853535-X},
   MRCLASS = {22E65 (58D15 81D15)},
  MRNUMBER = {900587},
MRREVIEWER = {Jouko Mickelsson},
}

@article {MR2531556,
    AUTHOR = {Guo, Zihua},
     TITLE = {Global well-posedness of {K}orteweg-de {V}ries equation in
              {$H^{-3/4}(\R)$}},
   JOURNAL = {J. Math. Pures Appl. (9)},
  FJOURNAL = {Journal de Math\'{e}matiques Pures et Appliqu\'{e}es. Neuvi\`eme S\'{e}rie},
    VOLUME = {91},
      YEAR = {2009},
    NUMBER = {6},
     PAGES = {583--597},
      ISSN = {0021-7824},
   MRCLASS = {35Q53 (35B30)},
  MRNUMBER = {2531556},
MRREVIEWER = {John Albert},
       DOI = {10.1016/j.matpur.2009.01.012},
       URL = {https://doi.org/10.1016/j.matpur.2009.01.012},
}

@article {MR4222602,
    AUTHOR = {Dubrovin, Boris and Yang, Di and Zagier, Don},
     TITLE = {On tau-functions for the {K}d{V} hierarchy},
   JOURNAL = {Selecta Math. (N.S.)},
  FJOURNAL = {Selecta Mathematica. New Series},
    VOLUME = {27},
      YEAR = {2021},
    NUMBER = {1},
     PAGES = {Paper No. 12, 47},
      ISSN = {1022-1824},
   MRCLASS = {37K10 (05A15 14N35 33E15 53D45)},
  MRNUMBER = {4222602},
       DOI = {10.1007/s00029-021-00620-x},
       URL = {https://doi.org/10.1007/s00029-021-00620-x},
}

@article {MR286402,
    AUTHOR = {Gardner, Clifford S.},
     TITLE = {Korteweg-de {V}ries equation and generalizations. {IV}. {T}he
              {K}orteweg-de {V}ries equation as a {H}amiltonian system},
   JOURNAL = {J. Mathematical Phys.},
  FJOURNAL = {Journal of Mathematical Physics},
    VOLUME = {12},
      YEAR = {1971},
     PAGES = {1548--1551},
 %     ISSN = {0022-2488},
   MRCLASS = {81.35 (35.00)},
  MRNUMBER = {286402},
MRREVIEWER = {Joel Smoller},
     %  DOI = {10.1063/1.1665772},
    %   URL = {https://doi.org/10.1063/1.1665772},
}

@book {MR1995460,
    AUTHOR = {Babelon, Olivier and Bernard, Denis and Talon, Michel},
     TITLE = {Introduction to classical integrable systems},
    SERIES = {Cambridge Monographs on Mathematical Physics},
 PUBLISHER = {Cambridge University Press, Cambridge},
      YEAR = {2003},
     PAGES = {xii+602},
 %     ISBN = {0-521-82267-X},
   MRCLASS = {37J35 (35Q51 35Q53 37K10 37K15 37K20 70H06)},
  MRNUMBER = {1995460},
MRREVIEWER = {	Fritz Gesztesy},
      % DOI = {10.1017/CBO9780511535024},
    %   URL = {https://doi.org/10.1017/CBO9780511535024},
}

@article {MR2326419,
    AUTHOR = {Constantin, A. and Kappeler, T. and Kolev, B. and Topalov, P.},
     TITLE = {On geodesic exponential maps of the {V}irasoro group},
   JOURNAL = {Ann. Global Anal. Geom.},
  FJOURNAL = {Annals of Global Analysis and Geometry},
    VOLUME = {31},
      YEAR = {2007},
    NUMBER = {2},
     PAGES = {155--180},
      ISSN = {0232-704X},
   MRCLASS = {58B25 (35Q35 37K65)},
  MRNUMBER = {2326419},
MRREVIEWER = {Valentin Ovsienko},
       DOI = {10.1007/s10455-006-9042-8},
       URL = {https://doi.org/10.1007/s10455-006-9042-8},
}

@article {MR2216268,
    AUTHOR = {Constantin, A. and Kolev, B.},
     TITLE = {Integrability of invariant metrics on the diffeomorphism group
              of the circle},
   JOURNAL = {J. Nonlinear Sci.},
  FJOURNAL = {Journal of Nonlinear Science},
    VOLUME = {16},
      YEAR = {2006},
    NUMBER = {2},
     PAGES = {109--122},
      ISSN = {0938-8974},
   MRCLASS = {37K65 (37K10 58B20)},
  MRNUMBER = {2216268},
MRREVIEWER = {Boris A. Khesin},
       DOI = {10.1007/s00332-005-0707-4},
       URL = {https://doi.org/10.1007/s00332-005-0707-4},
}

@book {MR2348643,
    AUTHOR = {Faddeev, Ludwig D. and Takhtajan, Leon A.},
     TITLE = {Hamiltonian methods in the theory of solitons},
    SERIES = {Classics in Mathematics},
   EDITION = {English},
      NOTE = {Translated from the 1986 Russian original by Alexey G. Reyman},
 PUBLISHER = {Springer, Berlin},
      YEAR = {2007},
     PAGES = {x+592},
   %   ISBN = {978-3-540-69843-2},
 %  MRCLASS = {37K10 (35P25 35Q51 35Q55 35R30 37J35 37N20 81R12)},
%  MRNUMBER = {2348643},
}

@book {MR1964513,
    AUTHOR = {Dickey, L. A.},
     TITLE = {Soliton equations and {H}amiltonian systems},
    SERIES = {Advanced Series in Mathematical Physics},
    VOLUME = {26},
   EDITION = {Second},
 PUBLISHER = {World Scientific Publishing Co., Inc., River Edge, NJ},
      YEAR = {2003},
     PAGES = {xii+408},
  %    ISBN = {981-238-173-2},
   MRCLASS = {37K10 (35F20 35Q51 35Q53 37K15)},
  MRNUMBER = {1964513},
MRREVIEWER = {Marco Pedroni},
   %    DOI = {10.1142/5108},
  %     URL = {https://doi.org/10.1142/5108},
}

@book {MR779467,
    AUTHOR = {Novikov, S. and Manakov, S. V. and Pitaevski\u{\i}, L. P. and
              Zakharov, V. E.},
     TITLE = {Theory of solitons},
    SERIES = {Contemporary Soviet Mathematics},
      NOTE = {The inverse scattering method,
              Translated from the Russian},
 PUBLISHER = {Consultants Bureau [Plenum], New York},
      YEAR = {1984},
     PAGES = {xi+276},
  %    ISBN = {0-306-10977-8},
   MRCLASS = {35Q20 (58F07 76B25)},
  MRNUMBER = {779467},
}

@article {MR3874652,
    AUTHOR = {Koch, Herbert and Tataru, Daniel},
     TITLE = {Conserved energies for the cubic nonlinear {S}chr\"{o}dinger
              equation in one dimension},
   JOURNAL = {Duke Math. J.},
  FJOURNAL = {Duke Mathematical Journal},
    VOLUME = {167},
      YEAR = {2018},
    NUMBER = {17},
     PAGES = {3207--3313},
%      ISSN = {0012-7094},
   MRCLASS = {35Q55 (35Q53 37K05 37K10)},
  MRNUMBER = {3874652},
MRREVIEWER = {Thierry Cazenave},
     %  DOI = {10.1215/00127094-2018-0033},
     %  URL = {https://doi.org/10.1215/00127094-2018-0033},
}

@book {MR874343,
    AUTHOR = {Schuur, Peter Cornelis},
     TITLE = {Asymptotic analysis of soliton problems},
    SERIES = {Lecture Notes in Mathematics},
    VOLUME = {1232},
      NOTE = {An inverse scattering approach},
 PUBLISHER = {Springer-Verlag, Berlin},
      YEAR = {1986},
     PAGES = {viii+180},
      ISBN = {3-540-17203-3},
   MRCLASS = {35-02 (15A60 35Q20 45M05)},
  MRNUMBER = {874343},
MRREVIEWER = {J. J. C. Nimmo},
       DOI = {10.1007/BFb0073054},
       URL = {https://doi.org/10.1007/BFb0073054},
}

@misc{harropgriffiths2020sharp,
      title={Sharp well-posedness for the cubic NLS and mKdV in $H^s(\mathbb R)$}, 
      author={Benjamin Harrop-Griffiths and Rowan Killip and Monica Visan},
      year={2020},
      eprint={2003.05011},
      archivePrefix={arXiv},
     % primaryClass={math.AP}
}

@article {MR2830706,
    AUTHOR = {Molinet, Luc},
     TITLE = {A note on ill posedness for the {K}d{V} equation},
   JOURNAL = {Differential Integral Equations},
  FJOURNAL = {Differential and Integral Equations. An International Journal
              for Theory \& Applications},
    VOLUME = {24},
      YEAR = {2011},
    NUMBER = {7-8},
     PAGES = {759--765},
      ISSN = {0893-4983},
   MRCLASS = {35Q53 (35B30)},
  MRNUMBER = {2830706},
MRREVIEWER = {Luiz Gustavo Farah},
}

@article {MR2927357,
    AUTHOR = {Molinet, Luc},
     TITLE = {Sharp ill-posedness results for the {K}d{V} and m{K}d{V}
              equations on the torus},
   JOURNAL = {Adv. Math.},
  FJOURNAL = {Advances in Mathematics},
    VOLUME = {230},
      YEAR = {2012},
    NUMBER = {4-6},
     PAGES = {1895--1930},
      ISSN = {0001-8708},
   MRCLASS = {35Q53 (35A01 35B30 35B45 35D30 35R25)},
  MRNUMBER = {2927357},
MRREVIEWER = {Corentin Audiard},
       DOI = {10.1016/j.aim.2012.03.026},
       URL = {https://doi.org/10.1016/j.aim.2012.03.026},
}

@article {MR535697,
    AUTHOR = {Kato, Tosio},
     TITLE = {On the {K}orteweg-de\thinspace {V}ries equation},
   JOURNAL = {Manuscripta Math.},
  FJOURNAL = {Manuscripta Mathematica},
    VOLUME = {28},
      YEAR = {1979},
    NUMBER = {1-3},
     PAGES = {89--99},
      ISSN = {0025-2611},
   MRCLASS = {35Q20},
  MRNUMBER = {535697},
       DOI = {10.1007/BF01647967},
       URL = {https://doi.org/10.1007/BF01647967},
}

@article {MR1329387,
    AUTHOR = {Kenig, Carlos E. and Ponce, Gustavo and Vega, Luis},
     TITLE = {A bilinear estimate with applications to the {K}d{V} equation},
   JOURNAL = {J. Amer. Math. Soc.},
  FJOURNAL = {Journal of the American Mathematical Society},
    VOLUME = {9},
      YEAR = {1996},
    NUMBER = {2},
     PAGES = {573--603},
      ISSN = {0894-0347},
   MRCLASS = {35Q53 (35Bxx)},
  MRNUMBER = {1329387},
MRREVIEWER = {F. Pempinelli},
       DOI = {10.1090/S0894-0347-96-00200-7},
       URL = {https://doi.org/10.1090/S0894-0347-96-00200-7},
}

@book {MR2456522,
    AUTHOR = {Khesin, Boris and Wendt, Robert},
     TITLE = {The geometry of infinite-dimensional groups},
    SERIES = {Ergebnisse der Mathematik und ihrer Grenzgebiete. 3. Folge. A
              Series of Modern Surveys in Mathematics [Results in
              Mathematics and Related Areas. 3rd Series. A Series of Modern
              Surveys in Mathematics]},
    VOLUME = {51},
 PUBLISHER = {Springer-Verlag, Berlin},
      YEAR = {2009},
     PAGES = {xii+304},
      ISBN = {978-3-540-77262-0},
   MRCLASS = {58B25 (22E65 37K10 37K30 37K65 58D30)},
  MRNUMBER = {2456522},
MRREVIEWER = {Daniel Belti\c{t}\u{a}},
}

@article {MR3990604,
    AUTHOR = {Killip, Rowan and Vi\c{s}an, Monica},
     TITLE = {Kd{V} is well-posed in {$H^{-1}$}},
   JOURNAL = {Ann. of Math. (2)},
  FJOURNAL = {Annals of Mathematics. Second Series},
    VOLUME = {190},
      YEAR = {2019},
    NUMBER = {1},
     PAGES = {249--305},
      ISSN = {0003-486X},
   MRCLASS = {35Q53 (35B30 37K10)},
  MRNUMBER = {3990604},
MRREVIEWER = {John Albert},
       DOI = {10.4007/annals.2019.190.1.4},
       URL = {https://doi.org/10.4007/annals.2019.190.1.4},
}

@article {MR3820439,
    AUTHOR = {Killip, Rowan and Vi\c{s}an, Monica and Zhang, Xiaoyi},
     TITLE = {Low regularity conservation laws for integrable {PDE}},
   JOURNAL = {Geom. Funct. Anal.},
  FJOURNAL = {Geometric and Functional Analysis},
    VOLUME = {28},
      YEAR = {2018},
    NUMBER = {4},
     PAGES = {1062--1090},
      ISSN = {1016-443X},
   MRCLASS = {35Q53 (35B45 35Q55 37K05)},
  MRNUMBER = {3820439},
MRREVIEWER = {Thierry Cazenave},
       DOI = {10.1007/s00039-018-0444-0},
       URL = {https://doi.org/10.1007/s00039-018-0444-0},
}

@article {MR4186010,
    AUTHOR = {Koch, Herbert and Liao, Xian},
     TITLE = {Conserved energies for the one dimensional
              {G}ross-{P}itaevskii equation},
   JOURNAL = {Adv. Math.},
  FJOURNAL = {Advances in Mathematics},
    VOLUME = {377},
      YEAR = {2021},
     PAGES = {107467, 83},
      ISSN = {0001-8708},
   MRCLASS = {35Q55 (37K10)},
  MRNUMBER = {4186010},
       DOI = {10.1016/j.aim.2020.107467},
       URL = {https://doi.org/10.1016/j.aim.2020.107467},
}

@article {MR3292346,
    AUTHOR = {Liu, Baoping},
     TITLE = {A priori bounds for {K}d{V} equation below {$H^{-\frac34}$}},
   JOURNAL = {J. Funct. Anal.},
  FJOURNAL = {Journal of Functional Analysis},
    VOLUME = {268},
      YEAR = {2015},
    NUMBER = {3},
     PAGES = {501--554},
      ISSN = {0022-1236},
   MRCLASS = {35Q53 (35B45 35D30)},
  MRNUMBER = {3292346},
MRREVIEWER = {Peter E. Zhidkov},
       DOI = {10.1016/j.jfa.2014.06.020},
       URL = {https://doi.org/10.1016/j.jfa.2014.06.020},
}

@article {MR2267286,
    AUTHOR = {Kappeler, T. and Topalov, P.},
     TITLE = {Global wellposedness of {K}d{V} in {$H^{-1}(\mathbb T,\R)$}},
   JOURNAL = {Duke Math. J.},
  FJOURNAL = {Duke Mathematical Journal},
    VOLUME = {135},
      YEAR = {2006},
    NUMBER = {2},
     PAGES = {327--360},
      ISSN = {0012-7094},
   MRCLASS = {35Q53 (35B30)},
  MRNUMBER = {2267286},
MRREVIEWER = {Pierre A. Vuillermot},
       DOI = {10.1215/S0012-7094-06-13524-X},
       URL = {https://doi.org/10.1215/S0012-7094-06-13524-X},
}

@article {MR2501679,
    AUTHOR = {Kishimoto, Nobu},
     TITLE = {Well-posedness of the {C}auchy problem for the {K}orteweg-de
              {V}ries equation at the critical regularity},
   JOURNAL = {Differential Integral Equations},
  FJOURNAL = {Differential and Integral Equations. An International Journal
              for Theory \& Applications},
    VOLUME = {22},
      YEAR = {2009},
    NUMBER = {5-6},
     PAGES = {447--464},
      ISSN = {0893-4983},
   MRCLASS = {35Q53 (35B30)},
  MRNUMBER = {2501679},
MRREVIEWER = {Liana L. Fleming},
}

@book {MR1736222,
    AUTHOR = {Miwa, T. and Jimbo, M. and Date, E.},
     TITLE = {Solitons},
    SERIES = {Cambridge Tracts in Mathematics},
    VOLUME = {135},
      NOTE = {Differential equations, symmetries and infinite-dimensional
              algebras,
              Translated from the 1993 Japanese original by Miles Reid},
 PUBLISHER = {Cambridge University Press, Cambridge},
      YEAR = {2000},
     PAGES = {x+108},
      ISBN = {0-521-56161-2},
   MRCLASS = {37K10 (17B67 35Q51 35Q53 37K20 37K30 58B25)},
  MRNUMBER = {1736222},
MRREVIEWER = {Stanislav Z. Pakuliak},
}

@incollection {MR730247,
    AUTHOR = {Sato, Mikio and Sato, Yasuko},
     TITLE = {Soliton equations as dynamical systems on infinite-dimensional
              {G}rassmann manifold},
 BOOKTITLE = {Nonlinear partial differential equations in applied science
              ({T}okyo, 1982)},
    SERIES = {North-Holland Math. Stud.},
    VOLUME = {81},
     PAGES = {259--271},
 PUBLISHER = {North-Holland, Amsterdam},
      YEAR = {1983},
   MRCLASS = {58F07 (35Q20 58D25)},
  MRNUMBER = {730247},
MRREVIEWER = {Vassil Tsanov},
}

@book {MR2085332,
    AUTHOR = {Hirota, Ryogo},
     TITLE = {The direct method in soliton theory},
    SERIES = {Cambridge Tracts in Mathematics},
    VOLUME = {155},
      NOTE = {Translated from the 1992 Japanese original and edited by
              Atsushi Nagai, Jon Nimmo and Claire Gilson,
              With a foreword by Jarmo Hietarinta and Nimmo},
 PUBLISHER = {Cambridge University Press, Cambridge},
      YEAR = {2004},
     PAGES = {xii+200},
      ISBN = {0-521-83660-3},
   MRCLASS = {37K10 (35Q51 37K35 37K40)},
  MRNUMBER = {2085332},
MRREVIEWER = {Anatoly G. Meshkov},
       DOI = {10.1017/CBO9780511543043},
       URL = {https://doi.org/10.1017/CBO9780511543043},
}

@incollection {MR1723386,
    AUTHOR = {Segal, Graeme},
     TITLE = {Integrable systems and inverse scattering},
 BOOKTITLE = {Integrable systems ({O}xford, 1997)},
    SERIES = {Oxf. Grad. Texts Math.},
    VOLUME = {4},
     PAGES = {53--119},
 PUBLISHER = {Oxford Univ. Press, New York},
      YEAR = {1999},
   MRCLASS = {37Kxx (34A55 34L25 35Q53 37J35)},
  MRNUMBER = {1723386},
MRREVIEWER = {Piotr G. Grinevich},
}

@article {MR252826,
    AUTHOR = {Miura, Robert M. and Gardner, Clifford S. and Kruskal, Martin
              D.},
     TITLE = {Korteweg-de {V}ries equation and generalizations. {II}.
              {E}xistence of conservation laws and constants of motion},
   JOURNAL = {J. Mathematical Phys.},
  FJOURNAL = {Journal of Mathematical Physics},
    VOLUME = {9},
      YEAR = {1968},
     PAGES = {1204--1209},
      ISSN = {0022-2488},
   MRCLASS = {35.36 (82.00)},
  MRNUMBER = {252826},
MRREVIEWER = {Joel Smoller},
       DOI = {10.1063/1.1664701},
       URL = {https://doi.org/10.1063/1.1664701},
}

@book {Benyamini,
    AUTHOR = {Benyamini, Yoav and Lindenstrauss, Joram},
     TITLE = {Geometric nonlinear functional analysis. {V}ol. 1},
    SERIES = {American Mathematical Society Colloquium Publications},
    VOLUME = {48},
 PUBLISHER = {American Mathematical Society, Providence, RI},
      YEAR = {2000},
     PAGES = {xii+488},
      ISBN = {0-8218-0835-4},
   MRCLASS = {46-02 (46Bxx 46T99 47-02)},
  MRNUMBER = {1727673},
MRREVIEWER = {Gilles Godefroy},
       DOI = {10.1090/coll/048},
       URL = {https://doi.org/10.1090/coll/048},
}

@article {Saut-KdV,
    AUTHOR = {Saut, J.-C.},
     TITLE = {Quelques g\'{e}n\'{e}ralisations de l'\'{e}quation de
              {K}orteweg-de\thinspace {V}ries. {II}},
   JOURNAL = {J. Differential Equations},
  FJOURNAL = {Journal of Differential Equations},
    VOLUME = {33},
      YEAR = {1979},
    NUMBER = {3},
     PAGES = {320--335},
      ISSN = {0022-0396,1090-2732},
   MRCLASS = {35Q20 (76B99)},
  MRNUMBER = {543702},
       DOI = {10.1016/0022-0396(79)90068-8},
       URL = {https://doi.org/10.1016/0022-0396(79)90068-8},
}

@article {KPV-higher1,
    AUTHOR = {Kenig, Carlos E. and Ponce, Gustavo and Vega, Luis},
     TITLE = {Higher-order nonlinear dispersive equations},
   JOURNAL = {Proc. Amer. Math. Soc.},
  FJOURNAL = {Proceedings of the American Mathematical Society},
    VOLUME = {122},
      YEAR = {1994},
    NUMBER = {1},
     PAGES = {157--166},
      ISSN = {0002-9939,1088-6826},
   MRCLASS = {35G25 (35Q53)},
  MRNUMBER = {1195480},
MRREVIEWER = {Michael\ Wiegner},
       DOI = {10.2307/2160855},
       URL = {https://doi.org/10.2307/2160855},
}

@incollection {KPV-higher2,
    AUTHOR = {Kenig, Carlos E. and Ponce, Gustavo and Vega, Luis},
     TITLE = {On the hierarchy of the generalized {K}d{V} equations},
 BOOKTITLE = {Singular limits of dispersive waves ({L}yon, 1991)},
    SERIES = {NATO Adv. Sci. Inst. Ser. B: Phys.},
    VOLUME = {320},
     PAGES = {347--356},
 PUBLISHER = {Plenum, New York},
      YEAR = {1994},
      ISBN = {0-306-44628-6},
   MRCLASS = {35Q53},
  MRNUMBER = {1321214},
MRREVIEWER = {Daniel\ B\"{a}ttig},
}

@article{Pilod,
    AUTHOR = {Pilod, Didier},
     TITLE = {On the {C}auchy problem for higher-order nonlinear dispersive
              equations},
   JOURNAL = {J. Differential Equations},
  FJOURNAL = {Journal of Differential Equations},
    VOLUME = {245},
      YEAR = {2008},
    NUMBER = {8},
     PAGES = {2055--2077},
      ISSN = {0022-0396,1090-2732},
   MRCLASS = {35Q53 (35B30)},
  MRNUMBER = {2446185},
MRREVIEWER = {Sebastian\ Herr},
       DOI = {10.1016/j.jde.2008.07.017},
       URL = {https://doi.org/10.1016/j.jde.2008.07.017},
}

@article{Grunrock,
    AUTHOR = {Gr\"{u}nrock, Axel},
     TITLE = {On the hierarchies of higher order m{K}d{V} and {K}d{V}
              equations},
   JOURNAL = {Cent. Eur. J. Math.},
  FJOURNAL = {Central European Journal of Mathematics},
    VOLUME = {8},
      YEAR = {2010},
    NUMBER = {3},
     PAGES = {500--536},
      ISSN = {1895-1074,1644-3616},
   MRCLASS = {35Q53 (37K10)},
  MRNUMBER = {2653659},
       DOI = {10.2478/s11533-010-0024-5},
       URL = {https://doi.org/10.2478/s11533-010-0024-5},
}

@article {Kenig-Pilod,
    AUTHOR = {Kenig, Carlos E. and Pilod, Didier},
     TITLE = {Well-posedness for the fifth-order {K}d{V} equation in the
              energy space},
   JOURNAL = {Trans. Amer. Math. Soc.},
  FJOURNAL = {Transactions of the American Mathematical Society},
    VOLUME = {367},
      YEAR = {2015},
    NUMBER = {4},
     PAGES = {2551--2612},
      ISSN = {0002-9947,1088-6850},
   MRCLASS = {35Q53 (35B30 37K05)},
  MRNUMBER = {3301874},
MRREVIEWER = {Peter\ E.\ Zhidkov},
       DOI = {10.1090/S0002-9947-2014-05982-5},
       URL = {https://doi.org/10.1090/S0002-9947-2014-05982-5},
}

@article {GuoKwakKwon,
    AUTHOR = {Guo, Zihua and Kwak, Chulkwang and Kwon, Soonsik},
     TITLE = {Rough solutions of the fifth-order {K}d{V} equations},
   JOURNAL = {J. Funct. Anal.},
  FJOURNAL = {Journal of Functional Analysis},
    VOLUME = {265},
      YEAR = {2013},
    NUMBER = {11},
     PAGES = {2791--2829},
      ISSN = {0022-1236,1096-0783},
   MRCLASS = {35Q53 (35B30 35B45)},
  MRNUMBER = {3096990},
MRREVIEWER = {Andrei\ V.\ Faminskii},
       DOI = {10.1016/j.jfa.2013.08.010},
       URL = {https://doi.org/10.1016/j.jfa.2013.08.010},
}

@article {BKV,
    AUTHOR = {Bringmann, Bjoern and Killip, Rowan and Visan, Monica},
     TITLE = {Global well-posedness for the fifth-order {K}d{V} equation in
              {$H^{-1}(\mathbb R)$}},
   JOURNAL = {Ann. PDE},
  FJOURNAL = {Annals of PDE. Journal Dedicated to the Analysis of Problems
              from Physical Sciences},
    VOLUME = {7},
      YEAR = {2021},
    NUMBER = {2},
     PAGES = {Paper No. 21, 46},
      ISSN = {2524-5317,2199-2576},
   MRCLASS = {35Q53},
  MRNUMBER = {4304314},
       DOI = {10.1007/s40818-021-00111-4},
       URL = {https://doi.org/10.1007/s40818-021-00111-4},
}

@article {Kappeler-Molnar,
    AUTHOR = {Kappeler, Thomas and Molnar, Jan-Cornelius},
     TITLE = {On the wellposedness of the {K}d{V}/{K}d{V}2 equations and
              their frequency maps},
   JOURNAL = {Ann. Inst. H. Poincar\'{e} C Anal. Non Lin\'{e}aire},
  FJOURNAL = {Annales de l'Institut Henri Poincar\'{e} C. Analyse Non
              Lin\'{e}aire},
    VOLUME = {35},
      YEAR = {2018},
    NUMBER = {1},
     PAGES = {101--160},
      ISSN = {0294-1449,1873-1430},
   MRCLASS = {37K10 (35B30 35Q53 35R25)},
  MRNUMBER = {3739929},
MRREVIEWER = {Jens\ Wirth},
       DOI = {10.1016/j.anihpc.2017.03.003},
       URL = {https://doi.org/10.1016/j.anihpc.2017.03.003},
}

@misc{harropgriffiths2023global,
      title={Global well-posedness for the derivative nonlinear Schr\"odinger equation in $L^2(\mathbb{R})$}, 
      author={Benjamin Harrop-Griffiths and Rowan Killip and Maria Ntekoume and Monica Visan},
      year={2023},
      eprint={2204.12548},
      archivePrefix={arXiv},
      primaryClass={math.AP}
}

@article {Gassot,
    AUTHOR = {Gassot, Louise},
     TITLE = {The third order {B}enjamin-{O}no equation on the torus:
              well-posedness, traveling waves and stability},
   JOURNAL = {Ann. Inst. H. Poincar\'{e} C Anal. Non Lin\'{e}aire},
  FJOURNAL = {Annales de l'Institut Henri Poincar\'{e} C. Analyse Non
              Lin\'{e}aire},
    VOLUME = {38},
      YEAR = {2021},
    NUMBER = {3},
     PAGES = {815--840},
      ISSN = {0294-1449,1873-1430},
   MRCLASS = {37K10 (35Q53)},
  MRNUMBER = {4227053},
       DOI = {10.1016/j.anihpc.2020.09.004},
       URL = {https://doi.org/10.1016/j.anihpc.2020.09.004},
}

@article {Zhou,
    AUTHOR = {Zhou, Yi},
     TITLE = {Uniqueness of weak solution of the {K}d{V} equation},
   JOURNAL = {Internat. Math. Res. Notices},
  FJOURNAL = {International Mathematics Research Notices},
      YEAR = {1997},
    NUMBER = {6},
     PAGES = {271--283},
      ISSN = {1073-7928,1687-0247},
   MRCLASS = {35Q53 (35D05)},
  MRNUMBER = {1440304},
MRREVIEWER = {Vladislav\ G.\ Dubrovsky},
       DOI = {10.1155/S1073792897000202},
       URL = {https://doi.org/10.1155/S1073792897000202},
}

\end{document}